\documentclass[10pt,a4paper]{book}

\usepackage{amsthm}
\usepackage{epsfig}
\usepackage{latexsym}
\usepackage{amssymb}
\usepackage{amsfonts}
\usepackage{amsmath}
\usepackage{pb-diagram}
\usepackage{euscript}

\newcommand{\hook}{{\setlength{\unitlength}{11pt}
\begin{picture}(.833,.8)
\put(.15,.08){\line(1,0){.35}}
\put(.5,.08){\line(0,1){.5}}
\end{picture}}}

\newcommand{\GH}{\mathfrak{H}}
\newcommand{\gh}{\mathfrak{h}}
\newcommand{\gs}{\mathfrak{s}}
\newcommand{\C}{\mathbb{C}}

\newcommand{\R}{\mathbb{R}}
\newcommand{\N}{\mathbb{N}}
\renewcommand{\S}{\mathbb{S}}
\newcommand{\Z}{\mathbb{Z}}
\renewcommand{\d}{\mathrm{d}}

\newcommand{\notD}{{\rm D\hskip-0.65em / }}
\newcommand{\expnotD}{{\rm D\hskip-0.51em / }}

\newcommand{\ExtB}{{\cal B}_I}
\newcommand{\KB}{{\cal B}_I^{\scriptscriptstyle\mathrm{KBL}}}

\newcommand{\hor}{{\mathfrak H}}
\newcommand{\ts}{{t^*}} 
\newcommand{\phis}{{\varphi^*}}
\newcommand{\thetas}{{\theta^*}}
\newcommand{\st}{\mbox{}^*\hspace{-0.02in}t}
\newcommand{\sphi}{\mbox{}^* \hspace{-0.02in} \varphi}
\newcommand{\stheta}{\mbox{}^*\hspace{-0.02in} \theta}
\newcommand{\scri}{\EuScript{I}}
\newcommand{\supp}{\mathrm{supp}}

\newtheorem{theorem}{Theorem}[chapter]
\newtheorem{proposition}{Proposition}[chapter]
\newtheorem{corollary}{Corollary}[chapter]
\newtheorem{lemma}{Lemma}[chapter]
\newtheorem{remark}{Remark}[chapter]

\numberwithin{equation}{chapter}

\begin{document}

\title{Creation of fermions by rotating charged black-holes}
\author{Dietrich H\"AFNER\\
Universit\'e Bordeaux 1, Institut de
  Math\'ematiques de Bordeaux,\\ 351, cours de
  la Lib\'eration, 33405 Talence cedex, France. \\ e-mail address:
  Dietrich.Hafner@math.u-bordeaux1.fr}
\date{}  
\maketitle

{\bf Abstract~:}

This work is devoted to the mathematical study of the Hawking effect for fermions in the
setting of the collapse of a rotating charged star. We show that an observer who is located
far away from the star and at rest with respect to the Boyer Lindquist coordinates 
observes the emergence of a thermal state when his proper time goes to infinity.
We first introduce a model of the collapse of the star. We suppose that the space-time outside the star is given by the Kerr-Newman metric. The assumptions on the asymptotic behavior of the surface of the star are inspired by the asymptotic behavior of certain timelike geodesics in the Kerr-Newman metric. 
The Dirac equation is then written using coordinates and a Newman-Penrose tetrad which are adapted to the collapse.
This coordinate system and tetrad are based on the so called simple null geodesics.
The quantization of Dirac fields in a globally hyperbolic space-time is described.
We formulate and prove a theorem about the Hawking effect in this setting. The proof 
of the theorem contains a minimal velocity estimate for Dirac fields that is slightly stronger than the usual ones
and an existence and uniqueness result for solutions of a characteristic Cauchy problem for Dirac fields
in the Kerr-Newman space-time. In an appendix we construct explicitly a Penrose compactification of block $I$ 
of the Kerr-Newman space-time based on simple null geodesics.\\

{\bf R\'esum\'e~:}

Ce travail est d\'edi\'e \`a l'\'etude math\'ematique de l'effet Hawking pour des fermions dans le cadre de l'effondrement d'une \'etoile charg\'ee en rotation. On d\'emontre qu'un observateur localis\'e loin de l'\'etoile et au repos par rapport aux coordonn\'ees de Boyer-Lindquist observe l'\'emergence d'un \'etat thermal quand son temps propre tend vers l'infini. On introduit d'abord un mod\`ele de l'effondrement de l'\'etoile. On suppose que l'espace-temps \`a l'ext\'erieur de l'\'etoile est donn\'e par la m\'etrique de Kerr-Newman. Les hypoth\`eses sur le comportement asymptotique de la surface de l'\'etoile sont inspir\'ees par le comportement asymptotique de certaines g\'eod\'esiques de type temps dans la m\'etrique de Kerr-Newman. L'\'equation de Dirac est alors \'ecrite en utilisant des coordonn\'ees et une t\'etrade de Newman-Penrose adapt\'es \`a l'effondrement. Ce syst\`eme de coordonn\'ees et cette t\'etrade sont bas\'es sur des g\'eod\'esiques qu'on appelle des g\'eod\'esiques simples isotropes. La quantification des champs de Dirac dans un espace-temps globalement hyperbolique est d\'ecrite. On formule un th\'eor\`eme sur l'effet Hawking dans ce cadre. La preuve du th\'eor\`eme contient une estimation de vitesse minimale pour les champs de Dirac l\'eg\`erement plus forte que les estimations usuelles ainsi qu'un r\'esultat d'existence et d'unicit\'e pour les solutions d'un probl\`eme caract\'eristique pour les champs de Dirac dans l'espace-temps de Kerr-Newman. Dans un appendice, nous construisons explicitement la compactification de Penrose du bloc $I$ de l'espace-temps de Kerr-Newman qui est bas\'ee sur les g\'eod\'esiques simples isotropes.\\ 


83C60.\\


Hawking effect, Dirac equation, Scattering theory, Characteristic Cauchy problem.
\tableofcontents
\chapter{Introduction}
It was in 1975 that S. W. Hawking published his famous paper about the creation of particles 
by black holes (see \cite{Haw}). Later this effect was analyzed by other authors in more 
detail ( see e.g. \cite{Wa1}) and we can say that the effect was well understood from
a physical point of view at the end of the 1970's. 
From a mathematical point of view, however, fundamental questions linked to the Hawking radiation
such as scattering theory for field equations on black-hole space-times had 
not been addressed at that time.

In the early 1980's Dimock and Kay started a research programme concerning scattering 
theory on curved space-times. They obtained an asymptotic completeness result for classical
and quantum massless scalar fields on the Schwarzschild metric (see \cite{DK1}-\cite{DK3}).
Their work was pushed further by Alain Bachelot in the 1990's. He showed asymptotic completeness
for Maxwell and Klein-Gordon fields (see \cite{Ba1}, \cite{Ba2}) and gave a mathematically precise
description of the Hawking effect (see \cite{Ba3}-\cite{Ba5}) in the spherically symmetric case. Meanwhile other authors
contributed to the subject such as Nicolas (see \cite{Ni1}), Jin (see \cite{Ji1})
and Melnyk (see \cite{Me1}, \cite{Me2}). All these works deal with the spherically symmetric case.

The more realistic case of a rotating black hole is more difficult. In the spherically symmetric case, the
study of a field equation can be reduced to the study of a $1+1$ 
dimensional equation with potential. In the Kerr case this reduction is no longer possible
and the methods used in the papers cited so far do not apply.
A paper by De Bi\`evre, Hislop, Sigal using different methods appeared in 1992
(see \cite{DHS}). By means of a Mourre estimate they show asymptotic completeness for the wave equation on non-compact
Riemannian manifolds; possible applications are therefore static situations such as the Schwarzschild case,
which they treat, but the Kerr geometry is not even stationary. In this context we also mention the paper
of Daud\'e about the Dirac equation in the Reissner-Nordstr\"om metric (see \cite{Da2}). A complete scattering theory for
the wave equation on stationary, asymptotically flat space-times, was obtained
by the author in 2001 (see \cite{Hae1}). To our knowledge the first asymptotic 
completeness result in the Kerr case was obtained by the author in \cite{Hae2},
for the non superradiant modes of the Klein-Gordon field. The first complete scattering
theory for a field equation in the Kerr metric was obtained by Nicolas and the author in \cite{HN} for massless
Dirac fields. This result was generalized by Daud\'e in \cite{Da2} to the massive charged Dirac field in
the Kerr-Newman metric. All these papers use Mourre theory.

The aim of the present paper is to give a mathematically precise description of the Hawking 
effect for spin 1/2 fields in the setting of the collapse of a rotating charged star.
We show that an observer who is located far away from the black hole and at rest 
with respect to the Boyer-Lindquist coordinates observes the emergence of a thermal state 
when his proper time $t$ goes to infinity. Let us give an idea of the
theorem describing the effect. Let $r_*$ be the Regge-Wheeler coordinate. We suppose
that the boundary of the star is described by $r_*=z(t,\theta)$. The space-time is then given by
\[ {\cal M}_{col}=\bigcup_t\Sigma^{col}_t,\, \Sigma^{col}_t=\{(t,r_*,\omega)\in \R_t\times \R_{r_*}\times S^2;\, r_*\ge z(t,\theta)\}. \]
The typical asymptotic behavior of $z(t,\theta)$ is ($\kappa_+>0$) :
\[z(t,\theta)=-t-A(\theta)e^{-2\kappa_+ t}+B(\theta)+{\mathcal O}(e^{-4\kappa_+ t}),\, t\rightarrow\infty. \]
Let ${\cal H}_t=L^2((\Sigma^{col}_t,\mathrm{dVol});\C^4)$. The Dirac equation can be written as
\begin{eqnarray}
\label{1.1}
\partial_t\Psi=i\notD_t\Psi+\quad\mbox{boundary condition.}
\end{eqnarray}
We will put a MIT boundary condition on the surface of the star. The evolution of 
the Dirac field is then described by an isometric propagator $U(s,t):{\cal H}_s\rightarrow{\cal H}_t$.
The Dirac equation on the whole exterior Kerr-Newman space-time ${\cal M}_{BH}$ will be written as 
\begin{eqnarray*}
\partial_t\Psi=i\notD\Psi.
\end{eqnarray*}
Here $\notD$ is a selfadjoint operator on ${\cal H}=L^2((\R_{r_*}\times S^2,dr_*d\omega);\C^4)$.
There exists an asymptotic velocity operator $P^{\pm}$ s.t. for all continuous functions
$J$ with $\lim_{|x|\rightarrow\infty}J(x)=0$ we have 
\[ J(P^{\pm})=s-\lim_{t\rightarrow\pm\infty}e^{-it\expnotD}J\left(\frac{r_*}{t}\right)e^{it\expnotD}. \]
Let ${\cal U}_{col}({\cal M}_{col})$ resp. ${\cal U}_{BH}({\cal M}_{BH})$ be the algebras of observables outside the collapsing body resp. on 
the space-time describing the eternal black-hole generated by $\Psi^*_{col}(\Phi_1)\Psi_{col}(\Phi_2)$ resp. $\Psi^*_{BH}(\Phi_1)\Psi_{BH}(\Phi_2)$. Here $\Psi_{col}(\Phi)$ resp. $\Psi_{BH}(\Phi)$ are the quantum spin fields on ${\cal M}_{col}$ resp. ${\cal M}_{BH}$. Let $\omega_{col}$ be a vacuum state on ${\cal U}_{col}({\cal M}_{col})$; $\omega_{vac}$
a vacuum state on ${\cal U}_{BH}({\cal M}_{BH})$ and $\omega_{Haw}^{\eta,\sigma}$ be a KMS-state on ${\cal U}_{BH}({\cal M}_{BH})$ with inverse temperature $\sigma>0$ and chemical potential $\mu=e^{\sigma \eta}$ (see Chapter \ref{sec4} for details). For a function $\Phi\in C_0^{\infty}({\cal M}_{BH})$ we define :
\begin{eqnarray*}
\Phi^T(t,r_*,\omega)=\Phi(t-T,r_*,\omega).
\end{eqnarray*}
The theorem about the Hawking effect is the following : 
\begin{theorem}[Hawking effect]
Let 
\[ \Phi_j\in (C_0^{\infty}({\cal M}_{col}))^4,\, j=1,2.\]
Then we have
\begin{eqnarray}
\label{HQ}
\lefteqn{\lim_{T\rightarrow\infty}\omega_{col}(\Psi^*_{col}(\Phi^T_1)\Psi_{col}(\Phi^T_2))}\nonumber\\
&=&\omega_{Haw}^{\eta,\sigma}(\Psi^*_{BH}({\bf 1}_{\R^+}(P^-)\Phi_1)\Psi_{BH}({\bf 1}_{\R^+}(P^-)\Phi_2))\nonumber\\
&+&\omega_{vac}(\Psi^*_{BH}({\bf 1}_{\R^-}(P^-)\Phi_1)\Psi_{BH}({\bf 1}_{\R^-}(P^-)\Phi_2)),\\
T_{Haw}&=&1/\sigma=\kappa_+/2\pi,\quad \mu=e^{\sigma\eta},\, \eta=\frac{qQr_+}{r_+^2+a^2}+\frac{aD_{\varphi}}{r_+^2+a^2}.\nonumber
\end{eqnarray}
\end{theorem} 
Here $q$ is the charge of the field, $Q$ the charge of the black-hole, $a$ the angular momentum per unit mass of the black-hole, $r_+=M+\sqrt{M^2-(a^2+Q^2)}$ defines the outer event horizon and $\kappa_+$ is the surface gravity of this horizon. The interpretation of (\ref{HQ}) is the following. We start with a vacuum state which we evolve in the proper time of an observer at rest with respect to the Boyer Lindquist coordinates. The limit when the proper time of this observer goes to infinity is a thermal state coming from the event horizon in formation and a vacuum state coming from infinity as expressed on the R.H.S of (\ref{HQ}).
The Hawking effect is often interpreted in terms of particles, the antiparticle falling into the black hole and the particle escaping to infinity. From our point of view this interpretation is somewhat misleading. The effect really comes from an infinite Doppler effect and the mixing of positive and negative frequencies. To explain this a little bit more we describe the analytic problem behind the effect. 
 Let $f(r_*,\omega)\in C_0^{\infty}(\R\times S^2)$.
The key result about the Hawking effect is :
\begin{eqnarray}
\label{1.2}
\lim_{T\rightarrow\infty}||{\bf 1}_{{[}0,\infty)}(\notD_0)U(0,T)f||^2_0
&=&\langle {\bf 1}_{\R^+}(P^-)f,\mu e^{\sigma \expnotD}(1+\mu e^{\sigma\expnotD})^{-1}{\bf 1}_{\R^+}(P^-)f\rangle\nonumber\\
&+&||{\bf 1}_{{[}0,\infty)}(\notD){\bf 1}_{\R^-}(P^-)f||^2,
\end{eqnarray}
where $\mu,\, \eta,\, \sigma$ are as in the above theorem. Equation (\ref{1.2}) implies (\ref{HQ}).

\begin{figure}
\centering\epsfig{figure=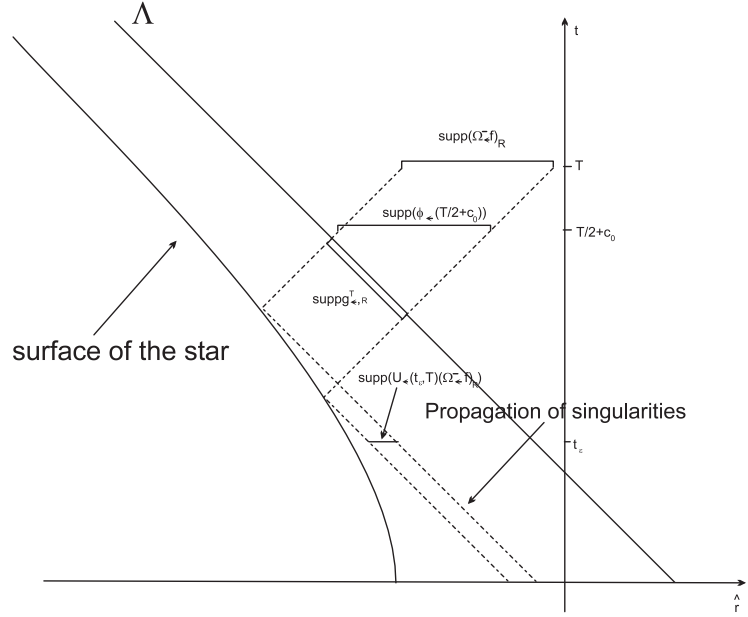,width=10cm}
\caption{The collapse of the star}
\label{figure1}
\end{figure}

The term on the L.H.S. comes from the vacuum state we consider. We have to project on the positive frequency solutions (see Chapter \ref{sec4} for details). Note that in (\ref{1.2}) we consider the time reversed evolution. This comes from
the quantization procedure. When time becomes large the solution hits the surface of the star
at a point closer and closer to the future event horizon. Figure \ref{figure1} shows the
situation for an asymptotic comparison dynamics, which satisfies Huygens' principle.
For this asymptotic comparison dynamics the support of the solution concentrates more and more
when time becomes large, which means that the frequency increases.
The consequence of the change in frequency is that the system does not stay in the vacuum state.

We conclude this introduction with some comments on the boson case which we do not treat 
in this paper. This case is more difficult because of the superradiance phenomenon. There exists no
positive conserved energy for the wave equation in block $I$ of the Kerr metric. This is 
linked to the fact that the Kerr metric is not stationary outside the black hole. 
Because of the difficulty linked to superradiance, there is at present no complete scattering theory for the wave equation on the Kerr metric, a necessary prerequisite for the mathematical
description of the Hawking effect. However some progress in this direction has been made by
Finster, Kamran, Smoller and Yau who obtained an integral representation for the propagator of the wave equation
on the Kerr metric (see \cite{FKSY}). We also refer to \cite{Ba7} for scattering results in 
a superradiant situation.
\vspace{0.5cm}

{\bf{\Large Notations}}

Let $({\cal M},g)$ be a smooth 4-manifold equipped with a lorentzian metric $g$ with signature 
$(+,-,-,-)$. We denote by $\nabla_a$ the Levi-Civita connection on $({\cal M},g)$. 

Many of our equations will be expressed using the two-component spinor notations and abstract index
formalism of R. Penrose and W. Rindler \cite{PR1}.

Abstract indices are
denoted by light face latin letters, capital for spinor indices and
lower case for tensor indices. Abstract indices are a notational
device for keeping track of the nature of objects in the course of
calculations, they do not imply any reference to a coordinate basis,
all expressions and calculations involving them are perfectly
intrinsic. For example, $g_{ab}$ will refer to the space-time metric
as an intrinsic symmetric tensor field of valence $\left
[ \begin{array}{cc} 0 \\ 2 \end{array} \right]$, i.e. a section of
$\mathrm{T}^* {\cal M} \odot \mathrm{T}^* {\cal M}$ and $g^{ab}$ will
refer to the inverse metric as an intrinsic symmetric tensor field of
valence $\left [ \begin{array}{cc} 2 \\ 0 \end{array} \right]$, i.e. a
section of $\mathrm{T} {\cal M} \odot \mathrm{T} {\cal M}$ (where
$\odot$ denotes the symmetric tensor product, $\mathrm{T} {\cal M}$
the tangent bundle to our space-time manifold $\cal M$ and
$\mathrm{T}^* {\cal M}$ its cotangent bundle).

Concrete indices defining components in reference to a basis are
represented by bold face latin letters. Concrete spinor indices,
denoted by bold face capital latin letters, take their values
in $\{ 0,1 \}$ while concrete tensor indices, denoted by bold face
lower case latin letters, take their values in $\{ 0,1,2,3
\}$. Consider for example a basis of $\mathrm{T} {\cal M}$, that is a
family of four smooth vector fields on $\cal M$~: ${\cal B} = \left\{
e_0 ,  e_1, e_2, e_3 \right\}$ such that at each point $p$ of $\cal
M$ the four vectors $e_0 (p) , e_1 (p), e_2 (p), e_3 (p)$ are linearly
independent, and the corresponding dual basis of $\mathrm{T}^* {\cal
M}$~: ${\cal B}^* = \left\{ e^0 , e^1, e^2, e^3 \right\}$ such that
$e^\mathbf{a} \left( e_\mathbf{b} \right) =
\delta^\mathbf{a}_\mathbf{b}$, $\delta^\mathbf{a}_\mathbf{b}$ denoting
the Kronecker symbol~; $g_\mathbf{ab}$ will refer to the components of
the metric $g_{ab}$ in the basis $\cal B$~: $g_\mathbf{ab} =
g(e_\mathbf{a} , e_\mathbf{b} )$ and $g^\mathbf{ab}$ will denote the
components of the inverse metric $g^{ab}$ in the dual basis ${\cal
B}^*$, i.e. the $4 \times 4$ real symmetric matrices $\left
( g_\mathbf{ab} \right)$ and $\left ( g^\mathbf{ab} \right)$ are the
inverse of one another. In the abstract index formalism, the basis
vectors $e_\mathbf{a}$, $\mathbf{a} = 0,1,2,3$, are denoted
${e_\mathbf{a}}^a$ or ${g_\mathbf{a}}^a$. In a coordinate basis, the
basis vectors $e_\mathbf{a}$ are coordinate vector fields and
will also be denoted by $\partial_\mathbf{a}$ or
$\frac{\partial}{\partial x^\mathbf{a}}$~; the dual basis covectors
$e^\mathbf{a}$ are coordinate 1-forms and will be denoted by $\d
x^\mathbf{a}$.

We adopt Einstein's convention for the same index appearing twice,
once up, once down, in the same term. For concrete indices, the sum is
taken over all the values of the index. In the case of abstract
indices, this signifies the contraction of the index, i.e. $f_a V^a$
denotes the action of the 1-form $f_a$ on the vector field $V^a$.

For a manifold
$Y$ we denote by $C_b^{\infty}(Y)$ the set of all $C^{\infty}$ functions
on $Y$, that are bounded together with all their derivatives. We
denote by $C_{\infty}(Y)$ the set of all continuous functions tending
to zero at infinity.
\vspace{0.5cm}

{\Large{\bf Acknowledgments}}

The author warmly thanks A. Bachelot, J.-F. Bony and J.-P. Nicolas for fruitful discussions.
This work was partially supported by the ANR project JC0546063 "Equations hyperboliques dans les espaces-temps de
la relativit\'e g\'en\'erale : Diffusion et r\'esonances."
\\

\chapter{Strategy of the proof and organization of the article}
\label{secstor}
\section{The analytic problem}
\begin{figure}
\centering\epsfig{figure=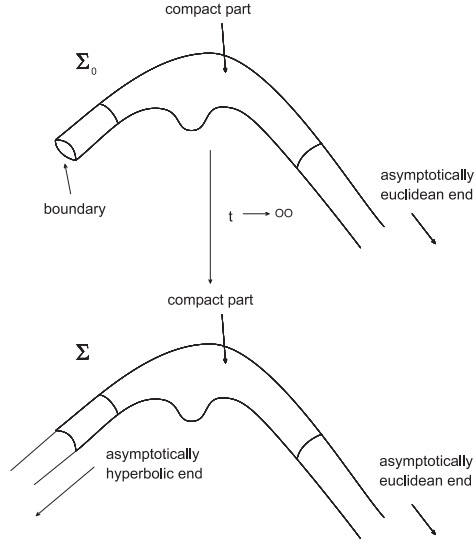,width=8cm}
\caption{The manifold at time $t=0$ $\Sigma_0$ and the limit manifold $\Sigma$.}
\label{movbound}
\end{figure}

Let us consider a model, where the eternal black-hole is described by a static space-time (although the Kerr-Newman
space-time is not even stationary, the problem will be essentially reduced to this kind of situation). Then the problem can be described as follows. Consider a riemannian manifold $\Sigma_0$ with one asymptotically euclidean end and a boundary. The boundary will move when $t$ becomes large asymptotically with the speed of light. The manifold at time $t$ is denoted $\Sigma_t$. The "limit" manifold $\Sigma$ is a manifold with two ends, one asymptotically euclidean and the other asymptotically 
hyperbolic (see Figure \ref{movbound}). The problem consists in evaluating the limit 
\begin{eqnarray*}
\lim_{T\rightarrow\infty}||{\bf 1}_{{[}0,\infty)}(\notD_0)U(0,T)f||_0,
\end{eqnarray*}
where $U(0,T)$ is the isometric propagator for the Dirac equation on the manifold with 
moving boundary and suitable boundary conditions. It is worth noting that the underlying scattering theory is not the scattering theory for the problem with moving boundary
but the scattering theory on the "limit" manifold. It is largely believed that the result does not depend on the boundary condition. We will show in this paper that it does not depend on the chiral angle in the MIT boundary condition. Note also that the boundary viewed in $\bigcup_t \{t\}\times\Sigma_t$ is only weakly timelike, a problem that has been rarely considered (but see \cite{Ba4}).  

One of the problems for the description of the Hawking effect is to derive a reasonable 
model for the collapse of the star. We will suppose that the metric outside the collapsing star
is always given by the Kerr-Newman metric. Whereas this is a genuine assumption in 
the rotational case, in the spherically symmetric case Birkhoffs theorem assures that the metric outside the star is the Reissner-Nordstr\"om
metric. We will suppose that a point on the surface of the star will move along a curve which behaves asymptotically like a
timelike geodesic with $L={\cal Q}=\tilde{E}=0$, where $L$ is the angular momentum,
$\tilde{E}$ the rotational energy and ${\cal Q}$ the Carter constant.  
The choice of geodesics is justified by the fact that the collapse creates the space-time, i.e. angular momenta 
and rotational energy should be zero with respect to the space-time.
We will need an additional asymptotic condition on the collapse. 
It turns out that there is a natural coordinate system $(t,\hat{r},\omega)$ associated to the collapse. In this coordinate
system the surface of the star is described by $\hat{r}=\hat{z}(t,\theta)$.
We need to assume the existence of a constant $C$ s.t.
\begin{eqnarray}
\label{1.3}
|\hat{z}(t,\theta)+t+C|\rightarrow 0,\, t\rightarrow \infty.
\end{eqnarray}
It can be checked that this asymptotic condition is fulfilled if we use the above geodesics for some appropriate initial condition. On the one hand we are not able to compute this initial condition explicitly, on the other hand it seems more natural to impose a (symmetric) asymptotic condition than an initial condition. If we would allow in (\ref{1.3}) a function $C(\theta)$ rather than a constant, the problem would become more difficult. Indeed one of the problems for treating the Hawking radiation in the rotational case is the high frequencies of the solution. In contrast with the spherically symmetric case, the difference between the Dirac operator and an operator with constant coefficients is near the horizon always a differential operator of order one
\footnote{In the spherically symmetric case we can diagonalize the operator. After diagonalization
the difference is just a potential.}. This explains that in the 
high energy regime we are interested in, the Dirac operator is not close to a constant
coefficient operator. Our method to prove (\ref{1.2}) is to use scattering arguments
to reduce the problem to a problem with a constant coefficient operator, for which we can compute 
the radiation explicitly. If we do not impose a condition of type (\ref{1.3}), then in all
coordinate systems the solution has high frequencies, in the radial as well as in the angular
directions. With condition (\ref{1.3}) these high frequencies only occur in the radial
direction. Our asymptotic comparison dynamics will differ from the real dynamics only 
by derivatives in angular directions and by potentials.

\section{Strategy of the proof}   
In this section we will give some ideas of the proof of (\ref{1.2}). We want to reduce the problem to the evaluation 
of a limit that can be explicitly computed. To do so, we use the asymptotic completeness
results obtained in \cite{HN} and \cite{Da2}. There exists a constant coefficient operator 
$\notD_{\leftarrow}$ s.t. the following limits exist :
\begin{eqnarray*}
W_{\leftarrow}^{\pm}&:=&s-\lim_{t\rightarrow\pm\infty}e^{-it\expnotD}e^{it\expnotD_{\leftarrow}}{\bf 1}_{\R^{\mp}}(P_{\leftarrow}^{\pm}),\\
\Omega_{\leftarrow}^{\pm}&:=&s-\lim_{t\rightarrow\pm\infty}e^{-it\expnotD_{\leftarrow}}e^{it\expnotD}{\bf 1}_{\R^{\mp}}(P^{\pm}).
\end{eqnarray*}
Here $P_{\leftarrow}^{\pm}$ is the asymptotic velocity operator associated to the dynamics 
$e^{it\expnotD_{\leftarrow}}$. Then the R.H.S. of (\ref{1.2}) equals :
\begin{eqnarray*}
||{\bf 1}_{{[}0,\infty)}(\notD){\bf 1}_{\R^-}(P^-)f||^2
+\langle \Omega_{\leftarrow}^-f,\mu e^{\sigma \expnotD_{\leftarrow}}(1+\mu e^{\sigma\expnotD_{\leftarrow}})^{-1}\Omega_{\leftarrow}^-f\rangle.
\end{eqnarray*}
The aim is to show that the incoming part is :
\[ \lim_{T\rightarrow\infty}||{\bf 1}_{{[}0,\infty)}(D_{\leftarrow,0})U_{\leftarrow}(0,T)\Omega_{\leftarrow}^-f||^2_0
=\langle \Omega_{\leftarrow}^-f,\mu e^{\sigma \expnotD_{\leftarrow}}(1+\mu e^{\sigma\expnotD_{\leftarrow}})^{-1}\Omega_{\leftarrow}^-f\rangle, \]
where the equality can be shown by explicit calculation. Here $\notD_{\leftarrow,t}$ and $U_{\leftarrow}(s,t)$
are the asymptotic operator with boundary condition and the associated propagator. 
The outgoing part is easy to treat. 

As already mentioned, we have to consider the solution in a high frequency regime. 
Using the Regge-Wheeler variable as a position variable and, say, the Newman-Penrose 
tetrad used in \cite{HN} we find that the modulus of the local velocity 
\[ {[}ir_*,\notD{]}=h^2(r_*,\omega)\Gamma^1 \]
is not equal to 1, whereas the asymptotic dynamics must have constant
local velocity. Here $h$ is a continuous function and $\Gamma^1$ a constant matrix. Whereas the $(r_*,\omega)$ coordinate system and the tetrad used in \cite{HN}
were well adapted to the time dependent scattering theory developed in \cite{HN}, they
are no longer well adapted when we consider large times and high frequencies. We are therefore looking 
for a variable $\hat{r}$ s.t.
\[ \{(t,\hat{r},\omega); \hat{r}\pm t = const. \} \]
are characteristic surfaces. By a separation of variables Ansatz we find a family of such variables
and we choose the one which is well adapted to the collapse of the star in the sense
that along an incoming null geodesic with $L={\cal Q}=0$ we have :
\[ \frac{\partial\hat{r}}{\partial t}=-1. \]
This variable turns out to be a generalized Bondi-Sachs variable. The null geodesics 
with $L={\cal Q}=0$ are generated by null vector fields $N^{\pm}$ that we choose to
be $l$ and $n$ in the Newman-Penrose tetrad. If we write down the hamiltonian
for the Dirac equation with this choice of coordinates and tetrad we find that the local
velocity now has modulus $1$ everywhere and our initial problem disappears. The new hamiltonian is
again denoted $\notD$. Let $\notD_{\leftarrow}$ be an asymptotic comparison dynamics near the
horizon with constant coefficients. Note that (\ref{1.2}) is of course independent of the choice 
of the coordinate system and the tetrad, i.e. both sides of (\ref{1.2}) are independent
of these choices. We now proceed as follows :
\begin{enumerate}
\item
We decouple the problem at infinity from the problem near the horizon by cut-off 
functions. The problem at infinity is easy to treat.  
\item We consider $U(t,T)f$ on a characteristic hypersurface $\Lambda$.
The resulting characteristic data is denoted $g^T$.  
We will approximate $\Omega_{\leftarrow}^-f$ by a function $(\Omega_{\leftarrow}^-f)_R$ 
with compact support and higher regularity in the angular derivatives. Let $U_{\leftarrow}(s,t)$
be the isometric propagator associated to the asymptotic hamiltonian $\notD_{\leftarrow}$ 
with MIT boundary conditions. We also consider $U_{\leftarrow}(t,T)(\Omega_{\leftarrow}^-f)_R$ on $\Lambda$.
The resulting characteristic data is denoted $g_{\leftarrow,R}^T$. The situation for the 
asymptotic comparison dynamics is shown in Figure \ref{figure1}.
\item We solve a characteristic Cauchy problem for the Dirac equation with data $g_{\leftarrow,R}^T$.
The solution at time zero can be written in a region near the boundary as 
\[ G(g^T_{\leftarrow,R})=U(0,T/2+c_0)\Phi(T/2+c_0), \]
where $\Phi$ is the solution of a characteristic Cauchy problem in the whole space (without the
star). 
The solutions of the characteristic problems for the asymptotic hamiltonian are written
in a similar way and denoted respectively $G_{\leftarrow}(g^T_{\leftarrow,R})$ and $\Phi_{\leftarrow}$. 
\item Using the asymptotic completeness result we show that $g^T-g^T_{\leftarrow,R}
\rightarrow 0$ when $T,R\rightarrow\infty$. By continuous dependence on the characteristic data
we see that :
\[ G(g^T)-G(g_{\leftarrow,R}^T)\rightarrow 0, T,R\rightarrow \infty. \]
\item We write 
\begin{eqnarray*}
G(g_{\leftarrow,R}^T)-G_{\leftarrow}(g^T_{\leftarrow,R})&=&U(0,T/2+c_0)(\Phi(T/2+c_0)-\Phi_{\leftarrow}(T/2+c_0))\\
&+&(U(0,T/2+c_0)-U_{\leftarrow}(0,T/2+c_0))\Phi_{\leftarrow}(T/2+c_0). 
\end{eqnarray*}
The first term becomes small near the boundary when $T$ becomes large. We then note that for all $\epsilon>0$ 
there exists $t_{\epsilon}>0$ s.t. 
\[ ||(U(t_{\epsilon},T/2+c_0)-U_{\leftarrow}(t_{\epsilon},T/2+c_0))\Phi_{\leftarrow}(T/2+c_0)||<\epsilon \]
uniformly in $T$ large. The function $U_{\leftarrow}(t_{\epsilon},T/2+c_0)\Phi_{\leftarrow}(T/2+c_0)$
will be replaced by a geometric optics approximation $F^T_{t_{\epsilon}}$ which has the following properties :
\begin{eqnarray}
\label{1.4}
\supp F^T_{t_{\epsilon}}\subset(-t_{\epsilon}-|{\cal O}(e^{-\kappa_+ T})|,-t_{\epsilon}),\\
\label{1.5}
F^T_{t_{\epsilon}}\rightharpoonup 0,\, T\rightarrow \infty,\\
\label{1.6}
\forall \lambda>0\quad Op(\chi(\langle \xi \rangle\le \lambda \langle q \rangle))F^T_{t_{\epsilon}}\rightarrow 0,\, T\rightarrow \infty.
\end{eqnarray}
Here $\xi$ and $q$ are the dual coordinates to $\hat{r},\theta$ respectively.
\item We show that for $\lambda$ sufficiently large possible singularities of  $Op(\chi(\langle \xi \rangle\ge \lambda \langle q \rangle))F^T_{t_{\epsilon}}$
are transported by the group $e^{-it_{\epsilon}\expnotD}$ in such a way that they always stay away from
the surface of the star. 
\item From the points 1. to 5. follows :
\[ \lim_{T\rightarrow\infty}||{\bf 1}_{{[}0,\infty)}(\notD_0)j_-U(0,T)f||_0^2=\lim_{T\rightarrow\infty}
||{\bf 1}_{{[}0,\infty)}(\notD_0)U(0,t_{\epsilon})F^T_{t_{\epsilon}}||_0^2, \]
where $j_-$ is a smooth cut-off which equals $1$ near the boundary and $0$ at infinity.
Let $\phi_{\delta}$ be a cut-off outside the surface of the star at time $0$. 
If $\phi_{\delta}=1$ sufficiently close to the surface of the star at time $0$ we see
by the previous point that  
\begin{eqnarray}
\label{1.7b}
(1-\phi_{\delta})e^{-it_{\epsilon}\expnotD}F^T_{t_{\epsilon}}\rightarrow 0,\, T\rightarrow\infty. 
\end{eqnarray}
Using (\ref{1.7b}) we show that (modulo a small error term): 
\[ (U(0,t_{\epsilon})-\phi_{\delta}e^{-it_{\epsilon}\expnotD})F^T_{t_{\epsilon}}\rightarrow 0, \, T\rightarrow \infty. \]
Therefore it remains to consider :
\[ \lim_{T\rightarrow\infty}||{\bf 1}_{{[}0,\infty)}(\notD_0)\phi_{\delta}e^{-it_{\epsilon}\expnotD}F^T_{t_{\epsilon}}||_0. \]
\item We show that we can replace ${\bf 1}_{{[0},\infty)}(\notD_0)$ by ${\bf 1}_{{[0},\infty)}(\notD)$. 
This will essentially allow to commute the energy cut-off and the group. We then show that we can replace
the energy cut-off by  ${\bf 1}_{{[0},\infty)}(\notD_{\leftarrow})$. We end up with :
\begin{eqnarray}
\label{1.7}
\lim_{T\rightarrow\infty}||{\bf 1}_{{[0},\infty)}(\notD_{\leftarrow})e^{-it_{\epsilon}\expnotD_{\leftarrow}}F^T_{t_{\epsilon}}||.
\end{eqnarray}
\item We compute the limit in (\ref{1.7}) explicitly.
\end{enumerate}
\section{Organization of the article}
The paper is organized as follows :
\begin{itemize}
\item In Chapter \ref{sec2} we present the model of the collapsing star. We first analyze 
the geodesics in the Kerr-Newman space-time and explain how the Carter constant can be
understood in terms of the hamiltonian flow. We construct the variable $\hat{r}$ and show that
\[ \frac{\partial{\hat{r}}}{\partial t}=\pm 1 \,\mbox{along null geodesics with}\, L={\cal Q}=0. \]
We then show that in the $(t,\hat{r},\omega)$ coordinate system we have along incoming timelike geodesics with 
$L={\cal Q}=\tilde{E}=0$ :
\begin{eqnarray}
\label{GAS}
\hat{r}=-t-\hat{A}(\theta,r_0)e^{-2\kappa_+ t}+\hat{B}(\theta,r_0)+{\cal O}(e^{-4\kappa_+ t})
\end{eqnarray}
with $\hat{A}(\theta,r_0)>0$. Our assumption will be that a point on the surface behaves asymptotically like (\ref{GAS}) with $\hat{B}(\theta_0,r_0(\theta_0))=const.$ Here $r_0(\theta_0)$ is a function defining the surface at time $t=0$.
\item In Chapter \ref{sec3} we describe classical Dirac fields. We introduce a new Newman-Penrose tetrad
and compute the new expression of the equation. New asymptotic hamiltonians are introduced 
and classical scattering results are obtained from scattering results in \cite{HN} and
\cite{Da2}. The MIT boundary condition is discussed in detail.  
\item Dirac quantum fields are discussed in Chapter \ref{sec4}. We first present the second quantization
of Dirac fields and then describe the quantization in a globally hyperbolic space-time.
The theorem about the Hawking effect is formulated and discussed in Section \ref{sec4.3}.
\item In Chapter \ref{sec7} we show additional scattering results that we will need later.
A minimal velocity estimate slightly stronger than the usual ones is established.
\item In Chapter \ref{sec8} we solve the characteristic problem for the Dirac equation. 
We approximate the characteristic surface by smooth spacelike hypersurfaces and recover
the solution in the limit. This method is close to that used by H\"ormander in \cite{Hoe2} for the wave equation.
\item Chapter \ref{Red} contains several reductions of the problem. We show that (\ref{1.2}) 
implies the theorem about the Hawking effect. We use the axial symmetry to fix the angular
momentum. Several technical results are collected.
\item Chapter \ref{sec9} is devoted to the comparison of the dynamics on the interval 
${[}t_{\epsilon},T{]}$. 
\item In Chapter \ref{sec10} we study the propagation of singularities for the Dirac equation 
in the Kerr-Newman metric. We show that "outgoing" singularities located in 
\[ \{(\hat{r},\omega,\xi,q);\, \hat{r}\ge -t_{\epsilon}-C^{-1},\, |\xi|\ge C |q| \} \]
stay away from the surface of the star for $C$ large. 
\item The main theorem is proven in Chapter \ref{mainth}.
\item Appendix \ref{AppA} contains the proof of the existence and uniqueness 
of solutions of the Dirac equation in the space-time describing the collapsing star.
\item In Appendix \ref{AppB} we show that we can compactify the block $I$ of the Kerr-Newman
space-time using null geodesics with $L={\cal Q}=0$ instead of principal null 
geodesics.
\end{itemize}

\chapter{The model of the collapsing star}
\label{sec2}
The purpose of this chapter is to describe the model of the collapsing star. We will suppose
that the metric outside the star is given by the Kerr-Newman metric, which is discussed in Section \ref{sec2.1}.
Geodesics are discussed in Section \ref{sec2.1.2}. We give a description of the Carter constant in terms of
the associated hamiltonian flow. A new position variable is introduced. In Section 
\ref{sec2.2} we give the precise asymptotic behavior of the boundary of the star using this new position variable.
We require that a point on the surface behaves asymptotically like incoming timelike geodesics with $L={\cal Q}=\tilde{E}=0$,
which are studied in Section \ref{secE.3.2.1}. The precise assumptions are given in Section \ref{secE.3.2.2}.
\\

\section{The Kerr-Newman metric}
\label{sec2.1}
We give a brief description of the Kerr-Newman metric, 
which describes an eternal rotating charged black-hole. 
A detailed description can be found e.g. in \cite{Wa}.
\subsection{Boyer-Lindquist coordinates}
\label{sec2.1.1}
In Boyer-Lindquist coordinates, a Kerr-Newman black-hole is described by 
a smooth 4-dimensional lorentzian manifold ${\cal M}_{BH} =
 \R_t \times \R_r \times S^2_\omega$, whose space-time metric $g$ and 
electromagnetic vector potential $\Phi_a$ are given by :

\begin{eqnarray}
\label{KerrNewBL}
g &=& \left( 1 + \frac{Q^2-2Mr}{\rho^2} \right)
\d t^2 + \frac{2a\sin^2\theta(2Mr-Q^2)}{\rho^2}
\, \d t \d\varphi - \frac{\rho^2}{\Delta} \d r^2  - \rho^2 \d\theta^2\nonumber\\
&-& \frac{\sigma^2}{\rho^2} \sin^2 \theta \, \d\varphi^2 ,
\end{eqnarray}
\begin{eqnarray*}
\rho^2 &=& r^2 + a^2 \cos^2 \theta \, ,~ \Delta = r^2 -2Mr +a^2+Q^2,\nonumber\\
\sigma^2 &=& \left( r^2 +a^2 \right) \rho^2 +(2Mr-Q^2)a^2 \sin^2\theta
=(r^2+a^2)^2-a^2\Delta\sin^2\theta,\nonumber\\
\Phi_a dx^a&=&-\frac{Qr}{\rho^2}(dt-a\sin^2\theta d\varphi).\nonumber
\end{eqnarray*}

Here $M$ is the mass of the black hole, $a$ its angular
momentum per unit mass and $Q$ the charge of the black-hole. If $Q=0$, $g$ reduces 
to the Kerr metric, and if $Q=a=0$ we recover the Schwarzschild metric. 
The expression (\ref{KerrNewBL}) of the Kerr metric has two types of
singularities. While the set of points $\{ \rho^2 =0\}$ (the equatorial
ring $\{ r=0 \, ,~\theta = \pi/2 \}$ of the $\{ r=0 \}$ sphere) is a
true curvature singularity, the spheres where $\Delta$ vanishes,
called horizons, are mere coordinate singularities. We will consider in this 
paper subextremal Kerr-Newman space-times, that is we suppose $Q^2+a^2<M^2$. In this 
case $\Delta$ has two real roots:
\begin{equation}
\label{2.3}
r_\pm = M \pm \sqrt{M^2 - (a^2+Q^2)}.
\end{equation}

The spheres $\{ r=r_- \}$ and $\{
r=r_+ \}$ are called event horizons.
The two horizons separate ${\cal M}_{BH}$ into three connected
components called Boyer-Lindquist blocks~: $B_I,B_{II},B_{III} 
(r_+<r, r_-<r<r_+, r<r_-)$.
No Boyer-Lindquist block is stationary, that is to say there exists no globally defined timelike
Killing vector field on any given block. In particular, block $I$
contains a toroidal region, called the ergosphere, surrounding the
horizon,
\begin{eqnarray}
\label{2.4b}
 {\cal E} = \left\{ (t,r,\theta ,\varphi )\, ;~ r_+ <r < M
  +\sqrt{M^2 -Q^2-a^2 \cos^2 \theta } \right\}, 
\end{eqnarray}
where the vector $\partial /\partial t$ is spacelike.

An important feature of the Kerr-Newman space-time is that it has Petrov type D (see 
e.g. \cite{ONe}). This means that the Weyl tensor has two double roots at each point.
These roots, referred to as the principal null directions of the Weyl tensor, are given
by the two vector fields 
\begin{eqnarray*}
V^{\pm}=\frac{r^2+a^2}{\Delta}\partial_t\pm\partial_r+\frac{a}{\Delta}\partial_{\varphi}.
\end{eqnarray*}
Since $V^+$ and $V^-$ are twice repeated null directions of the Weyl tensor, by the
Goldberg-Sachs theorem (see for example \cite[Theorem 5.10.1]{ONe}) their integral curves define
shear-free null geodesic congruences. We shall refer to the integral curves of $V^+$ (respectively $V^-$)
as the outgoing (respectively incoming) principal null geodesics and write from now on
PNG for principal null geodesic. The plane determined at each point by the two prinipal null directions 
is called the principal plane. 

We will often use a Regge-Wheeler type coordinate $r_*$ in $B_I$ 
instead of $r$ (see e.g. \cite{Ch}), which is given by
\begin{eqnarray}
\label{2.5b}
r_*=r+\frac{1}{2\kappa_+}\ln|r-r_+|-\frac{1}{2\kappa_-}\ln|r-r_-|+R_0,
\end{eqnarray}
where $R_0$ is any constant of integration and
\begin{eqnarray}
\label{surfgr}
\kappa_{\pm}=\frac{r_+-r_-}{2(r_{\pm}^2+a^2)}
\end{eqnarray}
are the surface gravities at the outer and inner horizons. 
The variable $r_*$ satisfies :
\begin{eqnarray}
\label{2.6b}
\frac{dr_*}{dr}=\frac{r^2+a^2}{\Delta}.
\end{eqnarray}
When $r$ runs from $r_+$ to $\infty$, $r_*$ runs from $-\infty$ to $\infty$. We put:
\begin{eqnarray}
\label{2.7b}
\Sigma:=\R_{r_*}\times S^2.
\end{eqnarray} 
We conclude this section with a useful identity on the coefficients of the metric :
\begin{eqnarray}
\label{usefulid}
1+\frac{Q^2-2Mr}{\rho^2}+\frac{a^2\sin^2\theta(2Mr-Q^2)^2}{\rho^2\sigma^2}=\frac{\rho^2\Delta}{\sigma^2}.
\end{eqnarray}
\subsection{Some remarks about geodesics in the Kerr-Newman space-time}
\label{sec2.1.2}
It is one of the most remarkable facts about the Kerr-Newman metric that there exist
four first integrals for the geodesic equations. If $\gamma$ is a geodesic in the Kerr-Newman
space-time, then $p:=\langle \gamma',\gamma'\rangle$ is conserved. The two Killing vector fields
$\partial_t,\, \partial_{\varphi}$ give two first integrals, the energy $E:=\langle \gamma', \partial_t \rangle$
and the angular momentum $L:=-\langle \gamma',\partial_{\varphi}\rangle$. There exists a 
fourth constant of motion, the so-called Carter constant ${\cal K}$ (see e.g. \cite{Ca}). Even if these
facts are well known we shall prove them here. The explicit form of the Carter constant in terms
of the hamiltonian flow appearing in the proof will be useful in the following. We will also use 
the Carter constant ${\cal Q}={\cal K}-(L-aE)^2$, which has a somewhat more geometrical meaning,
but gives in general more complicated formulas. Let 
\begin{eqnarray}
\label{DefDP}
{\bf P}:=(r^2+a^2)E-aL,\, {\bf D}:=L-aE\sin^2\theta.
\end{eqnarray}
We will consider the hamiltonian flow of
the principal symbol of $\frac{1}{2}\Box_g$ and then use the fact that a geodesic 
can be understood as the projection of the hamiltonian flow on ${\cal M}_{BH}$.
The d'Alembert operator associated to the Kerr-Newman metric is given by :
\begin{eqnarray}
\label{2.8}
\Box_g&=&\frac{\sigma^2}{\rho^2\Delta}\partial_t^2-\frac{2a(Q^2-2Mr)}{\rho^2\Delta}\partial_{\varphi}
\partial_t-\frac{\Delta-a^2\sin^2\theta}{\rho^2\Delta \sin^2\theta}\partial_{\varphi}^2\nonumber\\
&-&\frac{1}{\rho^2}\partial_r\Delta\partial_r-\frac{1}{\rho^2}\frac{1}{\sin\theta}\partial_{\theta}\sin\theta\partial_{\theta}.
\end{eqnarray}
The principal symbol of $\frac{1}{2}\Box_g$ is :
\begin{eqnarray}
\label{2.9}
P:=\frac{1}{2\rho^2}\left(\frac{\sigma^2}{\Delta}\tau^2-\frac{2a(Q^2-2Mr)}{\Delta}q_{\varphi}\tau
-\frac{\Delta-a^2\sin^2\theta}{\Delta \sin^2\theta}q_{\varphi}^2
-\Delta|\xi|^2-q_{\theta}^2\right).
\end{eqnarray}
Let 
\[{\cal C}_p:=\left\{(t,r,\theta,\varphi;\tau,\xi,q_{\theta},q_{\varphi});P(t,r,\theta,\varphi;\tau,\xi,q_{\theta},q_{\varphi})=\frac{1}{2}p\right\}. \]
Here $(\tau,\xi,q_{\theta},q_{\varphi})$ is dual to $(t,r,\theta,\varphi)$.
We have the following :
\begin{theorem}
\label{prop2.2}
$(i)$ Let $x_0=(t_0,r_0,\varphi_0,\theta_0,\tau_0,\xi_0,q_{\theta_0},q_{\varphi_0})\in {\cal C}_p$
and\\
$x(s)=(t(s),r(s),\theta(s),\varphi(s);\tau(s),\xi(s),q_{\theta}(s),q_{\varphi}(s))$
be the associated hamiltonian flow line. Then we have the following constants of motion :
\begin{eqnarray}
\label{2.10}
p=2P,E=\tau,\, L=-q_{\varphi},\,{\cal K}=q_{\theta}^2+\frac{{\bf D}^2}{\sin^2\theta}+pa^2\cos^2\theta=\frac{{\bf P}^2}{\Delta}-\Delta|\xi|^2-pr^2,
\end{eqnarray}
where ${\bf D},\, {\bf P}$ are defined in (\ref{DefDP}). 
\end{theorem}
It follows
\begin{corollary}
\label{th2.1}
Let $\gamma$ with $\gamma'=t'\partial_t+r'\partial_r+\theta'\partial_{\theta}
+\varphi'\partial_{\varphi}$ be a geodesic in the Kerr-Newman space-time. 
Then there exists a constant ${\cal K}={\cal K}_{\gamma}$ such that 
\begin{eqnarray}
\label{2.4}
\rho^4(r')^2&=&R(r)=\Delta(-pr^2-{\cal K})+{\bf P}^2,\\
\label{2.5}
\rho^4(\theta')^2&=&\Theta(\theta)={\cal K}-pa^2\cos^2\theta-\frac{{\bf D}^2}{\sin^2\theta},\\
\label{2.6}
\rho^2\varphi'&=&\frac{{\bf D}}{\sin^2\theta}+\frac{a{\bf P}}{\Delta},\\
\label{2.7}
\rho^2t'&=&a{\bf D}+(r^2+a^2)\frac{{\bf P}}{\Delta}.
\end{eqnarray}
\end{corollary} 

\begin{remark}
\label{rem2.3}
Theorem \ref{prop2.2} explains the link between the Carter constant ${\cal K}$
and the separability of the wave equation. Looking for $f$ in the form 
\[f(t,r,\theta,\varphi)=e^{i\omega t}e^{in\varphi}f_{r}(r)f_{\theta}(\theta)\]
we find :
\begin{eqnarray*}
\lefteqn{\Box_g f=0}\\
&\Leftrightarrow &\left(-\frac{(r^2+a^2)^2}{\Delta}\omega^2+\frac{2a(Q^2-2Mr)}{\Delta}n\omega
+\frac{r^2+a^2}{\Delta}D_{r_*}(r^2+a^2)D_{r_*}-\frac{a^2n^2}{\Delta}\right)f\\
&+&\left(a^2\sin^2\theta\omega^2+\frac{n^2}{\sin^2\theta}+\frac{1}{\sin\theta}D_{\theta}\sin\theta D_{\theta}\right)f=0.
\end{eqnarray*}
For fixed $\omega$ $P_{S^2}^{\omega}:=a^2\sin^2\theta\omega^2+\frac{D_{\varphi}^2}{\sin^2\theta}+
\frac{1}{\sin\theta} D_{\theta}\sin\theta D_{\theta}$ is a positive elliptic operator
on the sphere with eigenfunctions of the form $e^{in\varphi}f_{\theta}(\theta)$.
This gives the separability of the equation. The Carter constant is the analogue
of the eigenvalue of $P_{S^2}^{\omega}$ in classical mechanics.
\end{remark}
{\bf Proof of Theorem \ref{prop2.2}.}

The hamiltonian equations are the following :
\begin{eqnarray}
\label{2.11}
\dot{t}&=&\left(\frac{\sigma^2}{\rho^2\Delta}\tau-a\frac{Q^2-2Mr}{\rho^2\Delta}q_{\varphi}\right),\\
\label{2.12}
\dot{\tau}&=&0,\\
\label{2.13}
\dot{r}&=&-\frac{\Delta}{\rho^2}\xi,\\
\label{2.14}
\dot{\xi}&=&-\partial_rP,\\
\label{2.15}
\dot{\theta}&=&-\frac{q_{\theta}}{\rho^2},\\
\label{2.16}
\dot{q}_{\theta}&=&\frac{a^2\sin\theta\cos\theta}{\rho^2}\tau^2+\frac{1}{2\rho^2}\partial_{\theta}\frac{q_{\varphi}^2}{\sin^2{\theta}}
-\frac{a^2\cos\theta\sin\theta}{\rho^2}p,\\
\label{2.17}
\dot{\varphi}&=&\left(\frac{a(2Mr-Q^2)}{\rho^2\Delta}\tau-\frac{\Delta-a^2\sin^2\theta}{\rho^2\Delta\sin^2\theta}q_{\varphi}\right),\\
\label{2.18}
\dot{q}_{\varphi}&=&0,
\end{eqnarray}
where we have used that $(t,r,\theta,\varphi,\tau,\xi,\,q_{\theta},q_{\varphi})$ stays in
${\cal C}_p$.
The first two constants of motion follow from (\ref{2.12}), (\ref{2.18}). We multiply 
(\ref{2.16}) with $q_{\theta}$ given by (\ref{2.15}) and obtain :
\begin{eqnarray*}
\dot{q}_{\theta}q_{\theta}&=&-a^2\sin\theta\cos\theta\tau^2\dot{\theta}-\frac{1}{2}\frac{d}{dt}
\frac{q_{\varphi}^2}{\sin^2\theta}+a^2\cos\theta\sin\theta\dot{\theta} p\\
&=&-\frac{1}{2}\frac{d}{dt}\left(a^2\sin^2\theta\tau^2+\frac{q_{\varphi}^2}{\sin^2\theta}+a^2\cos^2\theta p\right)\\
&\Rightarrow& q_{\theta}^2+a^2\sin^2\theta\tau^2+\frac{q_{\varphi}^2}{\sin^2\theta}+a^2\cos^2\theta p=\hat{\cal K}=const.
\end{eqnarray*}
To obtain the second expression for ${\cal K}$ we use the fact that the flow stays 
in ${\cal C}_p$.
\qed
\\

The case $L=0$ is of particular interest. Let $\gamma$ be a null geodesic 
with energy $E>0$, Carter constant ${\cal K}$, angular momentum $L=0$ and given signs of 
$r_0',\, \theta_0'$. We can associate a hamiltonian flow line using (\ref{2.10}) to define the
initial data $\tau_0,\xi_0,q_{\theta_0},q_{\varphi_0}$ given $t_0,r_0,\theta_0,\varphi_0$. The signs of $q_{\theta_0}$ and $\xi_0$ are fixed by $sign q_{\theta_0}=-sign\theta_0',\, sign \xi_0=-sign r_0'$. From (\ref{2.10}) we infer conditions under which $\xi,\, q_{\theta}$ do not change their signs.
\begin{eqnarray*}
{\cal K}<\min_{r\in (r_+,\infty)}\frac{(r^2+a^2)^2E^2}{\Delta}&\Rightarrow&\xi\,\mbox{does not change its sign},\\
{\cal Q}\ge 0 &\Rightarrow&q_{\theta}\,\mbox{does not change its sign}.
\end{eqnarray*}
Note that in the case ${\cal Q}=0$ $\gamma$ is either in the equatorial plane or it does not cross it. 
Under the above conditions $\xi$ resp. $q_{\theta}$ can be understood as a function of $r$ resp. $\theta$ alone. 
In this case let $k_{{\cal K},E}$ and $l_{{\cal K},E}$ s.t.
\begin{eqnarray}
\label{C1a}
\frac{dk_{{\cal K},E}(r)}{dr}=\frac{\xi(r)}{E},\quad l'_{{\cal K},E}=\frac{q_{\theta}(\theta)}{E},
\quad\hat{r}_{{\cal K},E}:=k_{{\cal K},E}(r)+l_{{\cal K},E}(\theta).
\end{eqnarray} 
It is easy to check that $(t,\hat{r}_{{\cal K},E},\omega)$ is a coordinate system on block $I$.
We note that by (\ref{2.7}) we have $\partial_s t=\frac{E\sigma^2}{\Delta}>0$, thus 
$r,\theta,\varphi$ can be understood as functions of $t$ along $\gamma$.
\begin{lemma}
\label{lem2.4}
We have :
\begin{eqnarray}
\label{2.23}
\frac{\partial\hat{r}_{{\cal K},E}}{\partial t}=-1\quad\mbox{along}\quad \gamma,
\end{eqnarray}
where $t$ is the Boyer-Lindquist time.
\end{lemma}
{\bf Proof.}

This is an explicit calculation using equations (\ref{2.4})-(\ref{2.7}), (\ref{2.11})-(\ref{2.18}) :
\begin{eqnarray*}
\frac{\partial\hat{r}_{{\cal K},E}}{\partial t}&=&\frac{\partial \hat{r}_{{\cal K},E}}{\partial r}
\frac{\partial r}{\partial s}\frac{\partial s}{\partial t}+
\frac{\partial \hat{r}_{{\cal K},E}}{\partial \theta}\frac{\partial\theta}{\partial s}\frac{\partial s}{\partial t}\\
&=&-(|\xi|^2\Delta+|q_{\theta}|^2)\frac{1}{E}\left(a{\bf D}+(r^2+a^2)\frac{\bf P}{\Delta}\right)^{-1}\\
&=&-\frac{1}{E}\left(\frac{{\bf P}^2}{\Delta}-\frac{{\bf D}^2}{\sin^2\theta}\right)\left(a{\bf D}+(r^2+a^2)\frac{{\bf P}}{\Delta}\right)^{-1}=-1.
\end{eqnarray*}
\qed

We will suppose from now on $r_0'<0$, i.e. our construction is based on incoming null
geodesics. 
\begin{remark}
\label{rem2.1a}

Using the axial symmetry of the Kerr-Newman space-time we can for many studies of field 
equations in this background fix the angular momentum $\partial_{\varphi}=in$ in the
expression of the operator. The principal symbol of the new operator is the principal 
symbol of the old one with $q_{\varphi}=0$. This explains the importance of Lemma 
\ref{lem2.4}.
\end{remark}

We will often use the $r_*$ variable and its dual variable $\xi^*$. In this case we have to replace $\xi(r)$ by 
$\frac{r^2+a^2}{\Delta}\xi^*(r_*)$. The function $k_{{\cal K},E}$ is then a function 
of $r_*$ satisfying :
\begin{eqnarray}
\label{C1b}
k'_{{\cal K},E}(r_*)=\frac{\xi^*}{E},
\end{eqnarray}
where the prime denotes derivation with respect to $r_*$. Using the explicit form of the
Carter constant in Theorem \ref{prop2.2} we find :
\begin{eqnarray}
\label{2.19}
(k'_{{\cal K},E})^2&=&1-\frac{\Delta{\cal K}}{(r^2+a^2)^2E^2},\\
\label{2.20}
(l'_{{\cal K},E})^2&=&\frac{\cal K}{E^2}-a^2\sin^2\theta.
\end{eqnarray}
In particular we have :
\begin{eqnarray}
\label{2.22}
(k'_{{\cal K},E})^2\frac{(r^2+a^2)^2}{\sigma^2}+(l'_{{\cal K},E})^2\frac{\Delta}{\sigma^2}=1.
\end{eqnarray}
We will often consider the case ${\cal Q}=0$ and write in this case simply $k,l$ instead of $k_{a^2E^2,E},\, l_{a^2E^2,E}$.
\begin{remark}
\label{rem2.5}
The incoming null geodesics with Carter constant ${\cal K}$, angular momentum $L=0$, energy $E>0$
and given sign of $\theta_0'$ are the integral curves of the following vector fields :
\begin{eqnarray}
\label{2.24}
M^a_{{\cal K},E}=\frac{E}{\rho^2}\left(\frac{\sigma^2}{\Delta}\partial_t-(r^2+a^2)k'_{{\cal K},E}(r_*(r))\partial_r-
l'_{{\cal K},E}(\theta)\partial_{\theta}+\frac{a(2Mr-Q^2)}{\Delta}\partial_{\varphi}\right).
\end{eqnarray}
Let us put 
\begin{eqnarray*}
V^a_1\partial_a&=&\sqrt{\frac{\sigma^2}{\rho^2\Delta}}\left(\partial_t+\frac{a(2Mr-Q^2)}{\sigma^2}\partial_{\varphi}\right),\\
V^a_2\partial_a&=&\sqrt{\frac{\Delta}{\rho^2\sigma^2}}\left((r^2+a^2)k'_{{\cal K},E}(r_*(r))\partial_r+
l'_{{\cal K},E}(\theta)\partial_{\theta}\right),\\
W^a_1\partial_a&=&\sqrt{\frac{\rho^2}{\sigma^2}}\frac{1}{\sin\theta}\partial_{\varphi},\\
W^a_2\partial_a&=&\frac{r^2+a^2}{\sqrt{\rho^2\sigma^2}}\left(\frac{\Delta}{r^2+a^2}l'_{{\cal K},E}\partial_r-k'_{{\cal K},E}\partial_{\theta}\right).
\end{eqnarray*}
Note that 
\begin{eqnarray*}
V_1^aV_{1a}&=&1,\, V_2^aV_{2a}=W_1^aW_{1a}=W_2^aW_{2a}=-1,\\
L^aN_a&=&0\,\forall N^a, L^a \in \{V_1^a,V_2^a,W_1^a,W_2^a\}\, L^a\neq N^a.
\end{eqnarray*}
Clearly the considered null geodesics lie in $\Pi=span\{ V_1^a,V_2^a\}$. Using the 
Frobenius theorem (see e.g. \cite[Theorem 1.7.4]{ONe}) we see that, in contrast to
the PNG case, in our case the distribution of planes $\Pi^{\bot}=span\{W_1^a,W_2^a\}$ is integrable.
\end{remark} 
\begin{corollary}
\label{col2.6}
For given Carter constant ${\cal K}$, energy $E>0$ and sign of $\theta_0'$ the surfaces 
\begin{eqnarray*}
{\cal C}^{c,\pm}_{{\cal K},E}=\{(t,r_*,\theta,\varphi);\,\pm t=\hat{r}_{{\cal K},E}(r_*,\theta)+c\}
\end{eqnarray*}
are characteristic.
\end{corollary}
{\bf Proof.}

By Lemma \ref{lem2.4} the incoming null geodesic $\gamma$ with Carter constant ${\cal K}$, energy $E$, angular momentum
$L=0$ and the correct sign of $\theta_0'$ lies entirely in ${\cal C}^{c,-}_{{\cal K},E}$ if the starting
point lies in it. The geodesic $\gamma(-s)$ lies in ${\cal C}^{c,+}_{{\cal K},E}$. 
\qed
\begin{remark}
\label{rem2.7}

$(i)$ The variable $\hat{r}_{{\cal K},E}$ is a Bondi-Sachs type coordinate. This coordinate system
is discussed in some detail in \cite{FL}. As in \cite{FL} we will call the null geodesics
with $L={\cal Q}=0$ simple null geodesics (SNG's). 

$(ii)$ A natural way of finding the variable $\hat{r}_{{\cal K},E}$
is to start with Corollary \ref{col2.6}. Look for functions $k_{{\cal K},E}(r_*)$ and $l_{{\cal K},E}(\theta)$
such that ${\cal C}_{{\cal K},E}^{c,\pm}=\{\pm t=k_{{\cal K},E}(r_*)+l_{{\cal K},E}(\theta)\}$ is characteristic.
The condition that the normal is null is equivalent to (\ref{2.22}). The curve 
generated by the normal lies entirely in ${\cal C}^{c,\pm}_{{\cal K},E}$.
\end{remark}
\begin{remark}
From the explicit form of the Carter constant in Theorem \ref{prop2.2} follows~:
\begin{eqnarray}
\label{rot}
q_{\theta}^2+(p-E^2)a^2\cos^2\theta+\frac{q_{\varphi}^2}{\sin^2\theta}={\cal Q}.
\end{eqnarray}
This is the equation of the $\theta$ motion and it is interpreted as conservation of the 
mechanical energy with $V(\theta)=(p-E^2)a^2\cos^2\theta+\frac{q_{\varphi}^2}{\sin^2\theta}$
as potential energy and $q_{\theta}^2$ in the role of kinetic energy. The quantity
$\tilde{E}=(E^2-p)a^2$ is usually called the rotational energy.
\end{remark}
\section{The model of the collapsing star}
\label{sec2.2}
Let ${\cal S}_0$ be the surface of the star at time $t=0$. We suppose that elements $x_0\in {\cal S}_0$ will
move along curves which behave asymptotically like certain incoming timelike geodesics $\gamma_p$. All these geodesics should have the same energy $E$, angular momentum $L$, Carter constant ${\cal K}$ 
(resp. ${\cal Q}={\cal K}-(L-aE)^2$) and "mass"
$p:=\langle \gamma'_p,\gamma'_p \rangle$. We will suppose :

(A) The angular momentum $L$ vanishes : $L=0$.

(B) The rotational energy vanishes : $\tilde{E}=a^2(E^2-p)=0$.

(C) The total angular momentum about the axis of symmetry vanishes : ${\cal Q}=0.$ 

The conditions (A)-(C) are imposed by the fact that the collapse itself creates the space-time, thus momenta 
and rotational energy should be zero with respect to the space-time. 
\subsection{Timelike geodesics with $L={\cal Q}=\tilde{E}=0$}
\label{secE.3.2.1}

We will study the above family of geodesics in the following. The starting point of the geodesic is denoted $(0,r_0,\theta_0,\varphi_0)$. Given a point in the space-time, the conditions (A)-(C) define a unique cotangent vector provided 
you add the condition that the corresponding tangent vector is incoming. The choice of $p$ is irrelevant because it just corresponds to a normalization of the proper time.
 
\begin{lemma}
\label{lem2.8}
Along the geodesic $\gamma_p$ we have :
\begin{eqnarray}
\label{2.25}
\frac{\partial \theta}{\partial t}=0,\\
\label{2.26}
\frac{\partial \varphi}{\partial t}=\frac{a(2Mr-Q^2)}{\sigma^2},
\end{eqnarray}
where $t$ is the Boyer-Lindquist time.
\end{lemma}
{\bf Proof.}

(\ref{2.25}) follows directly from (\ref{2.5}). We have 
\begin{eqnarray*}
\frac{\partial\varphi}{\partial t}&=&\frac{\partial\varphi}{\partial s}\frac{\partial s}{\partial t}
=\frac{1}{\rho^2}\left(-aE+\frac{a(r^2+a^2)E}{\Delta}\right)\rho^2\left(-a^2ES^2+\frac{(r^2+a^2)^2}{\Delta}E\right)^{-1}\\
&=&\frac{\Delta a}{\sigma^2\Delta}(r^2+a^2-\Delta)=\frac{a(2Mr-Q^2)}{\sigma^2}.
\end{eqnarray*}
\qed

The function $\frac{\partial\varphi}{\partial t}=\frac{a(2Mr-Q^2)}{\sigma^2}$
is usually called the {\it local angular velocity of the space-time}. Our next aim is to adapt 
our coordinate system to the collapse of the star. The most natural way of doing this is to choose
an incoming null geodesic $\gamma$ with $L=Q=0$ and then use the Bondi-Sachs type coordinate
as in the previous section. In addition we want that $k(r_*)$ behaves like $r_*$ when $r_*\rightarrow-\infty$.
We therefore put :
\begin{eqnarray}
\label{2.27}
k(r_*)&=&r_*+\int_{-\infty}^{r_*}\left(\sqrt{1-\frac{a^2 \Delta(s)}{(r(s)^2+a^2)^2}}-1\right)ds,\\
\label{2.28}
l(\theta)&=&a \sin\theta.
\end{eqnarray}
The choice of the sign of $l'$ is not important, the opposite sign would have been 
possible. Recall that $\cos \theta $ does not change its sign along a null geodesic with 
$L={\cal Q}=0$. We fix the notation for the null vector fields generating $\gamma$
and the corresponding outgoing vector field (see Remark \ref{rem2.5})~:
\begin{eqnarray}
\label{2.29}
N^{\pm,a}\partial_a=\frac{E\sigma^2}{\rho^2\Delta}\left(\partial_t\pm\frac{(r^2+a^2)^2}{\sigma^2}k'(r_*)\partial_{r_*}
\pm\frac{\Delta}{\sigma^2}a\cos\theta\partial_{\theta}+\frac{a(2Mr-Q^2)}{\sigma^2}\partial_{\varphi}\right).
\end{eqnarray}
These vector fields will be important for the construction of the Newman-Penrose tetrad. We put :
\begin{eqnarray}
\label{hatr}
\hat{r}=k(r_*)+l(\theta)
\end{eqnarray}
and by Lemma \ref{lem2.4} we have :
\begin{eqnarray}
\frac{\partial \hat{r}}{\partial t}=-1\quad\mbox{along}\, \gamma.
\end{eqnarray}
Note that in the $(t,\hat{r},\omega)$ coordinate system the metric is given by :
\begin{eqnarray}
\label{hatrg}
g &=& \left( 1 + \frac{Q^2-2Mr}{\rho^2} \right)
\d t^2 + \frac{2a\sin^2\theta(2Mr-Q^2)}{\rho^2}
\, \d t \d\varphi\nonumber\\
&-& \frac{\rho^2\Delta}{(r^2+a^2)^2}\frac{(d\hat{r}
-l'(\theta)d\theta)^2}{k'(r_*)^2}- \rho^2 \d\theta^2- \frac{\sigma^2}{\rho^2} \sin^2 \theta \, \d\varphi^2.
\end{eqnarray}
In order to describe the model of the collapsing star we have to evaluate $\frac{\partial \hat{r}}{\partial t}$ along
$\gamma_p$. We start by studying $r(t,\theta).$ Recall that $\theta(t)=\theta_0=const.$ along $\gamma_p$ and that $\kappa_+=\frac{r_+-r_-}{2(r_+^2+a^2)}$ is the surface gravity of the outer horizon.
In what follows a dot will denote derivation in $t$.
\begin{lemma}
\label{lem2.9}
There exist smooth functions $\hat{C}_1(\theta,r_0),\, \hat{C}_2(\theta,r_0)$ such that along
$\gamma_p$ we have uniformly in $\theta,\, r_0\in{[}r_1,r_2{]}\subset(r_+,\infty)$:
\begin{eqnarray*}
\frac{1}{2\kappa_+}\ln|r-r_+|&=&-t-\hat{C}_1(\theta,r_0)e^{-2\kappa_+ t}+\hat{C}_2(\theta,r_0)+{\cal O}(e^{-4\kappa_+t}),\quad t\rightarrow \infty.
\end{eqnarray*}
\end{lemma}

{\bf Proof.}

Note that $(r^2+a^2)>\Delta$ on ${[}r_+,\infty)$. Therefore $\frac{\partial r}{\partial t}$ cannot change its sign. 
As $\gamma_p$ is incoming, the minus sign has to be chosen.
From (\ref{2.4}), (\ref{2.7}) we find with $p=E^2$ :
\begin{eqnarray}
\label{2.32}
\frac{\partial r}{\partial t}=-((r^2+a^2)(2Mr-Q^2))^{1/2}\frac{\Delta}{\sigma^2}.
\end{eqnarray}
We can consider $\theta$ as a parameter. Equation (\ref{2.32}) gives :
\begin{eqnarray}
\label{2.33}
t&=&-\int_{r_0}^r\frac{\sigma^2}{((s^2+a^2)(2Ms-Q^2))^{1/2}\Delta}ds\nonumber\\
&=&-\frac{1}{2\kappa_+}\ln|r-r_+|+\hat{C}_2-\int_{r_+}^{r}P(s)ds
\end{eqnarray}
with
\begin{eqnarray*}
P(r)&:=&\left(\frac{\sigma^2}{((r^2+a^2)(2Mr-Q^2))^{1/2}(r-r_-)}
-\frac{1}{2\kappa_+}\right)\frac{1}{r-r_+}\\
&=&\left(\frac{P_1(r)}{P_2(r)}-\frac{1}{2\kappa_+}\right)\frac{1}{r-r_+}\\
&=&\frac{P_1^2(r)-\frac{1}{4\kappa_+^2}P_2^2(r)}{P_2(r)(P_1(r)+\frac{1}{2\kappa_+}P_2(r))}\frac{1}{r-r_+},\\
\hat{C}_2&:=&\frac{1}{2\kappa_+}\ln|r_0-r_+|+\int_{r_+}^{r_0}P(s)ds.
\end{eqnarray*}
Note that $P_1^2(r_+)-\frac{1}{4\kappa_+^2}P_2^2(r_+)=0$. 
As $P_1^2(r)-\frac{1}{4\kappa_+^2}P_2^2(r)$ is a polynomial 
we infer that $P(r)$ is smooth at $r_+$. Let 
\begin{eqnarray*}
F(r):=\frac{1}{r-r_+}\int_{r_+}^rP(s)ds.
\end{eqnarray*}
Clearly $F$ is smooth and $\lim_{r\rightarrow r_+}F(r)=P(r_+)=:F(r_+).$ From (\ref{2.33})
we infer :
\begin{eqnarray}
\label{2.34}
r-r_+&=&e^{-2 \kappa_+ t}e^{2 \kappa_+ (\hat{C}_2-(r-r_+)F(r))},\\
\label{2.35}
r-r_+&=&e^{-2 \kappa_+ t}e^{2 \kappa_+ \hat{C}_2}e^{f(t)}
\end{eqnarray}
with $f(t)={\cal O}(e^{-2 \kappa_+ t})$. Putting (\ref{2.35}) into (\ref{2.34}) we obtain :
\begin{eqnarray*}
r-r_+&=&e^{-2\kappa_+t}e^{2\kappa_+(\hat{C}_2-e^{-2\kappa_+t}e^{2\kappa_+\hat{C}_2}e^{f(t)}F(r))}\\
\Rightarrow \frac{1}{2 \kappa_+}\ln |r-r_+|&=&-t+\hat{C}_2-e^{-2 \kappa_+ t}
e^{2 \kappa_+ \hat{C}_2}e^{f(t)}F(r)\\
&=&-t+\hat{C}_2-e^{-2 \kappa_+ t}e^{2 \kappa_+ \hat{C}_2}F(r_+)\\
&+&e^{2 \kappa_+(\hat{C}_2-t)}
(F(r_+)-e^{f(t)}F(r))
\end{eqnarray*}
and it remains to show :
\begin{eqnarray}
\label{2.36}
F(r_+)-e^{f(t)}F(r)={\cal O}(e^{-2\kappa_+ t}).
\end{eqnarray}
We write :
\begin{eqnarray}
\label{2.37}
F(r_+)-e^{f(t)}F(r)=F(r_+)-F(r)+(1-e^{f(t)})F(r). 
\end{eqnarray}
Noting that 
\[ e^{f(t)}=1+{\cal O}(e^{-2 \kappa_+ t}) \]
we obtain the required estimate for the second term in (\ref{2.37}). To estimate 
the first term we write :
\begin{eqnarray*}
F(r_+)-F(r)=\int_{r_+}^r\frac{P(r_+)-P(s)}{r-r_+}ds\\
\end{eqnarray*}
As $|P(s)-P(r_+)|\lesssim |s-r_+|$ we have :
\begin{eqnarray*}
|F(r_+)-F(r)|\lesssim\int_{r_+}^r\left|\frac{s-r_+}{r-r_+}\right| ds
\lesssim |r-r_+|={\cal O}(e^{-2\kappa_+ t}).
\end{eqnarray*}
(\ref{2.37}) follows. From the explicit form of the equations it is clear that everything is uniform in $\theta,r_0\in{[}r_1,r_2{]}$. \qed
\begin{lemma}
\label{lem2.10}
$(i)$ There exist smooth functions $A(\theta,r_0)$, $B(\theta,r_0)$ such that along $\gamma_p$
we have uniformly in $\theta,\,r_0\in{[}r_1,r_2{]}\subset(r_+,\infty)$:
\begin{eqnarray*}
r_*&=&-t-A(\theta,r_0)e^{-2 \kappa_+ t}+B(\theta,r_0)+{\cal O}(e^{-4 \kappa_+ t}),\, t\rightarrow \infty.
\end{eqnarray*}
$(ii)$ There exist smooth functions $\hat{A}(\theta,r_0)> 0,\, \hat{B}(\theta,r_0)$ such that
along $\gamma_p$ we have uniformly in $\theta,\,r_0\in{[}r_1,r_2{]}\subset(r_+,\infty)$:
\begin{eqnarray}
\label{inftyhatr}
\hat{r}&=&-t-\hat{A}(\theta,r_0)e^{-2 \kappa_+ t}+\hat{B}(\theta,r_0)+{\cal O}(e^{-4 \kappa_+ t}),\, t\rightarrow \infty.
\end{eqnarray}
Furthermore there exists $k_0>0$ s.t. for all $t>0,\theta\in {[}0,\pi{]}$ we have :
\[ \left(\frac{(r^2+a^2)^2}{\sigma^2}k'^2-\dot{\hat{r}}^2\right)\ge k_0e^{-2 \kappa_+ t}. \]

\end{lemma}
{\bf Proof.}

$(i)$ Recall that :
\begin{eqnarray*}
r_*&=&r+\frac{1}{2\kappa_+}\ln|r-r_+|-\frac{1}{2\kappa_-}\ln|r-r_-|\\
&=&r-t-\hat{C}_1(\theta,r_0)e^{-2 \kappa_+ t}+\hat{C}_2(\theta,r_0)-\frac{1}{2\kappa_-}\ln|r-r_-|+{\cal O}(e^{-4\kappa_+ t})\\
&=&r_++e^{-2 \kappa_+ t}e^{-2 \kappa_+ \hat{C}_1(\theta,r_0)e^{-2\kappa_+ t}}e^{2\kappa_+ \hat{C}_2(\theta,r_0)}
e^{{\cal O}(e^{-4 \kappa_+ t})}\\
&-&\frac{1}{2\kappa_-}\ln\left|e^{-2\kappa_+ t}e^{-2 \kappa_+\hat{C}_1(\theta,r_0)e^{-2 \kappa_+ t}}
e^{2\kappa_+\hat{C}_2(\theta,r_0)}e^{{\cal O}(e^{-4\kappa_+ t})}+r_+-r_-\right|\\
&-&t-\hat{C}_1(\theta,r_0)e^{-2\kappa_+ t}+\hat{C}_2(\theta,r_0)+{\cal O}(e^{-4 \kappa_+ t})\\
&=&-t+e^{-2\kappa_+\hat{C}_1(\theta,r_0)}e^{2\kappa_+\hat{C}_2(\theta,r_0)}e^{-2\kappa_+ t}\\
&-&\frac{1}{2\kappa_-(r_+-r_-)}e^{-2\kappa_+\hat{C}_1(\theta,r_0)}e^{2\kappa_+\hat{C}_2(\theta,r_0)}
e^{-2\kappa_+ t}\\
&+&\hat{C}_2(\theta,r_0)-\hat{C}_1(\theta,r_0)e^{-2\kappa_+ t}+r_+-\frac{1}{2\kappa_-}\ln|r_+-r_-|+{\cal O}(e^{-4\kappa_+ t})\\
&=&-t-A(\theta,r_0)e^{-2\kappa_+t}+B(\theta,r_0)+{\cal O}(e^{-4\kappa_+ t}),
\end{eqnarray*}
where we have used the Taylor expansions of the functions $e^x,\,\ln(1+x)$.

$(ii)$ By part $(i)$ of the lemma we have :
\begin{eqnarray}
\label{2.38}
\frac{\partial \hat{r}}{\partial t}&=&\sqrt{1-\frac{a^2\Delta}{(r^2+a^2)^2}}\frac{\partial r_*}{\partial t}\nonumber\\
&=&\sqrt{1-\frac{a^2\Delta}{(r^2+a^2)^2}}\left(-1+2\kappa_+A(\theta,r_0)e^{-2\kappa_+ t}
+{\cal O}(e^{-4 \kappa_+ t})\right).
\end{eqnarray}
By Lemma \ref{lem2.9} we find : 
\begin{eqnarray}
\label{2.47A}
\frac{a^2\Delta}{(r^2+a^2)^2}=G(r_0,\theta)e^{-2\kappa_+t}+{\cal O}(e^{-4\kappa_+ t}), 
\end{eqnarray}
and thus 
\begin{eqnarray*}
\sqrt{1-\frac{a^2\Delta}{(r^2+a^2)^2}}&=&1-\frac{1}{2}\frac{a^2\Delta}{(r^2+a^2)^2}
+{\cal O}(e^{-4\kappa_+ t})\\
&=&1-\frac{1}{2}G(r_0,\theta)e^{-2\kappa_+t}+{\cal O}(e^{-4\kappa_+ t}).
\end{eqnarray*}
Putting this into (\ref{2.38}) gives :
\[ \frac{\partial \hat{r}}{\partial t}=-1+2\kappa_+\hat{A}(\theta,r_0)e^{-2\kappa_+ t}
+{\cal O}(e^{-4\kappa_+}),\,t\rightarrow \infty. \]
It remains to show that $\hat{A}(\theta,r_0)>0$.
The curve $t\mapsto (t,\hat{r}(t),\theta(t),\varphi(t))$ has to be timelike. Using (\ref{hatrg}) and (\ref{usefulid})
we find 
\begin{eqnarray*}
E^2=p=\langle \gamma'_p,\gamma'_p \rangle &=&\frac{\rho^2\Delta}{(r^2+a^2)^2k'^2}
\left(\frac{(r^2+a^2)^2}{\sigma^2}k'^2-\dot{\hat{r}}^2\right)(\partial_s t)^2\\
&=&\frac{E^2\sigma^4}{(r^2+a^2)^2k'^2\rho^2\Delta}\left(\frac{(r^2+a^2)^2}{\sigma^2}k'^2-\dot{\hat{r}}^2\right).
\end{eqnarray*}
It follows :
\[ \left(\frac{(r^2+a^2)^2}{\sigma^2}k'^2-\dot{\hat{r}}^2\right)\ge c_0\Delta\ge k_0e^{-2 \kappa_+ t}. \]
In particular we have $\hat{A}(\theta,r_0)>0$.
\qed
\\
\subsection{Precise assumptions}
\label{secE.3.2.2}
Let us now make the precise assumptions on the collapse. We will suppose that the surface at time 
$t=0$ is given in the $(t,\hat{r},\theta,\varphi)$ coordinate system by ${\cal S}_0=\{(\hat{r}_0(\theta_0),\theta_0,\varphi_0);\, (\theta_0,\varphi_0)\in S^2\},$
where $\hat{r}_0(\theta_0)$ is a smooth function. As $\hat{r}_0$ does not depend on $\varphi_0$, we will suppose that $\hat{z}(t,\theta_0,\varphi_0)$ will be independent of 
$\varphi_0: \hat{z}(t,\theta_0,\varphi_0)=\hat{z}(t,\theta_0)=\hat{z}(t,\theta)$ as this is the case for $\hat{r}(t)$ describing the geodesic. Thus the surface 
of the star is given by :
\begin{eqnarray}
\label{star1}
{\cal S}=\{(t,\hat{z}(t,\theta),\omega); t\in \R,\, \omega \in S^2 \}. 
\end{eqnarray}
The function $\hat{z}(t,\theta)$ satisfies
\begin{eqnarray}
\label{star4}
&\forall t\le 0,\theta\in{[}0,\pi{]}\quad \hat{z}(t,\theta)=\hat{z}(0,\theta)<0,\\
\label{star0}
&\forall t>0,\, \theta\in {[}0,\pi{]}\quad \dot{\hat{z}}(t,\theta)<0,\\
\label{star0b}
&\exists k_0>0\, \forall t>0,\theta\in {[}0,\pi{]}\nonumber\\
&\left(\left(\frac{(r^2+a^2)^2}{\sigma^2}k'^2\right)(\hat{z}(t,\theta),\theta)-\dot{\hat{z}}^2(t,\theta)\right)
\ge k_0e^{-2 \kappa_+ t},\\
\label{star2} 
&\exists \hat{A}\in C^{\infty}({[}0,\pi{]}),\,\hat{\xi}\in C^{\infty}(\R\times {[}0,\pi{]})\quad \hat{z}(t,\theta)=-t-\hat{A}(\theta)e^{-2 \kappa_+ t}+\hat{\xi}(t,\theta)\nonumber\\
&\hat{A}(\theta)>0,\quad
\forall 0\le \alpha,\beta\le 2\quad|\partial^{\alpha}_t\partial^{\beta}_{\theta}\hat{\xi}(t,\theta)|\le C_{\alpha,\beta}e^{-4 \kappa_+ t}\, \forall \theta\in {[}0,\pi{]},\, t>0.
\end{eqnarray}
As already explained these assumptions are motivated by the preceding analysis. We do not suppose that a point on the surface moves exactly on a geodesic.
Note that (\ref{star0}), (\ref{star0b}) imply :
\[ \forall t>0,\theta\in {[}0,\pi{]}\quad -1<\dot{\hat{z}}(t,\theta)<0. \]
Equations (\ref{star4})-(\ref{star2}) summarize our assumptions on the collapse.
The space-time of the collapsing star is given by :
\[ {\cal M}_{col}=\{(t,\hat{r},\theta,\varphi);\, \hat{r}\ge \hat{z}(t,\theta)\}. \]
We will also note :
\[ \Sigma^{col}_t=\{(\hat{r},\theta,\varphi);\, \hat{r}\ge \hat{z}(t,\theta) \}. \] 
Thus :
\[ {\cal M}_{col}=\bigcup_t\Sigma_t^{col}. \]
Note that in the $(t,r_*,\theta,\varphi)$ coordinate system ${\cal M}_{col}$ and 
$\Sigma_t^{col}$ are given by :
\[ {\cal M}_{col}=\{(t,r_*,\theta,\varphi);\, r_*\ge z(t,\theta)\},\, 
\Sigma^{col}_t=\{(r_*,\theta,\varphi);\, r_*\ge z(t,\theta) \} \]
with
\[ z(t,\theta)=-t-A(\theta)e^{-2\kappa_+ t}+B(\theta)+{\cal O}(e^{-4\kappa_+ t}) \] 
for some appropriate $A(\theta),\, B(\theta)$.
\begin{remark}
\label{remC2.4}
(i) Let us compare assumptions (\ref{star1})-(\ref{star2}) to the preceding discussion on geodesics. The assumption
(\ref{star2}) contains with respect to the previous discussion an additional
asymptotic assumption. Comparing to Lemma \ref{lem2.10} this condition can be expressed
as $\hat{B}(\theta,r_0(\theta))=const.\, (r_0(\theta)=r(\hat{r}_0(\theta),\theta))$. Using the freedom of the constant of integration in (\ref{2.5b})
we can suppose
\begin{eqnarray}
\label{star3}
\hat{B}(\theta,r_0(\theta))=0.
\end{eqnarray}
(ii) The Penrose compactification of block $I$ can be constructed based on the SNG's rather than on the principal
null geodesics (PNG's). This construction is explained in Appendix \ref{AppB}.
Starting from this compactification we could establish a model of the collapsing
star that is similar to the one established by Bachelot for the Schwarzschild case (see \cite{Ba6}). In this
model the function $\hat{z}$ would be independent of $\theta$. 
\end{remark}
We finish this chapter with a lemma which shows that the asymptotic form (\ref{star2}) can be accomplished by incoming timelike geodesics with $L={\cal Q}=\tilde{E}=0$.
\begin{lemma}
\label{lemasgeo}
There exists a smooth function $\hat{r}_0(\theta)$ with the following property. Let $\gamma$ be a timelike incoming geodesic with 
${\cal Q}=L=\tilde{E}=0$ and starting point $(0,\hat{r}_0(\theta_0),\theta_0,\varphi_0)$. Then we have along $\gamma$:
\[\hat{r}+t\rightarrow 0,\quad t\rightarrow \infty. \]
\end{lemma}
{\bf Proof.}

Let $M_{\theta}(r)=\hat{C}_2(\theta,r)$, $\hat{C}_2(\theta,r)$ as in Lemma \ref{lem2.9}.
We have (see the explicit form of $\hat{C}_2(\theta,r)$ in the proof of Lemma \ref{lem2.9}) :
\[\lim_{r\rightarrow r_+}M_{\theta}(r)=-\infty,\, \lim_{r\rightarrow \infty}M_{\theta}(r)=\infty,\, \frac{\partial M_{\theta}}{\partial r}\ge \epsilon>0\, \forall r\in (r_+,\infty),\, \theta\in {[}0,\pi{]}.\]
Therefore $M_{\theta}^{-1}$ exists and we put 
\[r_0(\theta)=M^{-1}_{\theta}(-a\sin \theta-r_++\frac{1}{2\kappa_-}\ln|r_+-r_-|),\, \hat{r}_0(\theta)=\hat{r}(r_0(\theta),\theta).\] 
Clearly $\hat{C}_2(\theta,r_0(\theta))=-a\sin\theta-r_++\frac{1}{2\kappa_-}\ln|r_+-r_-|$. Following the proof of Lemma \ref{lem2.10} we see that 
\[ B(\theta,r_0)=\hat{C}_2(\theta,r_0)+r_+-\frac{1}{2\kappa_-}\ln|r_+-r_-|=-a\sin\theta.\] 
Using (\ref{2.27}), (\ref{2.28}), (\ref{hatr}) and (\ref{inftyhatr}) we see that \[\hat{B}(\theta,r_0(\theta_0))=\lim_{t\rightarrow\infty}\hat{r}+t=B(\theta,r_0(\theta))+a\sin\theta=0.\] \qed
\chapter{Classical Dirac Fields}
\label{sec3}
In this chapter we describe classical Dirac fields on $B_I$ as well as on ${\cal M}_{col}$. The main results of this chapter are collected in Section \ref{sec3.1}. Sections \ref{sec3.2} and \ref{sec3.3} contain a discussion about spin structures and Dirac fields which is valid in general globally hyperbolic space-times. In Section \ref{sec3.4} we introduce a new Newman-Penrose tetrad which is adapted to our problem and we discuss scattering results as far as they are needed for the formulation and discussion of the main theorem. Other scattering results are collected in Chapter \ref{sec7}. The boundary condition is discussed in Section \ref{sec3.5}. The constructions in this section are crucial for what follows. However the reader who wishes to get a first idea of the main theorem can in a first reading accept the results of Section \ref{sec3.1} and skip the rest of this chapter before coming back to it later.
\section{Main results}
\label{sec3.1}
Let ${\cal H}=L^2((\R_{\hat{r}}\times S^2,d\hat{r}d\omega);\C^4),\quad \Gamma^1=Diag(1,-1,-1,1)$.
\begin{proposition}
There exists a Newman-Penrose tetrad such that the Dirac equation in the Kerr-Newman space-time can be written as 
\begin{eqnarray*}
\partial_t\psi=iH\psi;\quad H=\Gamma^1D_{\hat{r}}+P_{\omega}+W,
\end{eqnarray*}
where $W$ is a real potential and $P_{\omega}$ is a differential operator of order one with derivatives only in the angular directions. The operator $H$ is selfadjoint with domain $D(H)=\{v\in {\cal H};\, H v\in {\cal H}\}.$
\end{proposition}
\begin{proposition}
There exist selfadjoint operators $P^{\pm}$ s.t. for all $g\in C_{\infty}(\R)$:
\begin{eqnarray}
g(P^{\pm})=s-\lim_{t\rightarrow\pm\infty}e^{-it H}g\left(\frac{\hat{r}}{t}\right)e^{it H}.
\end{eqnarray}
\end{proposition}
Let
\begin{eqnarray*}
H_{\leftarrow}&=&\Gamma^1D_{\hat{r}}-\frac{a}{r_+^2+a^2}D_{\varphi}-\frac{qQr_+}{r_+^2+a^2},\\
{\cal H}^+&=&\{v=(v_1,v_2,v_3,v_4)\in {\cal H};\, v_1=v_4=0 \},\\
{\cal H}^-&=&\{v=(v_1,v_2,v_3,v_4)\in {\cal H};\, v_2=v_3=0 \}. 
\end{eqnarray*}
The operator $H_{\leftarrow}$ is selfadjoint on ${\cal H}$ with domain  
$D(H_{\leftarrow})=\{v\in{\cal H}; H_{\leftarrow} v\in{\cal H}\}$.
\begin{theorem}
The wave operators
\begin{eqnarray*}
W_{\leftarrow}^{\pm}&=&s-\lim_{t\rightarrow\pm\infty}e^{-itH}e^{itH_{\leftarrow}}P_{{\cal H}^{\mp}},\\
\Omega_{\leftarrow}^{\pm}&=&s-\lim_{t\rightarrow\pm\infty}e^{-itH_{\leftarrow}}e^{itH}{\bf 1}_{\R^\mp}(P^{\pm})
\end{eqnarray*}
exist.
\end{theorem}
There exist similar wave operators at infinity using a modified asymptotic dynamics $U_D(t)$. Using the above tetrad the Dirac equation with MIT boundary condition (chiral angle $\nu$) can be written in the following form:
\begin{eqnarray}
\label{4.63.2}
\left. \begin{array}{rcl} \partial_t\Psi&=&iH\Psi,\quad\hat{z}(t,\theta)<\hat{r},\\
(\sum_{\hat{\mu}\in\{t,\hat{r},\theta,\varphi\}}{\cal N}_{\hat{\mu}}\hat{\gamma}^{\hat{\mu}})\Psi(t,\hat{z}(t,\theta),\omega)
&=&-ie^{-i\nu \gamma^5}\Psi(t,\hat{z}(t,\theta),\omega), \\
\Psi(t=s,.)&=&\Psi_s(.).  \\
\end{array}\right\}
\end{eqnarray}
Here ${\cal N}_{\hat{\mu}}$ are the coordinates of the conormal, $\hat{\gamma}^{\hat{\mu}}$ are some appropriate Dirac matrices and $\gamma^5=Diag(1,1,-1,-1)$.
Let ${\cal H}_t=L^2((\{(\hat{r},\omega)\in \R\times S^2;\hat{r}\ge \hat{z}(t,\theta)\},d\hat{r}d\omega);\, \C^4)$.
\begin{proposition}
The equation (\ref{4.63.2}) can be solved by a unitary propagator $U(t,s): {\cal H}_s\rightarrow {\cal H}_t$.
\end{proposition} 
\section{Spin structures}
\label{sec3.2}
Let $({\cal M},g)$ be a smooth 4-manifold with a lorentzian metric $g$ with signature
$(+,-,-,-)$ which is assumed to be oriented, time oriented and globally hyperbolic.
Global hyperbolicity
implies :
\begin{enumerate} 
\item 
$({\cal M},g)$ admits a spin structure (see R.P. Geroch \cite{Ge1,Ge2,Ge3} and
E. Stiefel \cite{Sti}) and we choose one. We denote by $\S$ (or $\S^A$ in the
abstract index formalism) the spin bundle over ${\cal M}$ and $\bar{\S}$
(or $\S^{A'}$) the same bundle with the complex structure replaced by
its opposite. The dual bundles $\S^*$ and $\bar{\S}^*$ will be denoted
respectively $\S_A$ and $\S_{A'}$. The complexified tangent
bundle to ${\cal M}$ is recovered as the tensor product of $\S$ and
$\bar{\S}$, i.e.

\[ T {\cal M} \otimes \C = \S \otimes \bar{\S} ~\mathrm{or}~T^a {\cal M}
\otimes \C = \S^A \otimes \S^{A'} \]

and similarly

\[ T^* {\cal M} \otimes \C = \S^* \otimes \bar{\S}^* ~\mathrm{or}~T_a
{\cal M} \otimes \C = \S_A \otimes \S_{A'} \, . \]

An abstract tensor index $a$ is thus understood as an unprimed spinor
index $A$ and a primed spinor index $A'$ clumped together~:
$a=AA'$. The symplectic forms on $\S$ and $\bar{\S}$ are denoted $\epsilon_{AB},$ 
${\epsilon}_{A'B'}$ and are referred to as the Levi-Civita symbols. $\epsilon_{AB}$
can be seen as an isomorphism from $\S$ to $\S^*$ which to ${\kappa}^A$ associates
${\kappa}_A={\kappa}^B\epsilon_{BA}$. Similarly, $\epsilon_{A'B'}$ and the corresponding
$\epsilon^{A'B'}$ can be regarded as lowering and raising devices for primed indices.
The metric $g$ is expressed in terms of the Levi-Civita symbols as 
$g_{ab}=\epsilon_{AB}\epsilon_{A'B'}$. 
\item 
There exists a global time function $t$ on ${\cal M}$. The level surfaces $\Sigma_t,\, t\in \R$, 
of the function $t$ define a foliation of ${\cal M}$, all $\Sigma_t$ being Cauchy
surfaces and homemorphic to a given smooth $3-$manifold $\Sigma$ (see Geroch \cite{Ge3}). 
Geroch's theorem does not say anything about the regularity of the leaves $\Sigma_t$; the
time function is only proved to be continuous and they are thus simply understood as topological
submanifolds of ${\cal M}$. A regularization procedure for the time function can be found in \cite{BS1}, \cite{BS2}. 
In the concrete cases which we consider in this paper the time function is smooth and all the leaves are diffeomorphic 
to $\Sigma$. The function $t$ is then a smooth time coordinate on ${\cal M}$ and it is 
increasing along any non space-like future oriented curve. Its gradient $\nabla^a t$
is everywhere orthogonal to the level surfaces $\Sigma_t$ of $t$ and it is therefore 
everywhere timelike; it is also future oriented. We identify ${\cal M}$ with the smooth manifold
$\R\times \Sigma$ and consider $g$ as a tensor valued function on $\R\times\Sigma$.
\end{enumerate}
Let $T^a$ be the future-pointing timelike vector field normal to $\Sigma_t$, normalized for later 
convenience to satisfy :
\[ T^aT_a=2, \]
i.e.
\begin{eqnarray*}
T^a=\frac{\sqrt{2}}{|\nabla t|}\nabla^at,\,\mbox{where}\, |\nabla t|=(g_{ab}\nabla^a t\nabla^bt)^{1/2}.
\end{eqnarray*}  
  
\section{The Dirac equation and the Newman-Penrose formalism}
\label{sec3.3}
In terms of two component spinors 
(sections of the bundles $\S^A$, $\S_A$,
$\S^{A'}$ or $\S_{A'})$, the charged Dirac equation takes the form (see \cite{PR1}, page 418):

\begin{equation} \label{DirEqPen}
\left\{ \begin{array}{l} { (\nabla_{A'}^A-iq\Phi^A_{A'}) \phi_A = \mu \chi_{A'} , } \\
({\nabla^{A'}_A -iq\Phi^{A'}_A)\chi_{A'} = \mu \phi_{A} \, ,~\mu =
  \frac{m}{\sqrt{2}}\, , }
\end{array} \right.
\end{equation}

where $m\geq 0$ is the mass of the field. 
The Dirac equation (\ref{DirEqPen})
possesses a conserved current (see for example \cite{Ni3}) on
general curved space-times, defined by the future oriented
non-spacelike vector field, sum of two future oriented null vector
fields~:

\[ V^a = \phi^A \bar{\phi}^{A'} + \bar{\chi}^A \chi^{A'} \, .\]

The vector field $V^a$ is divergence free, i.e. $\nabla^aV_a=0.$ Consequently the 
3-form $\omega=*V_a dx^a$ is closed. Let $\Sigma$ be a spacelike or characteristic hypersurface, 
$d\Omega$ the volume form on ${\cal M}$ induced by the metric ($d\Omega=\rho^2dt\wedge dr\wedge d\omega$ for the
Kerr-Newman metric), ${\cal N}^a$ the (future pointing) normal to $\Sigma$ and ${\cal L}^a$
transverse to $\Sigma$ with ${\cal N}^a{\cal L}_a=1$. Then 
\begin{eqnarray*}
\lefteqn{\int_{\Sigma}*(\phi_A\bar{\phi}_{A'}dx^{AA'}+\bar{\chi}_A\chi_{A'}dx^{AA'})}\\
&=&\int_{\Sigma}{\cal N}_{BB'}{\cal L}^{BB'}
*(\phi_A\bar{\phi}_{A'}dx^{AA'}+\bar{\chi}_A\chi_{A'}dx^{AA'})\\
&=&\int_{\Sigma}{\cal L}^{CC'}\partial_{CC'}\hook({\cal N}_{BB'}dx^{BB'}\wedge*(\phi_A\bar{\phi}_{A'}dx^{AA'}+\bar{\chi}_A\chi_{A'}dx^{AA'}))\\
&=&\int_{\Sigma}{\cal N}^{AA'}(\phi_A\bar{\phi}_{A'}+\bar{\chi}_A\chi_{A'})({\cal L}^{BB'}\partial_{BB'})\hook d\Omega.
\end{eqnarray*}
If $\Sigma$ is spacelike we can take ${\cal L}^{AA'}={\cal N}^{AA'}$ and the integral 
defines a norm and by this norm the space $L^2(\Sigma;\S_A\oplus \S^{A'})$ as completion of
$C_0^{\infty}(\Sigma;\S_A\oplus \S^{A'})$. Note that if $\Sigma$ is characteristic 
$\int_{\Sigma}*\phi_A\bar{\phi}_{A'}dx^{AA'}=0$ does not entail $\phi_A=0$ on $\Sigma$ (see Remark
\ref{rem3.1}). If $\Sigma_t$ are the level surfaces of $t$, then we see by Stokes' theorem that the 
total charge
\begin{eqnarray}
\label{CC} 
C(t)=\frac{1}{\sqrt{2}}\int_{\Sigma_t}V_aT^ad\sigma_{\Sigma_t} 
\end{eqnarray}
is constant throughout time. Here $d\sigma_{\Sigma_t}=\frac{1}{\sqrt{2}}T^a\hook d\Omega.$

Using the Newman-Penrose formalism, equation (\ref{DirEqPen}) can be
expressed as a system of partial differential equations with respect
to a coordinate basis. This formalism is based on the choice of a null
tetrad, i.e. a set of four vector fields $l^a$, $n^a$, $m^a$ and
$\bar{m}^a$, the first two being real and future oriented, $\bar{m}^a$
being the complex conjugate of $m^a$, such that
all four vector fields are null and $m^a$ is orthogonal to $l^a$ and
$n^a$, that is to say

\begin{equation} \label{NP1}
l_a l^a = n_a n^a = m_a m^a = l_a m^a = n_a m^a = 0 \, .
\end{equation}

The tetrad is said to be normalized if in addition

\begin{equation} \label{NP2}
l_a n^a = 1 \, ,~ m_a \bar{m}^a =-1 \, .
\end{equation}
The vectors $l^a$ and $n^a$ usually describe "dynamic" or scattering directions, i.e.
directions along which light rays may escape towards infinity (or more generally
asymptotic regions corresponding to scattering channels). The vector $m^a$ tends to 
have, at least spatially, bounded integral curves, typically $m^a$ and $\bar{m}^a$ generate 
rotations. 
The principle of the Newman-Penrose formalism is to decompose the
covariant derivative into directional covariant derivatives along the
frame vectors. We introduce a spin-frame $\{o^A ,\iota^A \}$, defined
uniquely up to an overall sign factor by the requirements that 

\begin{equation} \label{NPdyad}
o^A \bar{o}^{A'} = l^a \, ,~\iota^A \bar{\iota}^{A'} = n^a \, ,~o^A
\bar{\iota}^{A'} = m^a \, ,~\iota^A \bar{o}^{A'} = \bar{m}^a \, ,~o_A
\iota^A = 1 \, .
\end{equation}
We will also denote the spin frame by $\{\epsilon_0^{\,\, A},\epsilon_1^{\,\, A}\}$.
The dual basis of $\S_A$ is $\{ \epsilon^{\,\,\, 0}_A,\epsilon^{\,\,\, 1}_A\}$, where
$\epsilon^{\,\,\, 0}_A=-\iota_A,\, \epsilon^{\,\,\, 1}_A=o_A.$  
Let $\phi_0$ and $\phi_1$ be the components of $\phi_A$ in
$\{o^A , \iota^A \}$, and $\chi_{0'}$ and $\chi_{1'}$ the components of
$\chi_{A'}$ in $(\bar{o}^{A'} , \bar{\iota}^{A'} )$~:

\[ \phi_0 = \phi_A o^A~,~~\phi_1 = \phi_A \iota^A ~,~~ \chi_{0'} =
\chi_{A'} \bar{o}^{A'} ~,~~ \chi_{1'} = \chi_{A'} \bar{\iota}^{A'} \,
.\]

The Dirac equation then takes the form (see for example \cite{Ch})
\begin{eqnarray}
\label{NewmanDirac}
 \left. \begin{array}{l}
{ n^\mathbf{a} (\partial_\mathbf{a}-iq\Phi_{\bf{a}}) \, \phi_0 - m^\mathbf{a}
(\partial_\mathbf{a}-iq\Phi_{\bf{a}})  \, \phi_1 + (\mu - \gamma )\phi_0 + (\tau - \beta )
\phi_1 = \frac{m}{\sqrt{2}} \chi_{1'} \, , } \\ \\

{ l^\mathbf{a} (\partial_\mathbf{a}-iq\Phi_{\bf{a}})  \, \phi_1 - \bar{m}^\mathbf{a}
(\partial_\mathbf{a}-iq\Phi_{\bf{a}})  \, \phi_0 + (\alpha - \pi )\phi_0 + (\varepsilon -
\tilde{\rho} ) \phi_1 = - \frac{m}{\sqrt{2}} \chi_{0'} \, , } \\ \\

{ n^\mathbf{a} (\partial_\mathbf{a}-iq\Phi_{\bf{a}})  \, \chi_{0'} - \bar{m}^\mathbf{a}
(\partial_\mathbf{a}-iq\Phi_{\bf{a}})  \, \chi_{1'} + (\bar{\mu} - \bar{\gamma} )\chi_{0'}
+ (\bar{\tau} - \bar{\beta} ) \chi_{1'} = \frac{m}{\sqrt{2}} \phi_{1} \,
, } \\ \\

{ l^\mathbf{a} (\partial_\mathbf{a}-iq\Phi_{\bf{a}})  \, \chi_{1'} - m^\mathbf{a}
(\partial_\mathbf{a}-iq\Phi_{\bf{a}})  \, \chi_{0'} + (\bar{\alpha} - \bar{\pi} )\chi_{0'}
+ (\bar{\varepsilon} - \bar{\tilde{\rho}} ) \chi_{1'} = - \frac{m}{\sqrt{2}}
\phi_{0}\, . } \end{array} \right\}
\end{eqnarray}

The $\mu,\gamma$ etc. are the so called spin coefficients, for example 
$\mu=-\bar{m}^a\delta n_a,\delta=m^a\nabla_a$. For the formulas of the spin 
coefficients and details about the Newman-Penrose formalism see e.g. \cite{PR1}.

It is often useful to allow simultaneous consideration of bases of $T^a{\cal M}$ and
$\S^A$, which are completely unrelated to one another. Let $\{e_0,e_1,e_2,e_3\}$ be such
a basis of $T^a{\cal M}$, which is not related to the Newman-Penrose tetrad.

We define the Infeld-Van der Waerden symbols as the 
spinor components of the frame vectors in the spin frame $\{\epsilon^{\,\, A}_0,\epsilon^{\,\, A}_1\}$:  
\begin{eqnarray*}
g_{\bf a}^{\, \, \bf AA'}=e_{\bf a}^{\, \, \bf AA'}=g_{\bf a}^{\,\, a}\epsilon_A^{\,\, \bf A}
\epsilon_{A'}^{\,\, \bf A'}
=\left(\begin{array}{cc} n_{\bf a} & -\bar{m}_{\bf a} \\ -m_{\bf a} & l_{\bf a} \end{array} \right)
\end{eqnarray*}
(recall that $g_{\bf a}^{\, a}=e_{\bf a}^{\, a}$ denotes the vector field $e_{\bf a}$).
We use these quantities to express (\ref{DirEqPen}) in terms of spinor components :
\begin{eqnarray}
\label{DirEqMat1}
\left. \begin{array}{rcl}-i\epsilon_{A'}^{\,\, {\bf A}'}(\nabla^{AA'}-iq\Phi^{AA'})\phi_{A}&=&-ig^{\bf aAA'}\epsilon_{{\bf A}}^{\,\, A}(\nabla_{\bf a}-iq\Phi_{\bf a})\phi_{A}=-i\mu\chi^{\bf A'},\\
-i\epsilon_{{\bf A}}^{\,\, A}(\nabla_{AA'}-iq\Phi_{AA'})\chi^{A'}
&=&-ig^{\bf a}_{\,\, \bf AA'}\epsilon_{A'}^{\,\, \bf A'}(\nabla_{\bf a}-iq\Phi_{\bf a})\chi^{A'}=i\mu\phi_{\bf A}, \end{array} \right\}
\end{eqnarray} 
where $\nabla_{\bf a}$ denotes $\nabla_{e_{\bf a}}$. For ${\bf a}=0,1,2,3,$ we introduce the 
$2\times 2$ matrices 
\begin{eqnarray*}
A^{\bf a}=\, ^t g^{\bf aAA'},\, B^{\bf a}=g^{\bf a}_{\, \, \bf AA'},
\end{eqnarray*}
and the $4\times 4$ matrices 
\begin{eqnarray}
\label{Diracmatrix*}
\gamma^{\bf a}=\left(\begin{array}{cc} 0 & i\sqrt{2}B^{\bf a} \\ -i \sqrt{2} A^{\bf a} & 0 \end{array} \right).
\end{eqnarray}
We find :
\begin{eqnarray}
\label{DiracMatrix}
\gamma^{\bf a}=i\sqrt{2}\left(\begin{array}{cccc} 0 & 0 & l^{\bf a} & m^{\bf a} \\
0 & 0 & \bar{m}^{\bf a} & n^{\bf a} \\
-n^{\bf a} & m^{\bf a} & 0 & 0 \\ \bar{m}^{\bf a} & -l^{\bf a} & 0 &0 \end{array} \right).
\end{eqnarray}
Putting $\Psi=\phi_A\oplus\chi^{A'}$, the components of $\Psi$ in the spin frame are
${\bf \Psi}=^t(\phi_0,\phi_1,\chi^{0'},\chi^{1'})$ and (\ref{DirEqMat1}) becomes
\begin{eqnarray}
\label{DirEqMat2}
\sum_{\bf a=0}^3\gamma^{\bf a}{\cal P}(\nabla_{e_{\bf a}}-iq\Phi_{\bf a})\Psi+im{\bf \Psi}=0,
\end{eqnarray}
where ${\cal P}$ is the mapping that to a Dirac spinor associates its components in the
spin frame~:
\[ \Psi=\phi_A\oplus\chi^{A'}\mapsto {\bf \Psi}=\phi_{\bf A}\oplus \chi^{\bf A'}. \]

\begin{remark}
\label{rem3.1}
We have 
\begin{eqnarray*}
\phi_A&=&-\phi_0 \iota_A+\phi_1o_A,\\
\chi_{A'}&=&-\chi_{0'}\bar{\iota}_{A'}+\chi_{1'}\bar{o}_{A'}.
\end{eqnarray*}
Thus :
\begin{eqnarray*}
\phi_A\bar{\phi}_{A'}&=&(-\phi_0\iota_A+\phi_1o_A)(-\bar{\phi}_0\bar{\iota}_{A'}+\bar{\phi}_1\bar{o}_{A'})\\
&=&|\phi_0|^2n_a-\phi_0\bar{\phi}_1\bar{m}_a-\phi_1\bar{\phi}_0m_a+|\phi_1|^2l_a,\\
\bar{\chi}_A\chi_{A'}&=&|\chi_{0'}|^2n_a-\bar{\chi}_{0'}\chi_{1'}\bar{m}_a-\bar{\chi}_{1'}\chi_{0'}m_a+|\chi_{1'}|^2l_a.
\end{eqnarray*}
Putting $\Psi=(\phi_0,\phi_1,\chi_{1'},-\chi_{0'})$ we obtain :
\begin{eqnarray*}
\phi_A\bar{\phi}_{A'}+\bar{\chi}_A\chi_{A'}&=&|\Psi_1|^2n_a-\Psi_1\bar{\Psi}_2\bar{m}_a-\Psi_2\bar{\Psi}_1m_a+|\Psi_2|^2l_a\\
&+&|\Psi_4|^2n_a+\bar{\Psi}_4\Psi_3\bar{m}_a+\bar{\Psi}_3\Psi_4m_a+|\Psi_3|^2l_a.
\end{eqnarray*}
Thus for a vector field $X^a$ we have :
\begin{eqnarray}
\label{XX}
X^{AA'}(\phi_A\bar{\phi}_{A'}+\bar{\chi}_A\chi_{A'})&=&\langle {\bf X}\Psi,\Psi \rangle_{\C^4},\\
{\bf X}&=&\left(\begin{array}{cccc} n_aX^a & -m_aX^a & 0 & 0 \\
-\bar{m}_aX^a & l_aX^a & 0 & 0 \\ 0 & 0 & l_aX^a & m_aX^a \\ 0 & 0 & \bar{m}_aX^a & n_aX^a \end{array}
\right).
\end{eqnarray}
If $\Sigma$ is a characteristic hypersurface with conormal $n_a$, then
\begin{eqnarray}
\label{ISC}
\lefteqn{\int_{\Sigma}*((\phi_A\bar{\phi}_{A'}+\bar{\chi}_A\chi_{A'})dx^{AA'})}\nonumber\\
&=&\int_{\Sigma}n^a(\phi_A\bar{\phi}_{A'}
+\bar{\chi}_A\chi_{A'})d\sigma_{\Sigma}=\int_{\Sigma}(|\Psi_2|^2+|\Psi_3|^2)d\sigma_{\Sigma},
\end{eqnarray}
where $d\sigma_{\Sigma}=(l^a\partial_a)\hook d\Omega.$ 
\end{remark}
\section{The Dirac equation on block $I$}
\label{sec3.4}
\subsection{A new Newman-Penrose tetrad}
\label{sec3.4.1}
The tetrad normally used to describe the Dirac equation on the Kerr-Newman background
is Kinnersley's tetrad (see \cite{Ki}). Kinnersley uses the type D structure of the space-time
and chooses the vectors $n^a,l^a$ to be $V^{\pm}$. In \cite{HN} it is argued that this tetrad is not adapted to
time dependent scattering problems. We introduce a new tetrad $L^a,N^a, M^a, \bar{M}^a$ which
is {\it adapted to the foliation} in the sense that
\[ L^a+N^a=T^a. \] 
One advantage of a tetrad adapted to the foliation is that
the conserved current simply reads :
\[ \frac{1}{\sqrt{2}}\int_{\Sigma_t}(|\phi_0|^2+|\phi_1|^2+|\chi_{0'}|^2+|\chi_{1'}|^2)d\sigma_{\Sigma_t}. \]
This follows from the formulas (\ref{CC}) and (\ref{XX}).
We choose $L^a$ and $N^a$ in the plane spanned by $T^a$ and $\partial_r$ and $L^a$
to be outgoing, $N^a$ to be incoming. The choice of $M^a$ is then imposed, except
for the freedom of a constant factor of modulus 1. In \cite{HN} this construction 
is done for the Kerr metric, the corresponding tetrad for the Kerr-Newman metric is calculated
in \cite{Da2} : 
\begin{eqnarray}
\label{3.2}
\left.\begin{array}{rcl}
L^a&=&\frac{1}{2}T^a+\sqrt{\frac{\Delta}{2\rho^2}}\partial_r
=\sqrt{\frac{\sigma^2}{2\Delta\rho^2}}(\partial_t+\frac{a(2Mr-Q^2)}{\sigma^2}\partial_{\varphi})+\sqrt{\frac{\Delta}{2\rho^2}}\partial_r,\\
N^a&=&\frac{1}{2}T^a-\sqrt{\frac{\Delta}{2\rho^2}}\partial_r
=\sqrt{\frac{\sigma^2}{2\Delta\rho^2}}(\partial_t+\frac{a(2Mr-Q^2)}{\sigma^2}\partial_{\varphi})-\sqrt{\frac{\Delta}{2\rho^2}}\partial_r,\\
M^a&=&\frac{1}{\sqrt{2\rho^2}}(\partial_{\theta}+\frac{\rho^2}{\sqrt{\sigma^2}}\frac{i}{\sin\theta}\partial_{\varphi})
.\end{array}\right\}
\end{eqnarray}
The Dirac equation in the Kerr-Newman metric is then described in the following way.
Let 
\begin{eqnarray*}
\mathbb{D}_{S^2}=\frac{1}{i}\left(\begin{array}{cc} 0 & \partial_{\theta}+\frac{\cot\theta}{2}+\frac{i}{\sin\theta}\partial_{\varphi} \\
\partial_{\theta}+\frac{\cot\theta}{2}-\frac{i}{\sin\theta}\partial_{\varphi} & 0 \end{array} \right),\quad 
\notD_{S^2}=\left(\begin{array}{cc} \mathbb{D}_{S^2} & 0 \\ 0 & -\mathbb{D}_{S^2} \end{array} \right).
\end{eqnarray*}
We define the Pauli matrices :
\[ \sigma^0=\left(\begin{array}{cc} 1 & 0 \\ 0 & 1 \end{array} \right),\quad 
\sigma^1=\left(\begin{array}{cc} 1 & 0 \\ 0 & -1 \end{array} \right),\quad
\sigma^2=\left(\begin{array}{cc} 0 & 1 \\ 1 & 0 \end{array} \right),\quad 
\sigma^3=i\left(\begin{array}{cc} 0 & -1 \\ 1 & 0 \end{array} \right) \]
and the Dirac matrices
\[ \gamma^0=i\left(\begin{array}{cc} 0 & \sigma^0 \\ -\sigma^0 & 0 \end{array} \right),\quad 
\gamma^k=i\left(\begin{array}{cc} 0 & \sigma^k \\ \sigma^k & 0 \end{array} \right),\quad k=1,2,3, \]
which satisfy the following anticommutation relations :
\[ \gamma^{\mu}\gamma^{\nu}+\gamma^{\nu}\gamma^{\mu}=2\eta^{\mu\nu}I_{\R^4},\quad \mu,\nu=0,...,3,\,\eta^{\mu\nu}=Diag(1,-1,-1,-1). \]
Let
\[ \gamma^5:=-i\gamma^0\gamma^1\gamma^2\gamma^3=Diag(1,1,-1,-1). \]
We will also need the matrices $\Gamma^j\,(j=1,...,4)$ defined by :
\[ 1\le k \le 3\, \Gamma^k=\left(\begin{array}{cc} \sigma^k & 0 \\ 0 & -\sigma^k \end{array} \right),\, 
\Gamma^4=i\left(\begin{array}{cc} 0 & -\sigma^0 \\ \sigma^0 & 0 \end{array} \right). \] 
Note that 
\[ \Gamma^1=Diag (1,-1,-1,1). \]
Let now $\Phi$ be the bi-spinor $^t(\phi_0,\phi_1,\chi_{1'},-\chi_{0'})$ and 
\[ \Psi=\sqrt{\frac{\sqrt{\Delta}\sigma\rho}{r^2+a^2}}\Phi. \]
Then the equation satisfied by $\Psi$ is (see \cite{Da2}) :
\begin{eqnarray}
\label{3.9a}
\partial_t\Psi&=&i\notD\Psi,\\
\label{3.10a}
\notD&=&h\notD_{\gs}h+V_{\varphi}D_{\varphi}+V_1,\\
\label{3.11a}
\notD_{\gs}&=&\Gamma^1D_{r_*}+a_0(r_*)\notD_{S^2}+b_0(r_*)\Gamma^4+c_1(r_*)+c_2^{\varphi}(r_*)D_{\varphi},\\
\label{3.12a}
a_0(r_*)&=&\frac{\sqrt{\Delta}}{r^2+a^2},\, b_0(r_*)=\frac{m\sqrt{\Delta}}{\sqrt{r^2+a^2}},\, 
c_1(r_*)=-\frac{qQr}{r^2+a^2},\\
c_2^{\varphi}(r_*)&=&-\frac{a(2Mr-Q^2)}{(r^2+a^2)^2},\\
h(r_*,\theta)&=&\sqrt{\frac{r^2+a^2}{\sigma}},\\
\label{3.13ab}
V_{\varphi}&=&-\frac{\sqrt{\Delta}}{\sigma\sin\theta}\left(\frac{\rho^2}{\sigma}-1\right)\Gamma^3-\frac{a(2Mr-Q^2)}{(r^2+a^2)\sigma^2}(r^2+a^2-\sigma),\\
V_1&=&\mathbb{V}_0+\frac{m\sqrt{\Delta}}{\sigma}(\rho-\sqrt{r^2+a^2})\Gamma^4-\frac{qQr}{\sigma^2}(r^2+a^2-\sigma),\\
\label{3.15a}
\mathbb{V}_0&=&\left(\begin{array}{cccc} v_0 & v_1 & 0 & 0\\ \bar{v}_1 & -v_0 & 0 & 0 \\
0 & 0 & v_0 & \bar{v}_1 \\ 0 & 0 & v_1 & -v_0 \end{array} \right),\\
\label{3.16a}
v_0&=&-\frac{a^3(2Mr-Q^2)\Delta\sin^2\theta\cos\theta}{2\rho^2\sigma^3},\\
\label{3.17a}
v_1&=&-\frac{aM\sqrt{\Delta}\sin\theta}{2\rho^2\sigma}\nonumber\\
&+&\frac{a\sqrt{\Delta}\sin\theta(2Mr-Q^2)(2r(r^2+a^2)-(r-M)a^2\sin^2\theta)}{2\rho^2\sigma^3}\nonumber\\
&-&i\frac{qQ\sqrt{\Delta}ra\sin\theta}{\sigma^2}.
\end{eqnarray}
Recall from \cite{Da2} that $\notD_{\gs},\, \notD$ are selfadjoint on 
\[ {\cal H}_*:=L^2((\R\times S^2;dr_*d\omega);\C^4) \]
with domain 
\[ {\cal H}^1_*:=D(\notD_{\gs})=D(\notD)=\{ u\in {\cal H}_*;\, \notD u \in {\cal H}_*\} \]
and that their spectrum is purely absolutely continuous.
We also define the asymptotic dynamics 
\begin{eqnarray}
\notD_H&=&\Gamma^1D_{r_*}-\frac{a}{r_+^2+a^2}D_{\varphi}-\frac{qQr_+}{r_+^2+a^2},
\, D(\notD_H)=\{u\in {\cal H}_*, \notD_Hu\in {\cal H}_*\},\\
\label{3.17b}
\notD_{\rightarrow}&=&\Gamma^1D_{r_*}+m\Gamma^4,\, D(\notD_{\rightarrow})=\{u\in {\cal H}_*;\,\notD_{\rightarrow} u\in {\cal H}_*\}.
\end{eqnarray}
Even if the above tetrad was successfully used for the proof of the asymptotic completeness 
result, it has a major drawback for the treatment of the Hawking effect. In fact in this 
representation and using the Regge-Wheeler type coordinate $r_*$ the modulus of local velocity
\[ v:={[}r_*,\notD{]}=h^2\Gamma^1\neq\Gamma^1 \]
is not equal to $1$. The consequence is that in the high 
energy regime which is characteristic of the Hawking effect the full dynamics $\notD$
and the free dynamics $\notD_H$ or $\notD_{\gs}$ are no longer close to each other.
Now recall that $\partial_t \hat{r}=-1$ along incoming null geodesics with the correct sign
of $\theta_0'$. This means that the observable $\hat{r}$ should increase (resp. decrease)
exactly like $t$ along the evolution if we focus on scattering directions on which
$\hat{r}$ increases (decreases) in this way. We therefore choose
\begin{eqnarray}
\label{3.4}
l^a=\lambda N^{+,a},\, n^a=\lambda N^{-,a}
\end{eqnarray}
for some normalization constant $\lambda$. The choice of $m^a$ is now imposed except for a factor 
of modulus $1$. We find :
\begin{eqnarray}
\label{3.5}
\left. \begin{array}{rcl} l^a&=&\sqrt{\frac{\sigma^2}{2\rho^2\Delta}}\left(\partial_t+
\frac{(r^2+a^2)^2}{\sigma^2}k'(r_*)\partial_{r_*}+\frac{\Delta}{\sigma^2}a\cos\theta\partial_{\theta}
+\frac{a(2Mr-Q^2)}{\sigma^2}\partial_{\varphi}\right)\\
n^a&=&\sqrt{\frac{\sigma^2}{2\rho^2\Delta}}\left(\partial_t-
\frac{(r^2+a^2)^2}{\sigma^2}k'(r_*)\partial_{r_*}-\frac{\Delta}{\sigma^2}a\cos\theta\partial_{\theta}
+\frac{a(2Mr-Q^2)}{\sigma^2}\partial_{\varphi}\right)\\
m^a&=&\sqrt{\frac{\rho^2}{2\sigma^2}}\left(i\frac{r^2+a^2}{\rho^2}(a\cos\theta\partial_{r_*}-k'(r_*)\partial_{\theta})
+\frac{1}{\sin\theta}\partial_{\varphi}\right). \end{array} \right\}
\end{eqnarray}
Note that in the $(t,\hat{r},\omega)$ coordinate system we have :
\begin{eqnarray*}
\left. \begin{array}{rcl} l^a&=&\sqrt{\frac{\sigma^2}{2\rho^2\Delta}}\left(\partial_t+
\partial_{\hat{r}}+\frac{\Delta}{\sigma^2}a\cos\theta\partial_{\theta}
+\frac{a(2Mr-Q^2)}{\sigma^2}\partial_{\varphi}\right)\\
n^a&=&\sqrt{\frac{\sigma^2}{2\rho^2\Delta}}\left(\partial_t-
\partial_{\hat{r}}-\frac{\Delta}{\sigma^2}a\cos\theta\partial_{\theta}
+\frac{a(2Mr-Q^2)}{\sigma^2}\partial_{\varphi}\right)\\
m^a&=&\sqrt{\frac{\rho^2}{2\sigma^2}}\left(-ik'\frac{r^2+a^2}{\rho^2}\partial_{\theta}
+\frac{1}{\sin\theta}\partial_{\varphi}\right). \end{array} \right\}
\end{eqnarray*}
and that the tetrad $l^a,n^a,m^a$ is adapted to the foliation.
\subsection{The new expression of the Dirac equation}
\label{sec3.4.2}
Let us put 
\begin{eqnarray*}
\left. \begin{array}{c} e_{\bf 1}^a=L^a,\, e_{\bf 2}^a=N^a,\, e_{\bf 3}^a=M^a,\, e_{\bf 4}^a=\bar{M}^a, \\
e^{\bf 1}_a=L_a,\, e^{\bf 2}_a=N_a,\, e^{\bf 3}_a=M_a,\, e^{\bf 4}_a=\bar{M}_a. \end{array} \right\}
\end{eqnarray*}

In order to find the new expression of the Dirac equation we express $l^a,n^a,m^a,\bar{m}^a$ in terms of $L^a,N^a,M^a,\bar{M}^a$ 
and find the Lorentz transformation :
\begin{eqnarray}
\label{3.6}
L^{\bf b}_{ \bf a}=L^b_ae_{\bf a}^{a}e_{b} ^{\bf b}=\frac{1}{2}\left(\begin{array}{cccc}
1+\alpha & 1-\alpha & \beta & \beta \\ 1-\alpha & 1+\alpha & -\beta & -\beta \\
 i \beta & -i\beta & -i(1+\alpha) & i(1-\alpha) \\
-i\beta & i\beta & -i(1-\alpha) & i(1+\alpha) \end{array}\right),
\end{eqnarray}
where $\alpha=k'h^2,\quad \beta=\sqrt{\frac{\Delta}{\sigma^2}}a\cos\theta$, in particular :
\begin{eqnarray}
\label{3.7}
\alpha^2+\beta^2=1.
\end{eqnarray}
We associate to the Newman-Penrose tetrads $L^a, N^a, M^a,\bar{M}^a$ resp. $l^a,n^a,m^a,\bar{m}^a$
the spin frames $\{\tilde{\epsilon}^{\,\, A}_0,\tilde{\epsilon}^{\,\, A}_1\}$ resp. $\{\epsilon_0^{\,\, A},\epsilon_1^{\,\, A}\}$
as explained in Section \ref{sec3.3}. 
The matrix $S^{\,\, \bf B}_{\bf A}$ of the corresponding spin transformation $S^{\,\, B}_A$ 
in the spin frame $\{\tilde{\epsilon}^{\,\, A}_0,\tilde{\epsilon}^{\,\, A}_1\}$ is uniquely determined, modulo sign, by
\begin{eqnarray*}
L^b_a=S^{\,\, B}_A\bar{S}^{\,\, B'}_{A'}\,\mbox{and}\,\det(S^{\,\, \bf B}_{\bf A})=1.
\end{eqnarray*}
The first condition can be expressed in terms of coordinates as 
\begin{eqnarray*}
L^{\bf b}_{\bf a}=\left(\begin{array}{cccc} |S_0^0|^2 & |S_0^1|^2 & S_0^0\bar{S}^{1'}_{0'} & S^1_0\bar{S}^{0'}_{0'} \\
|S_1^0|^2 & |S_1^1|^2 & S^0_1\bar{S}^{1'}_{1'} & S^1_1\bar{S}^{0'}_{1'} \\
S^0_0\bar{S}^{0'}_{1'} & S^1_0\bar{S}^{1'}_{1'} & S^0_0\bar{S}^{1'}_{1'} & S^1_0\bar{S}^{0'}_{1'} \\
S^0_1\bar{S}^{0'}_{0'} & S^1_1\bar{S}^{1'}_{0'} & S^0_1\bar{S}^{1'}_{0'} & S^1_1\bar{S}^{0'}_{0'} 
\end{array} \right).
\end{eqnarray*}

We find :
\begin{eqnarray}
\label{3.8}
S^{\,\, \bf B}_{\bf A}=\left(\begin{array}{cc} S^0_0 & S^1_0 \\ S^0_1 & S^1_1 \end{array} \right) 
=\frac{1}{\sqrt{2(1+\alpha)}}e^{-i\pi/4}\left(\begin{array}{cc} 1+\alpha & \beta \\ -i\beta & 
i(1+\alpha) \end{array} \right)=:U. 
\end{eqnarray}
Let $\phi_{\bf C}$ (resp. $\tilde{\phi}_{\bf C}$) be the components of $\phi_A$ in
the spin frame $\{{\epsilon}_0^{\,\, A},{\epsilon}_1^{\,\, A}\}$ (resp. $\{\tilde{\epsilon}_0^{\,\, A},\tilde{\epsilon}_1^{\,\, A}\}$) 
and $\chi^{\bf C'}$ (resp. $\tilde{\chi}^{\bf C'}$) the components of $\chi^{A'}$
in the spin frames $\{{\epsilon}_{0'}^{\,\, A'},{\epsilon}_{1'}^{\,\, A'}\}$ (resp. $\{\tilde{\epsilon}_{0'}^{\,\, A'},\tilde{\epsilon}_{1'}^{\,\, A'}\}$).
We have :
\begin{eqnarray*}
\phi_{\bf C}&=&\phi_A\epsilon_{\bf C}^{\,\, A}=\phi_AS_B^{\,\, A}\tilde{\epsilon}_{\bf C}^{\,\, B}=
\tilde{\phi}_{\bf A}\tilde{\epsilon}_A^{\,\, \bf A}S_B^{\,\, A}\tilde{\epsilon}_{\bf C}^{\,\, B}=
\tilde{\phi}_{\bf A}S_{\bf C}^{\,\, \bf A},\\ 
\tilde{\chi}^{\bf C'}&=&\chi^{A'}\tilde{\epsilon}_{A'}^{\,\, \bf C'}=\chi^{\bf A'}\epsilon_{\bf A'}^{\,\, A'}\tilde{\epsilon}_{A'}^{\,\, \bf C'}
=\chi^{\bf A'}\bar{S}_{B'}^{\,\, A'}\tilde{\epsilon}_{\bf A'}^{\,\, B'}\tilde{\epsilon}_{A'}^{\,\, \bf C'}
=\chi^{\bf A'}\bar{S}_{\bf A'}^{\,\, \bf C'}.
\end{eqnarray*}
Noting that the matrix $\bar{S}_{\bf A'}^{\,\, \bf C'}$ is the inverse of $S_{\bf C}^{\,\, \bf A}$
it follows \footnote{We make the usual convention that ${\bf A}={\bf A'}$ numerically.} :
\begin{eqnarray*}
\chi^{\bf C'}=\sum_{{\bf A}=0}^1 S_{\bf C}^{\, \, \bf A}\tilde{\chi}^{\bf A'}.
\end{eqnarray*}
We put :
\begin{eqnarray}
{\cal U}:=\left(\begin{array}{cc} U & 0 \\ 0 & U \end{array}\right).
\end{eqnarray}
Our aim is to calculate 
\begin{eqnarray}
\label{3.10}
\hat{\notD}:={\cal U}\notD\,{\cal U}^*.
\end{eqnarray}
We define :
\[ 1\le j\le 4 : \quad \hat{\Gamma}^j:={\cal U}\Gamma^j{\cal U}^* \]
and find :
\begin{eqnarray*}
& &1\le j \le 3 \quad \hat{\Gamma}^j=\left(\begin{array}{cc} \hat{\sigma}^j & 0 \\ 0 & -\hat{\sigma}^j \end{array} \right),
\quad \hat{\Gamma}^4=\Gamma^4,\\
\hat{\sigma}^1&:=&\left(\begin{array}{cc} \alpha & i \beta \\ -i\beta & -\alpha \end{array} \right),
\,\hat{\sigma}^2:=\left(\begin{array}{cc} \beta & -i \alpha \\ i\alpha & -\beta \end{array} \right),
\,\hat{\sigma}^3:=\left(\begin{array}{cc} 0 & -1 \\ -1 & 0 \end{array} \right).
\end{eqnarray*}
The Dirac operator $\hat{\notD}$ can be written in the following form :
\begin{eqnarray*}
\hat{\notD}&=&\left(\begin{array}{cc} M_{r_*} & 0 \\ 0 & -M_{r_*} \end{array}\right)+
\left(\begin{array}{cc} M_{\theta} & 0 \\ 0 & -M_{\theta} \end{array}\right)
+\frac{a_0h^2}{\sin\theta}\hat{\Gamma}^3D_{\varphi}\\
&+&h^2c_1+h^2c_2^{\varphi}D_{\varphi}+\hat{V}_{\varphi}D_{\varphi}+\hat{V}_1,\\
\hat{V}_{\varphi}&=&-\frac{\sqrt{\Delta}}{\sigma\sin\theta}\left(\frac{\rho^2}{\sigma}-1\right)\hat{\Gamma}^3-\frac{a(2Mr-Q^2)}{(r^2+a^2)\sigma^2}(r^2+a^2-\sigma),\\
\hat{V}_1&=&\hat{\mathbb{V}}_0+\frac{m\sqrt{\Delta}}{\sigma}(\rho-\sqrt{r^2+a^2})\hat{\Gamma}^4-\frac{qQr}{\sigma^2}(r^2+a^2-\sigma),\\
\hat{\mathbb{V}}_0&=&\left(\begin{array}{cc} \hat{\mathbb{V}}_1 & 0 \\ 0 & \bar{\hat{\mathbb{V}}}_1 
\end{array} \right),\, \hat{\mathbb{V}}_1=U\left(\begin{array}{cc} v_0 & v_1 \\ \bar{v}_1 & -v_0 \end{array} \right)U^*,\\
M_{r_*}&=&\frac{1}{2}\left(\begin{array}{cc} m^1_{r_*} & im^2_{r_*} \\
-im^2_{r_*} & -m^1_{r_*} \end{array}\right),\\
m^1_{r_*}&=&\sqrt{\alpha+1}hD_{r_*}h\sqrt{\alpha+1}-\frac{\beta h}{\sqrt{\alpha+1}}
D_{r_*}\frac{\beta h}{\sqrt{\alpha+1}},\\
m^2_{r_*}&=&\sqrt{\alpha+1}hD_{r_*}\frac{\beta h}{\sqrt{\alpha+1}}+\frac{\beta h}{\sqrt{\alpha+1}}D_{r_*}h\sqrt{\alpha+1},\\
 \\
M_{\theta}&=&\frac{1}{2}\left(\begin{array}{cc} m_{\theta}^1 & im_{\theta}^2 \\ 
-im_{\theta}^2 & -m_{\theta}^1 \end{array} \right),\\
m_{\theta}^1&=&\frac{\beta h\sqrt{a_0}}{\sqrt{\alpha+1}}\left(D_{\theta}+\frac{\cot\theta}{2i}\right)\sqrt{a_0}h\sqrt{\alpha+1}+\sqrt{\alpha+1}h\sqrt{a_0}\left(D_{\theta}+\frac{\cot\theta}{2i}\right)\frac{\beta h\sqrt{a_0}}{\sqrt{\alpha+1}},\\
m_{\theta}^2&=& \frac{\beta h \sqrt{a_0}}{\sqrt{\alpha+1}}\left(D_{\theta}+\frac{\cot\theta}{2i}\right)\frac{h\beta\sqrt{a_0}}{\sqrt{\alpha+1}}-\sqrt{\alpha+1}h\sqrt{a_0}\left(D_{\theta}+\frac{\cot\theta}{2i}\right)\sqrt{a_0}h\sqrt{\alpha+1}.
\end{eqnarray*}
We now want to use the variable introduced in (\ref{hatr}). Let ${\cal H}:=L^2((\R\times S^2,d\hat{r}d\omega);\C^4)$.
We define :
\begin{eqnarray*}
{\cal V} : \begin{array}{rcl} {\cal H}_*&\rightarrow&{\cal H} \\
u(r_*,\omega)&\mapsto&u(r_*(\hat{r},\theta),\omega)k'^{-1/2}(r_*(\hat{r},\theta)). \end{array}
\end{eqnarray*}
${\cal V}$ is a unitary transformation with inverse :
\begin{eqnarray*}
{\cal V}^* : \begin{array}{rcl} {\cal H}&\rightarrow&{\cal H}_* \\
v(\hat{r},\omega)&\mapsto&v(\hat{r}(r_*,\theta),\omega)k'^{1/2}(r_*). \end{array}
\end{eqnarray*}
The hamiltonian we want to work with is
\begin{eqnarray}
\label{3.11}
H:={\cal V}\hat{\notD}{\cal V}^*,
\end{eqnarray}
which acts on ${\cal H}.$ The operator $H$ is selfadjoint with domain
\[D(H)=\{u\in{\cal H};\, {\cal U}^*{\cal V}^*u\in D(\notD)\} \]
and its spectrum is purely absolutely continuous. This follows from the corresponding results for $\notD$. 
In order to calculate $H$ we first observe that
\[ {\cal V}D_{r_*}{\cal V}^*=(k')^{1/2}D_{\hat{r}}(k')^{1/2},\quad {\cal V}D_{\theta}{\cal V}^*=
a\cos\theta D_{\hat{r}}+D_{\theta}. \]
Now observe that :
\begin{eqnarray*}
\frac{1}{2}(h^2(\alpha+1)k'-\frac{\beta^2 h^2}{(\alpha+1)}k'+2\beta h^2a_0a\cos\theta)
=\alpha h^2k'+\beta h^2a_0a\cos\theta=\alpha^2+\beta^2=1,\\
i(2h^2\beta k'+\frac{\beta^2h^2}{\alpha+1}a_0a\cos\theta-(\alpha+1)h^2a_0a\cos\theta)=i(2\alpha\beta+(1-\alpha)\beta-(1+\alpha)\beta)=0
\end{eqnarray*}
and that
\begin{eqnarray*}
\lefteqn{\sqrt{\alpha+1}h(k')^{1/2}\partial_{\hat{r}}((k')^{1/2}h\sqrt{\alpha+1})-
\frac{\beta h}{\sqrt{\alpha+1}}(k')^{1/2}\partial_{\hat{r}}(\frac{\beta h}{\sqrt{\alpha+1}}(k')^{1/2})}\\
&+&\frac{\beta h \sqrt{a_0}}{\sqrt{\alpha+1}}a\cos\theta\partial_{\hat{r}}(h\sqrt{a_0}\sqrt{\alpha+1})
+\sqrt{\alpha+1}h\sqrt{a_0}a\cos\theta\partial_{\hat{r}}(\frac{\sqrt{a_0}\beta h}{\sqrt{\alpha+1}})\\
&=&\frac{1}{2}\partial_{\hat{r}}(k'h^2(\alpha+1))-\frac{1}{2}\partial_{\hat{r}}\left(k'\frac{\beta^2 h^2}{\alpha+1}\right)
+\partial_{\hat{r}}(\beta h^2 a_0a\cos\theta)\\
&=&\partial_{\hat{r}}(\alpha^2+\beta^2)=0,
\end{eqnarray*}
\begin{eqnarray*}
\lefteqn{\sqrt{\alpha+1}h(k')^{1/2}\partial_{\hat{r}}\left((k')^{1/2}\frac{\beta h}{\sqrt{\alpha+1}}\right)
+\frac{\beta h}{\sqrt{\alpha+1}}(k')^{1/2}\partial_{\hat{r}}((k')^{1/2}h\sqrt{\alpha+1})}\\
&+&\frac{\beta h\sqrt{a_0}}{\sqrt{\alpha+1}}\partial_{\hat{r}}\left(a\cos\theta\frac{\sqrt{a_0}h\beta}{\sqrt{\alpha+1}}\right)
-\sqrt{\alpha+1}h\sqrt{a_0}\partial_{\hat{r}}(a\cos\theta\sqrt{a_0}h\sqrt{\alpha+1})\\
&=&\partial_{\hat{r}}(\alpha\beta)+1/2\partial_{\hat{r}}\left(\frac{\beta^3}{\alpha+1}-\beta(\alpha+1)\right)=0.
\end{eqnarray*}
Therefore we obtain :
\begin{eqnarray}
\label{3.12}
H&=&\Gamma^1D_{\hat{r}}+\left(\begin{array}{cc} M_{\theta} & 0 \\ 0 & -M_{\theta} \end{array} \right)
+h^2a_0\hat{\Gamma}^3\frac{D_{\varphi}}{\sin\theta}+h^2c_1+h^2c_2^{\varphi}D_{\varphi}+\hat{V}_{\varphi}D_{\varphi}+\hat{V}_1\nonumber\\
&=&\Gamma^1D_{\hat{r}}+(m^1_{\theta}\Gamma^1-m^2_{\theta}\Gamma^3)-h^2a_0\Gamma^2\frac{D_{\varphi}}{\sin\theta}
+h^2c_1+h^2c_2^{\varphi}D_{\varphi}\nonumber\\
&+&\hat{V}_{\varphi}D_{\varphi}+\hat{V}_1.
\end{eqnarray}
Recalling that
\begin{eqnarray}
\label{POMEGA} 
P_{\omega}:=\left(\begin{array}{cc} M_{\theta} & 0 \\ 0 & -M_{\theta} \end{array} \right)
+\hat{\Gamma}^3h^2a_0\frac{D_{\varphi}}{\sin\theta}={\cal U}h\sqrt{a_0}\notD_{S^2}\sqrt{a_0}h{\cal U}^*,
\end{eqnarray}
where all functions have to be evaluated at $(r_*(\hat{r},\theta),\theta)$, we see that the angular
part is regular. Let us also define :
\begin{eqnarray}
\label{W}
W:=H-\Gamma^1D_{\hat{r}}-P_{\omega}.
\end{eqnarray} 
We put 
\[ {\cal H}^1:=D(H),\quad ||u||^2_{{\cal H}^1}=||Hu||^2+||u||^2. \]
\begin{lemma}
\label{lem6.5}
\begin{eqnarray}
\label{6.19}
(i) \qquad\forall u\in D({\notD})\quad ||D_{r_*}u||_{{\cal H}_*}&\lesssim&||{\notD}u||_{{\cal H}_*}+||u||_{{\cal H}_*},\\
\label{6.20}
 ||a_0{\notD}_{S^2}u||_{{\cal H}_*}&\lesssim& ||{\notD}u||_{{\cal H}_*}+||u||_{{\cal H}_*},\\
\label{6.21}
(ii) \qquad \forall v\in D(H)\quad ||D_{\hat{r}}v||_{\cal H}&\lesssim&||Hv||_{\cal H}+||v||_{\cal H},\\
\label{6.22}
 ||a_0\notD_{S^2}v||_{{\cal H}}&\lesssim& ||Hv||_{\cal H}+||v||_{\cal H}.
\end{eqnarray}
\end{lemma}
{\bf Proof.}

For part $(i)$ see \cite{Da2}. Let us show $(ii)$. From (\ref{6.19}) we infer for $v\in D(H)$:
\[ -\langle {\cal V}{\cal U}\partial_{r_*}^2{\cal U}^*{\cal V}^*v,v \rangle\lesssim (||Hv||+||v||)^2. \]
But we have :
\begin{eqnarray*}
-{\cal V}{\cal U}\partial_{r_*}^2{\cal U}^*{\cal V}^*
&=&-{[}{\cal U}(k')^{1/2},
\partial_{\hat{r}}{]}k'\partial_{\hat{r}}(k')^{1/2}{\cal U}^*\\
&-&\partial_{\hat{r}}
{\cal U}(k')^{3/2}{[}\partial_{\hat{r}},(k')^{1/2}{\cal U}^*{]}-\partial_{\hat{r}}(k')^2\partial_{\hat{r}}\\
&\ge& (\epsilon-\delta)\partial_{\hat{r}}^2-C_{\epsilon}
\end{eqnarray*}
with $\delta=\min(k'^2)>0$.
This gives (\ref{6.21}). Inequality (\ref{6.22}) can be established in a similar way.
We use :
\[ {\cal V}D_{\theta}^2{\cal V}^*=l'^2D_{\hat{r}}^2+2l'D_{\hat{r}}D_{\theta}+D_{\theta}^2 \]
as well as (\ref{6.21}).
\qed

\begin{remark}
\label{constants}
A precise analysis of the constants in \cite{Da2} shows that
\begin{eqnarray*}
\forall \epsilon>0,\, \exists C_{\epsilon}>0\, \forall u\in D(\notD)\quad
||D_{r_*}u||_{{\cal H}_*}&\le& (1+\epsilon)||\notD u||_{{\cal H}_*}+C_{\epsilon}||u||_{{\cal H}_*},\\
\forall \epsilon>0,\, \exists C_{\epsilon}>0\, \forall u\in D(H)\quad
||D_{\hat{r}}u||_{{\cal H}}&\le& (1+\epsilon)||H u||_{{\cal H}}+C_{\epsilon}||u||_{{\cal H}}.
\end{eqnarray*}
\end{remark}
\subsection{Scattering results}
\label{sec3.4.3}
A complete scattering theory for massless Dirac fields in the Kerr metric was obtained in \cite{HN}.
This result has been generalized to the case of massive charged Dirac fields in \cite{Da2}.
In both works the comparison dynamics are the dynamics which are natural in the $(t,r_*,\omega)$ coordinate
system and the $L^a,N^a,M^a$ Newman-Penrose tetrad. We will need comparison dynamics which are natural
with respect to the $(t,\hat{r},\theta,\varphi)$ coordinate system and the $l^a,n^a,m^a$ tetrad.
To this purpose we define :
\begin{eqnarray*}
H_{\leftarrow}&=&\Gamma^1D_{\hat{r}}-\frac{a}{r_+^2+a^2}D_{\varphi}-\frac{qQr_+}{r_+^2+a^2},\\
H_{\rightarrow}&=&{\cal V}{\cal U}\notD_{\rightarrow}{\cal U}^*{\cal V}^*,
\end{eqnarray*}
where $D_{\rightarrow}$ is defined in (\ref{3.17b}).
These operators are selfadjoint on ${\cal H}$ with domains  
$D(H_{\leftarrow})=\{v\in{\cal H}; H_{\leftarrow} v\in{\cal H}\}$,
$D(H_{\rightarrow})=\{v\in {\cal H};{\cal U}^*{\cal V}^*v\in D(\notD_{\rightarrow})\}$. 
Let 
\begin{eqnarray*}
{\cal H}^+&=&\{v=(v_1,v_2,v_3,v_4)\in {\cal H};\, v_1=v_4=0 \},\\
{\cal H}^-&=&\{v=(v_1,v_2,v_3,v_4)\in {\cal H};\, v_2=v_3=0 \}. 
\end{eqnarray*}
We note that 
\[ {\bf 1}_{\R^{\pm}}(-\Gamma^1)=P_{{\cal H}^{\pm}}, \]
where $P_{{\cal H}^{\pm}}$ is the projection from ${\cal H}$ to ${\cal H}^{\pm}$.
We also define the projections $P_{2,3}$ and $P_{1,4}$:
\begin{eqnarray*}
P_{2,3}: \left\{ \begin{array}{ccc} {\cal H} &\rightarrow& (L^2(\R\times S^2))^2 \\
(v_1,v_2,v_3,v_4) &\mapsto & (v_2,v_3) \end{array} \right. \\
P_{1,4}: \left\{ \begin{array}{ccc} {\cal H} &\rightarrow& (L^2(\R\times S^2))^2 \\
(v_1,v_2,v_3,v_4) &\mapsto & (v_1,v_4). \end{array} \right.
\end{eqnarray*}

Let 
\begin{eqnarray}
\label{cv}
\mathfrak{v}:=D_{r_*}\notD_{\rightarrow}^{-1},\, \hat{\mathfrak{v}}={\cal V}{\cal U}\mathfrak{v}{\cal U}^*
{\cal V}^*
\end{eqnarray}
be the "classical velocity operators" associated
to $\notD_{\rightarrow},\, H_{\rightarrow}$.
The following proposition gives the existence of the asymptotic velocity :
\begin{proposition}
\label{prop3.1}
There exist selfadjoint operators $P^{\pm}$ s.t. for all $g\in C_{\infty}(\R)$:
\begin{eqnarray}
\label{asympvel}
g(P^{\pm})=s-\lim_{t\rightarrow\pm\infty}e^{-it H}g\left(\frac{\hat{r}}{t}\right)e^{it H}.
\end{eqnarray}
The operators $P^{\pm}$ commute with $H$.
Furthermore we have :
\begin{eqnarray}
\label{3.13}
g(P^{\pm}){\bf 1}_{\R^{\pm}}(P^{\pm})&=&s-\lim_{t\rightarrow\pm\infty}e^{-itH}g(\hat{\mathfrak{v}})e^{itH}
{\bf 1}_{\R^\pm}(P^{\pm}),\\
\label{3.14}
g(P^{\pm}){\bf 1}_{\R^\mp}(P^{\pm})&=&s-\lim_{t\rightarrow\pm\infty}e^{-itH}g(-\Gamma^1)e^{itH}
{\bf 1}_{\R^\mp}(P^{\pm}),\\
\label{3.15}
\sigma(P^{\pm})&=&\{-1\}\cup{[}0,1{]},\\
\label{3.16}
\sigma_{pp}(H)&=&{\bf 1}_{\{0\}}(P^{\pm})=\emptyset.
\end{eqnarray}
\end{proposition}
\begin{remark}
a) For limits of the form (\ref{asympvel}) we will write the following :
\[ P^{\pm}=s-C_{\infty}-\lim_{t\rightarrow\pm \infty}e^{-itH}\frac{\hat{r}}{t}e^{itH}. \]
b) We can construct in the same way
\[ P^{\pm}_{\rightarrow}=s-C_{\infty}-\lim_{t\rightarrow\pm \infty}e^{-itH_{\rightarrow}}\frac{\hat{r}}{t}
e^{itH_{\rightarrow}}. \]
\end{remark}
{\bf Proof.}

Using the definition of $H$ we see that it is sufficient to show the existence
of 
\[ s-\lim_{t\rightarrow\pm\infty}e^{-it\expnotD}g\left(\frac{\hat{r}(r_*,\theta)}{t}\right)e^{it\expnotD}. \]
By \cite[Theorem 5.2]{Da2} we know that 
\begin{eqnarray}
\label{3.17}
g(\tilde{P}^{\pm})=s-\lim_{t\rightarrow\pm\infty}e^{-it\expnotD}g\left(\frac{r_*}{t}\right)e^{it\expnotD}\quad\mbox{exists.}
\end{eqnarray}
As $g$ is uniformly continuous and $|\hat{r}(r_*,\theta)-r_*|\lesssim 1$ we have :
\[ s-\lim_{t\rightarrow\pm\infty}e^{-it\expnotD}\left(g\left(\frac{\hat{r}(r_*,\theta)}{t}\right)-g\left(\frac{r_*}{t}\right)\right)e^{it\expnotD}=0. \]
This gives the existence. Equalities (\ref{3.13}),\, (\ref{3.15}), (\ref{3.16}) follow directly from the definitions and the corresponding equalities 
in \cite{Da2}, so does (\ref{3.14}) if we replace $\Gamma^1$ by $\hat{\Gamma}^1$. In order to
replace $\hat{\Gamma}^1$ again by $\Gamma^1$ it is sufficient by a density argument
to show 
\begin{eqnarray}
\label{C3.32}
\lefteqn{\forall \epsilon>0,\, \forall J\in C_{\infty}(\R),\, \supp J\subset (-\infty,-\epsilon)}\nonumber\\
&&s-\lim_{t\rightarrow\infty}e^{-it\expnotD}\left(g\left(-\hat{\Gamma}^1\right)
-g\left(-\Gamma^1\right)\right)J\left(\frac{r_*}{t}\right)e^{it\expnotD}=0.
\end{eqnarray}
We then note that it is sufficient to show (\ref{C3.32}) for smooth $g$. Indeed as
$g$ is uniformly continuous we can approximate it in $L^{\infty}$ norm by smooth functions.
Equation (\ref{C3.32}) follows for smooth $g$ from the observation that 
\[ \hat{\Gamma}^1-\Gamma^1=(\Gamma_{ij}),\quad \Gamma_{ij}={\cal O}(e^{\kappa_+ r_*}),\, r_*\rightarrow -\infty. \]
\qed
 
At infinity we define the Dollard modification :
\begin{eqnarray}
\label{Dollard}
U_D(t)=\left\{ \begin{array}{cc}  e^{itH_{\rightarrow}}T\left(e^{i\int_0^t((m-b_0)(\hat{s\mathfrak{v}})mH_{\rightarrow}^{-1}
+c_1(s\hat{\mathfrak{v}}))ds}\right) & m\neq 0,\\
e^{itH_{\rightarrow}} T\left(e^{i\int_0^tc_1(s\hat{\mathfrak{v}})ds}\right) & m=0, \end{array}
\right.
\end{eqnarray}
where $T$ denotes time ordering :
\[ T\left(e^{\int_s^tW(u)du}\right):=\sum_{n=0}^{\infty}\int_{t\ge u_n\ge...\ge u_1\ge s}
...\int W(u_n)...W(u_1)du_n...du_1 \]
and $\hat{\mathfrak{v}}$ is the classical velocity operator (see (\ref{cv})).
\begin{theorem}
\label{th3.2}
The wave operators
\begin{eqnarray}
\label{3.19}
W_{\rightarrow}^{\pm}&=&s-\lim_{t\rightarrow\pm\infty}e^{-itH}U_D(t){\bf 1}_{\R^\pm}(P^{\pm}_{\rightarrow}),\\
\label{3.20}
\Omega_{\rightarrow}^{\pm}&=&s-\lim_{t\rightarrow\pm\infty}U_D(-t)e^{itH}{\bf 1}_{\R^\pm}(P^{\pm}),\\
\label{3.21}
W_{\leftarrow}^{\pm}&=&s-\lim_{t\rightarrow\pm\infty}e^{-itH}e^{itH_{\leftarrow}}P_{{\cal H}^{\mp}},\\
\label{3.22}
\Omega_{\leftarrow}^{\pm}&=&s-\lim_{t\rightarrow\pm\infty}e^{-itH_{\leftarrow}}e^{itH}{\bf 1}_{\R^\mp}(P^{\pm})
\end{eqnarray}
exist and satisfy
\begin{eqnarray}
\label{3.23}
(W_{\leftarrow}^{\pm})^*=\Omega_{\leftarrow}^{\pm},\, (\Omega_{\leftarrow}^{\pm})^*=W_{\leftarrow}^{\pm},\, (W_{\rightarrow}^{\pm})^*=\Omega_{\rightarrow}^{\pm}
,\, (\Omega_{\rightarrow}^{\pm})^*=W_{\rightarrow}^{\pm}.
\end{eqnarray}
\end{theorem}
\begin{remark}
For the proof of the theorem about the Hawking effect we only need the asymptotic completeness
result near the horizon.
\end{remark}
 
{\bf Proof.}

Let
\[ \notD_{\leftarrow}:={\cal U}^*{\cal V}^*H_{\leftarrow}{\cal V}{\cal U},\quad \tilde{P}_{\pm}:=
{\cal U}^*{\cal V}^*P_{{\cal H}^{\pm}}{\cal V}{\cal U}. \]
It is sufficient to show that the following limits exist :
\begin{eqnarray}
\label{3.25}
\tilde{W}_{\rightarrow}^{\pm}&=&s-\lim_{t\rightarrow\pm\infty}e^{-it\expnotD}\tilde{U}_D(t){\bf 1}_{\R^\pm}(\tilde{P}^{\pm}_{\rightarrow}),\\
\label{3.26}
\tilde{\Omega}_{\rightarrow}^{\pm}&=&s-\lim_{t\rightarrow\pm\infty}\tilde{U}_D(-t)e^{it\expnotD}{\bf 1}_{\R^\pm}(\tilde{P}^{\pm}),\\
\label{3.27}
\tilde{W}_{\leftarrow}^{\pm}&=&s-\lim_{t\rightarrow\pm\infty}e^{-it\expnotD}e^{it\expnotD_{\leftarrow}}\tilde{P}_{\mp},\\
\label{3.28}
\tilde{\Omega}_{\leftarrow}^{\pm}&=&s-\lim_{t\rightarrow\pm\infty}e^{-it\expnotD_{\leftarrow}}e^{it\expnotD}{\bf 1}_{\R^\mp}(\tilde{P}^{\pm})
\end{eqnarray}
with 
\begin{eqnarray*}
\tilde{U}_D(t)&=&{\cal U}^*{\cal V}^*U_D(t){\cal V}{\cal U},\\ 
\tilde{P}^{\pm}&=&s-C_{\infty}-\lim_{t\rightarrow\pm\infty} e^{-it\expnotD}\frac{r_*}{t}e^{it\expnotD},
\tilde{P}_{\rightarrow}^{\pm}=s-C_{\infty}-\lim_{t\rightarrow\pm \infty}e^{-it\expnotD_{\rightarrow}}\frac{r_*}{t}e^{it\expnotD_{\rightarrow}}.
\end{eqnarray*}
The existence of the first two limits is contained in \cite[Theorem 5.5]{Da2}. 
We have 
\[ \notD_{\leftarrow}=\left(\begin{array}{cc} A & 0 \\ 0 & -A \end{array} \right)
-\frac{a}{r_+^2+a^2}D_{\varphi}-\frac{qQr_+}{r_+^2+a^2} \]
with 
\begin{eqnarray*} 
A&=&\frac{1}{2}\left(\begin{array}{cc} a_{11} & a_{12} \\ a_{21} & a_{22} \end{array}
\right),\\
a_{11}&=&\sqrt{\alpha+1}(k')^{-1/2}D_{r_*}(k')^{-1/2}\sqrt{\alpha+1}\\
&-&\frac{\beta}{\sqrt{\alpha+1}} (k')^{-1/2}D_{r_*}(k')^{-1/2}\frac{\beta}{\sqrt{\alpha+1}}=-a_{22}\\
a_{12}&=&\sqrt{\alpha+1}(k')^{-1/2}D_{r_*}(k')^{-1/2}\frac{\beta}{\sqrt{\alpha+1}}+hc=a_{21}.
\end{eqnarray*}
By \cite[Theorem 5.4]{Da2} we know that the limits 
\begin{eqnarray*}
\tilde{W}_H^{\pm}&=&s-\lim_{t\rightarrow\pm\infty}e^{-it\expnotD}e^{it\expnotD_H}P_{{\cal H}_*^{\mp}},\\
\tilde{\Omega}_H^{\pm}&=&s-\lim_{t\rightarrow\pm\infty}e^{-it\expnotD_H}e^{it\expnotD}{\bf 1}_{\R^\mp}(\tilde{P}^{\pm})
\end{eqnarray*}
exist. Here ${\cal H}_*^{\mp}$ denote
\begin{eqnarray*}
{\cal H}_*^+&=&\{u=(u_1,u_2,u_3,u_4);\, u_1=u_4=0\},\\ 
{\cal H}_*^-&=&\{u=(u_1,u_2,u_3,u_4);\, u_2=u_3=0\}.
\end{eqnarray*}
It is therefore sufficient to show the existence of
\begin{eqnarray*}
W_{c}^{\pm}&=&s-\lim_{t\rightarrow\pm\infty}e^{-it\expnotD_{\leftarrow}}e^{it\expnotD_H}P_{{\cal H}_*^{\mp}},\\
\Omega_{c}^{\pm}&=&s-\lim_{t\rightarrow\pm\infty}e^{-it\expnotD_H}e^{it\expnotD_{\leftarrow}}\tilde{P}_{\mp}.
\end{eqnarray*}
and that 
\begin{eqnarray}
\label{C3.44}
P_{{\cal H}_*^{\mp}}\Omega_c^{\pm}=\Omega_c^{\pm}.
\end{eqnarray}
The existence of the first limit follows from Cook's method and the fact that 
\begin{eqnarray*}
k'&=&1+{\cal O}(e^{\kappa_+ r_*}),\, r_*\rightarrow-\infty,\\
\alpha&=&1+{\cal O}(e^{\kappa_+r_*}),\, r_*\rightarrow -\infty,\\
\beta&=&{\cal O}(e^{\kappa_+ r_*}),\, r_*\rightarrow-\infty.
\end{eqnarray*}
The existence of $\Omega_{c}^{\pm}$ follows from the existence of 
\[ \hat{\Omega}^{\pm}_{c}=s-\lim_{t\rightarrow\pm\infty}e^{-itH_H}
e^{itH_{\leftarrow}}P_{{\cal H}^{\mp}}, \]
where $H_H={\cal V}{\cal U}\notD_H{\cal U}^*{\cal V}^*$. This allows us to apply 
Cook's method also in this case. We omit the details. It remains to check (\ref{C3.44}).
We note that $\tilde{P}_{\mp}={\bf 1}_{\R^\mp}(\tilde{P}_{\leftarrow}^{\pm})$,
where $\tilde{P}_{\leftarrow}^{\pm}$ is the asymptotic velocity associated to $\notD_{\leftarrow}$.
This follows from the argument used in the proof of Proposition \ref{prop3.1}.
Thus 
\[  \Omega^{\pm}_{c}=\Omega^{\pm}_{c}{\bf 1}_{\R^\mp}(\tilde{P}_{\leftarrow}^{\pm})
={\bf 1}_{\R^\mp}(\tilde{P}_H^{\pm})\Omega_c^{\pm}=P_{{\cal H}_*^{\mp}}\Omega_{c}^{\pm},\]
where $\tilde{P}_H^{\pm}=-\Gamma^1$ is the asymptotic velocity associated to $\notD_H$.
\qed
\section{The Dirac equation on ${\cal M}_{col}$}
\label{sec3.5}
We want to impose a boundary condition on the surface of the star such that the evolution
can be described by a unitary propagator. We will use the conserved current 
\[ V_a=\phi_A\bar{\phi}_{A'}+\bar{\chi}_A\chi_{A'}. \]
Integrating over the domain indicated in Figure \ref{DomInt}\footnote{The Penrose compactification of block $I$ shown in Figure \ref{DomInt} is performed in 
Appendix \ref{AppB}.} and supposing that the field is
0 in a neighborhood of $i^0$ gives by Stokes' theorem :
\begin{eqnarray*}
\int_{\Sigma_s}V_aT^a d\sigma_{\Sigma_s} -\int_{\Sigma_t} V_aT^a d\sigma_{\Sigma_t} +\int_{{\cal S}_I}V_a {\cal N}^a
d\sigma_{{\cal S}_I} =0,
\end{eqnarray*}
where ${\cal N}^a$ is the normal to the surface of the star. Therefore 
the necessary condition for charge conservation outside the collapsing body is 
\begin{eqnarray}
\label{charexst}
V_a{\cal N}^a=0\quad \mbox{on}\quad {\cal S}.
\end{eqnarray}
We will impose :
\begin{eqnarray}
\label{bound1}
\left.\begin{array}{rcl} {\cal N}^{AA'}\phi_A&=&\frac{1}{\sqrt{2}}e^{-i\nu}\chi^{A'},\\
{\cal N}^{AA'}\chi_{A'}&=&\frac{1}{\sqrt{2}}e^{i\nu}\phi^A. \end{array} \right\}\quad\mbox{on}\quad {\cal S}.
\end{eqnarray}
Here $\nu$ is the so called chiral angle.
We note that (\ref{bound1}) implies (\ref{charexst}) :
\begin{eqnarray*}
\sqrt{2}{\cal N}^aV_a&=&\sqrt{2}({\cal N}^{AA'}\phi_A\bar{\phi}_{A'}+{\cal N}^{AA'}\bar{\chi}_A\chi_{A'})\\
&=&e^{-i\nu}\chi^{A'}\bar{\phi}_{A'}+e^{-i\nu}\bar{\phi}^{A'}\chi_{A'}=0.
\end{eqnarray*}
\begin{figure}
\centering\epsfig{figure=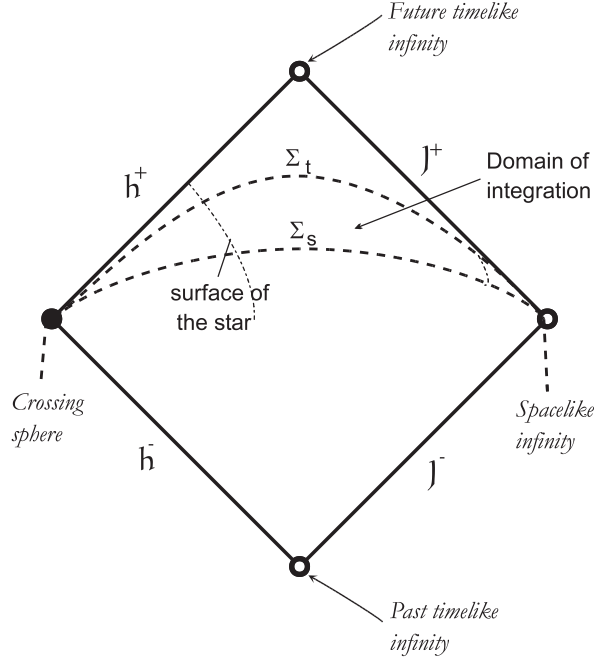,width=8cm}
\caption{The surface of the star and the domain of integration.}
\label{DomInt}
\end{figure}
From (\ref{bound1}) we obtain :
\begin{eqnarray*}
{\cal N}_{BA'}{\cal N}^{AA'}\phi_A=-\frac{1}{2}\phi_B.
\end{eqnarray*}
We have :
\begin{eqnarray*}
\epsilon_{AC}{\cal N}_{BA'}{\cal N}^{AA'}&=&{\cal N}_{BA'}{\cal N}_C^{\,\, A'}=:\kappa_{BC},\\
\kappa_{CB}&=&{\cal N}_{CA'}{\cal N}_B^{\, \, A'}=-\kappa_{BC}.
\end{eqnarray*}
From the antisymmetry of $\kappa_{BC}$ follows :
\[ \kappa_{BC}=\frac{1}{2}\kappa_A^{\,\, A}\epsilon_{BC}=\frac{1}{2}{\cal N}^a{\cal N}_a\epsilon_{BC}. \]
Thus (\ref{bound1}) implies :
\[ {\cal N}^a{\cal N}_a\phi_B=-\phi_B. \]
We therefore impose 
\[ {\cal N}^a{\cal N}_a=-1. \]
This avoids that the boundary condition imposes $(\phi_A,\chi^{A'})=0$ on ${\cal S}$.
We will from now on suppose ${\cal N}_a$ to be past directed, but the opposite
choice would of course be possible.
Let us now rewrite condition (\ref{bound1}) using a coordinate system and a spin frame.
We have :
\begin{eqnarray*}
{\cal N}^{\bf AA'}&=&g^{\bf aAA'}{\cal N}_{\bf a}=\,^tA^{\bf a}{\cal N}_{\bf a},\\
{\cal N}_{\bf AA'}&=&g^{\bf a}_{\, \, \bf AA'}{\cal N}_{\bf a}=B^{\bf a}{\cal N}_{\bf a}.
\end{eqnarray*}
The boundary condition (\ref{bound1}) implies
\begin{eqnarray}
\label{bound2}
& &\sqrt{2}\left(\begin{array}{cc} 0 & B^{\bf a} \\ A^{\bf a} & 0 \end{array} \right) 
{\cal N}_{\bf a}\left(\begin{array}{c} \phi_{\bf A} \\ \chi^{\bf A'} \end{array}\right)
=\left(\begin{array}{cc} -e^{i\nu} & 0 \\ 0 & e^{-i\nu}  \end{array}\right)
\left(\begin{array}{c} \phi_{\bf A} \\ \chi^{\bf A'} \end{array}\right)\nonumber\\
&\Leftrightarrow&\left(\begin{array}{cc} \frac{1}{i} & 0 \\ 0 & -\frac{1}{i} \end{array} \right) 
\gamma^{\bf a}{\cal N}_{\bf a}\left(\begin{array}{c} \phi_{\bf A} \\ \chi^{\bf A'} \end{array}\right)
=\left(\begin{array}{cc} -e^{i\nu} & 0 \\ 0 & e^{-i\nu} \end{array}\right)
\left(\begin{array}{c} \phi_{\bf A} \\ \chi^{\bf A'} \end{array}\right)\nonumber\\
&\Leftrightarrow& 
\gamma^{\bf a}{\cal N}_{\bf a}\left(\begin{array}{c} \phi_{\bf A} \\ \chi^{\bf A'} \end{array}\right)
=-i{\cal B}\left(\begin{array}{c} \phi_{\bf A} \\ \chi^{\bf A'} \end{array}\right).
\end{eqnarray}
where ${\cal B}=e^{-i\nu\gamma^5},\, \gamma^5=-i\gamma^0\gamma^1\gamma^2\gamma^3=Diag(1,1,-1,-1)$,
$\nu\in \R$. The boundary condition 
(\ref{bound2}) is usually called a MIT  boundary condition. Using formula (\ref{DiracMatrix})
we find\footnote{As we will see we do not need the explicit form of $\gamma^{\varphi},\,
\hat{\gamma}^{\varphi}$.} :
\begin{eqnarray*}
\gamma^t=\kappa_0\gamma^0,\, \gamma^{r_*}=\kappa_0h^2\gamma^1,\, \gamma^{\theta}=\kappa_0a_0h^2\gamma^2,\, \kappa_0=\frac{\sigma}{\rho\sqrt{\Delta}}.
\end{eqnarray*}
In the $(t,r_*,\theta,\varphi)$ coordinate system and the $L^a,\,N^a,\,M^a,\,\bar{M}^a$ tetrad
the boundary condition that we impose reads :
\begin{eqnarray}
\label{C3.48}
\sum_{\mu\in\{t,r_*,\theta,\varphi\}}{\cal N}_{\mu}\gamma^{\mu}\Phi=-i{\cal B}\Phi\quad \mbox{on}\quad{\cal S},
\end{eqnarray}
where ${\cal N}_{\mu}$ are the coordinates of the conormal of the surface of the star
in the $(t,r_*,\theta,\varphi)$ coordinate system. 
Let us now consider the $(t,\hat{r},\theta,\varphi)$ coordinate system and the $l^a,n^a,m^a,\bar{m}^a$
tetrad. We denote $\hat{\gamma}^t,\, \hat{\gamma}^{\hat{r}},\,\hat{\gamma}^{\theta},\,
\hat{\gamma}^{\varphi}$ the Dirac matrices with respect to these choices. We find :
\begin{eqnarray*}
\hat{\gamma}^t&=&\kappa_0\gamma^0,\,
\hat{\gamma}^{\hat{r}}=\kappa_0\gamma^1,\,
\hat{\gamma}^{\theta}=\kappa_0a_0h^2\hat{\gamma}^2,\,
\hat{\gamma}^2=i\left(\begin{array}{cc} 0 & \hat{\sigma}^2 \\ \hat{\sigma}^2 & 0 \end{array}
\right).
\end{eqnarray*}
Putting $\Psi={\cal V}{\cal U}\Phi$ we find the following boundary condition for $\Psi$ : 
\begin{eqnarray}
\label{boundhr}
\sum_{\hat{\mu}\in\{t,\hat{r},\theta,\varphi\}}{\cal N}_{\hat{\mu}}\hat{\gamma}^{\hat{\mu}}\Psi&=&-i{\cal B}\Psi\quad \mbox{on}\quad{\cal S}.
\end{eqnarray}
Here ${\cal N}_{\hat{\mu}}$ are the coordinates of the conormal in the $(t,\hat{r},\theta,\varphi)$
coordinate system. We will use in the following the $(t,\hat{r},\theta,\varphi)$
coordinate system.

We introduce the following Hilbert spaces:
\begin{eqnarray}
\label{3.32}
{\cal H}_t=((L^2(\Sigma^{col}_t,d\hat{r}d\omega))^4,||.||_t),
\end{eqnarray}
where the norm $||.||_t$ is defined by
\begin{eqnarray}
\label{3.33}
||\Psi||_t=||{[}\Psi{]}_L||,\quad {[}\Psi{]}_L(\hat{r},\omega)=\left\{\begin{array}{cc} \Psi(\hat{r},\omega) & \hat{r}\ge \hat{z}(t,\theta)\\
0 & \hat{r}\le \hat{z}(t,\theta). \end{array} \right.
\end{eqnarray}
Let
\begin{eqnarray*}
{\cal H}^1_t=\{u\in{\cal H}_t;\, Hu\in {\cal H}_t\},
\quad ||u||_{{\cal H}^1_t}^2=||u||_t^2+||Hu||_t^2.
\end{eqnarray*}
We also need an extension from ${\cal H}^1_t$ to ${\cal H}^1$. To this purpose we put for 
$\phi\in {\cal H}^1_t$ :
\begin{eqnarray*}
{[}\phi{]}_H(\hat{r},\omega)=\left\{\begin{array}{cc} \phi(\hat{r},\omega) & \hat{r}\ge \hat{z}(t,\theta)\\ 
\phi(2\hat{z}(t,\theta)-\hat{r},\omega) & \hat{r}\le \hat{z}(t,\theta) \end{array} \right..
\end{eqnarray*}
It is easy to check that ${[}\phi{]}_H$ is in ${\cal H}^1$. 
The operator $\notD_t$, the spaces ${\cal H}_{*t},{\cal H}^{1}_{*t}$ as well as the extension
${[}...{]}_H^*$ are defined 
in an analogous way using the $(t,r_*,\theta,\varphi)$ coordinate system and the $(L^a, N^a, M^a)$
tetrad.
On ${\cal M}_{col}$ we consider the following mixed problem~:
\begin{eqnarray}
\label{3.34}
\left. \begin{array}{rcl} \partial_t\Psi&=&iH_t\Psi, \quad  \hat{z}(t,\theta)<\hat{r},\\
(\sum_{\hat{\mu}\in\{t,\hat{r},\theta,\varphi\}}{\cal N}_{\hat{\mu}}\hat{\gamma}^{\hat{\mu}})\Psi(t,\hat{z}(t,\theta),\omega)
&=&-i{\cal B}\Psi(t,\hat{z}(t,\theta),\omega), \\
\Psi(t=s,.)&=&\Psi_s(.).  \\
\end{array}\right\}
\end{eqnarray}
Here the operator $H_t$ is given by 
\begin{eqnarray*}
H_t&=&H,\\
D(H_t)&=&\left\{\Psi\in {\cal H}^1_t;\, 
\left(\sum_{\hat{\mu}\in\{t,\hat{r},\theta,\varphi\}}{\cal N}_{\hat{\mu}}\hat{\gamma}^{\hat{\mu}}\Psi\right)(t,\hat{z}(t,\theta),\omega)
=-i{\cal B}\Psi(t,\hat{z}(t,\theta),\omega)\right\}.
\end{eqnarray*}
\begin{remark}[Explicit Calculation]
\label{rem3.4}

It will be helpful in the following to have a more explicit form of the boundary condition.
We choose the $(t,r_*,\theta,\varphi)$ coordinate system and the $L^a,N^a,M^a,\bar{M}^a$
tetrad. The conormal of the surface of the star is :
\[ {\cal N}_a=w(t,\theta)\rho a_0h^2(\dot{z}dt-dr_*+(\partial_{\theta}z)d\theta) \]
and the boundary condition reads :
\begin{eqnarray}
\label{EB1}
w(t,\theta)\left( \begin{array}{cccc} 0 & 0 & \dot{z}-h^2 & (\partial_{\theta}z)a_0h^2 \\
0 & 0 & (\partial_{\theta}z)a_0h^2 & \dot{z}+h^2 \\ -(\dot{z}+h^2) & a_0h^2(\partial_{\theta}z) & 0 & 0 \\
a_0h^2(\partial_{\theta}z) & -\dot{z}+h^2 & 0 & 0 \end{array} \right)\Phi&=&-{\cal B}\Phi\\
\Leftrightarrow w(t,\theta)\Gamma^4(-\dot{z}-h^2\Gamma^1+(\partial_{\theta}z)a_0h^2\Gamma^2)\Phi
&=&-i{\cal B}\Phi\nonumber.
\end{eqnarray}
Here $w(t,\theta)$ is a smooth function. We compute
\[ {\cal N}^a=\rho^{-1}a_0^{-1}h^{-2}w(\dot{z}\partial_t+h^4\partial_{r_*}
+\dot{z}\frac{a(2Mr-Q^2)}{\sigma^2}\partial_{\varphi}-a_0^2h^4(\partial_{\theta}z)\partial_{\theta}). \]
Normalization ${\cal N}^a{\cal N}_a=-1$ gives :
\[ w(t,\theta)=(h^4-\dot{z}^2+a_0^2h^4(\partial_{\theta}z)^2)^{-1/2}. \]
In an analogous way we find in the $(t,\hat{r},\theta,\varphi)$
coordinate system and using the\\
$(l^a,n^a,m^a,\bar{m}^a)$ tetrad :
\[ {\cal N}_a=\hat{w}\rho a_0h^2(\dot{\hat{z}}dt-d\hat{r}+(\partial_{\theta}\hat{z})d\theta). \]
Therefore the boundary condition reads :
\[ \hat{w}(t,\theta)\hat{\Gamma}^4(-\dot{\hat{z}}-\Gamma^1+(\partial_{\theta}\hat{z})a_0h^2\hat{\Gamma}^2)\Psi=-i{\cal B}\Psi. \]
We compute :
\[ {\cal N}^a=\rho^{-1}a_0^{-1}h^{-2}\hat{w}(\dot{\hat{z}}\partial_t+\dot{\hat{z}}
\frac{a(2Mr-Q^2)}{\sigma^2}\partial_{\varphi}+(1-(\partial_{\theta}\hat{z})\frac{\Delta}{\sigma^2}l')
\partial_{\hat{r}}+(l'-(\partial_{\theta}\hat{z}))\frac{\Delta}{\sigma^2}\partial_{\theta}). \]
Normalization gives :
\[ \hat{w}(t,\theta)=(h^4k'^2-\dot{\hat{z}}^2+a_0^2h^4((\partial_{\theta}\hat{z})-l')^2)^{-1/2}. \]
Note that by assumption (\ref{star0b}) $h^4k'^2>\dot{\hat{z}}^2$ and thus $h^4>\dot{z}^2$.
We also note that by the previous considerations we have :
\[ Rank\left(-\dot{\hat{z}}-\Gamma^1+(\partial_{\theta}\hat{z})a_0h^2\hat{\Gamma}^2+\frac{i}{\hat{w}}\hat{\Gamma}^4{\cal B}\right)=2. \]
\end{remark}
Let $\check{\notD}_t=\notD_t+\dot{z}D_{r_*},\, \check{H}_t=H_t+\dot{\hat{z}}D_{\hat{r}}$.
We argue that $\check{H}_tv\in {\cal H}\Leftrightarrow H_tv\in {\cal H}$. Let 
$\check{H}=H+\dot{\hat{z}}D_{\hat{r}}.$ From $\check{H}_tu\in {\cal H}_t$ we infer 
$\check{H}{[}u{]}_H\in {\cal H}$. But ($0<\delta<\epsilon_1,\, \epsilon_2>0$) :
\begin{eqnarray*}
||\check{H}{[}u{]}_H||&\ge& ||H{[}u{]}_H||-(1-\epsilon_1)||D_{\hat{r}}{[}u{]}_H||\\
&\ge&(\epsilon_1-\delta)||D_{\hat{r}}{[}u{]}_H||+\epsilon_2||a_0\notD_{S^2}{[}u{]}_H||
-C_1||{[}u{]}_H||\\
&\ge& \tilde{\epsilon}||H{[}u{]}_H||-\tilde{C}||{[}u{]}_H||,\, \tilde{\epsilon}>0.
\end{eqnarray*}
Here we have used Remark \ref{constants} and (\ref{star0b}). The implication
$H_tu\in {\cal H}\Rightarrow \check{H}_tu\in {\cal H}$ is shown in a similar
way. We therefore put : $D(\check{\notD}_t)=D(\notD_t),\, D(\check{H}_t)=D(H_t).$ Note that 
$\check{\notD}_0=\notD_0,\, \check{H}_0=H_0$.
\begin{lemma}
\label{Selfadj}
The operators $(\check{\notD}_t,D(\check{\notD}_t))$ and $(\check{H}_t,D(\check{H}_t))$
are selfadjoint.
\end{lemma}
{\bf Proof.}
The selfadjointness of $(\check{H}_t,D(\check{H}_t))$ follows from the selfadjointness
of $(\check{\notD}_t,D(\check{\notD}_t))$, which we show in the following. We calculate for
$u,\, v\in D(\check{\notD}_t)$
\begin{eqnarray*}
\langle \check{\notD}_tu,v\rangle&=&\langle u,\check{\notD}_tv\rangle 
+\frac{1}{i}\int_{S^2}\langle (-\dot{z}-\Gamma^1h^2+(\partial_{\theta}z)a_0h^2\Gamma^2)u,v\rangle (z(t,\theta),\omega)d\omega\\
&=&-i\int_{S^2}(-(h^2+\dot{z})u_1\bar{v}_1+(\partial_{\theta}z)a_0h^2u_2\bar{v}_1+(\partial_{\theta}z)a_0h^2u_1\bar{v}_2)(z(t,\theta),\omega)d\omega\\
&-&i\int_{S^2}((h^2-\dot{z})u_2\bar{v}_2+(h^2-\dot{z})u_3\bar{v}_3-a_0h^2(\partial_{\theta}z)u_4\bar{v}_3)(z(t,\theta),\omega)d\omega\\
&-&i\int_{S^2}(-a_0h^2(\partial_{\theta}z)u_3\bar{v}_4
-(h^2+\dot{z})u_4\bar{v}_4)(z(t,\theta),\omega)d\omega\\
&=&-i\int_{S^2}\frac{e^{i\nu}}{w}(-u_3\bar{v}_1-u_4\bar{v}_2+u_3\bar{v}_1+u_4\bar{v}_2)(z(t,\theta),\omega)d\omega=0.
\end{eqnarray*}
Therefore $\check{\notD}_t$ is symmetric. 
We have to show that
$D(\check{\notD}_t^*)=D(\check{\notD}_t)$. We have :
\begin{eqnarray}
\label{Selfadj1}
v\in D(\check{\notD}_t^*)\Leftrightarrow \forall u\in D(\check{\notD}_t)\quad
|\langle\check{\notD}_tu,v\rangle|\le C||u||.
\end{eqnarray}
Taking $u\in C_0^{\infty}(\Sigma_t^{col};\C^4)$ in (\ref{Selfadj1}) we find that
$\check{\notD}_tv\in {\cal H}_*$, thus $v(z(t,\theta),\omega)$ is well defined.
For $u\in D(\check{\notD}_t),\, v\in D(\check{\notD}_t^*)$ we compute :
\begin{eqnarray*}
\lefteqn{\langle \check{\notD}_tu,v\rangle}\\
&=&\langle u,\check{\notD}_tv\rangle +i\int_{S^2}\langle (-\dot{z}-h^2\Gamma^1
+(\partial_{\theta}z)a_0h^2\Gamma^2)u,v\rangle(z(t,\theta),\omega)d\omega\\
&=&-i\int_{S^2}(-e^{i\nu}u_3\bar{v}_1-e^{i\nu}u_4\bar{v}_2+(h^2-\dot{z})u_3\bar{v}_3-a_0h^2(\partial_{\theta}z)u_4\bar{v}_3d\omega\\
&+&i\int_{S^2}a_0h^2(\partial_{\theta}z)u_3\bar{v}_4-(h^2+\dot{z})u_4\bar{v}_4)(z(t,\theta),\omega)d\omega
+\langle u,\check{\notD}_tv\rangle\\ 
&=&-i\int_{S^2}u_3\overline{(-e^{-i\nu}v_1+(h^2-\dot{z})v_3-a_0h^2(\partial_{\theta}z)v_4)}
(z(t,\theta),\omega)d\omega\\
&-&i\int_{S^2}u_4\overline{(-e^{-i\nu}v_2-(\partial_{\theta}z)a_0h^2v_3-(h^2+\dot{z})v_4)}(z(t,\theta),\omega)
d\omega+\langle u,\check{\notD}_tv\rangle.
\end{eqnarray*}
It follows from (\ref{Selfadj1}) :
\begin{eqnarray}
\label{Selfadj2}
\left|\int_{S^2}u_3\overline{(-e^{-i\nu}v_1+(h^2-\dot{z})v_3-a_0h^2(\partial_{\theta}z)v_4)}
(z(t,\theta),\omega)d\omega\right.\nonumber\\
+\left.\int_{S^2}u_4\overline{(-e^{-i\nu}v_2-(\partial_{\theta}z)a_0h^2v_3-(h^2+\dot{z})v_4)}
(z(t,\theta),\omega)d\omega\right|\le C||u||.
\end{eqnarray}
Let $\phi\in C_0^{\infty}(\R),\, \phi(0)=1,\, u_{\epsilon}(r_*,\omega)=\phi\left(\frac{r_*-z(t,\theta)}{\epsilon}\right)
u(r_*,\omega)$. Clearly 
\[u_{\epsilon}(z(t,\theta),\omega)=u(z(t,\theta),\omega)\,\mbox{and}\, u_{\epsilon}\in D(\check{\notD}_t).\] 
We estimate :
\begin{eqnarray*}
||u||^2&=&\int_{S^2}\int_{z(t,\theta)}^{\infty}|u_{\epsilon}|^2dr_*d\omega\\
&\lesssim&\epsilon ||{[}u{]}_H^*||^2_{L^{\infty}(\R;(L^2(S^2))^4)}\\
&\lesssim& \epsilon ||{[}u{]}_H^*||_{H^1(\R;(L^2(S^2))^4)}
\lesssim \epsilon ||{[}u{]}_H^*||_{{\cal H}_*^1}^2.
\end{eqnarray*}
Thus if we replace $u$ by $u_{\epsilon}$ in (\ref{Selfadj2}), then the term on the R.H.S.
goes to zero when $\epsilon\rightarrow 0$, whereas the term on the L.H.S. remains
unchanged. It follows :
\begin{eqnarray}
\label{Selfadj3}
\forall u\in D(\check{\notD}_t)&\quad&
\int_{S^2}u_3\overline{(-e^{-i\nu}v_1+(h^2-\dot{z})v_3-a_0h^2(\partial_{\theta}z)v_4)}
(z(t,\theta),\omega)d\omega\nonumber\\
&+&\int_{S^2}u_4\overline{(-e^{-i\nu}v_2-(\partial_{\theta}z)a_0h^2v_3-(h^2+\dot{z})v_4)}
(z(t,\theta),\omega)d\omega=0.\nonumber\\
\end{eqnarray}
Let ${\cal R}=\{(z(t,\theta),\omega);\, \omega\in S^2\}.$ We claim:
\begin{eqnarray}
\label{Selfadj4}
\forall f\in L^2({\cal R};\C^2)\, \exists u\in D(\check{\notD}_t)\quad u_{3,4}|_{\cal R}=f.
\end{eqnarray}
From (\ref{Selfadj3}), (\ref{Selfadj4}) follows that $v$ satisfies the boundary condition.
Therefore it remains to show (\ref{Selfadj4}). Let 
\begin{eqnarray*}
Ker\left(w\Gamma^4\left(-\dot{z}-h^2\Gamma^1+(\partial_{\theta}z)a_0h^2\Gamma^2+\frac{i}{w}\Gamma^4{\cal B}\right)\right)
=span\{w_1,w_2\}
\end{eqnarray*}
with $w_1=\, ^t(w_{11},...,w_{14})$ and $w_2=\, ^t(w_{21},...,w_{24})$. The vectors $^t(w_{13},w_{14})$ and $^t(w_{23},w_{34})$ are linearly independent.
Indeed supposing 
\begin{eqnarray*}
 \left(\begin{array}{c} w_{13} \\ w_{14} \end{array} \right)=\alpha \left(\begin{array}{c} w_{23} \\ w_{24} \end{array} \right)\quad 
\mbox{we find}\quad M\left(\begin{array}{c} w_{11} \\ w_{12} \end{array} \right)=\alpha M
\left(\begin{array}{c} w_{21} \\ w_{22} \end{array} \right) 
\end{eqnarray*}
with
\[ M=\left(\begin{array}{cc} -(\dot{z}+h^2) & a_0h^2(\partial_{\theta}z) \\ a_0h^2(\partial_{\theta}z) & h^2-\dot{z} \end{array} \right). \]
The matrix $M$ being invertible we find :
\[ \left(\begin{array}{c} w_{11} \\ w_{12} \end{array} \right)=\alpha \left(\begin{array}{c} w_{21} \\ w_{22} \end{array} \right) \]
which is a contradiction.  
If $f\in L^2({\cal R};\C^2)$, define $g\in L^2({\cal R};\C^4)$ in the following way~:
\begin{eqnarray*}
\exists \alpha,\, \beta\quad \alpha \left(\begin{array}{c} w_{13} \\ w_{14} \end{array}\right)+\beta\left(\begin{array}{c} w_{23} \\ w_{24} \end{array} \right)=f.
\end{eqnarray*} 
We put $g:=\alpha w_1+\beta w_2.$ By the surjectivity of the trace operator there exists
$u\in {\cal H}^1_{*t}$ s.t. $u|_{\cal R}=g$. By construction of $g$, $u\in D(\tilde{\notD}_t)$.
\qed 
\\

The problem (\ref{3.34}) is solved by the following proposition.
\begin{proposition}
\label{prop3.3}
Let $\Psi_s\in D(H_s)$. Then there exists a unique solution 
\[{[}\Psi(.){]}_H={[}U(.,s)\Psi_s{]}_H\in C^1(\R_t;{\cal H})\cap C(\R_t;{\cal H}^1)\]
of (\ref{3.34}) s.t. 
for all $t\in\R$ $\Psi(t)\in D(H_t)$. Furthermore we have $||\Psi(t)||=||\Psi_s||$
and $U(t,s)$ possesses an extension to an isometric and strongly continuous propagator 
from ${\cal H}_s$ to ${\cal H}_t$ s.t. for all $\Psi_s\in D(H_s)$ we have:
\begin{eqnarray*} 
\frac{d}{dt}U(t,s)\Psi_s=iH_tU(t,s)\Psi_s.
\end{eqnarray*}
\end{proposition}
The proposition follows from Proposition \ref{prop6.3}, which is proven in the appendix.
\chapter{Dirac Quantum Fields}
\label{sec4}
We adopt the approach of Dirac quantum fields in the spirit of \cite{Di80} and \cite{Di82}.
This approach is explained in Section \ref{4.2}. In Section \ref{4.1} we recall the second 
quantization of Dirac fields (see e.g. \cite{BR} for a detailed discussion of the second
quantization procedure and \cite{Th} for the special case of the Dirac equation). In Section \ref{sec4.3} we present the theorem about the Hawking effect.

\section{Second Quantization of Dirac Fields}
\label{sec4.1}
We first explain the construction in the case of one kind of noninteracting fermions.
The one fermion space is a complex Hilbert space ${\GH}$ with scalar product $\langle.,.\rangle$ 
that we suppose linear with respect to the first argument. The space of $n$ fermions is 
the antisymmetric $n$-tensor product of $\mathfrak{H}$:
\begin{eqnarray*}
{\cal F}^{(0)}(\mathfrak{H}):=\C,\, 1\le n\Rightarrow {\cal F}^{(n)}(\mathfrak{H}):=\Lambda _{\nu=1}^n \mathfrak{H}.
\end{eqnarray*}
The Fermi-Fock space is defined as:
\begin{eqnarray*}
{\cal F}(\GH):=\oplus_{n=0}^{\infty}{\cal F}^{(n)}(\mathfrak{H}).
\end{eqnarray*}
For $f\in {\GH}$ we construct the {\it fermion annihilation operator} $a_{\GH}(f)$, and the 
{\it fermion creation operator} $a^*_{\GH}(f)$ by putting:
\begin{eqnarray*}
a_{\GH}(f)&:& {\cal F}^{(0)}(\GH)\rightarrow\{0\},\, 1\le n\quad a_{\GH}(f):{\cal F}^{(n)}(\GH)\rightarrow{\cal F}^{(n-1)}(\GH),\\
a_{\GH}(f)(f_1\wedge...\wedge f_n)&=&\frac{\sqrt{n}}{n!}\sum_{\sigma}\epsilon({\sigma})\langle f_{\sigma(1)},f\rangle f_{\sigma(2)}\wedge ... \wedge f_{\sigma(n)},\\
0\le n \quad a^*_{\GH}(f)&:&{\cal F}^{(n)}(\GH)\rightarrow{\cal F}^{(n+1)}(\GH),\\
a^*_{\GH}(f)(f_1\wedge...\wedge f_n)&=&\frac{\sqrt{n+1}}{n!}\sum_{\sigma}\epsilon({\sigma})f\wedge f_{\sigma(1)}\wedge ... \wedge f_{\sigma(n)},
\end{eqnarray*}
where the sum runs over all permutations $\sigma$ of $\{1,...,n\}$ and $\epsilon(\sigma)$
is one if $\sigma$ is even and $-1$ if $\sigma$ is odd. In contrast to the boson case,
$a_{\GH}(f)$ and $a^*_{\GH}(f)$ have bounded extensions on ${\cal F}(\GH)$, which we still denote 
$a_{\GH}(f), a^*_{\GH}(f)$. They satisfy:
\begin{eqnarray}
\label{4.1}
||a_{\GH}(f)||&=&||a_{\GH}^*(f)||=||f||,\\
\label{4.2}
a_{\GH}^*(f)&=&(a_{\GH}(f))^*.
\end{eqnarray}
Another important feature of the creation and annihilation operators is that they satisfy
the canonical anti-commutation relations (CAR's):
\begin{eqnarray}
\label{4.3}
a_{\GH}(f)a_{\GH}(g)+a_{\GH}(g)a_{\GH}(f)&=&0,\\
\label{4.4}
a^*_{\GH}(f)a^*_{\GH}(g)+a^*_{\GH}(g)a^*_{\GH}(f)&=&0,\\
\label{4.5}
a_{\GH}^*(f)a_{\GH}(g)+a_{\GH}(g)a_{\GH}^*(f)&=&\langle f,g\rangle {\bf 1}.
\end{eqnarray}
The CAR algebra on $\GH$ is the $C^*$-algebra ${\cal U}(\GH)$ generated by 
the identity and the $a_{\GH}(f),\,f\in\GH$. 

Let us now consider a situation where the 
classical fields obey the Schr\"odinger type equation ($H$ selfadjoint on $\GH$):
\begin{eqnarray*}
\partial_t\Psi=iH\Psi.
\end{eqnarray*}
A gauge invariant quasi-free state $\omega$ on ${\cal U}(\GH)$ satisfies the $(\beta,\mu)$-KMS
condition, $0<\beta$, if it is characterized by the two point function
\begin{eqnarray*}
\omega(a^*_{\GH}(f)a_{\GH}(f))=\langle ze^{\beta H}(1+ze^{\beta H})^{-1}f,g\rangle,
\end{eqnarray*}
where $z$ is the {\it activity} given by $z=e^{\beta\mu}$. This state is a model for the {\it ideal Fermi 
gas} with temperature $0<T=\beta^{-1}$ and {\it chemical potential $\mu$}.

In the case of charged spinor fields we have to consider both kinds of fermions, the particles
and the antiparticles. The space of the classical charged spinor fields is given by a complex
Hilbert space $\GH$ together with an anti-unitary operator $\Upsilon$ on $\GH$ (the {\it charge 
conjugation}). We assume $\GH$ is split into two orthogonal subspaces :
\begin{eqnarray*}
\GH=\GH_+\oplus\GH_-.
\end{eqnarray*}
We define the one particle space:
\begin{eqnarray*}
\gh_+=\GH_+
\end{eqnarray*}
and the one antiparticle space
\begin{eqnarray*}
\gh_-=\Upsilon\GH_-.
\end{eqnarray*}
The space of $n$ particles and $m$ antiparticles is given by the tensor product of the previous
spaces:
\begin{eqnarray*}
{\cal F}^{(n,m)}:={\cal F}^{(n)}(\gh_+)\otimes{\cal F}^{(m)}(\gh_-).
\end{eqnarray*}
The {\it Dirac-Fermi-Fock} space is given by :
\begin{eqnarray*}
{\cal F}(\GH):=\bigoplus_{n,m=0}^{\infty}{\cal F}^{(n,m)}.
\end{eqnarray*}
We will denote the elements $\Psi$ of ${\cal F}(\GH)$ by sequences
\begin{eqnarray*}
\Psi=(\Psi^{(n,m)})_{n,m\in \N},\, \Psi^{(n,m)}\in {\cal F}^{(n,m)},
\end{eqnarray*}
the {\it vacuum vector} is the vector $\Omega_{vac}$ defined by
\begin{eqnarray*}
\Omega_{vac}^{(0,0)}=1,\, (n,m)\neq (0,0)\Rightarrow\Omega_{vac}^{(n,m)}=0.
\end{eqnarray*}
For $\phi_{\pm}\in\gh_{\pm}$ we define the {\it particle annihilation operator}, $a(\phi_+)$, 
the {\it particle creation operator}, $a^*(\phi_+)$, the {\it antiparticle annihilation operator},
$b(\phi_-)$, and the {\it antiparticle creation operator}, $b^*(\phi_-)$, by putting for 
$\Psi_+^{(n)}\otimes\Psi_-^{(m)}\in {\cal F}^{(n,m)}$:
\begin{eqnarray}
\label{4.6}
a(\phi_+)(\Psi_+^{(n)}\otimes\Psi_-^{(m)})
=(a_{\gh_+}(\phi_+)(\Psi_+^{(n)}))\otimes\Psi_-^{(m)}\in {\cal F}^{(n-1,m)},\\
\label{4.7}
a^*(\phi_+)(\Psi_+^{(n)}\otimes\Psi_-^{(m)})
=(a^*_{\gh_+}(\phi_+)(\Psi_+^{(n)}))\otimes\Psi_-^{(m)}\in {\cal F}^{(n+1,m)},\\
\label{4.8}
b(\phi_-)(\Psi_+^{(n)}\otimes\Psi_-^{(m)})
=\Psi_+^{(n)}\otimes(b_{\gh_-}(\phi_-)(\Psi_-^{(m)}))\in {\cal F}^{(n,m-1)},\\
\label{4.9}
b^*(\phi_-)(\Psi_+^{(n)}\otimes\Psi_-^{(m)})
=\Psi_+^{(n)}\otimes(b^*_{\gh_-}(\phi_-)(\Psi_-^{(m)}))\in {\cal F}^{(n,m+1)}.
\end{eqnarray}
All these operators have bounded extensions on ${\cal F}(\GH)$ and satisfy the CAR's.
The main object of the theory is the {\it quantized Dirac field operator} $\Psi$:
\begin{eqnarray*}
f\in \GH\mapsto \Psi(f):=a(\Pi_+f)+b^*(\Upsilon \Pi_-f)\in {\cal L}({\cal F}(\GH)),
\end{eqnarray*}
where we have denoted by $\Pi_{\pm}$ the orthogonal projector from $\GH$ to $\GH_{\pm}$.
The mapping $f\mapsto\Psi(f)$ is antilinear and bounded:
\begin{eqnarray}
\label{4.10}
||\Psi(f)||=||f||.
\end{eqnarray}
Its adjoint denoted by $\Psi^*(f)$ is given by
\begin{eqnarray}
\label{4.11}
\Psi^*(f)=a^*(\Pi_+f)+b(\Upsilon \Pi_-f)
\end{eqnarray}
and the CAR's are satisfied:
\begin{eqnarray}
\label{4.12}
\Psi(f)\Psi(g)+\Psi(g)\Psi(f)&=&0,\\
\label{4.13}
  \Psi^*(f)\Psi^*(g)+\Psi^*(g)\Psi^*(f)&=&0,\\
\label{4.14}
\Psi^*(f)\Psi(g)+\Psi(g)\Psi^*(f)&=&\langle f,g \rangle {\bf 1}.
\end{eqnarray}
The {\it Field Algebra} is the $C^*$-algebra generated by ${\bf 1}$ and the $\Psi(f),\, f\in \GH$.
If we take $f$ only in $\GH_{+(-)}$ we get a subalgebra isometric to ${\cal U}(\gh_{+(-)})$.
The {\it vacuum state} $\omega_{vac}$ is defined by
\begin{eqnarray}
\label{4.15}
A\in {\cal U}(\GH),\, \omega_{vac}(A):=\langle A \Omega_{vac},\Omega_{vac} \rangle,
\end{eqnarray}
or by the two point function:
\begin{eqnarray}
\label{4.16}
\omega_{vac}(\Psi^*(f)\Psi(g))=\langle \Pi_-f, \Pi_-g \rangle .
\end{eqnarray}
Now assume that the classical fields satisfy a Dirac type equation :
\begin{eqnarray*}
\partial_t\Psi=iH\Psi,
\end{eqnarray*}
where $H$ is selfadjoint on $\GH$ and leaves $\GH_+,\GH_-$ invariant. Then
\begin{eqnarray}
\label{4.17}
H^+:=H|_{\gh_+},\, H^-:=-\Upsilon H|_{\GH_-}\Upsilon^{-1},
\end{eqnarray}
are respectively selfadjoint on $\gh_+$ and $\gh_-$, and the classical fields of one 
particle, $\phi_+$, and of one antiparticle, $\phi_-$, are solutions to a Schr\"odinger
type equation on $\gh_{+(-)}$:
\begin{eqnarray}
\label{4.18}
\partial_t\phi_{+(-)}=iH^{+(-)}\phi_{+(-)}.
\end{eqnarray}
A usual splitting of $\GH$ is the choice
\begin{eqnarray}
\label{4.19}
\GH_+={\bf 1}_{(-\infty,0)}(H),\quad {\GH}_-={\bf 1}_{{[}0,\infty)}(H).
\end{eqnarray}
We say that a state $\omega_{\beta,\mu}$ on ${\cal U}(\GH)$ satisfies the $(\beta,\mu)$-KMS
condition, $0<\beta$ if it is characterized by the two-point function
\begin{eqnarray}
\label{4.20}
\omega_{\beta,\mu}(\Psi^*(f)\Psi(g))=\langle ze^{\beta H}(1+ze^{\beta H})^{-1}f, g \rangle,\quad                                                             
z=e^{\beta \mu}.
\end{eqnarray}
We want to apply this procedure to several states at time $t=0$ and 
in the future. We first describe the quantization at time $t=0$. Let $\GH={\cal H}_0,\, H=H_0$.
We will emphasize the importance of the charge of the field by denoting the hamiltonian $H_0=H_0(q)$. A charge conjugation 
for $H_0(q)$ is given by :
\[ \Upsilon\phi=U_{\Upsilon}\bar{\phi}\quad\mbox{with}\quad U_{\Upsilon}=\gamma^3
=\left(\begin{array}{cccc} 0 & 0 & 0 & 1\\ 0 & 0 & -1 & 0\\ 0 & 1 & 0 & 0\\ -1 & 0
& 0 & 0 \end{array} \right). \]
We note that $H^-=H_0(-q)$ and that $\psi=\Upsilon\phi$ satisfies the boundary condition 
(\ref{boundhr}) if $\phi$ satisfies it.

As we do not know whether $H_0$ has the eigenvalue $0$ or not, there is a slight ambiguity
in the definition of particles and antiparticles. We will put:
\begin{eqnarray}
\label{4.23}
(\Pi_+,\Pi_-)
=({\bf 1}_{(-\infty,0)}(H_0),{\bf 1}_{{[}0,\infty)}(H_0)),
\end{eqnarray}
but the choice $({\bf 1}_{(-\infty,0{]}}(H_0),{\bf 1}_{(0,\infty)}(H_0))$
would also have been possible. We denote $\Psi_0$ the quantum field at time $t=0$ constructed 
in the previous way. We define the Boulware quantum state $\omega_0$ on the field algebra 
${\cal U}({\cal H}_0)$ as the vacuum state
\begin{eqnarray*}
\Phi^j_0\in {\cal H}_0\quad \omega_0(\Psi_0^*(\Phi^1_0)\Psi_0(\Phi^2_0))=\langle \Pi_- \Phi^1_0,\Pi_- \Phi^2_0 \rangle.
\end{eqnarray*}
At time $t=\infty$ we will take
\begin{eqnarray}
\label{4.24}
{\GH}={\cal H},
\end{eqnarray}
$H$ the Dirac hamiltonian in the Kerr-Newman space-time and the same charge conjugation.
We put 
\begin{eqnarray}
\label{4.26}
(\Pi_+,\Pi_-)
=({\bf 1}_{(-\infty,0)}(H),{\bf 1}_{{[}0,\infty)}(H)).
\end{eqnarray}
The fields obtained in this way are denoted $\Psi(f)$. From this we obtain the definitions 
of vacuum and KMS states with respect to $H$.
\section{Quantization in a globally hyperbolic space-time}
\label{sec4.2}
Following J. Dimock \cite{Di82} we construct the local algebra of observables in the space-time outside 
the collapsing star. This construction does not depend on the choice of the representation of the CAR's, or
on the spin structure of the Dirac field, or on the choice of the hypersurface. 
In particular we can consider the Fermi-Dirac-Fock representation and the following foliation of our space-time (see Section \ref{sec2.2}):
\begin{eqnarray*}
{\cal M}_{col}=\bigcup_{t\in \R}\Sigma_t^{col},\quad \Sigma_t^{col}=\{(t,\hat{r},\theta,\varphi);\hat{r}\ge \hat{z}(t,\theta)\}.
\end{eqnarray*}
We construct the Dirac field $\Psi_0$ and the $C^*$-algebra ${\cal U}({\cal H}_0)$ as explained in Section \ref{sec4.1}.
We define the operator:
\begin{eqnarray}
\label{4.24a}
S_{col}:\Phi\in (C_0^{\infty}({\cal M}_{col}))^4\mapsto S_{col}\Phi:=\int_{\R}U(0,t)\Phi(t)dt\in {\cal H}_0,
\end{eqnarray}
where $U(0,t)$ is the propagator defined in Proposition \ref{prop3.3}. The quantum spin field is defined by~:
\begin{eqnarray*}
\Psi_{col}:\Phi\in (C_0^{\infty}({\cal M}_{col}))^4\mapsto \Psi_{col}(\Phi):=\Psi_0(S_{col}\Phi)\in {\cal L}({\cal H}_0)
\end{eqnarray*}
and for an arbitrary set ${\cal O}\subset {\cal M}_{col}$, we introduce ${\cal U}_{col}({\cal O})$, the $C^*$-algebra
generated by $\psi_{col}^*(\Phi_1)\Psi_{col}(\Phi_2),$ $\supp\Phi_j\subset {\cal O},\, j=1,2.$
Eventually, we have:
\[ {\cal U}_{col}({\cal M}_{col})=\overline{\bigcup_{{\cal O}\subset{\cal M}_{col}}{\cal U}_{col}({\cal O})}. \]
Then we define the fundamental state on ${\cal U}_{col}({\cal M}_{col})$ as follows:
\[ \omega_{col}(\Psi^*_{col}(\Phi_1)\Psi_{col}(\Phi_2)):=\omega_{vac}(\Psi^*_0(S_{col}\Phi_1)\Psi_0(S_{col}\Phi_2))
=\langle {\bf 1}_{{[}0,\infty)}(H_0)S_{col}\Phi_1,S_{col}\Phi_2\rangle . \]
Let us now consider the future black-hole. We consider the space-time ${\cal M}_{BH}$
with the Dirac hamiltonian $H$ for a field with one particle. Let $\Psi(\Phi)$ be the Dirac field as constructed 
in Section \ref{sec4.1} and 
\[ S: \Phi\in (C_0^{\infty}({\cal M}_{BH}))^4\mapsto S\Phi:=\int_{\R}e^{-itH}\Phi(t)dt.\]
We also introduce :
\[ \Psi_{BH}: \Phi\in (C_0^{\infty}({\cal M}_{BH}))^4\mapsto \Psi_{BH}(\Phi):=\Psi(S\Phi)\]
and the $C^*$-algebra ${\cal U}_{BH}({\cal O})$ generated by 
$\Psi_{BH}(\Phi_1)\Psi^*_{BH}(\Phi_2),\,\Phi_1,\Phi_2\in (C_0^{\infty}({\cal O}))^4$,\\
${\cal O}\subset{\cal M}_{BH}$.
As before we put
\[ {\cal U}_{BH}({\cal M}_{BH})=\overline{\bigcup_{{\cal O}\subset{\cal M}_{BH}}{\cal U}_{BH}({\cal O})}. \]
We also define the thermal Hawking state:
\begin{eqnarray}
\label{Hawst}
\omega^{\eta,\sigma}_{Haw}(\Psi^*_{BH}(\Phi_1)\Psi_{BH}(\Phi_2))
&=&\langle \mu e^{\sigma H}(1+\mu e^{\sigma H})^{-1}S\Phi_1,S\Phi_2\rangle _{{\cal H}}\nonumber\\
&=:&\omega_{KMS}^{\eta,\sigma}(\Psi^*(S\Phi_1)\Psi(S\Phi_2))
\end{eqnarray}
with
\[ T_{Haw}=\sigma^{-1},\, \mu=e^{\sigma\eta},\, \sigma>0,\]
where $T_{Haw}$ is the Hawking temperature and $\mu$ is the chemical potential. 
We will also need a vacuum state which is given by :
\[ \omega_{vac}(\Psi^*_{BH}(\Phi_1)\Psi_{BH}(\phi_2))=\langle {\bf 1}_{{[}0,\infty)}(H)S\phi_1,S\phi_2\rangle. \]
\section{The Hawking effect}
\label{sec4.3}
In this section we formulate the main result of this paper.
Let $\Phi\in (C_0^{\infty}({\cal M}_{col}))^4$. We put
\begin{eqnarray}
\label{5.1}
\Phi^T(t,\hat{r},\omega)=\Phi(t-T,\hat{r},\omega).
\end{eqnarray}
\begin{theorem}[Hawking effect]
\label{th5.1}
Let 
\[ \Phi_j\in (C_0^{\infty}({\cal M}_{col}))^4,\, j=1,2.\]
Then we have
\begin{eqnarray}
\label{5.0}
\lefteqn{\lim_{T\rightarrow\infty}\omega_{col}(\Psi^*_{col}(\Phi^T_1)\Psi_{col}(\Phi^T_2))}\nonumber\\
&=&\omega_{Haw}^{\eta,\sigma}(\Psi^*_{BH}({\bf 1}_{\R^+}(P^-)\Phi_1)\Psi_{BH}({\bf 1}_{\R^+}(P^-)\Phi_2))\nonumber\\
&+&\omega_{vac}(\Psi^*_{BH}({\bf 1}_{\R^-}(P^-)\Phi_1)\Psi_{BH}({\bf 1}_{\R^-}(P^-)\Phi_2)),\\
\label{mueta}
T_{Haw}&=&1/\sigma=\kappa_+/2\pi,\quad \mu=e^{\sigma\eta},\, \eta=\frac{qQr_+}{r_+^2+a^2}+\frac{aD_{\varphi}}{r_+^2+a^2}.
\end{eqnarray}
\end{theorem}
In the above theorem $P^{\pm}$ is the asymptotic velocity introduced in Chapter \ref{sec3}. The projections ${\bf 1}_{\R^{\pm}}(P^{\pm})$ separate outgoing and incoming solutions.
\begin{remark}
The result is independent of the choices of coordinate system and tetrad, 
i.e. both sides of (\ref{5.0}) are independent
of these choices. Indeed a change of coordinate system or a change of tetrad
is equivalent to a conjugation of the operators by a unitary transformation. We also note that the result
is independent of the chiral angle $\nu$ in the boundary condition. Let 
\begin{eqnarray*}
S_{\leftarrow}\phi=\int_{\R}e^{-itH_{\leftarrow}}\phi(t)dt. 
\end{eqnarray*}
It is easy to check :
\begin{eqnarray*}
\lefteqn{\omega_{Haw}^{\eta,\sigma}(\Psi_{BH}^*({\bf 1}_{\R^+}(P^-)\phi_1)\Psi_{BH}({\bf 1}_{\R^+}(P^-)\phi_2))}\\
&=&\langle \mu e^{\sigma H_{\leftarrow}}(1+e^{\sigma H_{\leftarrow}})^{-1}S_{\leftarrow}\Omega_{\leftarrow}^-\phi_1, S_{\leftarrow}\Omega_{\leftarrow}^-\phi_2\rangle.
\end{eqnarray*}
We define $S_{\rightarrow}$ in an analogous way to $S_{\leftarrow}$ using the Dollard
modified dynamics :
\[ S_{\rightarrow}\phi=\int_{\R}U_D(-t)\phi(t)dt. \]
We find :
\[ \omega_{vac}(\Psi^*_{BH}({\bf 1}_{\R^-}(P^-)\phi_1)\Psi_{BH}({\bf 1}_{\R^-}(P^-)\phi_2)
=\langle {\bf 1}_{{[}0,\infty)}(H_{\rightarrow})S_{\rightarrow}\Omega_{\rightarrow}^-\phi_1,S_{\rightarrow}\Omega_{\rightarrow}^-\phi_2\rangle. \]
In particular our result coincides with the result of Melnyk (see \cite{Me2}) in the
case of a Reissner-Nordstr\"om black-hole.
\end{remark}
\chapter{Additional scattering results}
\label{sec7}
In this chapter we state some scattering results that we shall need in what follows. 
\section{Spin weighted spherical harmonics}
\label{sec6.4}
We will now introduce spin weighted spherical harmonics $Y^l_{sn}$ 
(for a complete definition, see e.g.
\cite{Ni2}). For each spinorial weight
$s$, $2s\in \Z$, the family $\{ Y^l_{sn}(\varphi,\theta) = e^{in\varphi}
  u^l_{sn}(\theta)\, ;~l-|s|\in \N \, ,~l-|n|\in \N\, \}$ forms a
Hilbert basis of $L^2(S^2_{\omega},d\omega)$ and we have the following
relations

\begin{eqnarray*}
\frac{du^{l}_{sn}}{d\theta} - \frac{n-s\cos\theta}{\sin\theta}
u^{l}_{sn} = -i \left[ (l+s)(l-s+1)\right]^{1/2} u^{l}_{s-1,n} \, , \\
\frac{du^{l}_{sn}}{d\theta} + \frac{n-s\cos\theta}{\sin\theta}
u^{l}_{sn} = -i \left[ (l+s+1)(l-s)\right]^{1/2} u^{l}_{s+1,n}\, .
\end{eqnarray*}
We define $\otimes_4$ as the following operation between two vectors
of $\C^4$
\begin{eqnarray*}
\forall v=(v_1,v_2,v_3,v_4),\, u=(u_1,u_2,u_3,u_4), \quad
v\otimes_4u=(u_1v_1,u_2v_2,u_3v_3,u_4v_4).
\end{eqnarray*}
Since the families
\begin{eqnarray*}
\{Y^l_{\frac{1}{2},n};(n,l)\in {\cal
  I}\},\quad\{Y^l_{-\frac{1}{2},n};(n,l)\in{\cal I}\},\quad {\cal
  I}=\{(n,l)/l-\frac{1}{2}\in \N, l-|n|\in \N\}
\end{eqnarray*}
form a Hilbert basis of $L^2(S^2_{\omega},d\omega)$, we express ${\cal
  H}_*$ as a direct sum
\begin{eqnarray*}
{\cal H}_*=\oplus_{(n,l)\in{\cal I}}{\cal H}^{nl}_*,\quad {\cal
  H}^{nl}_*=L^2\left( \left( \R;dr_*\right) ; \C^4 \right) \otimes_4
  Y_{nl},\\
Y_{nl}=(Y^l_{-\frac{1}{2},n},Y^l_{\frac{1}{2},n},Y^l_{-\frac{1}{2},n},Y^l_{\frac{1}{2},n}).
\end{eqnarray*}
We shall henceforth identify ${\cal H}^{nl}_*$ and $L^2 \left( \left( \R
    ;dr \right) ; \C^4 \right)$ as well as $\psi_{nl}\otimes_4 Y_{nl}$
and $\psi_{nl}$. We see that  
\begin{equation*}
{\notD}_{\gs} = \oplus_{nl}{\notD}_{\gs}^{nl} \quad\mbox{with} \quad 
{\notD}_{\gs}^{nl}:=\Gamma^1D_{r_*}+a_0(r_*)\Gamma^2(l+1/2)+b_0(r_*)\Gamma^{\nu}+c^n.
\end{equation*}
In a similar way we find the decompositions:
\[ {\notD}_H=\oplus_{(n,l)\in{\cal I}}{\notD}_H^{nl},\quad 
{\notD}_{\rightarrow}=\oplus_{(n,l)\in {\cal I}}{\notD}_{\rightarrow}^{nl}.
\]
Note that the operator ${\notD}^{nl}_{\gs}$ is selfadjoint on ${\cal H}^{nl}$ with domain
$(H^1(\R))^4\otimes_4Y_{nl}$. This leads to a useful characterization of the domain $D({\notD})$ 
(see\cite{Da2}) :
\[ {\cal H}_*^1=D({\notD})=D({\notD}_{\gs})=\{u=\sum_{nl}u_{nl};\,\forall n,l\, u_{nl}\in H^1(\R),\,\sum||{\notD}^{nl}_{\gs}u_{nl}||^2<\infty\}.\]
\begin{remark}
\label{smooth}
Note that $n\in \Z+1/2$. Indeed we are working with quantities of spin weight $1/2$.
Such quantities are multiplied by $e^{i\varphi/2}$ under rotation of angle $\varphi$. 
In particular $Y^l_{\frac{1}{2},n}$ are not smooth on $S^2_{\omega}$, but they are smooth 
on ${[}0,2\pi{]}_{\varphi}\times {[}0,\pi{]}_{\theta}$. Using the axial symmetry of our equations we will
often fix the angular momentum $D_{\varphi}=n$. When we do so we will always suppose 
$n\in \Z+1/2$.
\end{remark}
\section{Velocity estimates}
We start with the maximal velocity estimate :
\begin{lemma}
\label{lem7.1}
Let $J\in C^{\infty}_b(\R),\, \supp J\subset{]}-\infty,-1-\epsilon{]}\cup{[}1+\epsilon,\infty{[}$ for some 
$\epsilon>0$. Then we have :
\begin{eqnarray*}
(i) & \int_1^{\infty}||J\left(\frac{\hat{r}}{t}\right)e^{itH}\phi||^2\frac{dt}{t}\lesssim ||\phi||^2,\\
(ii) & s-\lim_{t\rightarrow\pm\infty}J\left(\frac{\hat{r}}{t}\right)e^{itH}=0.
\end{eqnarray*}
\end{lemma}
The lemma can be easily deduced from the equivalent statement for the dynamics $\notD$ in \cite[Proposition 4.4]{Da2}.
The minimal velocity estimate is given by the following lemma :
\begin{lemma}
\label{lem7.1a}
Let $\chi\in C_0^{\infty}(\R)$ such that $\supp\chi\subset \R\setminus \{-m,m\}$.
Then there exists a strictly positive constant $\epsilon_{\chi}$ such that we have :
\[ \int_1^{\infty}||{\bf 1}_{{[}0,\epsilon_{\chi}{]}}\left(\frac{|\hat{r}|}{t}\right)e^{itH}\chi(H)\Phi||^2\frac{dt}{t}\le C_{\chi}||\Phi||^2. \]
Furthermore 
\[ s-\lim_{t\rightarrow \infty}{\bf 1}_{{[}0,\epsilon_{\chi}{]}}\left(\frac{|\hat{r}|}{t}\right)e^{itH}\chi(H)=0. \]
\end{lemma}
The lemma follows from \cite[Proposition 4.3]{Da2}. The change of variables and tetrads is
treated in the usual way. It turns out that we need a stronger
minimal velocity estimate near the horizon~:
\begin{lemma}
\label{lem7.2}
Let $m\in C(\R^+),\,0<\epsilon<1,\, \lim_{t\rightarrow\infty}m(t)=\infty,\, 0\le m(|t|)\le\epsilon |t|$ for all $|t|\ge 1$. Then
\begin{eqnarray}
\label{7.28}
s-\lim_{t\rightarrow\pm\infty}{\bf 1}_{{[}0,1{]}}\left(\frac{|\hat{r}|}{|t|-m(|t|)}\right)e^{itH}{\bf 1}_{\R^\mp}(P^{\pm})=0.
\end{eqnarray}
An analogous result holds if we replace $e^{itH}{\bf 1}_{\R^\mp}(P^{\pm})$ by 
$e^{itH_{\leftarrow}}$.
\end{lemma}

{\bf Proof.}

By the asymptotic completeness result and a density argument it is sufficient to 
show:
\begin{eqnarray}
\label{7.29}
\forall f\in (C_0^{\infty}(\R\times S^2))^4\quad s-\lim_{t\rightarrow\pm \infty}{\bf 1}_{{[}0,1{]}}
\left(\frac{|\hat{r}|}{|t|-m(|t|)}\right)e^{itH_{\leftarrow}}f=0.
\end{eqnarray}
Let $\supp f\subset {[}R_1,R_2{]}\times S^2$. We define $\tilde{H}_{\leftarrow}=\Gamma^1D_{\hat{r}}.$
Using ${[}\tilde{H}_{\leftarrow},\frac{a}{r_+^2+a^2}D_{\varphi}+\frac{qQr_+}{r_+^2+a^2}{]}=0$
we see that it is sufficient to show :
\[ s-\lim_{t\rightarrow\pm \infty}{\bf 1}_{{[}0,1{]}}
\left(\frac{|\hat{r}|}{|t|-m(|t|)}\right)e^{it\tilde{H}_{\leftarrow}}f=0. \]
We only treat the case $t\rightarrow\infty$,
the case $t\rightarrow-\infty$ being analogous. We have 
\begin{eqnarray*}
(e^{it\tilde{H}_{\leftarrow}}f)(\hat{r},\omega)=\left(\begin{array}{c} f_1(\hat{r}+t,\omega)\\
f_2(\hat{r}-t,\omega)\\
f_3(\hat{r}-t,\omega)\\
f_4(\hat{r}+t,\omega) \end{array} \right).
\end{eqnarray*}
On $\supp \left({\bf 1}_{{[}0,1{]}}
\left(\frac{|\hat{r}|}{|t|- m(t)}\right)e^{itH_{\leftarrow}}f\right)_{2,3}$ we have
for $t\ge -R_1$ :
\[ t-m(t)\ge|\hat{r}|\ge t+R_1,\quad\mbox{which is impossible for $t$ sufficiently large}. \]
In the same way we find on $\supp \left({\bf 1}_{{[}0,1{]}}
\left(\frac{|\hat{r}|}{|t|- m(t)}\right)e^{itH_{\leftarrow}}f\right)_{1,4}$ 
for $t\ge R_2$ :
\[ t-m(t)\ge|\hat{r}|\ge t-R_2,\quad\mbox{which is impossible for $t$ sufficiently large}. \]
Thus for $t$ sufficiently large we have:
\begin{eqnarray*}
{\bf 1}_{{[}0,1{]}}
\left(\frac{|\hat{r}|}{|t|-m(t)}\right)e^{itH_{\leftarrow}}f=0.
\end{eqnarray*}
\qed
\section{Wave operators}
\label{sec7.2}
We shall need a characterization of the wave operators in terms of cut-off functions.
Let $J\in C^{\infty}(\R)$ s.t.
\begin{eqnarray}
\label{7.1}
\exists\, b,c\in \R,\, 0<b<c\quad J(\hat{r})=\left\{\begin{array}{cc} 1 & \hat{r}<b\\0 & \hat{r}>c \end{array}\right. 
\end{eqnarray}
\begin{proposition}
\label{prop7.3}
\begin{eqnarray*}
W_{\rightarrow}^{\pm}&=&s-\lim_{t\rightarrow\pm\infty}e^{-itH}(1-J)U_D(t),\\
\Omega_{\rightarrow}^{\pm}&=&s-\lim_{t\rightarrow\pm\infty}U_D(-t)(1-J)e^{itH},\\
W_{\leftarrow}^{\pm}&=&s-\lim_{t\rightarrow\pm\infty}e^{-itH}Je^{itH_{\leftarrow}},\\
\Omega_{\leftarrow}^{\pm}&=&s-\lim_{t\rightarrow\pm\infty}e^{-itH_{\leftarrow}}Je^{itH}.
\end{eqnarray*}
\end{proposition}
We refer to \cite{HN} for the link between the wave operators in terms of cut-off functions and the wave
operators using the asymptotic velocity. The operator ${\notD}$ is a short range perturbation of ${\notD}_{\gs}$. 
More precisely we have the following (\cite[Theorem 5.1]{Da2}):
\begin{proposition}
\label{prop7.2}
The wave operators 
\begin{eqnarray*}
\tilde{W}_{\gs}^{\pm}&=&s-\lim_{t\rightarrow\pm\infty}e^{-it{\expnotD}}e^{it{\expnotD}_{\gs}},\\
\tilde{\Omega}_{\gs}^{\pm}&=&s-\lim_{t\rightarrow\pm\infty}e^{-it{\expnotD}_{\gs}}e^{it{\expnotD}}
\end{eqnarray*}
exist and we have:
\[ \tilde{\Omega}_{\gs}^{\pm}=(\tilde{W}_{\gs}^{\pm})^*,\, \tilde{W}^{\pm}_{\gs}\tilde{\Omega}^{\pm}_{\gs}=\tilde{\Omega}_{\gs}^{\pm}\tilde{W}_{\gs}^{\pm}={\bf 1}_{{\cal H}_*}.\]
\end{proposition}
We define ${\cal S}^{\rho}(\R)$ as a subspace of $C^{\infty}(\R)$ by :
\[ f\in {\cal S}^{\rho}(\R)\Leftrightarrow \forall \alpha \in \N\, |f^{(\alpha)}(x)|\le C_{\alpha} \langle x\rangle^{\rho-\alpha}, \, 
\langle x \rangle :=\sqrt{1+x^2}. \]
We shall need (see \cite[Corollary 3.2]{Da2}):
\begin{lemma}
\label{lem7.3}
$(i)$ Let $\chi\in {\cal S}^0(\R)$. Then
\[ (\chi({\notD})-\chi({\notD}_{\gs}))({\notD}+i)^{-1}\quad\mbox{is compact.}\]
$(ii)$ If $\chi\in C_0^{\infty}(\R)$, then 
\[ (\notD-\notD_{\gs})\chi(\notD),\, (\notD-\notD_{\gs})\chi(\notD_{\gs})\quad\mbox{are compact.}\]
\end{lemma}
\begin{lemma}
\label{lem7.4}
The following wave operator exists:
\begin{eqnarray*}
W_0^{\pm}:=s-\lim_{t\rightarrow\pm\infty}e^{-it{H}_0}(1-J)e^{it{H}}.
\end{eqnarray*}
Furthermore we have 
\begin{eqnarray}
\label{a}
W_0^{\pm}{\bf 1}_{{[}0,\infty)}(H)={\bf 1}_{{[}0,\infty)}(H_0)W_0^{\pm},\\
\label{b}
\forall f \in {\cal H}\, ||W_0^{\pm}f||=||{\bf 1}_{\R^\pm}(P^{\pm})f||.
\end{eqnarray}
\end{lemma}
{\bf Proof.}

We only treat the case $t\rightarrow \infty$, the case $t\rightarrow -\infty$ being analogous.
By a density argument it is sufficient to show for $\chi\in C_0^{\infty}(\R),\, \supp
\chi\subset \R\setminus \{-m,m\}$ the existence of 
\[ s-\lim_{t\rightarrow \infty}e^{-itH_0}(1-J)e^{itH}\chi(H). \]
Let $\epsilon_{\chi}$ as in Lemma \ref{lem7.1a}. 
Let $J_0\in C_0^{\infty}(\R),\,J_0\ge 0,\, \supp J_0\subset(\frac{\epsilon_{\chi}}{4},1+2\epsilon_{\chi}),
\,J_0=1$ on ${[}\frac{\epsilon_{\chi}}{2},1+\epsilon_{\chi}{]}$.
We first show that
\begin{eqnarray}
\label{C12.1}
s-\lim_{t\rightarrow\infty}e^{-it{H}_0}(1-J)e^{it{H}}\chi(H)=s-\lim_{t\rightarrow\infty}e^{-it{H}_0}J_0\left(\frac{\hat{r}}{t}\right)e^{it{H}}\chi(H)
\end{eqnarray}
if the R.H.S. exists.
For this purpose let $J_0,J_{\pm}$ be a partition of unity with $J_{\pm}\ge 0$, 
$\supp J_-\subset (-\infty,\frac{\epsilon_{\chi}}{2}),\,\supp J_+\subset(1+\epsilon_{\chi},\infty),\, J_-+J_0+J_+=1.$
We have
\[ s-\lim_{t\rightarrow-\infty}e^{-it{H}_0}(1-J)J_+\left(\frac{\hat{r}}{t}\right)e^{it{H}}\chi(H)=0
\]
by Lemma \ref{lem7.1}. We write $J_-=J_-^1+J_-^2$, where 
$J_-^1,J_-^2\ge 0$ and $J_-^2$ equals one in a small neighborhood of $-1$ and is supported in a slightly larger
neighborhood $I$ of $-1$. We have
\[ s-\lim_{t\rightarrow\infty}e^{-it{H}_0}(1-J)J_-^1\left(\frac{\hat{r}}{t}\right)e^{it{H}}\chi(H)=0\]
by the maximal and minimal velocity estimates. Obviously we have for $I$ sufficiently small
and $t$ sufficiently large :
\[e^{-it{H}_0}(1-J)J^2_-\left(\frac{\hat{r}}{t}\right)e^{it{H}}\chi(H)=0.\]
Let
\[ J_0\left(\frac{\hat{r}}{t}\right)={\bf 1}_{{[}\frac{\epsilon_{\chi}}{2},1+\epsilon_{\chi}{]}}\left(\frac{\hat{r}}{t}\right)+J_0^1\left(\frac{\hat{r}}{t}\right).\]
Clearly
\[ (1-J)(\hat{r}){\bf 1}_{{[}\frac{\epsilon_{\chi}}{2},1+\epsilon_{\chi}{]}}\left(\frac{\hat{r}}{t}\right)
={\bf 1}_{{[}\frac{\epsilon_{\chi}}{2},1+\epsilon_{\chi}{]}}\left(\frac{\hat{r}}{t}\right)\]
for $t$ sufficiently large and 
\[s-\lim_{t\rightarrow\infty}e^{-it{H}_0}(1-J)J_0^1\left(\frac{\hat{r}}{t}\right)e^{it{H}}\chi(H)=0 \]
by the maximal and minimal velocity estimates. Again by the maximal and minimal velocity estimates
we obtain (\ref{C12.1}) :
\[ s-\lim_{t\rightarrow\infty}e^{-it{H}_0}{\bf 1}_{{[}\epsilon_{\chi},1+\epsilon_{\chi}{]}}\left(\frac{\hat{r}}{t}\right)e^{it{H}}\chi(H)
=s-\lim_{t\rightarrow\infty}e^{-it{H}_0}J_0\left(\frac{\hat{r}}{t}\right)e^{it{H}}\chi(H).\]
It remains to show that the limit on the R.H.S. of (\ref{C12.1}) exists :
\begin{eqnarray*}
\frac{d}{dt}e^{-it{H}_0}J_0\left(\frac{\hat{r}}{t}\right)e^{it{H}}\chi(H)&=&e^{-it{H}_0}({H}_0J_0\left(\frac{\hat{r}}{t}\right)-J_0\left(\frac{\hat{r}}{t}\right){H})e^{it{H}}\chi(H)\\
&=&e^{-it{H}_0}J'_0\left(\frac{\hat{r}}{t}\right)\frac{1}{t}\Gamma^1e^{it{H}}\chi(H)
\end{eqnarray*}
and this last expression is integrable in $t$ by the maximal and minimal velocity estimates. 
Let us now prove (\ref{a}). Let $\chi_n\in C_0^{\infty}(\R)$ with 
\begin{eqnarray*}
\chi_n(x)=\left\{ \begin{array}{cc} 1 & 0\le x\le n, \\ 0 & x\le -\frac{1}{n},\, x\ge n+1 \end{array} \right.
\end{eqnarray*}
We have 
\begin{eqnarray*}
\chi_n(H)=\int e^{-itH}\hat{\chi}_n(t)dt,
\end{eqnarray*}
where $\hat{\chi}_n$ stands for the Fourier transform of $\chi_n$. Using $W_0^{\pm}e^{-itH}=e^{-itH_0}W_0^{\pm}$
we see that 
\begin{eqnarray}
\label{a*}
W_0^{\pm}\chi_n(H)=\chi_n(H_0)W_0^{\pm}.
\end{eqnarray}
Taking the strong limit in (\ref{a*}) gives (\ref{a}). Let us now prove (\ref{b}). 
By a density argument we see that it is sufficient to establish (\ref{b}) for $\chi(H)f,\, \chi\in C_0^{\infty}(\R),\,
\supp\chi\subset \R\setminus \{-m,m\}$.
But by (\ref{C12.1}) we see that 
\begin{eqnarray*}
||W_0^+\chi(H)f||&=&\lim_{t\rightarrow \infty}||e^{-itH_0}J_0\left(\frac{\hat{r}}{t}\right)e^{itH}\chi(H)f||\\
&=&||J_0(P^+)\chi(H)f||=||{\bf 1}_{\R^+}(P^+)\chi(H)f||,
\end{eqnarray*}
where we have used the minimal and maximal velocity estimates.
\qed
\section{Regularity results}    
We first need a result on $\supp \Omega_{\leftarrow}^-f$:
\begin{lemma}
\label{lem7.7}
We have
\begin{eqnarray}
\label{7.16}
{\bf 1}_{(-\infty,R)}(\hat{r})\Omega_{\leftarrow}^-{\bf 1}_{{[}R,\infty)}(\hat{r})=0\quad \forall R.
\end{eqnarray}
\end{lemma}

{\bf Proof.}

We have :
\begin{eqnarray*}
\Omega_{\leftarrow}^-&=&s-\lim_{t\rightarrow-\infty}e^{-itH_{\leftarrow}}Je^{itH},\\
P_{{\cal H}^+}\Omega_{\leftarrow}^{-}&=&\Omega_{\leftarrow}^-.
\end{eqnarray*}
It follows
\begin{eqnarray}
\label{7.19}
\Omega_{\leftarrow}^-=s-\lim_{t\rightarrow-\infty}e^{-itH_{\leftarrow}}P_{{\cal H}^+}Je^{itH}.
\end{eqnarray}
Let $f\in{\cal H},\, \supp f\subset{[}R,\infty)\times S^2$. By the finite propagation speed 
we have for $t\le 0$ :
\begin{eqnarray*}
\supp e^{itH}f\subset{[}R+t,\infty)\times S^2.
\end{eqnarray*}
Therefore :
\begin{eqnarray}
\label{7.21}
\supp e^{-itH_{\leftarrow}}P_{{\cal H}^+}J e^{itH}f\subset {[}R,\infty)\times S^2.
\end{eqnarray}
The equations (\ref{7.19}), (\ref{7.21}) prove the lemma.
\qed

For $\Omega_{\leftarrow}^-$ we need the convergence in $H^1(\R;(L^2(S^2))^4)$:
\begin{lemma}
\label{lem7.8}
We have for all $f\in {\cal H}^1$:
\begin{eqnarray}
\label{7.22}
\lim_{t\rightarrow-\infty}||J(\hat{r})e^{itH}f-e^{itH_{\leftarrow}}\Omega_{\leftarrow}^-f||_{(H^1(\R;(L^2(S^2))^4)}=0.
\end{eqnarray}
\end{lemma}
{\bf Proof.}

The wave operators acting on 
${\cal H}_*$ will be denoted with a tilde. By conjugation with ${\cal V}{\cal U}$
we obtain wave operators acting on ${\cal H}$, e.g.
\begin{eqnarray*}
\Omega_{\gs}^-={\cal V}{\cal U}\tilde{\Omega}_{\gs}^-{\cal U}^*{\cal V}^*
=s-\lim_{t\rightarrow -\infty}e^{-itH_{\gs}}e^{itH}
\end{eqnarray*}
with $H_{\gs}={\cal V}{\cal U}\notD_{\gs}{\cal U}^*{\cal V}^*$. We note that
\begin{eqnarray*}
{\cal V}{\cal U}: {\cal H}_*^1\rightarrow{\cal H}^1;\quad {\cal V}{\cal U} : 
H^1(\R_{r_*};(L^2(S^2))^4)\rightarrow H^1(\R_{\hat{r}};(L^2(S^2))^4),\\
{\cal U}^*{\cal V}^*: {\cal H}^1\rightarrow {\cal H}^1_*;\quad {\cal U}^*{\cal V}^* : 
H^1(\R_{\hat{r}};(L^2(S^2))^4)\rightarrow H^1(\R_{r_*};(L^2(S^2))^4)
\end{eqnarray*}
are continuous. This follows from the definition of the spaces ${\cal H}^1,\, {\cal H}^1_*$ 
and from the estimate $0<\delta\le k'\le 1$.
We show:
\begin{eqnarray}
\label{7.23}
\forall f\in{\cal H}_*^1\, \lim_{t\rightarrow -\infty}||e^{it{\expnotD}}f-e^{it{\expnotD}_{\gs}}\tilde{\Omega}^-_{\gs}f||_{{\cal H}^1_*}=0,\\
\label{7.24}
\forall f\in{\cal H}_*^1\, \lim_{t\rightarrow -\infty}||J(\hat{r}(r_*,\theta))e^{it{\expnotD}_{\gs}}f-e^{it{\expnotD}_{H}}\tilde{\Omega}^-_{H,\gs}f||_{H^1(\R;(L^2(S^2))^4)}=0,
\end{eqnarray}
\begin{eqnarray}
\label{7.25}
\lefteqn{\forall f\in H^1(\R_{\hat{r}};(L^2(S^2))^4)}\nonumber\\
&&\lim_{t\rightarrow-\infty}||e^{itH_H}{\bf 1}_{\R^+}(P_H^-)f-e^{itH_{\leftarrow}}\Omega_{H,\leftarrow}^-f||_{(H^1(\R;(L^2(S^2))^4)}=0, 
\end{eqnarray}
where
\begin{eqnarray*}
\tilde{\Omega}^-_{H,\gs}&=&s-\lim_{t\rightarrow-\infty}e^{-it\expnotD_H}e^{it\expnotD_{\gs}}{\bf 1}_{\R^+}(\tilde{P}_{\gs}^-),\\
\tilde{P}^-_{\gs}&=&s-C_{\infty}-\lim_{t\rightarrow-\infty}e^{-it\expnotD_{\gs}}\frac{r_*}{t}e^{it\expnotD_{\gs}},\\
H_H&=&{\cal V}{\cal U}\notD_H{\cal U}^*{\cal V}^*,\\
\Omega_{H,\leftarrow}^-&=&s-\lim_{t\rightarrow-\infty}e^{-itH_{\leftarrow}}e^{itH_H}{\bf 1}_{\R^+}(P_H^-),\\
P_H^-&=&s-C_{\infty}-\lim_{t\rightarrow-\infty}e^{-itH_H}\frac{\hat{r}}{t}e^{itH_H}.
\end{eqnarray*}
We refer to \cite{Da1} for the existence of $\tilde{\Omega}^-_{H,\gs},\,\tilde{P}^-_{\gs}$
\footnote{Note that $\notD_{\gs}$ can be understood as the Dirac operator in a "Reissner-Nordstr\"om
type" space-time, see \cite{Da2} for details.}.
The existence of $P_H^-$ follows from the existence of 
\[ \tilde{P}_H^-=s-C_{\infty}-\lim_{t\rightarrow -\infty}e^{-it\expnotD_H}\frac{r_*}{t}
e^{it\expnotD_H} \]
and the existence of $\Omega_{H,\leftarrow}^-$ follows from the existence of 
$W_c^{\pm}$ (see proof of Theorem \ref{th3.2}).
Let us first argue that (\ref{7.23})-(\ref{7.25}) imply (\ref{7.22}).
Let 
\[ \Omega_H^-=s-\lim_{t\rightarrow-\infty}e^{-itH_H}e^{itH}{\bf 1}_{\R^+}(P^-). \]
The existence of $\Omega_H^-$ follows from the existence of $\tilde{\Omega}_H^-$
(see proof of Theorem \ref{th3.2}). We have :
\begin{eqnarray*}
\lefteqn{||(J(\hat{r})e^{itH}-e^{itH_{\leftarrow}}\Omega_{\leftarrow}^-)f||_{H^1(\R;(L^2(S^2))^4)}}\\
&\le&||J(\hat{r})(e^{itH}-e^{itH_{\gs}}\Omega_{\gs}^-)f||_{H^1(\R;(L^2(S^2))^4)}\\
&+&||(J(\hat{r})e^{itH_{\gs}}\Omega_{\gs}^--e^{itH_H}\Omega_H^-)f||_{H^1(\R;(L^2(S^2))^4)}\\
&+&||(e^{itH_H}\Omega_H^--e^{itH_{\leftarrow}}\Omega_{\leftarrow}^-)f||_{H^1(\R;(L^2(S^2))^4)}
=:I_1+I_2+I_3.
\end{eqnarray*}
We first estimate $I_1$. We have, thanks to (\ref{6.19}) and (\ref{7.23}) :
\begin{eqnarray*}
I_1&\le& ||(e^{it\expnotD}-e^{it\expnotD_{\gs}}\tilde{\Omega}^-_{\gs}){\cal V}^*{\cal U}^*f||_{H^1(\R;(L^2(S^2))^4)}\\
&\lesssim& ||(e^{it\expnotD}-e^{it\expnotD_{\gs}}\tilde{\Omega}^-_{\gs}){\cal V}^*{\cal U}^*f||_{{\cal H}_*^1}\rightarrow 0,\, t\rightarrow -\infty. 
\end{eqnarray*}
In order to estimate $I_2$ we observe that $\tilde{\Omega}_{\gs}^-: {\cal H}_*^1\rightarrow {\cal H}_*^1$ and that
$\tilde{\Omega}_H^-(\tilde{\Omega}^-_{\gs})^*=\tilde{\Omega}_{H,\gs}^-$. We obtain using (\ref{7.24}) :
\begin{eqnarray*}
I_2&\lesssim& ||(J(\hat{r}(r_*,\theta))e^{it\expnotD_{\gs}}\tilde{\Omega}_{\gs}^--e^{it\expnotD_H}\tilde{\Omega}_H^-){\cal V}^*{\cal U}^*f||_{H^1(\R;(L^2(S^2))^4)}\\
&=&||(J(\hat{r}(r_*,\theta))e^{it\expnotD_{\gs}}-e^{it\expnotD_H}\tilde{\Omega}_{H,\gs}^-)\tilde{\Omega}_{\gs}^-{\cal V}^*{\cal U}^*f||_{H^1(\R;(L^2(S^2))^4)}
\rightarrow 0,\, t\rightarrow -\infty.
\end{eqnarray*}
We now estimate $I_3$. We observe that $D(H_H)=H^1(\R;(L^2(S^2))^4)$ and that the graph
norm of $H_H$ is equivalent to the norm of $H^1(\R;(L^2(S^2))^4)$. This entails $\Omega_H^-f\in 
H^1(\R;(L^2(S^2))^4)$. Observe also that 
\[ \Omega_{\leftarrow}^-=\Omega_{H,\leftarrow}^-\Omega_H^-. \]
We obtain using (\ref{7.25}) :
\[ I_3=||(e^{itH_H}{\bf 1}_{\R^+}(P_H^-)-e^{itH_{\leftarrow}}\Omega_{H,\leftarrow}^-)\Omega_H^-f||_{H^1(\R;(L^2(S^2))^4)}
\rightarrow0,\, t\rightarrow-\infty. \]
The proof of (\ref{7.24}) is analogous to the proof of \cite[Lemma 6.3]{Me2} and we therefore omit it. 

Let us show (\ref{7.23}). Using the uniform
estimates 
\[||e^{it{\expnotD}}f||_{{\cal H}_*^1}\lesssim ||f||_{{\cal H}_*^1}\,\mbox{and}\, ||e^{it{\expnotD_{\gs}}}f||_{{\cal H}_*^1}\lesssim ||f||_{{\cal H}_*^1}\]
we argue that we can replace
$f$ by $\chi(\notD)f,\,\chi\in C_0^{\infty}(\R)$. We then write:
\begin{eqnarray*}
\notD(e^{it{\expnotD}}-e^{it{\expnotD}_{\gs}}\tilde{\Omega}^-_{\gs})\chi(\notD)f=(e^{it{\expnotD}}-e^{it{\expnotD}_{\gs}}\tilde{\Omega}^-_{\gs}){\notD}\chi(\notD)f
+({\notD}_{\gs}-\notD)\chi({\notD}_{\gs})e^{it\expnotD_{\gs}}\tilde{\Omega}^-_{\gs}f.
\end{eqnarray*}
When $t\rightarrow-\infty$ the first term goes to zero by Proposition \ref{prop7.2}
and the second term goes to zero because $({\notD}_{\gs}-\notD)\chi({\notD}_{\gs})$ is compact
by Lemma \ref{lem7.3} $(ii)$.

Let us show (\ref{7.25}). By a density argument we can replace $f$ by 
$\chi(H_H)f$ with $\chi\in C_0^{\infty}(\R)$. Then we have :
\begin{eqnarray}
\label{7.27}
\lefteqn{H_{\leftarrow}(e^{itH_H}{\bf 1}_{\R^+}(P_H^-)-e^{itH_{\leftarrow}}\Omega^-_{H,\leftarrow})\chi(H_H)f}\nonumber\\
&=&(H_{\leftarrow}-H_H)\chi(H_H)e^{itH_H}{\bf 1}_{\R^+}(P_H^-)f\nonumber\\
&+&(e^{itH_H}{\bf 1}_{\R^+}(P_H^-)-e^{itH_{\leftarrow}}\Omega^-_{H,\leftarrow})
H_H\chi(H_H)f.
\end{eqnarray}
The second term in (\ref{7.27}) goes to zero by definition of $\Omega^-_{H,\leftarrow}$.
In order to show
that the first term in (\ref{7.27}) goes to zero it is sufficient to show
\[ (\notD_{\leftarrow}-\notD_H)\chi(\notD_H)e^{it\expnotD_H}{\bf 1}_{\R^+}(\tilde{P}_H^-)f\rightarrow 0,\, t\rightarrow -\infty. \]
We have :
\begin{eqnarray}
\label{N6.11}
||(\notD_{\leftarrow}-\notD_H)\chi(\notD_H)e^{it\expnotD_H}{\bf 1}_{\R^+}(\tilde{P}_H^-)f||_{{\cal H}_*}\le 
||m(r_*)\notD_H\chi(\notD_H)e^{it\expnotD_H}{\bf 1}_{\R^+}(\tilde{P}_H^-)f||_{{\cal H}_*}
\end{eqnarray}
with $m(r_*)\rightarrow 0,\, |r_*|\rightarrow \infty$. We then use the spherical symmetry of the expression
on the R.H.S of (\ref{N6.11}) and the fact that $m(r_*)\notD_H^{nl}\chi(\notD^{nl}_H)$ is compact on 
${\cal H}_*^{nl}$ to conclude that the R.H.S. of (\ref{N6.11}) goes to zero when 
$t\rightarrow-\infty$. 
\qed
\chapter{The characteristic Cauchy problem}
\label{sec8}
The aim of this chapter is to solve a characteristic Cauchy problem in a space-time 
region near the collapsing star. The results of this chapter will be used later in the proof of the main theorem. 
We start by studying a characteristic Cauchy problem for the Dirac equation in the whole exterior Kerr-Newman space-time. In Section \ref{sec8.1} we formulate the main result of this chapter. Section \ref{sec8.2} is devoted to the 
usual Cauchy problem with data on a lipschitz space-like surface. The main theorem is proven 
in Section \ref{sec8.3}. In Section \ref{sec8.4} we use these results to solve 
the characteristic Cauchy problem near the collapsing star. Our strategy is similar
to that of H\"ormander in \cite{Hoe2} for the characteristic Cauchy problem for the wave
equation (see also \cite{Ni4} for weaker assumptions on the metric). A characteristic
Cauchy problem for the Dirac equation has been considered in \cite{MN} in a somewhat different
setting.
\section{Main results}
\label{sec8.1}
\begin{figure}
\centering\epsfig{figure=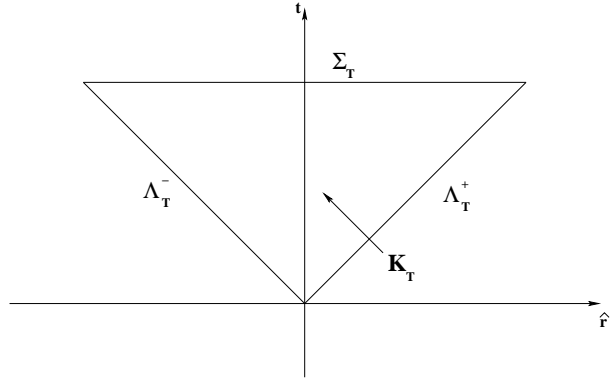,width=8cm}
\caption{The characteristic Cauchy problem}
\label{Chpbws}
\end{figure}

Let (see Figure \ref{Chpbws}) :
\begin{eqnarray*}
\Lambda^{\pm}_T&:=&\{(\pm\hat{r},\hat{r},\omega);\, 0\le\pm\hat{r}\le T,\, \omega\in S^2\},
\Lambda_T:=\Lambda_T^+\cup\Lambda_T^-,\\
K_T&:=&\{(t,\hat{r},\omega);\, |\hat{r}|\le T,\, |\hat{r}|\le t\le T,\, \omega\in S^2\},\\
\Sigma_T&:=&\{(T,\hat{r},\omega);\, |\hat{r}|\le T,\, \omega\in S^2\}.
\end{eqnarray*}
We need the spaces :
\begin{eqnarray*}
{\cal H}_T&:=&L^2(({[}-T,T{]}\times S^2,d\hat{r}d\omega);\C^4),\\
{\cal H}^1_T&:=&\{u\in {\cal H}_T;\, Hu\in {\cal H}_T\},\, ||u||_{{\cal H}^1_T}^2=||u||_{{\cal H}_T}^2+||Hu||_{{\cal H}_T}^2,\\
L^2_{T,-}&:=&L^2(({[}-T,0{]}\times S^2,d\hat{r}d\omega);\C^2),\\
L^2_{T,+}&:=&L^2(({[}0,T{]}\times S^2,d\hat{r}d\omega);\C^2).
\end{eqnarray*}
Let $\Phi_T\in C^{\infty}(\Sigma_T;\S_A\oplus \S^{A'})$. By the usual theorems 
for hyperbolic equations we can associate to $\Phi_T$ a smooth solution $\phi_A\oplus\chi^{A'}
\in C^{\infty}(K_T;\S_A\oplus\S^{A'})$ (see \cite{Ni3} for details). We use the $l^a,n^a,m^a$
tetrad and the $(t,\hat{r},\omega)$ coordinate system. Let 
$\Psi=\sqrt[4]{\frac{\rho^2\Delta\sigma^2}{(r^2+a^2)^2k'^2}}(\phi_0,\phi_1,\chi_{1'},-\chi_{0'})$
be the associated density spinor. The spinor fields $o^A$ and $\iota^A$ are smooth and non vanishing on $\Lambda^{\pm}_T$,
therefore we can associate to this solution the smooth trace of $\Psi$ :
\begin{eqnarray*}
\lefteqn{{\cal T}:\, \Phi_T\mapsto (\Psi_2,\Psi_3)(-\hat{r},\hat{r},\omega)\oplus
(\Psi_1,\Psi_4)(\hat{r},\hat{r},\omega)}\\
&\in& C^{\infty}({[}-T,0{]}\times S^2;\C^2)\oplus
C^{\infty}({[}0,T{]}\times S^2;\C^2).
\end{eqnarray*}
Using the conserved current we obtain by Stokes' theorem :
\begin{eqnarray}
\label{VII.1}
\int_{\Sigma_T}*(\phi_A\bar{\phi}_{A'}dx^{AA'}+\bar{\chi}_A\chi_{A'}dx^{AA'})
=\int_{\Lambda_T}*(\phi_A\bar{\phi}_{A'}dx^{AA'}
+\bar{\chi}_A\chi_{A'}dx^{AA'})
\end{eqnarray}
Let $\varphi(\hat{r})=|\hat{r}|.$ The normal to $\Lambda^-_T$ is $n^a$, the normal to $\Lambda^+_T$ is $l^a$.
We compute : 
\[ (l^a\partial_a\hook d\Omega)|_{\Lambda^+_T}=(n^a\partial_a\hook d\Omega)|_{\Lambda^-_T}
=\sqrt{\frac{2\rho^2\Delta\sigma^2}{(r^2+a^2)^2k'^2}}d\hat{r}\wedge d\omega. \]
Following (\ref{ISC}) we find :
\begin{eqnarray}
\label{VII.1b}
\lefteqn{\int_{\Lambda_T}*(\phi_A\bar{\phi}_{A'}dx^{AA'}+\bar{\chi}_A\chi_{A'}dx^{AA'})}\nonumber\\
&=&\sqrt{2}\int_{-T}^0\int_{S^2}(|\Psi_2|^2+|\Psi_3|^2)(-\hat{r},\hat{r},\omega)d\hat{r}d\omega\nonumber\\
&+&\sqrt{2}\int^{T}_0\int_{S^2}(|\Psi_1|^2+|\Psi_4|^2)(\hat{r},\hat{r},\omega)d\hat{r}d\omega,
\end{eqnarray}
where $\Psi=(\Psi_1,\Psi_2,\Psi_3,\Psi_4).$ It follows : 
\begin{eqnarray}
\label{VII.2}
\lefteqn{\int_{\Sigma_T}*(\phi_A\bar{\phi}_{A'}dx^{AA'}+\bar{\chi}_A\chi_{A'}dx^{AA'})}\nonumber\\
&=&\sqrt{2}\int_{-T}^0\int_{S^2}(|\Psi_2|^2+|\Psi_3|^2)(-\hat{r},\hat{r},\omega)d\hat{r}d\omega\nonumber\\
&+&\sqrt{2}\int^{T}_0\int_{S^2}(|\Psi_1|^2+|\Psi_4|^2)(\hat{r},\hat{r},\omega)d\hat{r}d\omega.
\end{eqnarray}
Therefore the operator ${\cal T}$ possesses an extension to a bounded operator 
\[ {\cal T}\in {\cal L}(L^2(\Sigma_T;\S_A\oplus\S^{A'});L^2({[}-T,0{]}\times S^2;\C^2)\oplus
L^2({[}0,T{]}\times S^2;\C^2)). \]
Our first result is  
\begin{theorem}
\label{thVII.1}
$\frac{1}{\sqrt[4]{2}}{\cal T}$ is an isometry.
\end{theorem}
The proof of the theorem will be given in Section \ref{sec8.3}.

The characteristic data ${\cal T}(\Phi_T)$ contains information only about 
$\phi_1,\chi_{1'}$ on $\Lambda^-_T$ and $\phi_0,\chi_{0'}$ on $\Lambda_T^+$.
The functions $\phi_0,\chi_{0'}$ on $\Lambda^-_T$ resp. $\phi_1,\chi_{1'}$ on $\Lambda^+_T$ are obtained 
from the given data by restriction of the equation to $\Lambda^{\pm}_T.$ On $\Lambda^-_T$ we have (see (\ref{NewmanDirac})) :
\begin{eqnarray*}
\left. \begin{array}{c}
n^{\bf{a}}(\partial_{\bf{a}}-iq\Phi_{\bf a})\phi_0-m^{\bf{a}}(\partial_{\bf{a}}-iq\Phi_{\bf a})\phi_1+(\mu-\gamma)\phi_0+(\tau-\beta)\phi_1=\frac{m}{\sqrt{2}}\chi_{1'},\\
n^{\bf{a}}(\partial_{\bf{a}}-iq\Phi_{\bf a})\chi_{0'}-\bar{m}^{\bf{a}}(\partial_{\bf{a}}-iq\Phi_{\bf a})\chi_{1'}+(\bar{\mu}-\bar{\gamma})\chi_{0'}+(\bar{\tau}-\bar{\beta})\chi_{1'}
=\frac{m}{\sqrt{2}}\phi_1,
\end{array} \right\}
\end{eqnarray*}  
where $\phi_1,\chi_{1'}$ have to be considered as source terms. Putting $g(\hat{r},\omega)=\Psi(-\hat{r},\hat{r},\omega),$\\
$\hat{g}(\hat{r},\omega)=\Psi(\hat{r},\hat{r},\omega)$
we find :
\begin{eqnarray}
\label{VII.3}
\left. \begin{array}{rcl} -D_{\hat{r}}g_{1,4}&=&((P_{\omega}+W)g)_{1,4},\\
g_{1,4}(0,\omega)&=&\hat{g}_{1,4}(0,\omega), \end{array} \right\},\\
\label{VII.4}
\left. \begin{array}{rcl} D_{\hat{r}}\hat{g}_{2,3}&=&((P_{\omega}+W)\hat{g})_{2,3},\\
\hat{g}_{2,3}(0,\omega)&=&g_{2,3}(0,\omega). \end{array} \right\}
\end{eqnarray}
Here $P_{\omega},\, W$ are defined in (\ref{POMEGA}), (\ref{W}).
We understand $\hat{r}$ as a time parameter which goes from $0$ to $-T$ for (\ref{VII.3})
and from $0$ to $T$ for (\ref{VII.4}). We write (\ref{VII.3}) as 
\begin{eqnarray}
\label{VII.5}
\left. \begin{array}{rcl} \partial_{\hat{r}}g_{1,4}&=&iA(\hat{r})g_{1,4}+S(\hat{r}),\\
g_{1,4}(0,\omega)&=&\hat{g}_{1,4}(0,\omega) \end{array}\right\}
\end{eqnarray}
with
\begin{eqnarray*}
A(\hat{r})&=&-\left\{\frac{m^1_{\theta}}{2}+\left(\begin{array}{cc} W_{11} & W_{14} \\ W_{41} & W_{44} \end{array} \right)\right\},\\
S(\hat{r})&=&-\left\{\left\{\frac{m^2_{\theta}}{2}+\frac{a_0h^2}{\sin\theta}D_{\varphi}\left(\begin{array}{cc} -1 & 0 \\ 0 & 1 \end{array} \right)
+\left(\begin{array}{cc} W_{12} & W_{13} \\ W_{42} & W_{43} \end{array} \right)
\right\}\left(\begin{array}{c} g_2 \\ g_3 \end{array} \right)\right\}.
\end{eqnarray*}
We want to show that (\ref{VII.5}) has a unique solution and to this purpose we must analyze $A(\hat{r})$.
We have :
\begin{eqnarray*}
\frac{m^1_{\theta}}{2}=\frac{1}{2}\frac{\beta h \sqrt{a_0}}{\sqrt{\alpha+1}}
\left(D_{\theta}+\frac{\cot\theta}{2i}\right)\sqrt{a_0}h\sqrt{\alpha+1}+hc.
\end{eqnarray*}
We want to show that the operator
\[ d_{\theta}:=D_{\theta}+\frac{\cot\theta}{2i} \]
is selfadjoint with some suitable domain. To this purpose we introduce the
following unitary transformation :
\begin{eqnarray*}
U\,:\,\left. \begin{array}{ccc} L^2(S^2,d\omega)&\rightarrow& L^2(S^2,d\theta d\varphi)\\
u(\theta,\varphi)&\mapsto& u(\theta,\varphi)\sqrt{\sin\theta}. \end{array} \right.
\end{eqnarray*}
Clearly
\[ \hat{d}_{\theta}=Ud_{\theta}U^*=D_{\theta}, \]
which is selfadjoint on $L^2(S^2,d\theta d\varphi)$ with domain 
\[ D(D_{\theta})=\{u\in L^2(S^2,d\theta d\varphi);\,D_{\theta}u\in L^2(S^2,d\theta d\varphi),u(0,.)=
u(\pi,.) \}. \]
Then $d_{\theta}$ is selfadjoint with domain $D(d_{\theta})=U^*D(D_{\theta})$.
It is easy to check that $\hat{r}\mapsto A(\hat{r})v$ is continuously differentiable for 
$v\in Y:=(D(d_{\theta}))^2$. We can apply \cite[Theorem 5.4.8]{Pa} 
and associate a unitary evolution system $V(\hat{r},\hat{r}')$.
For smooth $g_{2,3},\hat{g}_{1,4}$ we have $S(\hat{r})\in C({[}-T,0{]};Y)$
and $\hat{g}(0,\omega)\in Y$. By \cite[Theorem 5.5.2]{Pa}
(\ref{VII.5}) possesses a unique $Y$-valued solution given by :
\begin{eqnarray}
\label{VII.6}
g_{1,4}(\hat{r},\omega)=V(\hat{r},0)\hat{g}_{1,4}(0,\omega)+\int_0^{\hat{r}}V(\hat{r},\hat{r}')S(\hat{r}')d\hat{r}'.
\end{eqnarray}
For given $g_{2,3},\, \hat{g}_{1,4}$ we define $g_{1,4},\, \hat{g}_{2,3}$ as the solutions of the partial differential equations (\ref{VII.3}), (\ref{VII.4}) and put :
\begin{eqnarray}
\label{VII.7}
g^H_{2,3}(\hat{r},\omega)&:=&\frac{1}{2}(-D_{\hat{r}}g_{2,3}+((P_{\omega}+W)g)_{2,3}),\\
\label{VII.8}
\hat{g}^H_{1,4}(\hat{r},\omega)&:=&\frac{1}{2}(D_{\hat{r}}\hat{g}_{1,4}+((P_{\omega}+W)\hat{g})_{1,4}).
\end{eqnarray}
As $g_{1,4}$ is a $Y$-valued solution we have (the argument for $\hat{g}_{2,3}$ is analogous):
\begin{eqnarray*}
g^H_{2,3}\in L^2(({[}-T,0{]}\times S^2,d\hat{r}d\omega);\C^2),\, \hat{g}_{1,4}^H\in 
L^2(({[}0,T{]}\times S^2,d\hat{r}d\omega);\C^2).
\end{eqnarray*}
We define $\tilde{H}^1$ as the completion of $C^{\infty}({[}-T,0{]}\times S^2;\C^2)\oplus C^{\infty}({[}0,T{]}\times S^2;\C^2)$
in the norm :
\begin{eqnarray}
\label{VII.9}
||(g_{2,3},\hat{g}_{1,4})||_{\tilde{H}^1}^2=2(||g_{2,3}||^2_{L^2_{T,-}}+||\hat{g}_{1,4}||_{L^2_{T,+}}^2
+||g^H_{2,3}||_{L^2_{T,-}}^2+||\hat{g}^H_{1,4}||_{L^2_{T,+}}^2).
\end{eqnarray}
We now start with $\Psi_T\in C^{\infty}(\Sigma_T;\C^4).$ Then we can associate a classical
solution $\Psi\in C^{\infty}(K_T;\C^4)$ and the traces 
\[ g_{2,3}(\hat{r},\omega)=\Psi_{2,3}(-\hat{r},\hat{r},\omega),\, -T\le \hat{r}\le 0\quad\mbox{and}\quad \hat{g}_{1,4}(\hat{r},\omega)=\Psi_{1,4}(\hat{r},\hat{r},\omega),\, 0\le \hat{r}\le T\] 
are well defined. By the previous discussion 
\[ g_{1,4}(\hat{r},\omega)=\Psi_{1,4}(-\hat{r},\hat{r},\omega),\, -T\le \hat{r}\le 0\quad\mbox{and}\quad \hat{g}_{2,3}(\hat{r},\omega)=\Psi_{2,3}(\hat{r},\hat{r},\omega),\, 0\le\hat{r}\le T\]
are solutions of equations (\ref{VII.3}), (\ref{VII.4}).
As $\Psi$ is a solution of the Dirac equation, so is 
$H\Psi\in C^{\infty}(K_T;\C^4).$ We want to calculate $H\Psi|_{\Lambda^{\pm}_T}$ in terms of 
$g_{2,3},\, \hat{g}_{1,4}$ and to this purpose we introduce 
characteristic coordinates :
\begin{eqnarray*}
\left(\begin{array}{rcl} X&=&t-\hat{r} \\ T&=&t+\hat{r} \end{array} \right)\Leftrightarrow
\left(\begin{array}{rcl} t&=&\frac{X+T}{2} \\ \hat{r}&=&\frac{T-X}{2} \end{array} \right).
\end{eqnarray*}
Then we have 
\begin{eqnarray*}
(H\Psi)_{2,3}&=&(D_X-D_T)\Psi_{2,3}+((P_{\omega}+W)\Psi)_{2,3}\\
&=&D_X\Psi_{2,3}+\frac{1}{2}((P_{\omega}+W)\Psi)_{2,3}\\
\Rightarrow (H\Psi)_{2,3}(-\hat{r},\hat{r},\omega)&=&\frac{1}{2}(-D_{\hat{r}}g_{2,3}+((P_{\omega}+W)g)_{2,3}).
\end{eqnarray*}
In a similar manner we find :
\begin{eqnarray*}
(H\Psi)_{1,4}(\hat{r},\hat{r},\omega)=\frac{1}{2}(D_{\hat{r}}\hat{g}_{1,4}+((P_{\omega}+W)\hat{g})_{1,4}).
\end{eqnarray*}
Using the identity (\ref{VII.2}) we find :
\begin{eqnarray*}
||H\Psi_T||_{{\cal H}_T}^2=2(||g^H_{2,3}||^2_{L^2_{T,-}}+||\hat{g}^H_{1,4}||_{L^2_{T,+}}^2)
\end{eqnarray*}
and therefore :
\begin{eqnarray}
\label{VII.10}
||\Psi_T||^2_{{\cal H}^1_T}=||(g_{2,3},\hat{g}_{1,4})||_{\tilde{H}^1}^2.
\end{eqnarray}
This means that the trace operator
\begin{eqnarray*}
{\cal T}\, : \left.\begin{array}{ccc}  C^{\infty}(\Sigma_T;\C^4) &\rightarrow & C^{\infty}({[}-T,0{]}\times S^2;\C^2)\oplus C^{\infty}({[}0,T{]}\times S^2;\C^2)\\
\Psi_T &\mapsto& (\Psi_{2,3}(-\hat{r},\hat{r},\omega),\Psi_{1,4}(\hat{r},\hat{r},\omega))
\end{array}\right.
\end{eqnarray*}
extends to a bounded operator 
\[ {\cal T}_H\in {\cal L}({\cal H}_T^1;\tilde{H}^1). \]
Our second result is 
\begin{theorem}
\label{thVII.2}
${\cal T}_H$ is an isometry.
\end{theorem}
The proof of this theorem will be given in Section \ref{sec8.3}.
\begin{remark}
\label{RVII.0}
If $\Psi_T\in {\cal H}^1_T$, then the trace $(\Psi_{2,3}(-\hat{r},\hat{r},\omega),\Psi_{1,4}(\hat{r},\hat{r},\omega))$
exists in the usual sense and it is in $\tilde{H}^1$ by (\ref{VII.10}). This means that the operator ${\cal T}_H$
is defined as the usual trace and that ${\cal T}$ is an extension of ${\cal T}_H$.
\end{remark}
\section{The Cauchy problem with data on a lipschitz space-like hypersurface}
\label{sec8.2}
Let $\varphi : {[}-T_{\varphi},T_{\varphi}{]}\rightarrow \R$ be lipschitz continuous, $-T_{\varphi}\le -T,\, T_{\varphi}\ge T$,
\begin{eqnarray}
\label{VII.11} 
|\varphi'(\hat{r})|\le \tilde{\alpha}<1\quad\mbox{a.e.}
\end{eqnarray}
Thanks to (\ref{VII.11}) the hypersurface
$\Lambda_{\varphi}:=\{(\varphi(\hat{r}),\hat{r},\omega);\, |\hat{r}|\le T_{\varphi},\, \omega\in S^2 \}$
is space-like. Indeed we have 
\begin{eqnarray*}
\lefteqn{g(\varphi'(\hat{r})\partial_t+\partial_{\hat{r}},\varphi'(\hat{r})\partial_t+\partial_{\hat{r}})}\\
&=&\left(1+\frac{Q^2-2Mr}{\rho^2}\right)\varphi'^2-\frac{\rho^2\Delta}{(r^2+a^2)^2k'^2}
<\left(1+\frac{Q^2-2Mr}{\rho^2}\right)-\frac{\rho^2\Delta}{\sigma^2}\\
&=&-\frac{a^2\sin^2\theta(2Mr-Q^2)^2}{\rho^2\sigma^2}\le 0,\\
g(\partial_{\varphi},\partial_{\varphi})&<&0,\quad g(\partial_{\theta},\partial_{\theta})<0.
\end{eqnarray*}
Let for $0\le t \le T$
\begin{eqnarray*}
\Sigma^{\varphi}_t&:=&\{(t,\hat{r},\omega)\in \{t\}\times{[}-T_{\varphi},T_{\varphi}{]}\times S^2;\, t>\varphi(\hat{r})\},\\
K^{\varphi}_T&:=&\bigcup_{0\le t\le T}\Sigma^{\varphi}_t,\\
{\cal R}_t^{\varphi}&:=&\{\hat{r}\in (-T_{\varphi},T_{\varphi});\, t>\varphi(\hat{r})\}.
\end{eqnarray*}
We will suppose ${\cal R}^{\varphi}_T=(-T_{\varphi},T_{\varphi}),\, 
\varphi(-T_{\varphi})=\varphi(T_{\varphi})=T.$ We also define the spaces :
\begin{eqnarray*}
{\cal H}_{t,\varphi}&=&L^2((\Sigma^{\varphi}_t,d\hat{r}d\omega);\C^4),\\ 
{\cal H}_{t,\varphi}^1&=&\{u\in {\cal H}_{t,\varphi};\, Hu\in {\cal H}_{t,\varphi}\},\,
||u||_{{\cal H}^1_{t,\varphi}}^2=||u||_{{\cal H}_{t,\varphi}}^2+||Hu||_{{\cal H}_{t,\varphi}}^2.
\end{eqnarray*}
The aim of this section is to solve the following Cauchy problem :
\begin{eqnarray}
\label{VII.12}
\left. \begin{array}{rcl} \partial_t\Psi&=&iH\Psi\quad (t,\hat{r},\omega)\in K^{\varphi}_T,\\
\Psi(\varphi(\hat{r}),\hat{r},\omega)&=&g(\hat{r},\omega)\quad (\hat{r},\omega)\in {[}-T_{\varphi},T_{\varphi}{]}\times S^2.
\end{array}\right\}
\end{eqnarray}
We first define the space of data. For $g\in C^{\infty}({[}-T_{\varphi},T_{\varphi}{]};\C^4)$
we define :
\begin{eqnarray}
\label{VII.13}
\lefteqn{g^H_{\varphi}(\hat{r},\omega)}\nonumber\\
&:=&\left(Diag\left(\frac{1}{1+\varphi'},-\frac{1}{1-\varphi'},
-\frac{1}{1-\varphi'},\frac{1}{1+\varphi'}\right)D_{\hat{r}}g\right)(\hat{r},\omega)\nonumber\\
&+&\left(Diag\left(\frac{1}{1+\varphi'},\frac{1}{1-\varphi'},
\frac{1}{1-\varphi'},\frac{1}{1+\varphi'}\right)\left(P_{\omega}+W\right)g\right)(\hat{r},\omega).
\end{eqnarray}
We define $\tilde{H}^1_{\varphi}$ as the completion of $C^{\infty}({[}-T_{\varphi},T_{\varphi}{]}\times S^2;\C^4)$ 
in the norm :
\begin{eqnarray}
\label{VII.14}
||g||_{\tilde{H}_{\varphi}^1}^2&:=&||(1-\varphi')^{1/2}g_{2,3}||^2_{L^2({[}-T_{\varphi},T_{\varphi}{]}\times S^2;\C^2)}\nonumber\\
&+&||(1-\varphi')^{1/2}(g^H_{\varphi})_{2,3}||^2_{L^2({[}-T_{\varphi},T_{\varphi}{]}\times S^2;\C^2)}\nonumber\\
&+&||(1+\varphi')^{1/2}g_{1,4}||^2_{L^2({[}-T_{\varphi},T_{\varphi}{]}\times S^2;\C^2)}\nonumber\\
&+&||(1+\varphi')^{1/2}(g^H_{\varphi})_{1,4}||^2_{L^2({[}-T_{\varphi},T_{\varphi}{]}\times S^2;\C^2)}.
\end{eqnarray}
Let $\Psi_T\in C^{\infty}(\Sigma^{\varphi}_T;\C^4)$ and $\Psi\in C^{\infty}(K^{\varphi}_T;\C^4)$
be the associated solution of the Dirac equation. Then $\Psi(\varphi(\hat{r}),\hat{r},\omega)\in 
H^1({[}-T_{\varphi},T_{\varphi}{]};(H^{\infty}(S^2))^4)$ and as in Section \ref{sec8.1} we find using Stokes' theorem~:
\begin{eqnarray}
\label{VII.15}
\lefteqn{\int_{-T_{\varphi}}^{T_{\varphi}}\int_{S^2}|\Psi|^2d\hat{r}d\omega}\nonumber\\
&=&\int_{-T_{\varphi}}^{T_{\varphi}}\int_{S^2}
\left((1-\varphi')\left(|\Psi_2|^2+|\Psi_3|^2\right)+(1+\varphi')\left(|\Psi_1|^2+|\Psi_4|^4\right)\right)(\varphi(\hat{r}),\hat{r},\omega)
d\hat{r}d\omega.\nonumber\\ 
\end{eqnarray}
We want to estimate $\int_{-T_{\varphi}}^{T_{\varphi}}\int_{S^2}|H\Psi|^2d\hat{r}d\omega$.
To this purpose we introduce a sequence of smooth functions $\varphi_{\epsilon}:
{[}-T^1_{\varphi_{\epsilon}},T^2_{\varphi_{\epsilon}}{]}\rightarrow \R$ s.t.
\begin{eqnarray}
\label{VII.16}
{[}-T_{\varphi},T_{\varphi}{]}&\subseteq&{[}-T_{\varphi_{\epsilon}}^1,T^2_{\varphi_{\epsilon}}{]},\nonumber\\
|\varphi'_{\epsilon}(\hat{r})|&\le& \tilde{\tilde{\alpha}}<1,\nonumber\\
\varphi_{\epsilon}&\rightarrow& \varphi \quad L^{\infty}({[}-T_{\varphi},T_{\varphi}{]}),\nonumber\\
\varphi_{\epsilon}'&\rightarrow& \varphi'\quad\mbox{a.e.}\, {[}-T_{\varphi},T_{\varphi}{]},\nonumber\\
\varphi_{\epsilon}(\hat{r})&\le& \varphi(\hat{r})\quad \forall \hat{r}\in {[}-T_{\varphi},T_{\varphi}{]},\\
\varphi_{\epsilon}(-T^1_{\varphi_{\epsilon}})&=&T=\varphi_{\epsilon}(T^2_{\varphi_{\epsilon}}).\nonumber
\end{eqnarray}
This approximation can be achieved by convolution with smooth functions (see \cite[Lemma 3]{Hoe2}).
Note that we may have to replace the approximation by $\varphi_{\epsilon}-|\varphi-\varphi_{\epsilon}|_{L^{\infty}}$
to achieve (\ref{VII.16}). In order to compute $(H\Psi)(\varphi_{\epsilon}(\hat{r}),\hat{r},\omega)$ 
we introduce the following change of variables :
\begin{eqnarray*}
\left. \begin{array}{rcl} \tau&=&t-\varphi_{\epsilon}(\hat{r}),\\ x&=&\hat{r} \end{array}
\right\}\Rightarrow \partial_t=\partial_{\tau};\, \partial_{\hat{r}}=\partial_x-\varphi_{\epsilon}'(\hat{r})\partial_{\tau}.
\end{eqnarray*}
We have 
\begin{eqnarray}
\label{VII.18}
\partial_t\Psi&=&iH\Psi\nonumber\\
\Leftrightarrow \partial_{\tau}\Psi&=&\left(1+\Gamma^1\varphi'_{\epsilon}\right)^{-1}
\left(\Gamma^1\partial_x\Psi+i(P_{\omega}+W)\Psi\right).
\end{eqnarray}
Using (\ref{VII.18}) we calculate :
\begin{eqnarray*}
H\Psi&=&Diag\left(\frac{1}{1+\varphi_{\epsilon}'},-\frac{1}{1-\varphi_{\epsilon}'},
-\frac{1}{1-\varphi_{\epsilon}'},\frac{1}{1+\varphi_{\epsilon}'}\right)D_x\Psi\nonumber\\
&+&Diag\left(\frac{1}{1+\varphi_{\epsilon}'},\frac{1}{1-\varphi_{\epsilon}'},
\frac{1}{1-\varphi_{\epsilon}'},\frac{1}{1+\varphi_{\epsilon}'}\right)\left(P_{\omega}+W\right)\Psi
\end{eqnarray*}
Putting $g_{\epsilon}(\hat{r},\omega)=\Psi(\varphi_{\epsilon}(\hat{r}),\hat{r},\omega)$
we find :
\begin{eqnarray}
\label{VII.19}
(H\Psi)(\varphi_{\epsilon}(\hat{r}),\hat{r},\omega)&=&\left(Diag\left(\frac{1}{1+\varphi_{\epsilon}'},-\frac{1}{1-\varphi_{\epsilon}'},
-\frac{1}{1-\varphi_{\epsilon}'},\frac{1}{1+\varphi_{\epsilon}'}\right)D_{\hat{r}}g_{\epsilon}\right)(\hat{r},\omega)\nonumber\\
&+&\left(Diag\left(\frac{1}{1+\varphi_{\epsilon}'},\frac{1}{1-\varphi_{\epsilon}'},
\frac{1}{1-\varphi_{\epsilon}'},\frac{1}{1+\varphi_{\epsilon}'}\right)\left(P_{\omega}+W\right)g_{\epsilon}\right)(\hat{r},\omega)\nonumber\\
&=:&g_{\epsilon}^H.
\end{eqnarray}
Using Stokes' theorem we obtain :
\begin{eqnarray}
\label{VII.20}
\int_{-T^1_{\varphi_{\epsilon}}}^{T^2_{\varphi_{\epsilon}}}\int_{S^2}|H\Psi|^2d\hat{r}d\omega
&=&||(1-\varphi'_{\epsilon})^{1/2}(g_{\epsilon}^H)_{2,3}||^2_{L^2({[}-T^1_{\varphi_{\epsilon}},T^2_{\varphi_{\epsilon}}{]}\times S^2;\C^2)}\nonumber\\
&+&||(1+\varphi'_{\epsilon})^{1/2}(g_{\epsilon}^H)_{1,4}||^2_{L^2({[}-T^1_{\varphi_{\epsilon}},T^2_{\varphi_{\epsilon}}{]}\times S^2;\C^2)}.
\end{eqnarray}
Using $T^j_{\varphi_{\epsilon}}\rightarrow T_{\varphi},\, j=1,2$ as well as the fact that
$\Psi$ is smooth we can take the limit in (\ref{VII.20}) and find :
\begin{eqnarray}
\label{VII.21}
\int_{-T_{\varphi}}^{T_{\varphi}}\int_{S^2}|H\Psi|^2d\hat{r}d\omega
&=&||(1-\varphi')^{1/2}(g_{\varphi}^H)_{2,3}||^2_{L^2({[}-T_{\varphi},T_{\varphi}{]}\times S^2;\C^2)}\nonumber\\
&+&||(1+\varphi')^{1/2}(g_{\varphi}^H)_{1,4}||^2_{L^2({[}-T_{\varphi},T_{\varphi}{]}\times S^2;\C^2)}.
\end{eqnarray}
Putting (\ref{VII.15}) and (\ref{VII.21}) together we find :
\begin{eqnarray}
\label{VII.22}
||\Psi||_{{\cal H}_{T,\varphi}^1}^2=||g||^2_{\tilde{H}_{\varphi}^1}.
\end{eqnarray}
The equality (\ref{VII.22}) shows that the trace operator 
\[ \Psi_T\in C^{\infty}(\Sigma^{\varphi}_T;\C^4)\mapsto \Psi(\varphi(\hat{r}),\hat{r},\omega) \]
possesses an extension to a bounded operator 
\[ \hat{\cal T}\in {\cal L}({\cal H}^1_{T,\varphi},\tilde{H}_{\varphi}^1). \] 
The result of this section is :
\begin{theorem}
\label{thVII.3}
$\hat{\cal T}$ is an isometry.
\end{theorem}

{\bf Proof.}

Because of (\ref{VII.22}) we only need to prove surjectivity. We have to construct 
for $g\in \tilde{H}^1_{\varphi}$ a solution of 
\begin{eqnarray}
\label{VII.23}
\left. \begin{array}{rcl} \partial_t\Psi&=&iH\Psi,\\
\Psi(\varphi(\hat{r}),\hat{r},\omega)&=&g(\hat{r},\omega). \end{array} \right\}
\end{eqnarray} 
We first suppose $g\in C^{\infty}({[}-T_{\varphi},T_{\varphi}{]}\times S^2;\C^2)$ 
and consider the approximate problem~:
\begin{eqnarray}
\label{VII.24}
\left. \begin{array}{rcl} \partial_t\Psi^{\epsilon}&=&iH\Psi^{\epsilon},\\
\Psi^{\epsilon}(\varphi_{\epsilon}(\hat{r}),\hat{r},\omega)&=&g(\hat{r},\omega), \end{array} \right\}
\end{eqnarray}
where $\varphi_{\epsilon}$ is as before and $g(\hat{r},\omega)$ is identified with a smooth extension
on ${[}-T^1_{\varphi_{\epsilon}},T^2_{\varphi_{\epsilon}}{]}\times S^2.$ As $\varphi_{\epsilon}$
is smooth it is well known that (\ref{VII.24}) possesses a smooth solution and we have 
the estimate~:
\begin{eqnarray}
\label{VII.25}
||\Psi^{\epsilon}||^2_{{\cal H}^1_{T,\varphi}}\le ||\Psi^{\epsilon}||^2_{{\cal H}^1_{T,\varphi_{\epsilon}}}
=||g||^2_{\tilde{H}^1_{\varphi_{\epsilon}}}.
\end{eqnarray}
Similar estimates hold on ${\cal H}^1_{t,\varphi},\, \min \varphi\le t \le T=\max\varphi. $
As the R.H.S. of (\ref{VII.25}) is uniformly bounded we find uniform bounds (here we also use
$\varphi_{\epsilon}\le \varphi$) :
\begin{eqnarray*}
||\Psi^{\epsilon}||_{(H^1(K^{\varphi}_T))^4}\lesssim 1,\quad ||\Psi^{\epsilon}||_{{\cal H}^1_{T,\varphi}}\lesssim 1.
\end{eqnarray*}
We can therefore extract a subsequence which we denote again $\Psi^{\epsilon}$ s.t.
\begin{eqnarray*}
\Psi^{\epsilon}\rightharpoonup \Psi\quad {\cal H}_{T,\varphi}^1,\\
\Psi^{\epsilon}\rightharpoonup \Psi\quad (H^1(K^{\varphi}_T))^4,\\
\Psi^{\epsilon}\rightarrow \Psi\quad (H^{1/2}(K^{\varphi}_T))^4.
\end{eqnarray*}
The limit $\Psi$ is a solution of the Dirac equation. We have to check that
$\Psi(\varphi(\hat{r}),\hat{r},\omega)=g(\hat{r},\omega).$ 
To this aim we estimate for $\epsilon$ small :
\begin{eqnarray*}
\lefteqn{\int_{-T_{\varphi}}^{T_{\varphi}}\int_{S^2}|\Psi^{\epsilon}(\varphi(\hat{r}),\hat{r},\omega)
-g(\hat{r},\omega)|^2d\hat{r}d\omega}\\
&=&\int_{-T_{\varphi}}^{T_{\varphi}}\int_{S^2}|\Psi^{\epsilon}(\varphi(\hat{r}),\hat{r},\omega)
-\Psi^{\epsilon}(\varphi_{\epsilon}(\hat{r}),\hat{r},\omega)|^2d\hat{r}d\omega\\
&=&\int_{-T_{\varphi}}^{T_{\varphi}}\int_{S^2}\left|\int_{\varphi_{\epsilon}(\hat{r})}^{\varphi(\hat{r})}
\partial_t\Psi^{\epsilon}(t,\hat{r},\omega)\right|^2dtd\hat{r}d\omega\\
&\le&\int_{-T_{\varphi}}^{T_{\varphi}}\int_{S^2}\int_{\varphi_{\epsilon}(\hat{r})}^{\varphi(\hat{r})}
|\partial_t\Psi^{\epsilon}(t,\hat{r},\omega)|^2dt|\varphi(\hat{r})-\varphi_{\epsilon}(\hat{r})|d\hat{r}d\omega\\
&\le&|\varphi-\varphi_{\epsilon}|_{L^{\infty}}\int_{\frac{1}{2}\min \varphi}^{\max \varphi}
\int_{-T_{\varphi}}^{T_{\varphi}}\int_{S^2}|H\Psi^{\epsilon}(t,\hat{r},\omega)|^2d\hat{r}d\omega dt\\
&\lesssim& |\varphi_{\epsilon}-\varphi|_{L^{\infty}}||H\Psi^{\epsilon}||^2_{{\cal H}_{T,\varphi}}
\lesssim|\varphi_{\epsilon}-\varphi|_{L^{\infty}}\rightarrow0,\, \epsilon\rightarrow 0.
\end{eqnarray*}
Here we have used the Cauchy-Schwarz inequality. On the other hand we know that 
\[ ||\Psi^{\epsilon}(\varphi(\hat{r}),\hat{r},\omega)-\Psi(\varphi(\hat{r}),\hat{r},\omega)||_{L^2(({[}-T_{\varphi},T_{\varphi}{]}\times S^2,d\hat{r}d\omega);\C^4)}
\le ||\Psi^{\epsilon}-\Psi||_{(H^{1/2}(K^{\varphi}_T))^4}\rightarrow 0. \]
It follows that $\Psi(\varphi(\hat{r}),\hat{r},\omega)=g(\hat{r},\omega)$.
The solution $\Psi$ satisfies 
\begin{eqnarray}
\label{VII.26}
||\Psi||_{{\cal H}^1_{T,\varphi}}\le ||g||_{\tilde{H}^1_{\varphi}},\, 
||\Psi||_{(H^1(K^{\varphi}_T))^4}\lesssim ||g||_{\tilde{H}^1_{\varphi}}.
\end{eqnarray}
If $g\in \tilde{H}^1_{\varphi}$, then we approximate it by a sequence $g^n$ of $C^{\infty}$ 
functions. Then by (\ref{VII.26}) the associated sequence $\Psi^n$ of solutions converges 
to some $\Psi$ in the norms ${\cal H}^1_{T,\varphi},\, (H^1(K^{\varphi}_T))^4$.
As $\Psi\in (H^1(K_T^{\varphi}))^4$, the trace $\Psi|_{\Lambda_{\varphi}}$ exists 
and we have :
\begin{eqnarray*}
\int_{-T_{\varphi}}^{T_{\varphi}}\int_{S^2}|\Psi(\varphi(\hat{r}),\hat{r},\omega)-g(\hat{r},\omega)|^2d\hat{r}d\omega&\le& ||\Psi-\Psi^n||^2_{(H^{1/2}(K^{\varphi}_T))^4}\\
&+&\int_{-T_{\varphi}}^{T_{\varphi}}
\int_{S^2}|g^n-g|^2d\hat{r}d\omega\rightarrow 0
\end{eqnarray*}
and thus $\hat{\cal T}\Psi=g$. This concludes the proof of the theorem.
\qed 
\section{Proof of Theorems \ref{thVII.1} and \ref{thVII.2}}
\label{sec8.3}
Because of (\ref{VII.2}) and (\ref{VII.10}) we only have to show surjectiviy. Let 
$\lambda<1$ and $\varphi_{\lambda}=\lambda|\hat{r}|.$ Then the hypersurface 
$\Lambda_{\varphi_{\lambda}}=\{(\varphi_{\lambda}(\hat{r}),\hat{r},\omega);\, |\hat{r}|\le T,\, \omega\in S^2\}$ 
is a lipschitz space-like hypersurface. We first suppose 
\[g_{2,3}\in C^{\infty}({[}-T,0{]}\times S^2;\C^2)\quad\mbox{and}\quad\hat{g}_{1,4}\in C^{\infty}({[}0,T{]}\times S^2;\C^2)\]
(which we extend to smooth functions on ${[}-\frac{T}{\lambda},0{]}\times S^2$ resp. ${[}0,\frac{T}{\lambda}{]}\times S^2$). Let $g_{1,4}(\hat{r},\omega)$
and $\hat{g}_{2,3}(\hat{r},\omega)$ be the solutions of (\ref{VII.3}) and (\ref{VII.4}).
We consider the approximate problem :
\begin{eqnarray}
\label{VII.27}
\left.\begin{array}{rcl} \partial_t\Psi^{\lambda}&=&iH\Psi^{\lambda},\\
\Psi^{\lambda}(-\lambda\hat{r},\hat{r},\omega)&=&g(\hat{r},\omega)\quad -\frac{T}{\lambda}\le\hat{r}\le 0,\\
\Psi^{\lambda}(\lambda\hat{r},\hat{r},\omega)&=&\hat{g}(\hat{r},\omega)\quad 0\le \hat{r}\le \frac{T}{\lambda}.
\end{array}\right\}
\end{eqnarray}
We put 
\begin{eqnarray*}
\tilde{g}(\hat{r},\omega):=\left\{\begin{array}{cc} g(\hat{r},\omega)\quad -\frac{T}{\lambda}\le \hat{r}\le 0,\\
\hat{g}(\hat{r},\omega)\quad 0\le \hat{r}\le \frac{T}{\lambda}.\end{array}\right.
\end{eqnarray*}
Starting with $\tilde{g}$ we define $\tilde{g}^H_{\varphi_{\lambda}}$ as in (\ref{VII.13}).
For $-\frac{T}{\lambda}\le\hat{r}\le 0$ we have :
\begin{eqnarray*}
\tilde{g}^H_{\varphi_{\lambda}}(\hat{r},\omega)&:=&\left(Diag\left(\frac{1}{1-\lambda},-\frac{1}{1+\lambda},
-\frac{1}{1+\lambda},\frac{1}{1-\lambda}\right)D_{\hat{r}}g\right)(\hat{r},\omega)\nonumber\\
&+&\left(Diag\left(\frac{1}{1-\lambda},\frac{1}{1+\lambda},
\frac{1}{1+\lambda},\frac{1}{1-\lambda}\right)\left(P_{\omega}+W\right)g\right)(\hat{r},\omega).
\end{eqnarray*}
Note that the first and fourth components are zero because $g$ is a solution of (\ref{VII.3}).
In a similar way we find for $0\le \hat{r}\le \frac{T}{\lambda}$:
\begin{eqnarray*}
\tilde{g}^H_{\varphi_{\lambda}}(\hat{r},\omega)&:=&\left(Diag\left(\frac{1}{1+\lambda},-\frac{1}{1-\lambda},
-\frac{1}{1-\lambda},\frac{1}{1+\lambda}\right)D_{\hat{r}}\hat{g}\right)(\hat{r},\omega)\nonumber\\
&+&\left(Diag\left(\frac{1}{1+\lambda},\frac{1}{1-\lambda},
\frac{1}{1-\lambda},\frac{1}{1+\lambda}\right)\left(P_{\omega}+W\right)\hat{g}\right)(\hat{r},\omega).
\end{eqnarray*}
Here the second and third components are zero because $\hat{g}$ is a solution of (\ref{VII.4}).
As $(g_{2,3},\hat{g}_{1,4})\in \tilde{H}^1$ we see that $\tilde{g}\in \tilde{H}^1_{\varphi_{\lambda}}.$
Therefore (\ref{VII.27}) possesses by Theorem \ref{thVII.3} a unique solution $\Psi^{\lambda}$ satisfying 
the energy estimate :
\begin{eqnarray}
\label{VII.28}
||\Psi^{\lambda}||^2_{{\cal H}_T^1}\le ||\Psi^{\lambda}||^2_{{\cal H}_{\varphi_{\lambda},T}^1}
= ||g||^2_{\tilde{\cal H}_{\varphi_{\lambda}}^1}.
\end{eqnarray}
As
\begin{eqnarray*}
(\tilde{g}^H_{\varphi_{\lambda}})_{1,4}=0\quad \forall -\frac{T}{\lambda}\le \hat{r}\le 0,\\
(\tilde{g}^H_{\varphi_{\lambda}})_{2,3}=0\quad \forall 0\le \hat{r}\le \frac{T}{\lambda}
\end{eqnarray*}
we see that the R.H.S. of (\ref{VII.28}) is uniformly bounded. Repeating the argument 
for the spaces ${\cal H}^1_{\varphi_{\lambda},t},\, \min\varphi_{\lambda}\le t \le T=\max\varphi_{\lambda}$
we see that we can extract a subsequence, still denoted $\Psi^{\lambda}$, s.t.
\begin{eqnarray*}
\Psi^{\lambda}\rightharpoonup \Psi\quad {\cal H}_T^1,\\
\Psi^{\lambda}\rightharpoonup\Psi\quad (H^1(K_T))^4,\\
\Psi^{\lambda}\rightarrow\Psi\quad (H^{1/2}(K_T))^4.
\end{eqnarray*}
$\Psi$ is a solution of the Dirac equation and we have :
\begin{eqnarray}
\label{VII.16b}
||\Psi||_{{\cal H}^1_T}\le ||g||_{\tilde{H}^1},\, ||\Psi||_{(H^1(K_T))^4}\lesssim ||g||_{\tilde{H}^1}.
\end{eqnarray}
We want to check that 
\begin{eqnarray*}
\Psi_{2,3}(-\hat{r},\hat{r},\omega)&=&g_{2,3}(\hat{r},\omega)\quad \forall -T\le \hat{r}\le 0,\\
\Psi_{1,4}(\hat{r},\hat{r},\omega)&=&\hat{g}_{1,4}(\hat{r},\omega)\quad \forall 0\le \hat{r}\le T.
\end{eqnarray*}
In fact we can even show :
\[ \Psi(|\hat{r}|,\hat{r},\omega)=\tilde{g}(\hat{r},\omega). \]
As in Section \ref{sec8.2} we estimate :
\begin{eqnarray*}
\lefteqn{\int_{-T}^T\int_{S^2}|\tilde{g}(\hat{r},\omega)-\Psi^{\lambda}(|\hat{r}|,\hat{r},\omega)|^2d\hat{r}d\omega}\\
&=&\int_{-T}^T\int_{S^2}|\Psi^{\lambda}(\lambda|\hat{r}|,\hat{r},\omega)-\Psi^{\lambda}(|\hat{r}|,\hat{r},\omega)|^2d\hat{r}d\omega\\
&=&\int_{-T}^T\int_{S^2}\left|\int_{|\hat{r}|}^{\lambda|\hat{r}|}\partial_t\Psi^{\lambda}(t,\hat{r},\omega)dt\right|^2d\hat{r}d\omega dt\\
&\le&|\lambda-1|T\int_{-T}^T\int_{S^2}\int_{\lambda |\hat{r}|}^{|\hat{r}|}|\partial_t\Psi^{\lambda}(t,\hat{r},\omega)|^2d\hat{r}d\omega dt\\
&\lesssim&T|\lambda-1|||H\Psi^{\lambda}||_{{\cal H}_T}^2
\lesssim T|\lambda-1|\rightarrow 0.
\end{eqnarray*}
On the other hand :
\begin{eqnarray*}
\int_{-T}^T\int_{S^2}|\Psi^{\lambda}(|\hat{r}|,\hat{r},\omega)-\Psi(|\hat{r}|,\hat{r},\omega)|^2d\hat{r}d\omega
\le||\Psi^{\lambda}-\Psi||^2_{(H^{1/2}(K_T))^4}\rightarrow 0.
\end{eqnarray*}
Thus $\Psi(|\hat{r}|,\hat{r},\omega)=\tilde{g}(\hat{r},\omega)$. If $(g_{2,3},\hat{g}_{1,4})\in \tilde{H}^1,$
we approach it by a sequence $(g^n_{2,3},\hat{g}_{1,4}^n)$ of smooth data and the
corresponding solutions converge to a solution $\Psi$. The trace of $\Psi$ on $\Lambda_T$
exists and we show as in the proof of Theorem \ref{thVII.3} that 
\[ \Psi_{2,3}(-\hat{r},\hat{r},\omega)=g_{2,3}(\hat{r},\omega),\, \Psi_{1,4}(\hat{r},\omega)=g_{1,4}(\hat{r},\omega). \]
If $(g_{2,3},\hat{g}_{1,4})\in L^2_{T,-}\oplus L^2_{T,+}$ we again approach it by a sequence of
smooth data $(g^n_{2,3},\hat{g}^n_{1,4})$. The corresponding solutions are in ${\cal H}^1_T$ and converge to some $\Psi$ in
${\cal H}_T$. By definition of the extension we have ${\cal T}\Psi=g$ (see Remark \ref{RVII.0}).
This concludes the proofs of Theorem \ref{thVII.1} and Theorem \ref{thVII.2}.
\qed
\section{The characteristic Cauchy problem on ${\cal M}_{col}$}
\label{sec8.4}
\begin{figure}
\centering\epsfig{figure=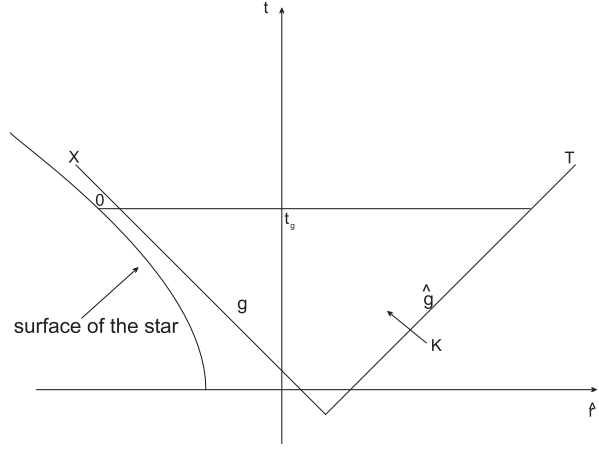,width=8cm}
\caption{The characteristic Cauchy problem outside the collapsing star.}
\label{Chpbcol}
\end{figure}

In this section we want to solve a characteristic Cauchy problem outside the collapsing star. 
The data will be a function $g_{2,3}(t,\omega)$ for which we suppose :
\begin{eqnarray}
\label{VII.29}
\exists t_g\quad \forall t\ge t_g\quad g_{2,3}(t,\omega)=0.
\end{eqnarray}
We want to solve the following characteristic Cauchy problem :
\begin{eqnarray}
\label{VII.28a}
\left.\begin{array}{rcl} \partial_t\Psi&=&iH\Psi\quad \hat{z}(t,\theta)\le \hat{r}\le -t+1,\, t\ge 0,\\
\Psi_{2,3}(t,-t+1,\omega)&=&g_{2,3}(t,\omega),\\
\sum_{\mu\in \{t,\hat{r},\theta,\varphi\}}({\cal N}_{\mu}\hat{\gamma}^{\mu}\Psi)(\hat{z}(t,\theta),\omega)&=&-i\Psi(\hat{z}(t,\theta),\omega),\\
t>t_g,\, \hat{r}\in{[}\hat{z}(t,\theta),-t+1{]}&\Rightarrow&\Psi(t,\hat{r},\omega)=0 \end{array}\right\}
\end{eqnarray}
and write the solution as
\[ \Psi(t)=U(t,t_g)\Psi_K(t_g), \]
where $\Psi_K(t_g)$ is the solution at time $t_g$ of a characteristic problem in $K$,
$K$ as in Figure \ref{Chpbcol}.
We first have to specify the regularity of the data. If 
$g_{1,4}(t,\omega)=\Psi_{1,4}(t,-t+1,\omega)$, then $g_{1,4}$ is solution of the equation 
\begin{eqnarray}
\label{VII.30}
\left.\begin{array}{rcl}
D_tg_{1,4}&=&((P_{\omega}+W)g)_{1,4},\\
g_{1,4}(t_g,\omega)&=&0, \end{array} \right\}
\end{eqnarray}
which we can write as
\begin{eqnarray*}
g_{1,4}(t,\omega)=\int_{t_g}^tV(t,s)S(s)ds
\end{eqnarray*}
with a propagator $V(t,s)$ and a source term $S(s)$ associated to (\ref{VII.30})
as in section \ref{sec8.1}. We then put 
\begin{eqnarray}
\label{VII.31}
g^H_{2,3}(t,\omega):=\frac{1}{2}(D_tg_{2,3}+((P_{\omega}+W)g)_{2,3}).
\end{eqnarray}
Let $\hat{H}^1({[}0,t_g{]}\times S^2;\C^2)$ be the completion of 
\[{\cal C}_0=\{u\in C^{\infty}({[}0,t_g{]}\times S^2;\C^2); u(t_g,\omega)=0\,\forall \omega\in S^2\}\]
in the norm :
\[ ||g_{2,3}||^2_{\hat{H}^1}=2(||g_{2,3}||^2_{L^2({[}0,t_g{]}\times S^2;\C^2)}
+||g^H_{2,3}||^2_{L^2({[}0,t_g{]}\times S^2;\C^2)}). \]
The result of this section is :
\begin{theorem}
\label{thVII.4}
Let $g_{2,3}\in \hat{H}^1$. Then (\ref{VII.28a}) possesses a unique solution $\Psi$ with
$\tilde{\Psi}\in C^1(\R;{\cal H})\cap C(\R;{\cal H}^1)$ and 
\[t\ge 0,\, \hat{r}\in{[}\hat{z}(t,\theta),-t+1{]}\Rightarrow \Psi(t,\hat{r},\omega)=\tilde{\Psi}(t,\hat{r},\omega). \]
Furthermore we have :
\begin{eqnarray}
\label{ECP}
\forall 0\le t_0\le t_g\quad \int_{S^2}\int_{\hat{z}(t_0,\theta)}^{-t_0+1}|\Psi|^2(t_0,\hat{r},\omega)d\hat{r}d\omega
=2\int_{t_0}^{t_g}\int_{S^2}|g_{2,3}|^2(t,\omega)dtd\omega. 
\end{eqnarray}
\end{theorem}
{\bf Proof.}
\begin{itemize}
\item We first show uniqueness. Let $0\le t_0\le t_g$,
\begin{eqnarray*}
B_1&:=&\{(t_0,\hat{r},\omega);\, \hat{z}(t_0,\omega)\le \hat{r}\le -t_0+1\},\\
B_2&:=&\{(t_g,\hat{r},\omega);\, \hat{z}(t_g,\theta)\le \hat{r}\le -t_g+1\},\\
B_3&:=&\{(t,\hat{z}(t,\theta),\omega);\, t_0\le t\le t_g,\, \omega\in S^2\},\\
B_4&:=&\{(t,-t+1,\omega);\, t_0\le t\le t_g,\, \omega\in S^2\},\, B:=B_1\cup B_2\cup B_3\cup B_4.
\end{eqnarray*}
By Stokes' theorem we have :
\begin{eqnarray}
\label{VII.32}
\lefteqn{\int_B*(\phi_A\bar{\phi}_{A'}dx^{AA'}+\bar{\chi}_A\chi_{A'}dx^{AA'})=0}\nonumber\\
&\Leftrightarrow&\int_{S^2}\int_{\hat{z}(t_0,\theta)}^{-t_0+1}|\Psi|^2(t_0,\hat{r},\omega)d\hat{r}d\omega
=2\int_{t_0}^{t_g}\int_{S^2}|g_{2,3}|^2(t,\omega)dtd\omega.
\end{eqnarray}
Indeed 
\[ \int_{B_2}*(\phi_A\bar{\phi}_{A'}dx^{AA'}+\bar{\chi}_A\chi_{A'}dx^{AA'})=0 \]
because the solution is zero on $B_2$ and 
\[ \int_{B_3}*(\phi_A\bar{\phi}_{A'}dx^{AA'}+\bar{\chi}_A\chi_{A'}dx^{AA'})=0 \]
because ${\cal N}^{AA'}(\phi_A\bar{\phi}_{A'}+\bar{\chi}_A\chi_{A'})=0.$
The equality (\ref{VII.32}) gives the uniqueness result.
\item Let us now prove existence. Let
\[ K:=\{(t,\hat{r},\omega);\, -1\le t\le t_g,\, -t+1\le \hat{r}\le t+3\}. \]
We put 
\begin{eqnarray}
\label{VII.33b}
\hat{g}_{1,4}(t,\omega)=\hat{g}_{1,4}(0,\omega)=g_{1,4}(0,\omega)=\int_{t_g}^0V(0,s)S(s)ds. 
\end{eqnarray}
The characteristic Cauchy problem 
\begin{eqnarray}
\left. \begin{array}{rcl} \partial_t\Psi&=&iH\Psi\quad (t,\hat{r},\omega)\in K,\\
\Psi_{2,3}(t,-t+1,\omega)&=&g_{2,3}(t,\omega),\\
\Psi_{1,4}(t,t+3,\omega)&=&\hat{g}_{1,4}(t,\omega) \end{array}\right\}
\end{eqnarray}
has, by the results of the previous sections, a unique solution $\Psi_{\hat{g}}$.
We now choose a smooth cut-off $\chi\in C^{\infty}(\R)$ with 
\begin{eqnarray*}
\chi(x)=\left\{ \begin{array}{cc} 1 & x\le t_g+2 \\ 0 & x\ge t_g+\frac{5}{2} \end{array} \right.
\end{eqnarray*}
and put 
\begin{eqnarray}
\label{8.53}
\Psi(t_g,\hat{r},\omega)=\left\{\begin{array}{cc} \chi\Psi_{\hat{g}}(t_g,\hat{r},\omega)& \hat{r}\ge -t_g+1,\\
0 & \hat{r}\le -t_g+1. \end{array} \right.
\end{eqnarray}
This defines $\Psi(t_g,\hat{r},\omega)$ for all $\hat{r}\ge \hat{z}(t_g,\theta)$. Note
that because of the finite propagation speed the solution in the domain we are interested
in is independent of $\Psi(t_g,\hat{r},\omega)|_{{[}t_g+2,\infty)\times S^2}.$
Because of (\ref{VII.33b}) we have $\Psi_{\hat{g}}(t_g,-t_g+1,\omega)=0$. Therefore $\Psi(t_g,\hat{r},\omega)\in {\cal H}^1.$
The restriction of 
\begin{eqnarray}
\label{8.53a}
\tilde{\Psi}(t)={[}U(t,t_g)\Psi(t_g){]}_H
\end{eqnarray}
to $\{(t,\hat{r},\omega);\, 0\le t\le t_g,\, \hat{z}(t,\theta)\le \hat{r}\le -t+1,\, \omega\in S^2\}$
solves the problem.
\end{itemize}
\qed
\begin{remark}
\label{remVII.1}
(a) We could of course permit data which do not vanish on $B_2$ and proceed as in the preceding sections.
However in the next sections, we shall need a description of the solutions as in (\ref{8.53}).

(b) Let $\Psi\in C^1(\R;{\cal H})\cap C(\R;{\cal H}^1)$ be a solution of the following characteristic 
problem :
\begin{eqnarray*}
\left. \begin{array}{rcl} \partial_t\Psi&=&iH\Psi,\\
\Psi_{2,3}(t,-t+1,\omega)&=&g_{2,3}(t,\omega). \end{array} \right\}
\end{eqnarray*}
Then we have the following energy estimate :
\begin{eqnarray}
\label{R7.1}
2\int_{\R}\int_{S^2}|g_{2,3}|^2(s,\omega)dsd\omega\le ||\Psi(t)||^2_{\cal H}\quad\forall t. 
\end{eqnarray}
Indeed Stokes' theorem gives for $T>0$ :
\begin{eqnarray*}
\lefteqn{\int_1^{2T+1}\int_{S^2}|\Psi|^2(0,\hat{r},\omega)d\hat{r}d\omega}\\
&=&2\left(\int_{-T}^0\int_{S^2}|g_{2,3}|^2(t,\omega)dtd\omega+\int_{-T}^0\int_{S^2}
|\Psi_{1,4}|^2(t,t+2T+1,\omega)dtd\omega\right)\\
&\Rightarrow& 2\int_{-\infty}^0\int_{S^2}|g_{2,3}|^2(t,\omega)dtd\omega\le\int_1^{\infty}\int_{S^2}|\Psi|^2(0,\hat{r},\omega)d\hat{r}d\omega.
\end{eqnarray*}
In the same way we can show :
\[ 2\int_0^{\infty}\int_{S^2}|g_{2,3}|^2(t,\omega)dtd\omega\le \int_{-\infty}^1\int_{S^2}|\Psi|^2(0,\hat{r},\omega)d\hat{r}d\omega. \]
Thanks to inequality (\ref{R7.1}) we can extend the trace operator $T : \Psi_T\mapsto \Psi(t,-t+1,\omega)$ to
a bounded operator $T\in {\cal L}({\cal H};L^2((\R\times S^2dtd\omega);\C^2))$.
\end{remark}
\chapter{Reductions}
\label{Red}
In this chapter we present the basic analytic problem that we have to solve in order
to prove the main theorem.
\section{The Key theorem}
\label{sec5.2}
Let $n\in \Z+1/2,\, \tilde{\Sigma}=\R\times {[}0,2\pi{]}_{\varphi}\times{[}0,\pi{]}_{\theta},
\,\eta^n=\frac{qQr_+}{r_+^2+a^2}+\frac{an}{r_+^2+a^2},\, \mu^n=e^{\sigma\eta^n}$
(see (\ref{mueta}) and Remark \ref{smooth}).
Theorem \ref{th5.1} will follow from
\begin{theorem}[Key theorem]
\label{th5.2}
Let $f(\hat{r},\omega)=e^{in\varphi}f^n(\hat{r},\theta)\in (C_0^{\infty}(\tilde{\Sigma}))^4$. Then:
\begin{eqnarray}
\label{5.2}
\lim_{T\rightarrow\infty}||{\bf 1}_{{[}0,\infty)}(H_0)U(0,T)f||^2_0&=&||{\bf 1}_{{[}0,\infty)}(H){\bf 1}_{\R^-}(P^-)f||^2\nonumber\\
&+&\langle \Omega_{\leftarrow}^{-}f,\mu^n e^{\sigma H_{\leftarrow}}(1+\mu^n e^{\sigma H_{\leftarrow}})^{-1}
\Omega_{\leftarrow}^{-}f\rangle.
\end{eqnarray}
\end{theorem}
{\bf Proof of Theorem \ref{th5.1} (using the result of Theorem \ref{th5.2})}

Using the axial symmetry of the problem it is clear that it is sufficient to show (\ref{5.0})
for $\Phi_j(t,\hat{r},\omega)=e^{in\varphi}\tilde{\Phi}_j(t,\hat{r},\theta)$ with $\eta$
replaced by $\eta^n$ and $\mu$ by $\mu^n$. We then use the polarization identity to
see that it is sufficient to evaluate for $\Phi(t,\hat{r},\omega)=e^{in\varphi}\tilde{\Phi}(t,\hat{r},\theta)\in (C_0^{\infty}({\cal M}_{col}))^4$:
\begin{eqnarray}
\label{5.3}
\lim_{T\rightarrow\infty}\omega_{col}(\Psi^*_{col}(\Phi^T)\Psi_{col}(\Phi^T))&=&\lim_{T\rightarrow\infty}
||{\bf 1}_{{[}0,\infty)}(H_0)S_{col}\Phi^T||_0^2\nonumber\\
&=&\lim_{T\rightarrow\infty}||{\bf 1}_{{[}0,\infty)}(H_0)U(0,T)S\Phi||_0^2.
\end{eqnarray}
Here we have used that for $T>0$ sufficiently large we have
\[S_{col}\Phi^T=U(0,T)S\Phi.\]
By a change of variables and using the compact support of $\Phi$ we have indeed :
\[ S_{col}\Phi^T=\int_a^bU(0,T)U(T,s+T)\Phi(s)ds,\quad -\infty<a<b<\infty.\]
There exist $\hat{r}_1,\, \hat{r}_2$ such that 
\[ \forall a\le s\le b \quad\supp e^{-isH}\Phi(s,\hat{r},\omega)\subset {[}\hat{r}_1,\hat{r}_2{]}\times S^2. \]
In order to replace $U(T,s+T)$ by $e^{-isH}$ it is sufficient to choose $T$ large enough such that 
\[\hat{z}(t,\theta)<\hat{r}_1\quad \forall t\ge\min\{a+T,T\},\,\theta\in {[}0,\pi{]}.\]  
Now using (\ref{5.3}) and Theorem \ref{th5.2} we obtain:
\begin{eqnarray*}
\lim_{T\rightarrow\infty}\omega_{col}(\Psi^*_{col}(\Phi^T)\Psi_{col}(\Phi^T))
&=&||{\bf 1}_{{[}0,\infty)}(H){\bf 1}_{\R^-}(P^-)S\Phi||^2\\
&+&\langle \Omega_{\leftarrow}^-S\Phi, \mu^n e^{\sigma H_{\leftarrow}}
(1+\mu^n e^{\sigma H_{\leftarrow}})^{-1}\Omega_{\leftarrow}^-S\Phi\rangle\\
&=&||{\bf 1}_{{[}0,\infty)}(H)S{\bf 1}_{\R^-}(P^-)\Phi||^2\\
&+&\langle S {\bf 1}_{\R^+}(P^-)\Phi, \mu^n e^{\sigma H}
(1+\mu^n e^{\sigma H})^{-1}S{\bf 1}_{\R^+}(P^-)\Phi\rangle.
\end{eqnarray*}
\qed

From now on we will always suppose $f(\hat{r},\omega)=e^{in\varphi}f^n(\hat{r},\theta).$
The proof of Theorem \ref{th5.2} resp. Theorem \ref{thkey2} below will be the purpose of the rest of this paper. It will be accomplished in Chapter \ref{mainth}.
\section{Fixing the angular momentum}
\label{sec6.1}
Thanks to the cylindrical symmetry of the Kerr-Newman space-time the angular momentum
of the solution is preserved. More precisely let for $n\in \Z+1/2$~:
\begin{eqnarray*}
{\cal H}^n:=\{e^{in\varphi}u\in {\cal H};\, u\in (L^2(\R\times {[}0,\pi{]};d\hat{r}sin\theta d\theta))^4\},\\ 
{\cal H}_*^n:=\{e^{in\varphi}u\in {\cal H}_*;\, u\in (L^2(\R\times {[}0,\pi{]};dr_*sin\theta d\theta))^4\}. 
\end{eqnarray*}
Then all the dynamics which were introduced so far preserve these spaces. Note that
\[ {\cal H}=\oplus_n{\cal H}^n. \]
We also define :
\[ {\cal H}^{n1}:={\cal H}^1\cap{\cal H}^n,\, {\cal H}^{n1}_*:={\cal H}^1_*\cap{\cal H}^n_*. \]
Let us put 
\begin{eqnarray}
\label{6.1}
{\notD}^n_{\gs}&:=&\Gamma^1 D_{r_*}+a_0(r_*){\notD}_{S^2}+b_0(r_*)\Gamma^4+c^n,\\
\label{6.2}
{\notD}^n&:=&h{\notD}_{\gs}^nh+V^n,\\
\label{6.3}
{\notD}_H^n&:=&\Gamma^1D_{r_*}-\frac{a}{r_+^2+a^2}n-\frac{qQr_+}{r_+^2+a^2},\\
\label{6.6}
{\notD}_{\rightarrow}^n&:=&\Gamma^1D_{r_*}+m\Gamma^4,\\
\label{6.7}
c^n&:=&nc_2^{\varphi}+c_1,\, V^n:=nV_{\varphi}+V_1,\, \eta^n:=\frac{an}{r_+^2+a^2}+\frac{qQr_+}{r_+^2+a^2},\\
V_0^n&:=&nV_{\varphi}+V_1=(V_{0ij}^n)_{ij},\\
W^n&:=&h^2c_1+h^2c_2^{\varphi}n+\hat{V}_{\varphi}n+\hat{V}_1.
\end{eqnarray}
The operators ${\notD}^n$ etc. will be understood as operators acting on ${\cal H}^n_*$ 
with domain $D({\notD}^n)=\{u\in {\cal H}^n_*;\, {\notD}^nu\in {\cal H}^n_*\}$ etc. 
They are selfadjoint with these domains (see \cite[Corollary 3.1]{Da2}), the graph norms of ${\notD}^n$
and ${\notD}_{\gs}^n$ are equivalent. We define the operators $\hat{\notD}^n,H^n$ by
\[ \hat{\notD}^n={\cal U}\notD^n{\cal U}^*,\, H^n={\cal V}\hat{\notD}^n{\cal V}^*. \]
Clearly $\hat{\notD}^n,\ H^n$ are selfadjoint with domains :
\[ D(\hat{\notD}^n)={\cal U}D({\notD}^n),\, D(H^n)={\cal V}D(\hat{\notD}^n)={\cal V}{\cal U}D(\notD^n). \]
We will also need the operator :
\begin{eqnarray*}
H_{\leftarrow}^n&:=&\Gamma^1D_{\hat{r}}-\eta^n.
\end{eqnarray*}
In order to describe the precise asymptotic behavior of all coefficients we
introduce the following symbol classes as subsets of $C^{\infty}(\Sigma)$:
\begin{eqnarray*}
f&\in& {\bf S}^{m,n}\quad\mbox{iff}\quad \forall \alpha,\beta\in
\N\quad \partial^{\alpha}_{r_*}\partial^{\beta}_{\theta}f\in
\left\{\begin{array}{cc} O(<r_*>^{m-\alpha}) & r_*\rightarrow
    +\infty \, , \\ O(e^{n\kappa_+ |r_*|}) & r_*\rightarrow - \infty \, ,
  \end{array} \right.\\
f&\in& {\bf S}^m\quad\mbox{iff}\quad \forall \alpha,\beta \in \N\quad
\partial^{\alpha}_{r_*}\partial^{\beta}_{\theta}f\in
O(<r_*>^{m-\alpha}),\, |r_*|\rightarrow\infty.
\end{eqnarray*}
We shall denote $f'$ the derivative of $f$ with respect to $r_*$ even for functions which depend also on $\omega$.
We recall \cite[Proposition 3.1, Lemma 3.2]{Da2} :
\begin{proposition}
\label{prop6.1}
We have 
\[ a_0\in {\bf S}^{-1,-1},\, b_0\in {\bf S}^{0,-1},\, b_0'\in{\bf S}^{-2,-1},\,
(c^n)'\in {\bf S}^{-2,-2},\, h^2-1\in {\bf S}^{-2,-2},\, V_{0ij}^n\in {\bf S}^{-2,-1}.\]
Furthermore there exist two constants $C_3,\epsilon>0$ s.t.
\begin{eqnarray}
\label{6.8}
(a_0(r_*)-\frac{1}{r_*})^{(i)}&\in&{\cal O}(\langle r_*\rangle^{-1-\epsilon-i}),\, r_*\rightarrow\infty\quad i=1,2\\
\label{6.9}
(a_0(r_*)-C_3e^{\kappa_+r_*})^{(i)}&\in&{\cal O}(e^{(\kappa_++\epsilon)r_*}),\, r_*\rightarrow -\infty\quad i=1,2\\
\label{6.10}
b_0-m&\in&{\cal O}(\langle r_*\rangle^{-1}),\, r_*\rightarrow \infty\\
\label{6.11}
c^n+\eta^n&\in&{\cal O}(e^{2\kappa_+r_*}),\, r_*\rightarrow -\infty.
\end{eqnarray}
\end{proposition}
\begin{remark}
\label{rem6.2}
$(i)$ Properties (\ref{6.8}), (\ref{6.9}) imply the existence of two constants\\
$R_0>0, C_0>0$ s.t.
\begin{eqnarray*}
\forall r_*\ge R_0\quad \frac{C_0^{-1}}{r_*}\le a_0(r_*)\le \frac{C_0}{r_*},\\
\forall r_*\le -R_0\quad C_0^{-1}e^{\kappa_+r_*}\le a_0(r_*)\le C_0e^{\kappa_+r_*}.
\end{eqnarray*}
$(ii)$ From the definition of $\hat{r}$ it is clear that we obtain equivalent statements
for $\hat{h}(\hat{r},\theta)=h(r_*(\hat{r},\theta),\theta)$ etc. if we define the symbol classes 
with respect to $\hat{r}$.
\end{remark}
\section{The basic problem}
\label{sec6.2}
For $\nu\in \R$ we put
\begin{eqnarray*}
\Gamma^{\nu}&:=&\left(\begin{array}{cc} 0 & \bar{a}_{\nu} \\ a_{\nu} & 0 \end{array} \right),\, 
a_{\nu}=ie^{i\nu}{\bf 1}_2,\\
{\notD}^{\nu,n}&:=&h{\notD}_{\gs}^{\nu,n}h+V^n,\quad 
{\notD}_{\gs}^{\nu,n}:=\Gamma^1D_{r_*}+a_0{\notD}_{S^2}+b_0\Gamma^{\nu}+c^n,\\
V^n&:=&V_1^{\nu}+nV_{\varphi},\,V_1^{\nu}=\mathbb{V}_0+\frac{m\sqrt{\Delta}}{\sigma}
(\rho-\sqrt{r^2+a^2})\Gamma^{\nu}-\frac{qQr}{\sigma^2}(r^2+a^2-\sigma),
\end{eqnarray*}
i.e. ${\notD}^{\nu,n}$ is obtained from ${\notD}^n$ by replacing $\Gamma^4$ by $\Gamma^{\nu}$.
In the same way we define $\notD_{\rightarrow}^{\nu,n}$. The operators $H^{\nu,n}$ and 
$H_{\rightarrow}^{\nu,n}$ are defined by $H^{\nu,n}={\cal V}{\cal U}\notD^{\nu,n}{\cal U}^*{\cal V}^*$,
$H^{\nu,n}_{\rightarrow}={\cal V}{\cal U}\notD_{\rightarrow}^{\nu,n}{\cal U}^*{\cal V}^*$.
We also define 
\[
{\cal H}_t^n:=\left\{u=e^{in\varphi}v\in {\cal H}_t;\, v\in (L^2(\Sigma^{col,\varphi}_t,d\hat{r}sin\theta d \theta))^4\right\}
\]
with $\Sigma^{col,\varphi}_t=\{(\hat{r},\theta)\in\R\times{[}0,\pi{]};\, \hat{r}\ge \hat{z}(t,\theta) \}$.
 
Let us consider the following problem
\begin{eqnarray}
\label{6.13}
\left.\begin{array}{rcl} \partial_t\Phi&=&iH_t^{\nu,n}\Phi,\qquad \hat{r}>\hat{z}(t,\theta),\\
\sum_{\mu\in\{t,\hat{r},\theta,\varphi\}}({\cal N}_{\mu}\hat{\gamma}^{\mu}\Phi)(t,\hat{z}(t,\theta),\omega)&=&-i\Phi(t,\hat{z}(t,\theta),\omega),\\
\Phi(t=s,.)&=&\Phi_s(.)\in D(H_s^{\nu,n}) \end{array}\right\}
\end{eqnarray}
with 
\begin{eqnarray*}
\lefteqn{D(H_s^{\nu,n})}\\
&=&\left\{u\in{\cal H};\, H_s^{\nu,n}u\in {\cal H}^n_s,\, 
\sum_{\mu\in\{t,\hat{r},\theta,\varphi\}}{\cal N}_{\mu}\hat{\gamma}^{\mu}u(t,\hat{z}(t,\theta),\omega)=-iu(t,\hat{z}(t,\theta),\omega)\right\}.
\end{eqnarray*}
\begin{proposition}
\label{prop6.3}
Let $\Psi_s\in D(H^{\nu,n}_s)$. Then there exists a unique solution 
\[{[}\Psi(.){]}_H={[}U^{\nu,n}(.,s)\Psi_s{]}_H\in C^1(\R_t;{\cal H}^n)\cap C(\R_t;{\cal H}^{n1})\]
of (\ref{6.13}) s.t. 
for all $t\in\R$ $\Psi(t)\in D(H^{\nu,n}_t)$. Furthermore we have $||\Psi(t)||=||\Psi_s||$
and $U^{\nu,n}(t,s)$ possesses an extension to an isometric and strongly continuous propagator 
from ${\cal H}^{n}_s$ to ${\cal H}^{n}_t$ s.t. for all $\Phi_s\in D(H_s)$ we have:
\begin{eqnarray*} 
\frac{d}{dt}U^{\nu,n}(t,s)\Phi_s=iH_tU^{\nu,n}(t,s)\Phi_s
\end{eqnarray*}
and if $R>\hat{z}(s,\theta)$ for all $\theta$ we have:
\begin{eqnarray*}
(\hat{r}>R\Rightarrow\Phi_s(\hat{r},\omega)=0)\Rightarrow (\hat{r}>R+|t-s|\Rightarrow (U^{\nu,n}(t,s)\Phi_s)(\hat{r},\omega)=0).
\end{eqnarray*}
\end{proposition}
Proposition \ref{prop6.3} is proven in Appendix \ref{AppA}.
Let us for the moment just note that Proposition \ref{prop6.3} implies Proposition \ref{prop3.3}.
Indeed if we define
\begin{eqnarray}
\label{6.14}
U(t,s)=\oplus_n e^{i\nu/2\gamma^5}U^{\nu,n}(t,s)e^{-i\nu/2\gamma^5},
\end{eqnarray}
then $U(t,s)$ has the required properties.
Let us now consider $H^{\nu,n}_{\eta^n}=H^{\nu,n}+\eta^n,$\\ 
$H_{\eta^n,t}^{\nu,n}:=H^{\nu,n}_t+\eta^n$.
Clearly 
\[{\bf 1}_{{[}0,\infty)}(H_0^{\nu,n})={\bf 1}_{{[}\eta^n,\infty)}(H_{\eta^n,0}^{\nu,n})
\]
and an equivalent equation for $H^{\nu,n}_{\eta^n}$. If $U^{\nu,n}_{\eta^n}(t,s)$ is the evolution system associated to $H_{\eta^n,t}^{\nu,n}$,
then we have the relation
\[
U^{\nu,n}_{\eta^n}(t,s)=e^{i(t-s)\eta^n}U^{\nu,n}(t,s).\]
Let also 
\[
H_{\leftarrow,\eta^n}=H_{\leftarrow}^n+\eta^n,\quad H^{\nu,n}_{\rightarrow,\eta^n}=H_{\rightarrow}^{\nu,n}+\eta^n,\,
W^n_{\eta^n}=W^n+\eta^n.
\]
Thus if $f(\hat{r},\omega)=\sum e^{in\varphi}f^n(\hat{r},\theta)$, then 
\[ ||{\bf 1}_{{[}0,\infty)}(H_0)U(0,T)f||_{\cal H}^2=\sum_n||{\bf 1}_{{[}\eta^n,\infty)}(H^{\nu,n}_{\eta^n,0})U^{\nu,n}_{\eta^n}(0,T)e^{i\nu/2\gamma^5}f^n||_{{\cal H}^n}^2\]
and therefore the key theorem will follow from the following
\begin{theorem}[Key theorem 2]
\label{thkey2}
Let $f(\hat{r},\omega)=e^{in\varphi}f^n(\hat{r},\theta)\in (C_0^{\infty}(\tilde{\Sigma}))^4$. Then:
\begin{eqnarray}
\label{16*}
\lefteqn{\lim_{T\rightarrow\infty}||{\bf 1}_{{[}\eta^n,\infty)}(H^{\nu,n}_{\eta^n,0})U^{\nu,n}_{\eta^n}(0,T)f||^2_0}\nonumber\\
&=&||{\bf 1}_{{[}\eta^n,\infty)}(H^{\nu,n}_{\eta^n})
{\bf 1}_{\R^-}(P^{\nu-}_n)f||^2\nonumber\\
&+&\langle \Omega_{\leftarrow}^{-,\nu,n}f, e^{\sigma{H}^n_{\leftarrow,\eta^n}}(1+e^{\sigma{H}^n_{\leftarrow,\eta^n}})^{-1}
\Omega_{\leftarrow}^{-,\nu,n}f\rangle,
\end{eqnarray}
where $\sigma$ is as in Theorem \ref{th5.1}.
\end{theorem}
Here $\Omega_{\leftarrow}^{-,\nu,n},\,P^{\nu-}_n$ denote~: 
\begin{eqnarray*}
\Omega_{\leftarrow}^{-,\nu,n}&=&s-\lim_{t\rightarrow-\infty}e^{-itH_{\leftarrow}^{\nu,n}}e^{itH^{\nu,n}}{\bf 1}_{\R^+}(P_n^{\nu,-}),\\
P^{\nu,-}_n&=&s-C_{\infty}-\lim_{t\rightarrow-\infty}e^{-itH^{\nu,n}}\frac{\hat{r}}{t}
e^{itH^{\nu,n}}.
\end{eqnarray*}
\section{The mixed problem for the asymptotic dynamics}
\label{sec6.5}
In this section we give an explicit formula for the mixed problem for the asymptotic 
dynamics near the horizon. We consider the following problem :
\begin{eqnarray}
\label{6.26}
\left.\begin{array}{rcl}\partial_t\Phi&=&iH_{\leftarrow,\eta^n,t}\Phi\qquad (\hat{r},\omega)\in \Sigma^{col}_t,\\
\Phi_2(t,\hat{z}(t,\theta),\omega)&=&-\hat{Z}(t,\theta)\Phi_4(t,\hat{z}(t,\theta),\omega)\\
\Phi_3(t,\hat{z}(t,\theta),\omega)&=&\hat{Z}(t,\theta)\Phi_1(t,\hat{z}(t,\theta),\omega)\\
\Phi(t=s,.)&=&\Phi_s(.), \end{array}\right\}
\end{eqnarray}
where $\hat{Z}(t,\theta)=\sqrt{\frac{1+\dot{\hat{z}}(t,\theta)}{1-\dot{\hat{z}}(t,\theta)}}$.
We note that the boundary condition is the MIT condition with $a_0=0$ and put :
\[ \hat{w}_0(t,\theta):=(1-\dot{\hat{z}}^2)^{-1/2}. \]
For $0>x_0>\hat{z}(0,\theta)$ we define $\hat{\tau}(x_0,\theta)$ by :
\[ \hat{z}(\hat{\tau}(x_0,\theta),\theta)+\hat{\tau}(x_0,\theta)=x_0. \]
We obtain:
\begin{eqnarray}
\label{6*}
\hat{\tau}(x_0,\theta)&=&-\frac{1}{2\kappa_+}\ln(-x_0)+\frac{1}{2\kappa_+}\ln(\hat{A}(\theta))+{\cal O}(x_0),\quad x_0\rightarrow0^-,\\
\label{6**}
1+\dot{\hat{z}}(\hat{\tau}(x_0,\theta))&=&-2\kappa_+x_0+{\cal O}(x_0^2),\quad x_0\rightarrow 0^-.
\end{eqnarray}
We denote $U_{\leftarrow}(t,s)$ the isometric propagator associated to (\ref{6.26}).
\begin{lemma}
\label{lem6.9}
For $t\le s$, given $f\in {\cal H}_s,\, u(t)=U_{\leftarrow}(t,s)f$ is given by :
  \begin{eqnarray*}
\lefteqn{\hat{r}>\hat{z}(t,\theta)}\\
&\Rightarrow& u_2(t,\hat{r},\omega)=f_2(\hat{r}-t+s,\omega),
\,u_3(t,\hat{r},\omega)=f_3(\hat{r}-t+s,\omega),\\
\lefteqn{\hat{r}>s+\hat{z}(s,\theta)-t}\\
&\Rightarrow&u_1(t,\hat{r},\omega)=f_1(\hat{r}+t-s,\omega),
\,u_4(t,\hat{r},\omega)=f_4(\hat{r}+t-s,\omega),\\
\lefteqn{\hat{z}(t,\theta)<\hat{r}<s+\hat{z}(s,\theta)-t}\\
&\Rightarrow&u_1(t,\hat{r},\omega)
=\hat{Z}^{-1}(\hat{\tau}(\hat{r}+t,\theta),\theta)f_3(\hat{r}+t+s-2\hat{\tau}(\hat{r}+t,\theta),\omega),\\
&&u_4(t,\hat{r},\omega)
=-\hat{Z}^{-1}(\hat{\tau}(\hat{r}+t,\theta),\theta)f_2(\hat{r}+t+s-2\hat{\tau}(\hat{r}+t,\theta),\omega).
\end{eqnarray*}
\end{lemma}
\section{The new hamiltonians}
In the remaining chapters we consider the operators $H^{\nu,n}_{\eta^n}$ etc. acting
on the Hilbert space ${\cal H}^n$. It is clear that all the results of the preceding chapters
hold also for these operators. We define the angular part $P^{\nu,n}_{\omega}$ of $H^{\nu,n}_{\eta^n}$
by 
\[ P^{\nu,n}_{\omega}:=H^{\nu,n}_{\eta^n}-\Gamma^1D_{\hat{r}}-W^n_{\eta^n}. \]
The indices $\nu,n,\eta^n$ will be suppressed from now on. In particular we have a new
hamiltonian which is slightly different from the hamiltonian considered in Chapters \ref{sec3}-\ref{Red}.
\chapter{Comparison of the dynamics}
\label{sec9}
\begin{figure}
\centering\epsfig{figure=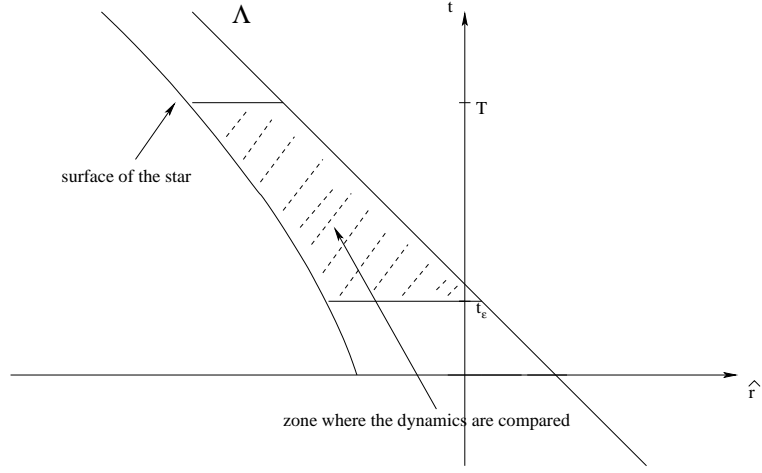,width=10cm}
\caption{Comparison of the dynamics}
\label{prop8}
\end{figure}

Let ${\cal J}\in C_b^{\infty}(\R),\, 0<a<b<1$ and 
\begin{eqnarray*}
{\cal J}(\hat{r})=\left\{ \begin{array}{cc} 1 & \hat{r}\le a \\ 0 & \hat{r}>b \end{array} \right.
\end{eqnarray*}
The aim of this chapter is to prove the following 
\begin{proposition}
\label{prop9.1}
Let $f(\hat{r},\omega)=e^{in\varphi}f^n(\hat{r},\theta)\in (C_0^{\infty}(\tilde{\Sigma}))^4,\, n\in \Z+1/2$.
Then we have 
\begin{eqnarray*}
\forall \epsilon>0\, \exists t_0>0\,\forall t_{\epsilon}\ge t_0\, \exists T_0>0\, \forall T\ge T_0 \, 
||{\cal J}(\hat{r}+t_{\epsilon})(U(t_{\epsilon},T)f-U_{\leftarrow}(t_{\epsilon},T)\Omega_{\leftarrow}^-)f||<\epsilon.
\end{eqnarray*}
\end{proposition}
Proposition \ref{prop9.1} compares the dynamics $U(t_{\epsilon},T)$ and $U_{\leftarrow}(t_{\epsilon},T)\Omega_{\leftarrow}$ (see Figure \ref{prop8}). The function ${\cal J}(\hat{r}+t_{\epsilon})$ is a cut-off in the region we are interested in. The proposition states that in this region the above dynamics are close to each other when $t_{\epsilon}$ becomes large and this uniformly in $T$. To prove Proposition \ref{prop9.1} we understand both $U(t_{\epsilon},T)f$ and
$U_{\leftarrow}(t_{\epsilon},T)\Omega_{\leftarrow}^{-}f$ as solutions of a characteristic problem.
In Section \ref{sec9.1} we compare the characteristic data, in Section \ref{sec9.2}
we compare the solutions of the characteristic problems for the operators $H$ 
and $H_{\leftarrow}$. Proposition \ref{prop9.1} is proven in Section \ref{sec9.3}.
We suppose for the whole chapter that $f(\hat{r},\omega)=e^{in\varphi}f^n(\hat{r},\theta)\in (C_0^{\infty}(\tilde{\Sigma}))^4,\, \supp f\subset{[}R_1,R_2{]}\times {[}0,2\pi{]}\times{[}0,\pi{]}$.
\section{Comparison of the characteristic data}
\label{sec9.1}
Let 
\begin{eqnarray*}
g^T(t,\omega)&:=&(P_{2,3}U(t,T)f)(-t+1,\omega),\\
g^T_{\leftarrow}(t,\omega)&:=&(P_{2,3}U_{\leftarrow}(t,T)\Omega_{\leftarrow}^-f)(-t+1,\omega).
\end{eqnarray*}
Note that it is a-priori not clear that $U_{\leftarrow}(t,T)\Omega_{\leftarrow}^-f$
is regular enough to take the trace 
\[(U_{\leftarrow}(t,T)\Omega_{\leftarrow}^-f)(-t+1,\omega)\]
in the usual sense, but it can be taken in the sense of Remark \ref{remVII.1} (b).
\begin{lemma}
\label{lem9.1}
We have :
\[ \int_0^{\infty}\int_{S^2}|g^T(t,\omega)-g_{\leftarrow}^T(t,\omega)|^2dt d\omega \rightarrow 0,\, T\rightarrow \infty. \]
\end{lemma}
{\bf Proof.}

First observe that :
\begin{eqnarray}
\label{9.1}
g^T(t,\omega)&=&P_{2,3}(e^{i(t-T)H}f)(-t+1,\omega),\\
\label{9.2}
g^T_{\leftarrow}(t,\omega)&=&P_{2,3}(e^{i(t-T)H_{\leftarrow}}\Omega_{\leftarrow}^-f)(-t+1,\omega).
\end{eqnarray}
Using Lemma \ref{lem7.7} we see that 
\begin{eqnarray}
\label{9.3}
\supp f,\supp\Omega_{\leftarrow}^-f\subset {[}R_1,\infty)\times {[}0,2\pi{]}\times{[}0,\pi{]}.
\end{eqnarray}
By the finite propagation speed this entails :
\begin{eqnarray*}
\supp \left(e^{i(t-T)H}f\right),\supp \left(e^{i(t-T)H_{\leftarrow}}\Omega_{\leftarrow}^-f\right)\subset{[}R_1-|T-t|,\infty)\times {[}0,2\pi{]}\times{[}0,\pi{]}.
\end{eqnarray*}
If $t>T$, the condition $-t+1\ge R_1-|T-t|$ implies $1\ge R_1+T$ and if $t<T$ the same condition
implies $t\le \frac{1+T-R_1}{2}$. Let $m(T)$ satisfy the conditions of Lemma \ref{lem7.2}. 
For $T$ sufficiently large we have :
\begin{eqnarray*}
I&:=&\int_0^{\infty}\int_{S^2}\left|g^T(t,\omega)-g_{\leftarrow}^T(t,\omega)\right|^2dtd\omega\\
&\lesssim&\int_0^{\frac{1+T-R_1}{2}}\int_{S^2}\left|P_{2,3}e^{itH}{\bf 1}_{(-1,\infty)}\left(\frac{\hat{r}}{T-m(T)}\right)e^{-iTH}f\right|^2(-t+1,\omega)dtd\omega\\
&+&\int_0^{\frac{1+T-R_1}{2}}\int_{S^2}\left|P_{2,3}e^{itH_{\leftarrow}}{\bf 1}_{(-1,\infty)}\left(\frac{\hat{r}}{T-m(T)}\right)e^{-iTH_{\leftarrow}}\Omega_{\leftarrow}^-f\right|^2(-t+1,\omega)dtd\omega\\
&+&\int_0^{\frac{1+T-R_1}{2}}\int_{S^2}\left|P_{2,3}\left(e^{itH}{\bf 1}_{(-\infty,-1)}\left(\frac{\hat{r}}{T-m(T)}\right)e^{-iTH}\right.\right.\\
&-&\left.\left.e^{itH_{\leftarrow}}{\bf 1}_{(-\infty,-1)}\left(\frac{\hat{r}}{T-m(T)}\right)e^{-iTH_{\leftarrow}}\Omega_{\leftarrow}^-\right)
f\right|^2(-t+1,\omega)dtd\omega\\
&=:&I_1+I_2+I_3.
\end{eqnarray*}
We want to show that 
\begin{eqnarray}
\label{VII.1z}
\lim_{T\rightarrow\infty}I_1=0.
\end{eqnarray}
By the energy estimate (\ref{R7.1}) we see that we can replace $f$ by $\chi(H)f,\,
\chi\in C_0^{\infty}(\R),$\\
$\supp\chi\subset \R\setminus \{-m+\eta^n,m+\eta^n\}$. Let $0<\epsilon_{\chi}<1$ be as in Lemma
\ref{lem7.1a}. Then
\begin{eqnarray*}
I_1&\lesssim&\int_0^{\frac{1+T-R_1}{2}}\int_{S^2}\left|P_{2,3}e^{itH}{\bf 1}_{(-1,\frac{\epsilon_{\chi}}{2})}\left(\frac{\hat{r}}{T-m(T)}\right)
e^{-iTH}\chi(H)f\right|^2(-t+1,\omega)dtd\omega\\
&+&\int_0^{\frac{1+T-R_1}{2}}\int_{S^2}\left|P_{2,3}e^{itH_{\leftarrow}}{\bf 1}_{(\frac{\epsilon_{\chi}}{2},\infty)}\left(\frac{\hat{r}}{T-m(T)}\right)
e^{-iTH}\chi(H)f\right|^2(-t+1,\omega)dtd\omega\\
&=:&I_{11}+I_{12}
\end{eqnarray*}
Let us first estimate $I_{12}$. We have 
\begin{eqnarray*}
\lefteqn{\supp \left(e^{itH}{\bf 1}_{(\frac{\epsilon_{\chi}}{2},\infty)}\left(\frac{\hat{r}}{T-m(T)}\right)e^{-iTH}\chi(H)f\right)}\\
&\subset& \left(\frac{\epsilon_{\chi}}{2}\left(T-m(T)\right)-t,\infty\right)\times {[}0,2\pi{]}\times{[}0,\pi{]}.
\end{eqnarray*}
But $-t+1\ge \frac{\epsilon_{\chi}}{2}\left(T-m(T)\right)-t$ implies $1\ge\frac{\epsilon_{\chi}}{2}\left(T-m(T)\right).$
Thus $I_{12}=0$ for $T$ sufficiently large. We now estimate $I_{11}.$ By the energy estimate (\ref{R7.1}) 
we have :
\begin{eqnarray*}
I_{11}&\le&||{\bf 1}_{(-1,\frac{\epsilon_{\chi}}{2})}\left(\frac{\hat{r}}{T-m(T)}\right)e^{-iTH}{\bf 1}_{\R^+}(P^-)\chi(H)f||\\
&+&||{\bf 1}_{(-1,\frac{\epsilon_{\chi}}{2})}\left(\frac{\hat{r}}{T-m(T)}\right)e^{-iTH}{\bf 1}_{\R^-}(P^-)\chi(H)f||\\
&=&I_{11}^a+I_{11}^b.
\end{eqnarray*}
We have :
\[ \lim_{T\rightarrow\infty}I_{11}^b=\lim_{T\rightarrow\infty}||{\bf 1}_{(-1,\frac{\epsilon_{\chi}}{2})}
\left(\frac{\hat{r}}{T-m(T)}\right){\bf 1}_{(\frac{3\epsilon_{\chi}}{4},\infty)}\left(\frac{\hat{r}}{T}\right)
e^{-iTH}\chi(H)f||^2=0. \]
By Lemma \ref{lem7.2} we have $\lim_{T\rightarrow\infty}I_{11}^a=0.$
(\ref{VII.1z}) follows. In the same manner we can show :
\[ \lim_{T\rightarrow\infty}I_2=0. \]
Let us now estimate $I_3$. We have ($t\ge 0$) : 
{\small
\begin{eqnarray*}
\lefteqn{\supp \left(e^{itH_{\leftarrow}}{\bf 1}_{(-\infty,-1)}\left(\frac{\hat{r}}{T-m(T)}\right)e^{-iTH_{\leftarrow}}\Omega_{\leftarrow}^-f\right),}\\
&&\supp \left(e^{itH}{\bf 1}_{(-\infty,-1)}\left(\frac{\hat{r}}{T-m(T)}\right)e^{-iTH}f\right)\subset (-\infty, m(T)-T+t)\times {[}0,2\pi{]}\times{[}0,\pi{]}.
\end{eqnarray*}}
The condition $-t+1\le m(T)-T+t$ implies $t\ge \frac{T-m(T)+1}{2}$. In particular\\ 
$1-{\cal J}(-t+1)=0$ for $T$ sufficiently large. Therefore we obtain :
{\small
\begin{eqnarray*}
I_3&\le&\int_{\frac{T-m(T)+1}{2}}^{\frac{1+T-R_1}{2}}\int_{S^2}\left|P_{2,3}\left({\cal J}(\hat{r})
e^{itH}{\bf 1}_{(-\infty,-1)}\left(\frac{\hat{r}}{T-m(T)}\right)e^{-iTH}\right.\right.\\
&-&\left.\left.e^{itH_{\leftarrow}}{\bf 1}_{(-\infty,-1)}\left(\frac{\hat{r}}{T-m(T)}\right)
e^{-iTH_{\leftarrow}}\Omega_{\leftarrow}^-\right)f\right|^2(-t+1,\omega)dtd\omega\\
&\lesssim& \int_{\frac{T-m(T)+1}{2}}^{\frac{1+T-R_1}{2}}\int_{S^2}\left|P_{2,3}e^{itH}{\bf 1}_{(-1,\infty)}\left(\frac{\hat{r}}{T-m(T)}\right)
e^{-iTH}f\right|^2(-t+1,\omega)dtd\omega\\
&+&\int_{\frac{T-m(T)+1}{2}}^{\frac{1+T-R_1}{2}}\int_{S^2}\left|P_{2,3}e^{itH_{\leftarrow}}{\bf 1}_{(-1,\infty)}\left(\frac{\hat{r}}{T-m(T)}\right)
e^{-iTH_{\leftarrow}}\Omega_{\leftarrow}^-f\right|^2(-t+1,\omega)dtd\omega\\
&+&\int_{\frac{T-m(T)+1}{2}}^{\frac{1+T-R_1}{2}}\int_{S^2}\left|P_{2,3}\left({\cal J}(\hat{r})e^{i(t-T)H}f-e^{i(t-T)H_{\leftarrow}}\Omega_{\leftarrow}^-\right)f\right|^2(-t+1,\omega)dtd\omega\\
&=:&L_1+L_2+L_3.
\end{eqnarray*}}
We have for $j=1,2$ 
\[ 0\le L_j\le I_j\rightarrow 0,\, T\rightarrow \infty. \]
We estimate :
\begin{eqnarray}
\label{9.6}
L_3&\le&\frac{m(T)-R_1}{2}\sup_{\sigma\le\frac{1-(T+R_1)}{2}}||{\cal J}(\hat{r})e^{i\sigma H}
f-e^{i\sigma H_{\leftarrow}}\Omega_{\leftarrow}^-f||^2_{L^{\infty}(\R;(L^2(S^2))^4)}\nonumber\\
&\lesssim&\frac{m(T)-R_1}{2}\sup_{\sigma\le\frac{1-(T+R_1)}{2}}||{\cal J}(\hat{r})e^{i\sigma H}
f-e^{i\sigma H_{\leftarrow}}\Omega_{\leftarrow}^-f||^2_{H^1(\R;(L^2(S^2))^4)}.
\end{eqnarray}
We can choose 
\[ m(T)=\min \left(\frac{T}{2},\left(\sup_{\sigma\le\frac{1-(T+R_1)}{2}}||{\cal J}(\hat{r})e^{i\sigma H}
f-e^{i\sigma H_{\leftarrow}}\Omega_{\leftarrow}^-f||^2_{H^1(\R;(L^2(S^2))^4)}\right)^{-1/2}\right). \]
Then by Lemma \ref{lem7.8} the R.H.S. of (\ref{9.6}) goes to zero when $T\rightarrow\infty.$
\qed
\section{Comparison with the asymptotic dynamics}
\label{sec9.2}
In the region $\hat{z}(t_{\epsilon},\theta)\le\hat{r}\le a-t_{\epsilon}$ we understand ${\cal J}(\hat{r}+t_{\epsilon})U(t_{\epsilon},T)f$ as the solution of the 
characteristic problem on ${\cal M}_{col}$ for the operator $H$ with data $g^T$ and
\[{\cal J}(\hat{r}+t_{\epsilon})U_{\leftarrow}(t_{\epsilon},T)\Omega_{\leftarrow}^-f\]
as solution of the characteristic problem for the operator $H_{\leftarrow}$ with data $g^T_{\leftarrow}$. We would like to solve 
the characteristic problem for $H$ with data $g^T_{\leftarrow}$ and write the solution as 
\[ G_{t_{\epsilon}}(g_{\leftarrow}^T)=U(t,T/2+c_0)\phi(T/2+c_0), \]
where $\phi$ is constructed as in (\ref{8.53}) with $t_g=T/2+c_0$ for some $c_0>0$. Unfortunately
$g_{\leftarrow}^T$ will in general not be regular enough to assure that $G_{t_{\epsilon}}(g_{\leftarrow}^T)$ 
is a strong solution. We shall therefore regularize $\Omega_{\leftarrow}^-f$. Let $\chi_R\in C_0^{\infty}(\R)$ with
\begin{eqnarray*}
\chi_R=\left\{\begin{array}{cc} 1 & R\ge \hat{r}\ge R_1,\\ 0 & \hat{r}\le R_1-1,\hat{r}\ge R+1. \end{array}
\right.
\end{eqnarray*}
Let 
\[ {\cal I}_n=\{l;\, l-\frac{1}{2}\in \N,\, l-|n|\in \N \},\, {\cal I}^N_n=\{ l\in {\cal I}_n;\, |(l,n)|\le N\}. \]
and
\[ \Omega_{\leftarrow}^-f=\sum_{l\in {\cal I}_n}(\Omega_{\leftarrow}^-f)^{nl},\, (\Omega_{\leftarrow}^-f)^{nl}\in {\cal H}^{nl}
=L^2((\R,d\hat{r});\C^4)\otimes_4 Y_{nl}\, \forall l. \] 
We put 
\[ (\Omega_{\leftarrow}^-f)^N=\sum_{l\in {\cal I}^N_n}(\Omega_{\leftarrow}^-f)^{nl},\quad
(\Omega_{\leftarrow}^-f)^N_R=\chi_R(\Omega_{\leftarrow}^-f)^N. \]
Let 
\[ \tilde{f}:=(\Omega_{\leftarrow}^-f)^N_R. \]
Clearly we have 
\[ \forall \epsilon>0\exists N_{\epsilon}>0,R_{\epsilon}>0\forall R\ge R_{\epsilon},N\ge N_{\epsilon}\quad 
||\Omega_{\leftarrow}^-f-(\Omega_{\leftarrow}^-f)^N_R||<\epsilon. \]
We put 
\[ g_{\leftarrow,R}^{T,N}=(P_{23}U_{\leftarrow}(t,T)(\Omega_{\leftarrow}^-f)^N_R)(t,-t+1,\omega). \]
The functions $(g_{\leftarrow,R}^{T,N})_{23}$ are compactly supported. The necessary regularity of $g_{\leftarrow,R}^{T,N}$ 
follows from the regularity of $e^{i(t-T)H_{\leftarrow}}(\Omega_{\leftarrow}^-f)^N_R$ 
by classical trace theorems \footnote{Recall that $(U_{\leftarrow}(t,T)(\Omega_{\leftarrow}^-f)^N_R)(-t+1,\omega)
=(e^{i(t-T)H_{\leftarrow}}(\Omega_{\leftarrow}^-f)^N_R)(-t+1,\omega)$.}. 
We put $c_0:=\frac{1-R_0}{2},\, R_0=R_1-1$. Let $\Phi^{R,N}(T/2+c_0)$ 
be the solution of the characteristic problem in the whole exterior Kerr-Newman space-time 
with data $g_{\leftarrow,R}^{T,N}$ on 
$\{(t,\hat{r},\omega);\, 0\le t,\, \hat{r}=-t+1\}$ as constructed in (\ref{8.53})
and 
\[ \Phi^{R,N}_{\leftarrow}(T/2+c_0)=e^{-i(T/2-c_0)H_{\leftarrow}}(\Omega_{\leftarrow}^-f)^N_R. \]
\begin{lemma}
\label{lem9.3}
We have uniformly in $t_{\epsilon}\ge 0$ :
\begin{eqnarray*}
||{\cal J}(\hat{r}+t_{\epsilon})U(t_{\epsilon},T/2+c_0)(\Phi^{R,N}(T/2+c_0,.)-\Phi_{\leftarrow}^{R,N}(T/2+c_0,.))||_{{\cal H}_{t_{\epsilon}}}\rightarrow 0,\, T\rightarrow \infty.
\end{eqnarray*}
\end{lemma}
{\bf Proof.}

Let 
\[ I:=||{\cal J}(\hat{r}+t_{\epsilon})U(t_{\epsilon},T/2+c_0)(\Phi^{R,N}(T/2+c_0,.)-\Phi_{\leftarrow}^{R,N}(T/2+c_0,.))||^2_{{\cal H}_{t_{\epsilon}}}. \]
First observe that by (\ref{ECP}) we have :
\[ I\le 2\int_{t_{\epsilon}}^{T/2+c_0}\int_{S^2}|\tilde{g}^{T,N}_R-g_{\leftarrow,R}^{T,N}|^2dtd\omega \]
with 
\[ \tilde{g}^{T,N}_R=(P_{2,3}e^{i(t-(T/2+c_0))H}
e^{-i(T/2-c_0)H_{\leftarrow}}\tilde{f})(-t+1,\omega). \]
We proceed as in the proof of Lemma \ref{lem9.1}. Let $m(T)$ satisfy the conditions of 
Lemma \ref{lem7.2}. Then we have :
{\footnotesize
\begin{eqnarray*}
I&\le&\int_{t_{\epsilon}}^{T/2+c_0}\int_{S^2}\left|P_{2,3}e^{itH}{\bf 1}_{(-1,\infty)}\left(\frac{\hat{r}}{T-m(T)}\right)e^{-i(T/2+c_0)H}e^{-i(T/2-c_0)H_{\leftarrow}}\tilde{f}\right|^2(-t+1,\omega)dtd\omega\\
&+&\int_{t_{\epsilon}}^{T/2+c_0}\int_{S^2}\left|P_{2,3}e^{itH_{\leftarrow}}{\bf 1}_{(-1,\infty)}\left(\frac{\hat{r}}{T-m(T)}\right)e^{-iTH_{\leftarrow}}\tilde{f}\right|^2(-t+1,\omega)dtd\omega\\
&+&\int_{t_{\epsilon}}^{T/2+c_0}\int_{S^2}\left|P_{2,3}\left(e^{itH}{\bf 1}_{(-\infty,-1)}\left(\frac{\hat{r}}{T-m(T)}\right)e^{-i(T/2+c_0)H}e^{-i(T/2-c_0)H_{\leftarrow}}\right.\right.\\
&-&\left.\left.e^{itH_{\leftarrow}}{\bf 1}_{(-\infty,-1)}\left(\frac{\hat{r}}{T-m(T)}\right)e^{-iTH_{\leftarrow}}\right)\tilde{f}\right|^2(-t+1,\omega)dtd\omega\\
&=:&I_1+I_2+I_3.
\end{eqnarray*}}
Let us first estimate $I_1$. We have :
{\footnotesize
\begin{eqnarray*}
I_1&\le&\int_{t_{\epsilon}}^{T/2+c_0}\int_{S^2}\left|P_{2,3}e^{itH}{\bf 1}_{(-1,1)}
\left(\frac{\hat{r}}{T-m(T)}\right)e^{-i\left(\frac{T}{2}+c_0\right)H}e^{-i\left(\frac{T}{2}-c_0\right)H_{\leftarrow}}\tilde{f}\right|^2(-t+1,\omega) dtd\omega\\
&+&\int_{t_{\epsilon}}^{T/2+c_0}\int_{S^2}\left|P_{2,3}e^{itH}{\bf 1}_{(1,\infty)}\left(\frac{\hat{r}}{T-m(T)}\right)
e^{-i\left(\frac{T}{2}+c_0\right)H}e^{-i\left(\frac{T}{2}-c_0\right)H_{\leftarrow}}\tilde{f}\right|^2(-t+1,\omega)dtd\omega\\
&=:&I_{11}+I_{12}.
\end{eqnarray*}}
Using a finite propagation speed argument we see that $I_{12}=0$ for $T$ sufficiently
large.
We estimate :
\begin{eqnarray*}
I_{11}&\le&\left|\left|{\bf 1}_{(-1,1)}\left(\frac{\hat{r}}{T-m(T)}\right)e^{-i\left(\frac{T}{2}+c_0\right)H}e^{-i\left(\frac{T}{2}-c_0\right)H_{\leftarrow}}\tilde{f}\right|\right|^2\\
&\le&\left|\left|{\bf 1}_{(-1,1)}\left(\frac{\hat{r}}{T-m(T)}\right)e^{-iTH}{\bf 1}_{\R^+}(P^-)W_{\leftarrow}^-\tilde{f}\right|\right|^2\\
&+&\left|\left|\left(e^{i\left(\frac{T}{2}-c_0\right)H}e^{-i\left(\frac{T}{2}+c_0\right)H_{\leftarrow}}P_{{\cal H}^+}
-W_{\leftarrow}^-\right)\tilde{f}\right|\right|^2\rightarrow 0,\, T\rightarrow \infty,
\end{eqnarray*} 
where we have used Lemma \ref{lem7.2}. In the same way we can show :
\[ \lim_{T\rightarrow\infty}I_2=0. \]
By the same arguments as in the proof of Lemma \ref{lem9.1} we see that 
{\footnotesize
\begin{eqnarray}
\label{I_3}
I_3&\le&\int_{\frac{T-m(T)+1}{2}}^{T/2+c_0}\int_{S^2}\left|\left(e^{i(t-(T/2+c_0)H}-e^{i(t-(T/2+c_0))H_{\leftarrow}}\right)
e^{-i(T/2-c_0)H_{\leftarrow}}\tilde{f}\right|^2(-t+1,\omega)dtd\omega\nonumber\\
&\lesssim&\left(c_0-\frac{1}{2}+\frac{m(T)}{2}\right)\sup_{\frac{1-m(T)}{2}-c_0\le \sigma\le 0}
||\left(e^{i\sigma H}- e^{i\sigma H_{\leftarrow}}\right)e^{-i(T/2-c_0)H_{\leftarrow}}\tilde{f}||_{H^1(\R,(L^2(S^2))^4)}.\nonumber\\
\end{eqnarray}}
We estimate :
\begin{eqnarray*}
\lefteqn{||(e^{i\sigma H}-e^{i\sigma H_{\leftarrow}})e^{-i(T/2-c_0)H_{\leftarrow}}\tilde{f}||_{H^1(\R,(L^2(S^2))^4)}}\\
&\le&\int_{\sigma}^0||e^{isH}(H-H_{\leftarrow})e^{i(\sigma-s-T/2+c_0)H_{\leftarrow}}\tilde{f}||_{H^1(\R,(L^2(S^2))^4)}ds=:I_{\sigma}.
\end{eqnarray*}
We have :
\begin{eqnarray}
\label{II.1}
I_{\sigma}&\le&\int_{\sigma}^0(||D_{\hat{r}}(P_{\omega}+W)e^{i(\sigma-s-T/2+c_0)H_{\leftarrow}}\tilde{f}||\nonumber\\
&+&||(P_{\omega}+W)e^{i(\sigma-s-T/2+c_0)H_{\leftarrow}}\tilde{f}||)ds.
\end{eqnarray}
We have :
\begin{eqnarray}
\label{II.2}
\lefteqn{\Gamma^1D_{\hat{r}}(P_{\omega}+W)e^{i(\sigma-s-T/2+c_0)H_{\leftarrow}}\tilde{f}}\nonumber\\
&=&{[}\Gamma^1D_{\hat{r}},P_{\omega}+W{]}e^{i(\sigma-s-T/2+c_0)H_{\leftarrow}}\tilde{f}\nonumber\\
&+&(P_{\omega}+W)e^{i(\sigma-s-T/2+c_0)H_{\leftarrow}}\Gamma^1D_{\hat{r}}\tilde{f}.
\end{eqnarray}
The last term equals :
\begin{eqnarray}
\label{II.3}
\lefteqn{(P_{\omega}+W)e^{i(\sigma-s-T/2+c_0)H_{\leftarrow}}\Gamma^1(\chi_R)'(\Omega_{\leftarrow}^-f)^N}\nonumber\\
&+&(P_{\omega}+W)e^{i(\sigma-s-T/2+c_0)H_{\leftarrow}}(\Omega^-_{\leftarrow}Hf)^N_R.
\end{eqnarray}
The first term can be treated in a similar manner.
Using (\ref{II.2}), (\ref{II.3}) as well as 
\[ \supp e^{i(\sigma-s-T/2+c_0)H_{\leftarrow}}\tilde{f}\subset(-\infty,R+\sigma-s-T/2+c_0)\times {[}0,2\pi{]}\times{[}0,\pi{]}. \]
we can push further the estimate (\ref{II.1}) ($0\ge \sigma\ge \frac{1-m(T)}{2}-c_0$) :
\begin{eqnarray*}
I_{\sigma}&\lesssim&\int_{\sigma}^0e^{\kappa_+(R+\sigma-s-T/2)}(N+1)(||\Omega_{\leftarrow}^-f||+||\Omega_{\leftarrow}^-Hf||)ds\\
&\le&\frac{1}{\kappa_+}e^{\kappa_+(R+c_0-\frac{T}{2})}(N+1)(||\Omega_{\leftarrow}^-f||+||\Omega_{\leftarrow}^-Hf||)
=:R(T)\rightarrow0,\, T\rightarrow\infty.
\end{eqnarray*}
Choosing $m(T)=\min(\frac{1}{2}T,R(T)^{-1/2})$ in (\ref{I_3}) we find $I_3\rightarrow 0,\, 
T\rightarrow\infty.$ This concludes the proof of the lemma.
\qed
\section{Proof of Proposition \ref{prop9.1}}
\label{sec9.3}
We start with the following lemma which analyzes the frequencies in $\notD_{S^2}$ and $D_{\hat{r}}$ 
of $U_{\leftarrow}(t,s)f,\, f\in {\cal H}^1$~:
\begin{lemma}
\label{lem9.4}
Let $f\in {\cal H}^1,\, \supp f\subset{[}R_1,R_2{]}\times {[}0,2\pi{]}\times{[}0,\pi{]}$. 
Then we have uniformly in $0\le t\le s$ :
\begin{eqnarray*}
(i)\quad ||\notD_{S^2}U_{\leftarrow}(t,s)f||_{{\cal H}_t}&\le& C(R_1,R_2) ||f||_{{\cal H}^1},\\
(ii)\quad ||D_{\hat{r}}U_{\leftarrow}(t,s)f||_{{\cal H}_t}&\le&C(R_1,R_2)e^{\kappa_+s}||f||_{{\cal H}^1}.
\end{eqnarray*} 

\end{lemma}
{\bf Proof.}

Let $u=U_{\leftarrow}(t,s)f$. Recall the explicit formula for $u$ from Section 
\ref{sec6.5}. We first show $(i)$. We note :
\begin{eqnarray}
\label{9.11}
1-\dot{\hat{z}}(t,\theta))\ge 1,\, 1+\dot{\hat{z}}(t,\theta))\gtrsim e^{-2\kappa_+t}.
\end{eqnarray}
This follows from (\ref{star0b}). We then claim that :
\begin{eqnarray}
\label{9.20}
\left|\frac{\partial\hat{\tau}}{\partial\theta}(\hat{r},\theta)\right|\lesssim 1\quad\mbox{uniformly in}\quad \hat{r},\theta.
\end{eqnarray}
Indeed from 
\begin{eqnarray*}
\hat{z}(\hat{\tau}(\hat{r},\theta),\theta)+\hat{\tau}(\hat{r},\theta)&=&\hat{r}\quad\mbox{follows with (\ref{9.11})}\\
\left|\frac{\partial\hat{\tau}}{\partial\theta}(\hat{r},\theta)\right|=\left|\frac{\frac{\partial\hat{z}}{\partial\theta}}{1+\dot{\hat{z}}}\right|
&\lesssim&\left|e^{2\kappa_+\hat{\tau}(\hat{r},\theta)}\frac{\partial\hat{z}}{\partial\theta}\right|\lesssim 1.
\end{eqnarray*}
It follows $(j=2,3)$ :
\begin{eqnarray}
\label{9.21}
\lefteqn{|\partial_{\theta}f_j(\hat{r}+t+s-2\hat{\tau}(\hat{r}+t,\theta),\omega)|}\nonumber\\
&\lesssim&|(\partial_{\theta }f_j)(\hat{r}+t+s-2\hat{\tau}(\hat{r}+t,\theta),\omega)|\nonumber\\
&+&|(\partial_{\hat{r}}f_j)(\hat{r}+t+s-2\hat{\tau}(\hat{r}+t,\theta),\omega)|.
\end{eqnarray}
We next claim :
\[ \left|\partial_{\theta}\sqrt{\frac{1-\dot{\hat{z}}(\hat{\tau}(\hat{r}+t,\theta),\theta)}
{1+\dot{\hat{z}}(\hat{\tau}(\hat{r}+t,\theta),\theta)}}\right|\lesssim \sqrt{\frac{1-\dot{\hat{z}}(\hat{\tau}(\hat{r}+t,\theta),\theta)}
{1+\dot{\hat{z}}(\hat{\tau}(\hat{r}+t,\theta),\theta)}}. \]
Indeed :
\begin{eqnarray*}
\left|\partial_{\theta}\sqrt{\frac{1-\dot{\hat{z}}(\hat{\tau}(\hat{r}+t,\theta),\theta)}
{1+\dot{\hat{z}}(\hat{\tau}(\hat{r}+t,\theta),\theta)}}\right|&=&\left|\frac{1}
{1-\dot{\hat{z}}}\frac{(\Ddot{\hat{z}}
\frac{\partial\hat{\tau}}{\partial\theta}+\frac{\partial\dot{\hat{z}}}{\partial\theta})}{1+\dot{\hat{z}}}
\sqrt{\frac{1-\dot{\hat{z}}(\hat{\tau}(\hat{r}+t,\theta),\theta)}{1+\dot{\hat{z}}(\hat{\tau}(\hat{r}+t,\theta),\theta)}}\right|\\
&\lesssim&\sqrt{\frac{1-\dot{\hat{z}}(\hat{\tau}(\hat{r}+t,\theta),\theta)}
{1+\dot{\hat{z}}(\hat{\tau}(\hat{r}+t,\theta),\theta)}},
\end{eqnarray*}
where we have used (\ref{9.11}) and the uniform boundedness of $\frac{\Ddot{\hat{z}}}{1+\dot{\hat{z}}},\frac{\partial\hat{\tau}}{\partial\theta},\frac{\frac{\partial\dot{\hat{z}}}{\partial\theta}}{1+\dot{\hat{z}}}$.
It follows :
\begin{eqnarray*}
||\notD_{S^2}U_{\leftarrow}(t,s)f||_{{\cal H}_t}&\lesssim&||U_{\leftarrow}(t,s)(\partial_{\theta}f)||+||U_{\leftarrow}(t,s)\frac{1}{\sin\theta}(\partial_{\varphi}f)||\\
&+&||U_{\leftarrow}(t,s)(\partial_{\hat{r}}f)||+||f||\\
&\le& C(R_1,R_2)||f||_{{\cal H}^1}.
\end{eqnarray*}
Let us now show $(ii)$. We first claim that on $\supp f_j(\hat{r}+t+s-2\hat{\tau}(\hat{r}+t,\theta),\omega)$
we have for $s$ sufficiently large :
\[ \left|\frac{\partial\hat{\tau}}{\partial\hat{r}}(\hat{r}+t,\theta)\right|\lesssim \frac{1}{|\hat{r}+t|}. \]
Indeed from 
\[ \hat{z}(\hat{\tau}(\hat{r}+t,\theta),\theta)+\hat{\tau}(\hat{r}+t,\theta)=\hat{r}+t\]
follows
\[ \left|\frac{\partial\hat{\tau}}{\partial\hat{r}}\right|\le \frac{1}{|\dot{\hat{z}}(\hat{\tau}(\hat{r}+t,\theta),\theta)+1|}\lesssim \frac{1}{|\hat{r}+t|}, \]
where we have used (\ref{6**}) and the fact that 
\[\hat{r}+t\rightarrow 0,s\rightarrow \infty\quad\mbox{on}\quad \supp f_j(\hat{r}+t+s-2\hat{\tau}(\hat{r}+t,\theta),\omega).\]
We next note that on $\supp f_j(\hat{r}+t+s-2\hat{\tau}(\hat{r}+t,\theta),\omega)$
we have :
\[ |\hat{r}+t|\gtrsim e^{-\kappa_+s}. \]
This follows from (\ref{6*}). We now estimate on $\supp f_j(\hat{r}+t+s-2\hat{\tau}(\hat{r}+t,\theta),\omega)$: 
\begin{eqnarray}
\label{8.12w1}
\lefteqn{\left|\partial_{\hat{r}}\left(\left(\frac{1-\dot{\hat{z}}}{1+\dot{\hat{z}}}\right)^{1/2}(\hat{\tau}(\hat{r}+t,\theta),\theta)\right)\right|}\nonumber\\
&=&\left|\left(\frac{1-\dot{\hat{z}}}{1+\dot{\hat{z}}}(\hat{\tau}(\hat{r}+t,\theta),\theta)\right)^{1/2}
\frac{\ddot{\hat{z}}(\hat{\tau}(\hat{r}+t,\theta),\theta)\partial_{\hat{r}}\hat{\tau}(\hat{r}+t,\theta)}{(1+\dot{\hat{z}})(1-\dot{\hat{z}})}\right|\nonumber\\
&\lesssim& e^{\kappa_+ s}\left(\frac{1-\dot{\hat{z}}}{1+\dot{\hat{z}}}\right)^{1/2}(\hat{\tau}(\hat{r}+t,\theta),\theta).
\end{eqnarray}
In the same way we estimate :
\begin{eqnarray}
\label{8.12w2}
\lefteqn{|\partial_{\hat{r}}f_j(\hat{r}+t+s-2\hat{\tau}(\hat{r}+t,\theta),\omega)|}\nonumber\\
&\le&|(\partial_{\hat{r}}f_j)(\hat{r}+t+s-2\hat{\tau}(\hat{r}+t,\theta),\omega)(1-2\partial_{\hat{r}}\hat{\tau})|\nonumber\\
&\lesssim& e^{\kappa_+ s}|(\partial_{\hat{r}}f_j)(\hat{r}+t+s-2\hat{\tau}(\hat{r}+t,\theta),\omega)|.
\end{eqnarray}
The estimates (\ref{8.12w1}) and (\ref{8.12w2}) give $(ii)$.
\qed
\\

\underline{Proof of Proposition \ref{prop9.1}}

Let $\epsilon>0$. We first note that 
\begin{eqnarray*}
{\cal J}(\hat{r}+t_{\epsilon})U(t_{\epsilon},T)f={\cal J}(\hat{r}+t_{\epsilon})G_{t_{\epsilon}}(g^T), 
\end{eqnarray*}
where $G_{t_{\epsilon}}(g^T)$ is the solution at time $t_{\epsilon}$ of the characteristic problem (\ref{VII.28a})
with data $g^T$. In the same way we denote $G_{t_{\epsilon}}(g_{\leftarrow,R}^{T,N})$ the solution 
at time $t_{\epsilon}$ of the characteristic Cauchy problem (\ref{VII.28a}) with data $g_{\leftarrow,R}^{T,N}$.
 We estimate :
\begin{eqnarray}
\label{9.12}
\lefteqn{||{\cal J}(\hat{r}+t_{\epsilon})(U(t_{\epsilon},T)f-U_{\leftarrow}(t_{\epsilon},T)\Omega_{\leftarrow}^-f)||_{{\cal H}_{t_{\epsilon}}}}\nonumber\\
&\le&||{\cal J}(\hat{r}+t_{\epsilon})(G_{t_{\epsilon}}(g^T)-G_{t_{\epsilon}}(g_{\leftarrow,R}^{T,N}))||_{{\cal H}_{t_{\epsilon}}}\nonumber\\
&+&||{\cal J}(\hat{r}+t_{\epsilon})(G_{t_{\epsilon}}(g_{\leftarrow,R}^{T,N})-U_{\leftarrow}(t_{\epsilon},T)(\Omega_{\leftarrow}^-f)^N_R)||_{{\cal H}_{t_{\epsilon}}}+||(\Omega_{\leftarrow}^-f)^N_R-\Omega_{\leftarrow}^-f||\nonumber\\
&\le& \left(\int_0^{\infty}\int_{S^2}|g^T-g_{\leftarrow}^T|^2dtd\omega\right)^{1/2}
+\left(\int_0^{\infty}\int_{S^2}|g_{\leftarrow}^T-g_{\leftarrow,R}^{T,N}|^2dtd\omega\right)^{1/2}\nonumber\\
&+&||{\cal J}(\hat{r}+t_{\epsilon})U(t_{\epsilon},T/2+c_0)(\Phi^{R,N}(T/2+c_0,.)-\Phi^{R,N}_{\leftarrow}(T/2+c_0,.))||_{{\cal H}_{t_{\epsilon}}}\nonumber\\
&+&||(U(t_{\epsilon},T/2+c_0)-U_{\leftarrow}(T/2+c_0,.))\Phi_{\leftarrow}^{R,N}(T/2+c_0,.)||_{{\cal H}_{t_{\epsilon}}}\nonumber\\
&+&||(\Omega_{\leftarrow}^-f)^N_R-\Omega_{\leftarrow}^-f||\nonumber\\
&\le&\left(\int_0^{\infty}\int_{S^2}|g^T-g_{\leftarrow}^T|^2dtd\omega\right)^{1/2}
+2||(\Omega_{\leftarrow}^-f)^N_R-\Omega_{\leftarrow}^-f||\nonumber\\
&+&||{\cal J}(\hat{r}+t_{\epsilon})U(t_{\epsilon},T/2+c_0)(\Phi^{R,N}(T/2+c_0,.)-\Phi^{R,N}_{\leftarrow}(T/2+c_0,.))||_{{\cal H}_{t_{\epsilon}}}\nonumber\\
&+&||(U(t_{\epsilon},T/2+c_0)-U_{\leftarrow}(t_{\epsilon},T/2+c_0))\Phi_{\leftarrow}^{R,N}(T/2+c_0,.)||_{{\cal H}_{t_{\epsilon}}}.
\end{eqnarray}
\begin{figure}
\centering\epsfig{figure=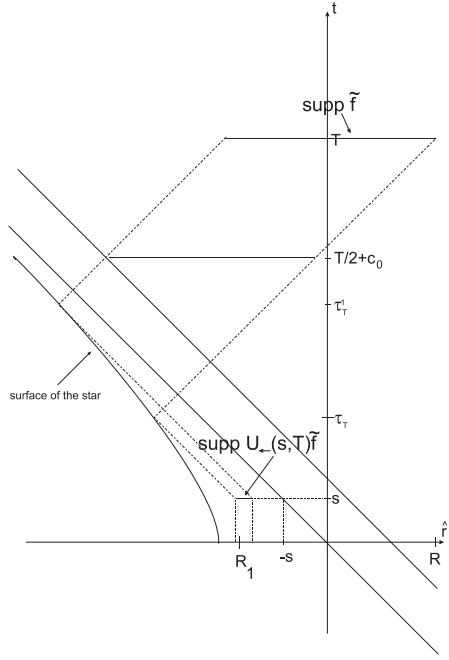,width=6cm}
\caption{The support of $(\Omega_{\leftarrow}^-f)^N_R$ transported by the asymptotic propagator $U_{\leftarrow}(s,T).$}
\label{CompDyn}
\end{figure}
We now fix $R,N$ s.t.
\[ 2||(\Omega_{\leftarrow}^-f)^N_R-\Omega_{\leftarrow}^-f||<\epsilon/4. \]
In order to estimate the last term in (\ref{9.12}) we want to use the Duhamel formula.
Let $\tau_T,\, \tau_T^1$ be defined by
\begin{eqnarray*}
\tau_T-T+R=\hat{z}(\tau_T,\theta)\Rightarrow\tau_T=\frac{T}{2}-\frac{R}{2}+{\cal O}(e^{-\kappa_+T})
\end{eqnarray*}
resp. an analogous definition for $\tau^1_T$ with $R$ replaced by $R_1$. For the 
above implication we have used (\ref{star2}).
We now observe that (see Figure \ref{CompDyn}) 
\begin{eqnarray*}
\lefteqn{s\ge \tau^1_T}\\
&\Rightarrow&\supp U_{\leftarrow}(s,T)(\Omega_{\leftarrow}^-f)_R^N\subset{[}\hat{z}(s,\theta),s+R-T{]},\\
\lefteqn{\tau^1_T\ge s\ge \tau_T}\\
&\Rightarrow&\supp U_{\leftarrow}(s,T)(\Omega_{\leftarrow}^-f)_R^N\subset{[}\hat{z}(s,\theta),\max(s+R-T,-s+2\tau^1_T-T+R_1){]},\\
\lefteqn{\tau_T\ge s}\\
&\Rightarrow&\supp U_{\leftarrow}(s,T)(\Omega_{\leftarrow}^-f)_R^N\subset{[}-s+2\tau_T-T+R,-s+2\tau^1_T-T+R_1{]}.
\end{eqnarray*}
We have to distinguish the cases $s\in{[}\tau_T,\tau^1_T{]}$ and $s\in {[}0,T{]}\setminus
{[}\tau_T,\tau^1_T{]}$. If $s\in {[}0,T{]}\setminus {[}\tau_T,\tau^1_T{]}$, then 
$U_{\leftarrow}(s,T)(\Omega_{\leftarrow}^-f)^N_R$ is zero on the boundary, in particular
\begin{eqnarray}
\label{8.13a}
U_{\leftarrow}(s,T)(\Omega_{\leftarrow}^-f)^N_R\in D(H_s)\quad \forall s\in {[}0,T{]}\setminus 
{[}\tau_T,\tau^1_T{]}.
\end{eqnarray}
Here we have used Lemma \ref{lem9.4} to establish the necessary regularity of $U_{\leftarrow}(s,T)(\Omega_{\leftarrow}^-f)^N_R.$
But we have :
\begin{eqnarray}
\label{8.13b}
U_{\leftarrow}(s,T)(\Omega_{\leftarrow}^-f)^N_R\notin D(H_s)\quad \forall s\in {[}\tau_T,\tau^1_T{]}.
\end{eqnarray}
Let 
\[ M(t,\theta)=\hat{w}(t,\theta)\hat{w}_0(t,\theta)\hat{\Gamma}^4\left(\hat{w}^{-1}\hat{\Gamma}^4+\frac{Z_1}{i}-
\frac{Z_2}{i}\hat{\Gamma}^2\right)\hat{\Gamma}^4\left(\hat{w}_0^{-1}\hat{\Gamma}^4-\frac{Z_1}{i}\right)\]
with $Z_1=\dot{\hat{z}}(t,\theta)+1,\, Z_2=(\partial_{\theta}\hat{z})(t,\theta)a_0(\hat{z}(t,\theta),\theta)h^2(\hat{z}(t,\theta),\theta)$.
The coefficients in $\hat{\Gamma}^2$ have to be evaluated at $(\hat{z}(t,\theta),\theta)$.
We first note that a matrix of type 
\[ V=\hat{\Gamma}^4\left(\hat{w}^{-1}\hat{\Gamma}^4-\frac{Z_1}{i}+\frac{Z_2}{i}\hat{\Gamma}^2\right)\]
is invertible. Indeed an elementary calculation using the anticommutation relations for
Dirac matrices gives :
\[ V\hat{\Gamma}^4\left(\hat{w}^{-1}\hat{\Gamma}^4+\frac{Z_1}{i}-\frac{Z_2}{i}\hat{\Gamma}^2\right)=\hat{w}^{-2}+Z_1^2-Z_2^2. \]
Let 
\[ A=\hat{\Gamma}^4(-\dot{\hat{z}}-\Gamma^1+i\hat{w}^{-1}\hat{\Gamma}^4+Z_2\hat{\Gamma}^2). \]
Then $V$ is an isomorphism from $Ker A$ to ${\cal K}=\{(\Psi_1,0,0,\Psi_4)\in \C^4\}$.
For dimensional reasons we only have to show $V Ker A\subset {\cal K},$
which follows from 
\[ P_{2,3}V\Psi=P_{2,3}\hat{\Gamma}^4\left(\hat{w}^{-1}\hat{\Gamma}^4-\frac{\dot{\hat{z}}}{i}-\frac{\Gamma^1}{i}
+\frac{Z_2}{i}\hat{\Gamma}^2\right)\Psi=0 \]
for $\Psi \in Ker A$. This shows that $M(s,\theta)U_{\leftarrow}(s,T)(\Omega_{\leftarrow}^-f)^N_R$ fulfills the boundary conditions for $H_s$.
We have 
\begin{eqnarray}
\label{8.13b1}
M=1+\frac{Z_2\hat{w}}{i}\hat{\Gamma}^2\hat{\Gamma}^4+Z_2Z_1\hat{w}\hat{w}_0\hat{\Gamma}^2+\frac{Z_1}{i}(\hat{w}-\hat{w}_0)\hat{\Gamma}^4+Z_1^2\hat{w}\hat{w}_0.
\end{eqnarray}
By (\ref{star0b}) we have :
\begin{eqnarray*}
\hat{w}&\lesssim& e^{\kappa_+ t},\, \hat{w}_0\lesssim e^{\kappa_+ t}.
\end{eqnarray*}
We estimate :
\begin{eqnarray*}
|\hat{w}-\hat{w}_0|\le \hat{w}\hat{w}_0(1-\dot{\hat{z}}^2)^{1/2}
\left|\left(1+\frac{(\partial_{\theta}\hat{z})^2a_0^2h^4-2(\partial_{\theta}\hat{z})l'a_0^2h^4}{1-\dot{\hat{z}}^2}\right)^{1/2}-1\right|.
\end{eqnarray*}
As 
\[ \left|\frac{(\partial_{\theta}\hat{z})^2a_0^2h^4-2(\partial_{\theta}\hat{z})l'a_0^2h^4}{1-\dot{\hat{z}}^2}\right|
\lesssim e^{-2\kappa_+ t} \]
we obtain :
\[ |\hat{w}-\hat{w}_0|\lesssim e^{-\kappa_+ t}. \]
This entails :
\begin{eqnarray}
\label{8.13b2}
M(t)=1+{\cal O}(e^{-2\kappa_+ t}),\quad \frac{d}{dt}M(t)={\cal O}(e^{-2\kappa_+ t}).
\end{eqnarray}
We write :
\begin{eqnarray*}
I&:=&\left(U\left(t_{\epsilon},\frac{T}{2}+c_0\right)-U_{\leftarrow}\left(t_{\epsilon},\frac{T}{2}+c_0\right)\right)\Phi^{R,N}_{\leftarrow}\left(\frac{T}{2}+c_0\right)\\
&=&U\left(t_{\epsilon},\tau_T\right)\left(U\left(\tau_T,\frac{T}{2}+c_0\right)-U_{\leftarrow}\left(\tau_T,\frac{T}{2}+c_0\right)\right)\Phi^{R,N}_{\leftarrow}\left(\frac{T}{2}+c_0\right)\\
&+&\left(U\left(t_{\epsilon},\tau_T\right)-U_{\leftarrow}\left(t_{\epsilon},\tau_T\right)\right)U_{\leftarrow}\left(\tau_T,\frac{T}{2}+c_0\right)\Phi^{R,N}_{\leftarrow}\left(\frac{T}{2}+c_0\right)=:I_1+I_2.
\end{eqnarray*}
Recalling that $U_{\leftarrow}(s,\frac{T}{2}+c_0)\Phi^{R,N}_{\leftarrow}(\frac{T}{2}+c_0)=
U_{\leftarrow}(s,T)(\Omega_{\leftarrow}^-f)^N_R$ and using (\ref{8.13a}) we can estimate 
the second term using the Duhamel formula :
\begin{eqnarray}
\label{8.13c}
||I_2||&\le& \int_{t_{\epsilon}}^{\tau_T}||(P_{\omega}+W)U_{\leftarrow}(s,T)(\Omega_{\leftarrow}^-f)^N_R||_{{\cal H}_s}ds\nonumber\\
&\le& C(R_1,R)\int_{t_{\epsilon}}^{\tau_T}e^{-\kappa_+ s}ds ||(\Omega_{\leftarrow}^-f)^N_R||_{{\cal H}^1}\nonumber\\
&\le&C(R_1,R)(N+1)(||\Omega_{\leftarrow}^-f||+||\Omega_{\leftarrow}^-Hf||)\frac{1}{\kappa_+}(e^{-\kappa_+ t_{\epsilon}}-e^{-\kappa_+\tau_T}).
\end{eqnarray}
Here we have used Lemma \ref{lem9.4}. To estimate the first term we write :
\begin{eqnarray*}
\lefteqn{\left(U\left(\tau_T,\frac{T}{2}+c_0\right)-U_{\leftarrow}\left(\tau_T,\frac{T}{2}+c_0\right)\right)\Phi^{R,N}_{\leftarrow}\left(\frac{T}{2}+c_0\right)}\\
&=&\left(U\left(\tau_T,\frac{T}{2}+c_0\right)M(T)-M(\tau_T)U_{\leftarrow}\left(\tau_T,\frac{T}{2}+c_0\right)\right)\Phi^{R,N}_{\leftarrow}\left(\frac{T}{2}+c_0\right)\\
&+&{\cal O}(e^{-\kappa_+ T}),
\end{eqnarray*}
where we have used (\ref{8.13b2}). Now we can use the Duhamel formula :
\begin{eqnarray*}
\lefteqn{\left|\left|\left(U\left(\tau_T,\frac{T}{2}+c_0\right)M(T)-M(\tau_T)U_{\leftarrow}\left(\tau_T,\frac{T}{2}+c_0\right)\right)\Phi^{R,N}_{\leftarrow}(\frac{T}{2}+c_0)\right|\right|}\\
&\le&\int_{\tau_T}^{\frac{T}{2}+c_0}||U(\tau_T,s)(HM(s)-M(s)H_{\leftarrow})U_{\leftarrow}(s,T)(\Omega_{\leftarrow}^-f)^N_R||ds\\
&+&\int_{\tau_T}^{\frac{T}{2}+c_0}||U(\tau_T,s)\frac{dM(s)}{ds}U_{\leftarrow}(s,T)(\Omega_{\leftarrow}^-f)^N_R||ds=:I^1_2+I_2^2
\end{eqnarray*}
To be more precise we should have used the operators $H_s$ and $H_{\leftarrow,s}$
in the above formula. But as $U_{\leftarrow}(s,T)(\Omega_{\leftarrow}^-f)^N_R$ is a smooth
function we have :
\begin{eqnarray*}
M(s)H_{\leftarrow,s}U_{\leftarrow,s}(s,T)(\Omega_{\leftarrow}^-f)^N_R&=&\Gamma^1D_{\hat{r}}MU_{\leftarrow}(s,T)(\Omega_{\leftarrow}^-f)^N_R\\
&+&{[}M,\Gamma^1D_{\hat{r}}{]}U_{\leftarrow}(s,T)(\Omega_{\leftarrow}^-f)^N_R.
\end{eqnarray*}
This will be used below. By (\ref{8.13b2}) we have :
\begin{eqnarray}
\label{8.13c1}
I^2_2\lesssim e^{-\kappa_+ T}||\Omega_{\leftarrow}^-f||.
\end{eqnarray}
Let us now estimate $I_2^1$ :
\begin{eqnarray*}
I^1_2&\le&\int_{\tau_T}^{\frac{T}{2}+c_0}||U(\tau_T,s)(H-H_{\leftarrow})MU_{\leftarrow}(s,T)(\Omega_{\leftarrow}^-f)^N_R||ds\\
&+&\int_{\tau_T}^{\frac{T}{2}+c_0}||U(\tau_T,s){[}M,H_{\leftarrow}{]}U_{\leftarrow}(s,T)(\Omega_{\leftarrow}^-f)^N_R||ds=:I_a+I_b.
\end{eqnarray*}
We have ${[}M,H_{\leftarrow}{]}={[}M,\Gamma^1{]}D_{\hat{r}}$. Recalling that $\Gamma^1=(\alpha\hat{\Gamma}^1+\beta\hat{\Gamma}^2)(\hat{z}(s,\theta),\theta)$ we find :
\[ {[}Z_2\hat{w}\hat{\Gamma}^2\hat{\Gamma}^4,\Gamma^1{]}=Z_2\hat{w}\beta{[}\hat{\Gamma}^2\hat{\Gamma}^4,\hat{\Gamma}^2{]}. \]
It follows :
\[ {[}M(s),\Gamma^1{]}={\cal O}(e^{-3\kappa_+ s}). \]
Therefore we can estimate $I_b$ using Lemma \ref{lem9.4} $(ii)$ :
\begin{eqnarray}
\label{8.13c2}
I_b&\le&C(R_1,R)\int_{\tau_T}^{\frac{T}{2}+c_0}e^{-3\kappa_+ s}e^{\kappa_+ T} (N+1)(||\Omega_{\leftarrow}^-f||
+||\Omega_{\leftarrow}^-Hf||)ds\nonumber\\
&\le&C(R_1,R)e^{-\frac{\kappa_+}{2}T}(N+1)(||\Omega_{\leftarrow}^-f||+||\Omega_{\leftarrow}^-Hf||).
\end{eqnarray}
We now estimate $I_a$ using Lemma \ref{lem9.4} $(i)$ :
\begin{eqnarray}
\label{8.14b}
I_a&\le&\int_{\tau_T}^{\frac{T}{2}+c_0}||(P_{\omega}+W)MU_{\leftarrow}(s,T)(\Omega_{\leftarrow}^-f)^N_R||ds\nonumber\\
&\le&C(R_1,R)\int_{\tau_T}^{\frac{T}{2}+c_0}e^{-\kappa_+s}||(\Omega_{\leftarrow}^-f)^N_R||_{{\cal H}^1}ds\nonumber\\
&\le&C(R_1,R)e^{-\kappa_+\frac{T}{2}}(N+1)(||\Omega_{\leftarrow}^-f||+||\Omega_{\leftarrow}^-Hf||).
\end{eqnarray}
Here we have used that $D_{\theta}M={\cal O}(e^{-2\kappa_+ t})$. Putting everything together
we find :
\begin{eqnarray}
\label{8.14c}
||I||\le C(R_1,R) e^{-\kappa_+ t_{\epsilon}}(N+1)(||\Omega_{\leftarrow}^-f||+||\Omega_{\leftarrow}^-Hf||)\quad
\mbox{uniformly in $T$ large.}
\end{eqnarray}
We fix $t_{\epsilon}$ large enough s.t. the term on the R.H.S. of (\ref{8.14c}) is controlled
by $\epsilon/4$ uniformly in $T$ large.
For $T$ sufficiently large we can estimate the first and the third terms in (\ref{9.12}) using Lemmas \ref{lem9.1} and \ref{lem9.3} :  
\begin{eqnarray*}
\left(\int_0^{\infty}\int_{S^2}|g^T-g_{\leftarrow}^T|^2dtd\omega\right)^{1/2}
<\epsilon/4,\\
||{\cal J}(\hat{r}+t_{\epsilon})U(t_{\epsilon},T/2+c_0)(\Phi^{R,N}(T/2+c_0,.)-\Phi^{R,N}_{\leftarrow}(T/2+c_0,.))||_{{\cal H}_{t_{\epsilon}}}<\epsilon/4.
\end{eqnarray*}
This concludes the proof of the proposition.
\qed
\chapter{Propagation of singularities}
\label{sec10}
So far we have compared the full dynamics $U(s,T)f$ to the dynamics $U_{\leftarrow}(s,t)\Omega_{\leftarrow}^-f$ on the 
interval ${[}t_{\epsilon},T{]}$. We will now replace the dynamics 
$U_{\leftarrow}(t_{\epsilon},T)\Omega_{\leftarrow}^-f$ by the so called geometric optics approximation. 
We suppose for the whole chapter that 
\[f(\hat{r},\omega)=e^{in\varphi}f^n(\hat{r},\theta),\, f^n\in 
C_0^{\infty}(\R\times {[}0,\pi{]}),\, n\in \Z+1/2.\] 
Let
\begin{eqnarray}
\label{10.1}
F_{t_0}^T(\hat{r},\omega)&:=&\frac{1}{\sqrt{-\kappa_+(\hat{r}+t_0)}}(\tilde{f}_3,0,0,-\tilde{f}_2)
\left(T+\frac{1}{\kappa_+}\ln(-(\hat{r}+t_0))-\frac{1}{\kappa_+}\ln \hat{A}(\theta),\omega\right),\nonumber\\
\end{eqnarray}
where $\tilde{f}=(\Omega^-_{\leftarrow}f)_R^N$ (see Section \ref{sec9.2}). Note that
\begin{eqnarray}
\label{10.2}
\supp F_{t_0}^T\subset (-t_0-|{\cal O}(e^{-\kappa_+T})|,-t_0)\times {[}0,2\pi{]}\times{[}0,\pi{]}
\end{eqnarray}
and that $F^T_{t_0}$ depends on $N,R$. All functions involved have fixed angular 
momentum, e.g. $F^T_{t_0}(\hat{r},\omega)=F^T_{t_0,n}(\hat{r},\theta)e^{in\varphi}$.
The functions $F^T_{t_0}$ and $F^T_{t_0,n}$ will often be identified. 
We therefore fix now the angular momentum $\partial_{\varphi}=in$
everywhere in the expression of $H$:
\[ H=\Gamma^1D_{\hat{r}}+\left(\begin{array}{cc} M_{\theta} & 0 \\ 0 & -M_{\theta} \end{array} \right)
+\frac{h^2a_0}{\sin\theta}\hat{\Gamma}^3n+h^2c_1+h^2c_2^{\varphi}n+\hat{V}_{\varphi}n+\hat{V}_1^{\nu}+\frac{qQr_+}{r_+^2+a^2}
+\frac{an}{r_+^2+a^2}.\]
Here $\hat{V}^{\nu}_1$ is obtained from $\hat{V}_1$ by replacing $\hat{\Gamma}^4$ by $\Gamma^{\nu}$.
Recall that 
\begin{eqnarray*}
\left(\begin{array}{cc} M_{\theta} & 0 \\ 0 & -M_{\theta} \end{array} \right)
+\frac{h^2a_0}{\sin\theta}\hat{\Gamma}^3n={\cal U}h\sqrt{a_0}\notD^n_{S^2}\sqrt{a_0}h{\cal U}^*, 
\end{eqnarray*}
where $\notD_{S^2}^n$ is the restriction of $\notD_{S^2}$ to $\{u=e^{in\varphi}\tilde{u}(\theta);\,
\tilde{u}(\theta)\in L^2(({[}0,\pi{]};\sin\theta d\theta);\, \C^4)\}$. Therefore $H$ is a regular operator and the
singularities in the expression of $H$ are coordinate singularities. 
We put 
\[ \hat{\cal H}=(L^2(\R\times {[}0,\pi{]};d\hat{r}\sin\theta d\theta))^4. \]
Clearly 
\[ \hat{\cal H}=\oplus_l\hat{\cal H}^l;\quad \hat{\cal H}^l=(L^2(\R))^4\otimes_4Y_l,\quad
Y_l=(span\{e^{2il\theta}\})^4. \]

For $\delta>0$ let $\phi_{\delta}\in C^{\infty}(\R)$ with
\begin{eqnarray*}
\phi_{\delta}(\hat{r})=\left\{\begin{array}{cc} 1 & \hat{r}\ge \delta,\\
0 & \hat{r}\le \delta/2. \end{array} \right.
\end{eqnarray*}
The aim of this chapter is to prove the following
\begin{proposition}
\label{prop10.1}
We have :
\begin{eqnarray*}
&\forall \epsilon>0\,\exists N_0>0,R_0>0\, \forall N\ge N_0,R\ge R_0\, \exists t_0>0\\
&\, \forall t_{\epsilon}\ge t_0\,\exists \delta=\delta(t_{\epsilon},N,R),\,T_0=T_0(t_{\epsilon},\delta,N,R)\,\forall T\ge T_0 \\
&||{\cal J}U(0,T)f-\phi_{\delta}(.-\hat{z}(0,\theta))e^{-it_{\epsilon}H}F_{t_{\epsilon}}^T||\le \epsilon.
\end{eqnarray*}
\end{proposition}
Note that we can consider ${\cal J}U(0,T)f-\phi_{\delta}(.-\hat{z}(0,\theta))e^{-it_{\epsilon}H}F_{t_{\epsilon}}^T$ as an element
of $\hat{\cal H}$ and it is sufficient to show :
\[ ||{\cal J}U(0,T)f-\phi_{\delta}(.-\hat{z}(0,\theta))e^{-it_{\epsilon}H}F_{t_{\epsilon}}^T||_{\hat{\cal H}}\le \epsilon. \]
We will use the pseudodifferential calculus on $\Sigma_1=\R_{\hat{r}}\times (0,\pi)$.
We note $\xi$ the dual variable to $\hat{r}$ and $q$ the dual variable to
$\theta$. Let $S^m(\Sigma_1)$ be the space of symbols of order $m$ and $\Psi^m(\Sigma_1)$
be the space of pseudodifferential operators of order $m$ (see \cite[Chapter XVIII]{Hoe1}):
\begin{eqnarray*}
\lefteqn{a(\hat{r},\theta,\xi,q)\in S^m(\Sigma_1)}\\
&\Leftrightarrow& \forall \alpha=(\alpha_1,\alpha_2),\, 
\beta=(\beta_1,\beta_2)\quad |\partial^{\alpha_1}_{\xi}\partial^{\alpha_2}_q\partial^{\beta_1}_{\hat{r}}\partial^{\beta_2}_{\theta}
a(\hat{r},\theta,\xi,q)|\le C_{\alpha,\beta}\langle(\xi,q)\rangle ^{m-|\alpha|}. 
\end{eqnarray*}
For a matrix $M=(m_{ij})$ of operators we shall write 
\[M\in \Psi^m(\Sigma_1)\Leftrightarrow\forall i,j\quad m_{ij}\in \Psi^m(\Sigma_1).\]
We use an analogous notation for a matrix of symbols. The matrix $Op(M)$ is the matrix of operators
$(Op(m_{ij}))$. 
Recall that ${\cal S}^{\rho}$ is defined as a subspace of $C^{\infty}(\R)$ by
\[ f\in S^{\rho}\Leftrightarrow \forall \alpha\in \N\, |f^{(\alpha)}(x)|\le C_{\alpha}\langle x \rangle^{\rho-\alpha}. \] 
We will study $U(0,t_{\epsilon})F_{t_{\epsilon}}^T$ microlocally. We first observe that $F_{t_{\epsilon}}^T$ has high frequencies 
in $\xi$. We show that for $L_0>0$ we have (see Lemma \ref{lem10.5}):
\begin{eqnarray}
\label{10.3}
Op\left(\chi\left(\frac{\langle \xi \rangle }{\langle q \rangle}
\le L_0\right)\right)F_{t_{\epsilon}}^T\rightarrow 0,\, T\rightarrow \infty,\\
\label{10.4}
F^T_{t_{\epsilon}}\rightharpoonup 0,\, T\rightarrow \infty.
\end{eqnarray}
We then study the propagation of singularities of $e^{-it_{\epsilon} H}$. Because of 
(\ref{10.2}), (\ref{10.3}) we are interested in the propagation of "outgoing" singularities 
located in 
\[\{(\hat{r},\omega;\xi,q);\hat{r}\ge -t_{\epsilon}-|{\cal O}(e^{-\kappa_+ T})|,\langle \xi \rangle \ge L_0 \langle q \rangle\}.\]
We will show that these singularities stay away from the surface of the star. Because of
(\ref{10.4}) it follows (modulo a small error term) : 
\begin{eqnarray}
\label{C10.5}
(1-\phi_{\delta})(.-\hat{z}(0,\theta))e^{-it_{\epsilon} H}F^T_{t_{\epsilon}}\rightarrow 0
\end{eqnarray}
for an appropriate choice of $\delta>0$. Using (\ref{C10.5}) we show that 
\[ ||(U(0,t_{\epsilon})-\phi_{\delta}(.-\hat{z}(0,\theta))e^{-it_{\epsilon}H}F^T_{t_{\epsilon}}|| \]
is small for $t_{\epsilon},\, T$ large.

This will prove Proposition \ref{prop10.1}.
\section{The geometric optics approximation and its properties}
We need the following lemmas :
\begin{lemma}
\label{lem10.2}
We have :
\[ ||F_{t_0}^T||_{L^1(\R;(L^2(S^2))^4)}\rightarrow 0,\, T\rightarrow \infty. \]
\end{lemma}
{\bf Proof.}
The lemma follows from the following calculation:
\begin{eqnarray*}
\lefteqn{||F_{t_0}^T||_{L^1(\R;(L^2(S^2))^4)}}\\
&=&\int_{-t_0-|{\cal O}(e^{-\kappa_+T})|}^{-t_0}\frac{1}{\sqrt{-\kappa_+(\hat{r}+t_0)}}\\
&\times&\left(\int_{S^2}(|\tilde{f}_3|^2+|\tilde{f}_2|^2)\left(T+\frac{1}{\kappa_+}\ln(-(\hat{r}+t_0))-\frac{1}{\kappa_+}\ln \hat{A}(\theta),\omega\right)d\omega\right)^{1/2}d\hat{r}\\
&=&\int_I-\sqrt{\kappa_+}e^{\kappa_+(y-T)/2}\left(\int_{S^2}\left(|\tilde{f}_3|^2+|\tilde{f}_2|^2\right)\left(y-\frac{1}{\kappa_+}\ln \hat{A}(\theta),\omega\right)d\omega\right)^{1/2}
dy\rightarrow 0.
\end{eqnarray*}
Here $I$ is a compact interval depending on the support of $\tilde{f}$. 
\qed
\begin{lemma}
\label{lem10.3}
We have 
\[ ||U_{\leftarrow}(t_0,T)\tilde{f}-F_{t_0}^T||\rightarrow 0,\, T\rightarrow \infty. \]
\end{lemma}
{\bf Proof.}

Let $u=U_{\leftarrow}(t_0,T)\tilde{f}$. Recall that for $T$ sufficiently large we have :
\begin{eqnarray*}
u_1(\hat{r},\omega)&=&\sqrt{\frac{1-\dot{\hat{z}}}{1+\dot{\hat{z}}}}(\hat{\tau}(\hat{r}+t_0,\theta),\theta)\tilde{f}_3
(\hat{r}+t_0+T-2\hat{\tau}(\hat{r}+t_0,\theta),\omega),\\
u_4(\hat{r},\omega)&=&-\sqrt{\frac{1-\dot{\hat{z}}}{1+\dot{\hat{z}}}}(\hat{\tau}(\hat{r}+t_0,\theta),\theta)\tilde{f}_2
(\hat{r}+t_0+T-2\hat{\tau}(\hat{r}+t_0,\theta),\omega),\\
u_2&=&u_3=0.
\end{eqnarray*}
We will also use:
\begin{eqnarray*}
\hat{\tau}(x_0,\theta)&=&-\frac{1}{2\kappa_+}\ln(-x_0)+\frac{1}{2\kappa_+}\ln\hat{A}(\theta)+{\cal O}(x_0),\, x_0\rightarrow 0^-,\\
1+\dot{\hat{z}}(\hat{\tau}(x_0,\theta))&=&-2\kappa_+x_0+{\cal O}(x_0^2),\, x_0\rightarrow 0^-.
\end{eqnarray*}
We have :
\[ \hat{l}(\hat{r},\theta):=\sqrt{-\kappa_+\hat{r}}\sqrt{\frac{1-\dot{\hat{z}}}{1+\dot{\hat{z}}}}(\hat{\tau}(\hat{r},\theta),\theta)
=\sqrt{\frac{-2\kappa_+\hat{r}+{\cal O}'(\hat{r}^2)}{-2\kappa_+\hat{r}+{\cal O}(\hat{r}^2)}},\, \hat{r}\rightarrow 0^-. \]
We calculate :
\begin{eqnarray*}
\lefteqn{||(U_{\leftarrow}(t_0,T)\tilde{f}-F_{t_0}^T)_1||^2}\\
&=&\int_{S^2}\int_{\R}|(u_1-F^T_1)|^2d\hat{r}d\omega\\
&=&\int_{S^2}\int_{I}
|\hat{l}(-\hat{A}(\theta)e^{\kappa_+(y-T)})\tilde{f}_3(y+{\cal O}(-\hat{A}(\theta)e^{\kappa_+(y-T)}),\omega)-\tilde{f}_3(y,\omega)|^2dyd\omega\\
&&\rightarrow 0,\, T\rightarrow \infty.
\end{eqnarray*}
Here $I$ is a compact interval depending on the support of $\tilde{f}$.
\qed
\\

Let $G^T=e^{-it_0H}F^T_{t_0}.$
\begin{corollary}
\label{cor10.4}
We have :
\[ F_{t_0}^T\rightharpoonup 0,\,G^T\rightharpoonup 0,\, U(t_0,T)\tilde{f}\rightharpoonup 0. \]
\end{corollary}
\begin{lemma}
\label{lem10.5}
Let $\chi\in {\cal S}^{-\rho},\, \rho>0$. Then we have for 
all $M>0$ and uniformly in $t_{\epsilon}$:
\begin{eqnarray}
\label{C10.9}
\quad Op\left(\chi\left(\frac{\langle \xi \rangle}{\langle q \rangle}\right)
\langle q \rangle^M\right) F_{t_{\epsilon},n}^T\rightarrow 0, \, T\rightarrow \infty\quad \hat{\cal H},\\
\label{C10.10}
\quad Op(\chi(\langle (\xi,q) \rangle\langle q \rangle^M) F_{t_{\epsilon},n}^T\rightarrow 0,\, T\rightarrow \infty\quad\hat{\cal H}.
\end{eqnarray}
\end{lemma}
{\bf Proof.}

We only show (\ref{C10.9}), the proof of (\ref{C10.10}) being analogous. 
Let us write 
\begin{eqnarray*}
F_{0,n}^T&=&K_T\hat{g},\\
\hat{g}(\hat{r},\theta)&=&\frac{1}{\sqrt{-\kappa_+\hat{r}}}(\tilde{f}_3,0,0,-\tilde{f}_2)\left(\frac{1}{\kappa_+}\ln(-\hat{r})-\frac{1}{\kappa_+}
\ln \hat{A}(\theta),\theta\right),\\
(K_T\hat{g})(\hat{r},\theta)&=&e^{\frac{\kappa_+}{2} T}\hat{g}((\hat{r}+t_{\epsilon})e^{\kappa_+T},\theta)
\end{eqnarray*}
and we shall also consider $K_T$ as an operator on $\hat{\cal H}^{l}$ (rather than on $\hat{\cal H}$).
Let us write 
\begin{eqnarray*}
\hat{g}&=&\sum_{l}\hat{g}^{l},\quad\hat{g}^l\in\hat{\cal H}^l\quad\forall l.\quad\mbox{Thus :}\\
\chi\left(\frac{\langle D_{\hat{r}} \rangle }{\langle D_{\theta} \rangle } \right)
\langle D_{\theta} \rangle^M K_T\hat{g}&=&\sum_{l}\chi\left(\frac{\langle D_{\hat{r}} \rangle }{\langle l \rangle } \right)
\langle l \rangle^M K_T\hat{g}^{l}.
\end{eqnarray*}
We have 
\[ {\cal F}(K_T\hat{g}^{l})(\xi)={\cal F}(\hat{g}^{l})(\xi e^{-\kappa_+ T})e^{-\frac{\kappa_+}{2} T}e^{it_{\epsilon}\xi}, \]
where ${\cal F}$ denotes the Fourier transform in $\hat{r}$. 
Note that 
\begin{eqnarray*}
\left|\left|\chi\left(\frac{\langle D_{\hat{r}} \rangle}{\langle D_{\theta} \rangle}\right)
\langle D_{\theta} \rangle^MK_T\hat{g}\right|\right|^2=\sum_l\left|\left|\chi\left(\frac{\langle \xi \rangle}{\langle l \rangle }\right)
\langle l \rangle^M{\cal F}(\hat{g}^l)(\xi e^{-\kappa_+ T})e^{-\frac{\kappa_+T}{2}}\right|\right|^2,\\
\forall T\quad \left|\left|\chi\left(\frac{\langle \xi \rangle}{\langle l \rangle}\right)\langle l \rangle^M
{\cal F}(\hat{g}^l)(\xi e^{-\kappa_+ T})e^{-\frac{\kappa_+ T}{2}}\right|\right|^2\le \langle l \rangle ^{2M}\left|\left|\hat{g}^l\right|\right|^2,\quad
\sum_l\langle l \rangle^{2M}||\hat{g}^l||^2<\infty.
\end{eqnarray*}
It is therefore sufficient to show :
\[ \forall l \quad \left|\left|\chi\left(\frac{\langle \xi \rangle}{\langle l \rangle}\right)\langle l \rangle^M
{\cal F}(\hat{g}^l)(\xi e^{-\kappa_+ T})e^{-\frac{\kappa_+ T}{2}}\right|\right|\rightarrow 0,\, T\rightarrow\infty. \]
But,
\begin{eqnarray*}
\left|\left|\chi\left(\frac{\langle \xi \rangle}{\langle l \rangle}\right)\langle l \rangle^M
{\cal F}(\hat{g}^l)(\xi e^{-\kappa_+ T})e^{-\frac{\kappa_+ T}{2}}\right|\right|^2
=\int\left|\chi\left(\frac{\langle \xi e^{\kappa_+ T}\rangle}{\langle l \rangle}\right)\langle l \rangle^M
{\cal F}(\hat{g}^l)(\xi)\right|^2d\xi\rightarrow 0
\end{eqnarray*}
by the Lebesgue Theorem.
This proves (\ref{C10.9}). 
\qed
\section{Diagonalization}
\label{sec9Diag}
Let $\nu_1>0,\, \nu_2>0,\, j_{\pm},\chi\in C^{\infty}(\R),\, \supp\chi\subset\R\setminus {[}-\nu_2,\nu_2{]},\, 
\chi\equiv 1\,\mbox{on}\, \R\setminus{[}-2\nu_2,2\nu_2{]},$
\[ j_+(x)=\left\{\begin{array}{cc} 1 & x\ge 2\nu_1 \\ 0 & x\le \nu_1 \end{array} \right.,\,
j_-(x)=\left\{\begin{array}{cc} 1 & x\le -2\nu_1 \\ 0 & x\ge -\nu_1 \end{array} \right.,\, 
j^2(x)=j^2_-(x)+j^2_+(x). \]
We put 
\begin{eqnarray}
\label{IX.1}
\tilde{W}_{\pm}(x)&=&\tilde{c}_{\pm}\left(\begin{array}{cc} a_0(x\pm l') & \pm k'-\sqrt{k'^2+a_0^2(x\pm l')^2} \\
\sqrt{k'^2+a_0^2(x\pm l')^2}\mp k' & a_0(x\pm l') \end{array} \right),\nonumber\\
\tilde{c}_{\pm}&=&\frac{e_{\pm}(x+l')}{\sqrt{2}\sqrt[4]{k'^2+a_0^2(x\pm l')^2}\sqrt{\sqrt{k'^2+a_0^2(x\pm l')^2}\mp k'}},\nonumber\\
e_+(x)&=&sign x,\, e_-(x)=1,\\
W^0_{\pm}&=&\tilde{W}_{\pm}\left(\frac{q}{|\xi|}\right)j_{\pm}\left(\frac{\xi}{|q|}\right)
\chi\left(h^2\sqrt{k'^2\xi^2+a_0^2(l'\xi+q)^2}\right),\nonumber\\
W_{\pm}&=&{\cal U}\left(\begin{array}{cc} W^0_{\pm} & 0 \\ 0 & W_{\pm}^0 \end{array} \right),\nonumber\\
W(\hat{r},\theta,\xi,q)&=&W_+(\hat{r},\theta,\xi,q)+W_-(\hat{r},\theta,\xi,q).
\end{eqnarray}
By (\ref{IX.1}) $W(\hat{r},\theta,\xi,q)$ is only defined for $q\neq 0$ and $\xi\neq 0$.
We note that for $\epsilon_q>0$ small enough and $|q|<\epsilon_q$ we have 
\begin{eqnarray}
\label{IX.2}
W(\hat{r},\theta,\xi,q)&=&{\cal U}\left[\left(\begin{array}{cc} \tilde{W}_+\left(\frac{q}{|\xi|}\right) & 0 \\
0 & \tilde{W}_+\left(\frac{q}{|\xi|}\right) \end{array} \right)+\left(\begin{array}{cc} \tilde{W}_-\left(\frac{q}{|\xi|}\right) & 0 \\
0 & \tilde{W}_-\left(\frac{q}{|\xi|}\right) \end{array} \right)\right]\nonumber\\
&\times&\chi\left(h^2\sqrt{k'^2\xi^2+a_0^2(l'\xi+q)^2}\right)
\end{eqnarray}
Indeed if $|q|<\epsilon_q$ we have on $\supp\chi\left(h^2\sqrt{k'^2\xi^2+a_0^2(l'\xi+q)^2}\right)$
\begin{eqnarray}
\label{IX.3}
C(|\xi|^2+|q|^2)\ge\nu_2^2\Rightarrow |\xi|^2\ge\frac{1}{C}\nu_2^2-\epsilon_q^2.
\end{eqnarray}
But if $|\xi|^2\ge 4\nu_1^2|q|^2$ we have $j_{\pm}\left(\frac{\xi}{|q|}\right)=1.$ Using
(\ref{IX.3}) we see that this is fulfilled if $\epsilon_q$ is small enough. Therefore we define
$W(\hat{r},\theta,\xi,q)$ for $q=0$ by (\ref{IX.2}). Similarly there exists $\epsilon_{\xi}>0$ 
s.t. if $|\xi|<\epsilon_{\xi}$ then 
\begin{eqnarray}
\label{IX.4}
W(\hat{r},\theta,\xi,q)=0.
\end{eqnarray}
Indeed on $\supp j_{\pm}\left(\frac{\xi}{|q|}\right)$ we have $|q|\le \frac{|\xi|}{\nu_1}$.
If $|\xi|\le \epsilon_{\xi}$, then
\[ h^2\sqrt{k'^2\xi^2+a_0^2(l'\xi+q)^2}\le C\epsilon_{\xi}< \nu_2 \]
for $\epsilon_{\xi}$ small enough and thus $\chi\left(h^2\sqrt{k'^2\xi^2+a_0^2(l'\xi+q)^2}\right)=0$.
Therefore we define $W(\hat{r},\theta,0,q)=0.$ $W(\hat{r},\theta,\xi,q)$ with these definitions is a matrix of smooth
functions.
We want to check that $W(\hat{r},\theta;\xi,q)\in S^0(\Sigma_1)$. To this purpose we apply
the symplectic change of coordinates 
\begin{eqnarray}
\label{NSC}
r_*=r_*(\hat{r},\theta),\, \theta^*=\theta,\, q^*=l'\xi+q,\, \xi^*=k'\xi. 
\end{eqnarray}
Under this change of coordinates we obtain the symbol :
\begin{eqnarray*}
\hat{W}(r_*,\theta^*,\xi^*,q^*))&=&{\cal U}\left[\left(\begin{array}{cc} \hat{W}_+\left(\frac{q^*}{|\xi^*|}\right) & 0 \\ 0 & \hat{W}_+\left(\frac{q^*}{|\xi^*|}\right) \end{array}
\right)j_+\left(\frac{\xi^*}{|k'q^*-l'\xi^*|}\right)\right.\\
&+&\left.\left(\begin{array}{cc} \hat{W}_-\left(\frac{q^*}{|\xi^*|}\right) & 0 \\ 0 & \hat{W}_-\left(\frac{q^*}{|\xi^*|}\right) \end{array} \right)
j_-\left(\frac{\xi^*}{|k'q^*-l'\xi^*|}\right)\right]\\
&\times&\chi\left(h^2\sqrt{|\xi^*|^2+a_0^2|q^*|^2}\right),\\
\hat{W}_{\pm}(x)&=&\hat{c}_{\pm}\left(\begin{array}{cc} a_0 x & \pm 1 -\sqrt{1+a_0^2x^2} \\
\sqrt{1+a_0^2x^2}\mp 1 & a_0 x \end{array} \right),\\
\hat{c}_{\pm}&=&\frac{e_{\pm}(x)}{\sqrt{2}\sqrt[4]{1+a_0^2x^2}\sqrt{\sqrt{1+a_0^2x^2}\mp 1}}.
\end{eqnarray*}
It is sufficient to check
\[ \hat{W}(r_*,\theta^*,\xi^*,q^*)\in S^0. \]
On $\supp j_{\pm}\left(\frac{\xi^*}{|k'q^*-l'\xi^*|}\right),$ 
$\frac{|q^*|}{|\xi^*|}$ remains bounded. The functions 
\[ f^{\pm}(x)=\frac{|x|}{\sqrt{2}\sqrt[4]{x^2+1}\sqrt{\sqrt{x^2+1}\mp 1}},\quad 
g^{\pm}=\pm \frac{\sqrt{\sqrt{x^2+1}\mp 1}}{\sqrt{2}\sqrt[4]{x^2+1}} \]
are $C^{\infty}(\R)$ functions with 
\begin{eqnarray}
\label{DerC}
|\partial^{\alpha}f^{\pm}(x)|\le C_{\alpha},\, |\partial^{\alpha}g^{\pm}(x)|\le C_{\alpha}. 
\end{eqnarray}
In order to see that the estimate holds for $f^+$ we note that $\check{f}(x)=\frac{\sqrt{\sqrt{x^2+1}-1}}{|x|}$
can be extended to an analytic function in a neighborhood of zero with $\check{f}(0)=\frac{1}{\sqrt{2}}$. We obtain :
\[ \left|\partial_{q^*}\hat{W}_{\pm}\left(\frac{q^*}{\xi^*}\right)\right|\le C \frac{1}{|\xi^*|},\quad
\left|\partial_{\xi^*}\hat{W}_{\pm}\left(\frac{q^*}{\xi^*}\right)\right|\le C \frac{|q^*|}{|\xi^*|^2}. \]
But on $\supp j_{\pm}\left(\frac{\xi^*}{|k'q^*-l'\xi^*|}\right)\cap
\supp \chi(h^2\sqrt{|\xi^*|^2+a_0^2|q^*|^2})$ we have :
\begin{eqnarray}
\label{xiq}
\frac{1}{|\xi^*|}\lesssim \frac{1}{\langle \xi^* \rangle}, \quad \frac{1}{|\xi^*|}\lesssim \frac{1}{\langle q^* \rangle}. 
\end{eqnarray}
We now have to estimate derivatives on $j_{\pm}\left(\frac{\xi^*}{|k'q^*-l'\xi^*|}\right)
\chi\left(h^2\sqrt{|\xi^*|^2+a_0^2|q^*|^2}\right)$. In the region $k'q^*-l'\xi^*>0$
we find :
\begin{eqnarray*}
\left|\partial_{\xi^*}j_{\pm}\left(\frac{\xi^*}{|k'q^*-l'\xi^*|}\right)\right|
&=&\left|j'_{\pm}\left(\frac{\xi^*}{|k'q^*-l'\xi^*|}\right)\left(\frac{1}{k'q^*-l'\xi^*}+\frac{\xi^*l'}{(k'q^*-l'\xi^*)^2}\right)\right|\\
&\lesssim&\frac{1}{|\xi^*|},\\
\left|\partial_{q^*}j_{\pm}\left(\frac{\xi^*}{|k'q^*-l'\xi^*|}\right)\right|&\lesssim&\frac{1}{|\xi^*|}.
\end{eqnarray*}
We then use (\ref{xiq}). Derivatives in $r_*,\theta^*$, derivatives on $\chi\left(h^2\sqrt{|\xi^*|^2+a_0^2|q^*|^2}\right),\, {\cal U}$
as well as higher order derivatives can be controlled in a similar way.
Let $\tilde{\chi}\in C^{\infty}(\R),\, \supp \tilde{\chi}\subset \R\setminus {[}-\frac{\nu_2}{2},\frac{\nu_2}{2}{]}$
with $\tilde{\chi}\chi=\chi.$ We next put :
\begin{eqnarray*}
\lambda&=&h^2\sqrt{k'^2\xi^2+a_0^2(l'\xi+q)^2}\tilde{\chi}(h^2\sqrt{k'^2\xi^2+a_0^2(l'\xi+q)^2}),\\ \check{H}_d&=&\lambda\Gamma^1,\, \hat{H}_d=Op(\check{H}_d).
\end{eqnarray*} 
We have :
\begin{eqnarray}
\label{IX.5}
HOp(W_{\pm})&=&Op(W_{\pm})\hat{H}_d+Op(R_{\pm})+Op(\tilde{R}_{\pm})\quad\mbox{with}\\
R_{\pm}&=&R_{\pm}^1+R_{\pm}^2,\, R_{\pm}^k=(r_{ij}^{k\pm}),\\
\supp r^{1\pm}_{ij}&\subset& \supp j_{\pm}\left(\frac{\xi}{|q|}\right)\cap \supp \chi\left(h^2\sqrt{k'^2\xi^2+a_0^2|l'\xi+q|^2}\right),\\
r^{1\pm}_{ij}&\in& S^0(\Sigma_1),\, r^{2\pm}_{ij}\in S^{-\infty}(\Sigma_1),\\
\left(\begin{array}{cc} r^{k\pm}_{13} & r^{k\pm}_{14} \\ r^{k\pm}_{23} & r^{k\pm}_{24} \end{array} \right)
&=&\left(\begin{array}{cc} r^{k\pm}_{31} & r^{k\pm}_{32} \\ r^{k\pm}_{41} & r^{k\pm}_{42} \end{array} \right)=0,\\
\tilde{R}_{\pm}&=&\tilde{R}^1_{\pm}+\tilde{R}^2_{\pm},\, \tilde{R}^2_{\pm}\in \Psi^{-1}(\Sigma_1),\\
\tilde{R}^1_{\pm}&=&{\cal U}p\Gamma^4\left(\begin{array}{cc} W^0_{\pm} & 0 \\ 0 & 
W^0_{\pm} \end{array} \right),\\
p(\hat{r},\theta)&=&\frac{m\sqrt{\Delta}}{\sigma}(\rho-\sqrt{r^2+a^2})+b_0.
\end{eqnarray}
We need a better estimate on the remainder.
\begin{lemma}
\label{lem10.6}
There exists $M=(m_{ij})$ s.t. for all $j$ $m_{jj}=1$ and for all $i\neq j$ $m_{ij}\in \Psi^{-1}(\Sigma_1)$ 
as well as $r_j\in S^0(\Sigma_1),\, j=1,...4$ s.t. for 
\[ H_d=\hat{H}_d+Op(Diag(r_1,r_2,r_3,r_4)) \]
we have :
\begin{eqnarray}
\label{IX.6}
HOp(W)M-Op(W)MH_d\in \Psi^{-1}(\Sigma_1).
\end{eqnarray}
\end{lemma}

{\bf Proof.}

We can construct independently $M_{\pm}$ and $r^{\pm}_j$ s.t. (\ref{IX.6}) is fulfilled
for $M$ replaced by $M_{\pm}$, $r_j$ by $r_j^{\pm}$ and $W$ by $W_{\pm}$. We then put
$M=\frac{1}{2}(M_++M_-)$, $r_j=\frac{1}{2}(r_j^++r_j^-)$. We only consider the $+$ case and drop the index $+$. 
We are looking for $M$ in the form 
\begin{eqnarray*}
M&=&\left(\begin{array}{cc} A & B \\ C & D \end{array} \right);\quad
B,C\in \Psi^{-1}(\Sigma_1),\\
A&=&\left(\begin{array}{cc} 1 & Op(\alpha_1) \\ Op(\alpha_2) & 1 \end{array} \right),\quad
D=\left(\begin{array}{cc} 1 & Op(\delta_1) \\ Op(\delta_2) & 1\end{array}\right),\\
\alpha_j,\delta_j&\in& S^{-1}(\Sigma_1),\, j=1,2.
\end{eqnarray*}
If $M$ is of this form it is sufficient that
\begin{eqnarray}
\label{IX.7}
Op(R)+Op(\tilde{R})-Op(W){[}M,\hat{H}_d{]}&=&Op(W)Op(Diag(r_1,...,r_4))\nonumber\\
&+&\hat{R},\, \hat{R}\in \Psi^{-1}(\Sigma_1).
\end{eqnarray}
Here we have used that 
\begin{eqnarray}
\label{10.9}
M=I_d+R_1,\, R_1\in \Psi^{-1}(\Sigma_1).
\end{eqnarray}
Therefore 
\[ MOp(Diag(r_1,...,r_4))=Op(Diag(r_1,...,r_4))+R_2,\, R_2\in \Psi^{-1}(\Sigma_1). \]
Recalling that $\hat{H}_d=\Lambda Diag(1_d,-1_d),\, \Lambda:=Op(\lambda)Diag(1,-1)$ we find:
\[ {[}M,\hat{H}_d{]}=\left(\begin{array}{cc} {[}A,\Lambda{]} & - \{B,\Lambda\} \\
\{C,\Lambda\} & -{[}D,\Lambda{]} \end{array} \right). \]
If ${\cal U}^*R^1=\left(\begin{array}{cc} R_{11} & 0 \\ 0 & R_{22} \end{array} \right)$
and ${\cal U}^*\tilde{R}^1=\left(\begin{array}{cc} 0 & R_{12} \\ R_{21} & 0 \end{array} \right)$
we have to find $A,B,C,D,r_j$ s.t. 
\begin{eqnarray}
\label{10.10}
-W^0{[}A,\Lambda{]}+R_{11}&=&W^0Diag(r_1,r_2)+\hat{R}_{11},\, \hat{R}_{11}\in \Psi^{-1}(\Sigma_1), \\
\label{10.11}
W^0{[}D,\Lambda{]}+R_{22}&=&W^0Diag(r_3,r_4)+\hat{R}_{22},\, \hat{R}_{22}\in \Psi^{-1}(\Sigma_1), \\
\label{10.12}
W^0\{B,\Lambda\}+R_{12}&\in&\Psi^{-1}(\Sigma_1),\\
\label{10.13}
-W^0\{C,\Lambda\}+R_{21}&\in&\Psi^{-1}(\Sigma_1).
\end{eqnarray}
We consider equations (\ref{10.10}), (\ref{10.13}). 
On 
\[\supp j_+\left(\frac{\xi}{|q|}\right)
\cap \supp\chi\left(h^2\sqrt{k'^2\xi^2+a_0^2(l'\xi+q)^2}\right)\]
the matrix $W^0$ is invertible. Let
\begin{eqnarray*}
(W^0)^{-1}R_{11}=\left(\begin{array}{cc} \check{r}_{11} & \check{r}_{12} \\ 
\check{r}_{21} & \check{r}_{22} \end{array} \right). 
\end{eqnarray*}
As 
\begin{eqnarray*}
{[}A,\Lambda{]}=\left(\begin{array}{cc} 0 & -2Op(\alpha_1\lambda) \\ 2Op(\alpha_2\lambda) & 0 \end{array}\right)+R_3,\,
R_3\in S^{-1}(\Sigma_1) 
\end{eqnarray*}
we have to solve on $\supp \left(j_+\left(\frac{\xi}{|q|}\right)\right)
\cap \supp\left(\chi\left(h^2\sqrt{k'^2\xi^2+a_0^2(l'\xi+q)^2}\right)\right)$:
\begin{eqnarray}
\label{10.14}
2\alpha_1\lambda+\check{r}_{12}&\in& S^{-1}(\Sigma_1),\\
\label{10.14ca}
-2\alpha_2\lambda+\check{r}_{21}&\in& S^{-1}(\Sigma_1),
\end{eqnarray}
which can be achieved by
\[ \alpha_1=-\frac{\check{r}_{12}}{2\lambda},\, \alpha_2=\frac{\check{r}_{21}}{2\lambda}. \]
In order to solve (\ref{10.13}) we have to use the special structure of $R_{21}$. 
Indeed we have $R_{21}=ipW^0$. We try
\[ C=i\left(\begin{array}{cc} Op(\gamma_1) & 0 \\ 0 & Op(\gamma_2) \end{array} \right);\quad \gamma_1,
\gamma_2\in S^{-1}(\Sigma_1). \]
Then :
\[ -W^0\{C,\Lambda\}+R_{21}=iW^0\left(\begin{array}{cc} -2Op(\gamma_1\lambda) & 0 \\ 0 & 2Op(\gamma_2\lambda) \end{array} \right)+iW^0p \]
and therefore we can take on $\supp j_+\left(\frac{\xi}{|q|}\right)
\cap \supp\chi\left(h^2\sqrt{k'^2\xi^2+a_0^2(l'\xi+q)^2}\right)$:
\[ \gamma_1=-\gamma_2=\frac{p}{2\lambda}. \]
This concludes the proof of the lemma.
\qed
\begin{lemma}
\label{lem10.10}
We have for all $t_{\epsilon}>0$:
\[\forall 0\le s \le t_{\epsilon}\quad (e^{-isH}-Op(W)Me^{-isH_d}M^{-1}Op(W^*))F_{t_{\epsilon}}^T\rightarrow 0,\quad T\rightarrow \infty. \]
\end{lemma}
{\bf Proof.}

We have 
\begin{eqnarray*}
\lefteqn{(e^{-isH}-Op(W)Me^{-isH_d}M^{-1}Op(W^*))F^T_{t_{\epsilon}}}\\
&=&\int_0^{s}e^{-i\tau H}(H Op(W)M-Op(W)MH_d)e^{-i(s-\tau)H_d}M^{-1}Op(W^*)
F^T_{t_{\epsilon}}d\tau\\
&+&e^{-isH}(1-Op(W)Op(W^*))F^T_{t_{\epsilon}}\\
&=:&I_1(s,T)+I_2(s,T).
\end{eqnarray*}
By Lemma \ref{lem10.5} we have :
\[ I_2(s,T)\rightarrow0,\, T\rightarrow\infty. \]
Using Lemmas \ref{lem10.6} and \ref{lem10.5} we see that the first term can be estimated by :
\[ ||I_1(s,T)||\lesssim ||\langle H_d\rangle^{-1}M^{-1}Op(W^*)F^T_{t_{\epsilon}}||\rightarrow 0,
\, T\rightarrow \infty. \]
\qed 

Let 
\begin{eqnarray*}
\tilde{N}_{\pm}&=&\tilde{W}_{\pm}(0),\\
N&=&{\cal U}\left(\left(\begin{array}{cc} \tilde{N}_+ & 0 \\ 0 & \tilde{N}_+ \end{array} \right) j_+\left(\frac{\xi}{|q|}\right)
+\left(\begin{array}{cc} \tilde{N}_- & 0 \\ 0 & \tilde{N}_- \end{array} \right)j_-\left(\frac{\xi}{|q|}\right)\right)\\
&\times&\chi\left(h^2\sqrt{k'^2\xi^2+a_0^2(l'\xi+q)^2}\right). 
\end{eqnarray*}
We have :
\begin{eqnarray*}
\lefteqn{\tilde{N}_+}\\
&=&\frac{1}{\sqrt{2}}\left(\begin{array}{cc} \sqrt{1+\alpha} & -\sqrt{1-\alpha}sign(\cos\theta) \\
\sqrt{1-\alpha}sign(\cos\theta) & \sqrt{1+\alpha} \end{array} \right)j_+\left(\frac{\xi}{|q|}\right)\\
&\times&\chi\left(h^2\sqrt{k'^2\xi^2+a_0^2(l'\xi+q)^2}\right),\\
\lefteqn{\tilde{N}_-}\\
&=&\frac{1}{\sqrt{2}}\left(\begin{array}{cc} -sign(\cos\theta)\sqrt{1-\alpha} & -\sqrt{1+\alpha} \\
\sqrt{1+\alpha} & -sign(\cos\theta)\sqrt{1-\alpha} \end{array} \right)j_-\left(\frac{\xi}{|q|}\right)\\
&\times&\chi\left(h^2\sqrt{k'^2\xi^2+a_0^2(l'\xi+q)^2}\right) 
\end{eqnarray*}
with $\alpha$ as in Chapter \ref{sec3}. Note that the matrices $\tilde{N}_{\pm}$ are smooth
($\alpha(\hat{r},\frac{\pi}{2})=1$). From the lipschitz continuity of $f^{\pm},g^{\pm}$ we infer
\begin{eqnarray}
\label{ApW}
(W-N)\frac{\langle\xi\rangle }{\langle q \rangle}\in S^0(\Sigma_1).
\end{eqnarray}
Therefore by Lemma \ref{lem10.5} :
\[ (Op(W^*)-Op(N^*))F^T_{t_{\epsilon}}\rightarrow 0. \]
Let us put 
\begin{eqnarray*}
N_{\leftarrow}&=&\left(j_+\left(\frac{\xi}{|q|}\right)+\left(\begin{array}{cccc} 0 & -1 & 0 & 0 \\
1 & 0 & 0 & 0 \\ 0 & 0 & 0 & -1 \\ 0 & 0 & 1 & 0 \end{array} \right)j_-\left(\frac{\xi}{|q|}\right)\right)
\chi\left(h^2\sqrt{k'^2\xi^2+a_0^2(l'\xi+q)^2}\right)\\
&=&\left(j_+\left(\frac{\xi}{|q|}\right)+N_{\leftarrow}^-j_-\left(\frac{\xi}{|q|}\right)\right)
\chi\left(h^2\sqrt{k'^2\xi^2+a_0^2(l'\xi+q)^2}\right).
\end{eqnarray*}
Note that we have uniformly in $T$ large :
\begin{eqnarray}
\label{DiagError}
Op(N^*)F^T_{t_{\epsilon}}=Op(N^*_{\leftarrow})F_{t_{\epsilon}}^T+{\cal O}_{N,R}(e^{-\kappa_+t_{\epsilon}}). 
\end{eqnarray}
Here we have used that the commutators \[\left[\sqrt{1-\alpha},Op\left(j_+\left(\frac{\xi}{|q|}\right)\chi(h^2\sqrt{k'^2\xi^2+a_0^2q^2})\right)\right]\]
etc. are all in $\Psi^{-1}(\Sigma_1)$.
\section{Study of the hamiltonian flow}
In this section we study the hamiltonian flow of
\begin{eqnarray}
\label{10.30}
P=h^2\sqrt{k'^2\xi^2+a_0^2(l'\xi+q)^2}.
\end{eqnarray}
We denote $\phi_t$ the hamiltonian flow of $P$. Let for $L>0$
\begin{eqnarray*}
{\cal E}_L&:=&\left\{(\hat{r},\theta;\xi,q);\, 
\frac{\xi}{|q|}\ge L\right\},\\
{\cal I}^{t_0}_{L}&:=&\left\{(\hat{r},\theta;\xi,q);\, \hat{r}\ge -t_0-L^{-1}\right\}.
\end{eqnarray*}
\begin{lemma}
\label{lem10.7}
For $t_0>0$ sufficiently large there exist $\delta>0,\ L_0>0$ s.t. for all $L\ge L_0$ we have:
\[ \forall 0\le s\le t_0\quad \phi_{s}({\cal I}^{t_0}_L\cap {\cal E}_L)\subset
\{(\hat{r},\theta;\xi,q);\, \hat{r}\ge \hat{z}(t_0-s,\theta)+\delta\}. \]
\end{lemma}
{\bf Proof.}

We use the coordinates $(r_*,\theta^*,\xi^*,q^*)$ given by (\ref{NSC})
and drop the star for $\theta$ : $\theta=\theta^*$.
Under this change of coordinates the hamiltonian becomes :
\[ P^*=h^2\sqrt{|\xi^*|^2+a^2_0|q^*|^2}=E=const. \]
The hamiltonian equations are :
\begin{eqnarray}
\label{10.31}
\dot{r}_*&=&\frac{h^4\xi^*}{E},\\
\label{10.32}
\dot{\xi^*}&=&-\partial_{r_*}P^*,\\
\label{10.33}
\dot{\theta}&=&\frac{h^4a_0^2q^*}{E},\\
\label{10.34}
\dot{q}^*&=&-(\partial_{\theta}h^2)h^{-2}E.
\end{eqnarray}
Multiplying (\ref{10.34}) by $q^*$ given by (\ref{10.33}) we obtain :
\begin{eqnarray}
\label{10.38}
\frac{1}{2}\frac{d}{dt}(q^*)^2&=&q^*\dot{q^*}=-(\partial_{\theta}h^2)h^{-6}a_0^{-2}E^2\dot{\theta}
=-\frac{1}{2}(\partial_{\theta}h^4)h^{-8}a_0^{-2}E^2\dot{\theta}\nonumber\\
&=&\frac{1}{2}\frac{(r^2+a^2)^2}{\sigma^4}(\partial_{\theta}a^2\Delta\cos^2\theta)
\frac{(r^2+a^2)^2}{\Delta}\frac{\sigma^4}{(r^2+a^2)^4}E^2\dot{\theta}\nonumber\\
&=&\frac{1}{2}\frac{d}{dt}a^2E^2\cos^2\theta\nonumber\\
\Rightarrow|q^*|^2&=&|q_0^*|^2+a^2E^2(\cos^2\theta-\cos^2\theta_0),
\end{eqnarray}
in particular ${\cal K}=|q^*|^2+a^2E^2\sin^2\theta=const.$ We have 
\begin{eqnarray*}
E^2&=&h^4\left(\xi^2h^{-4}+2a_0^2l'\xi q+a_0^2q^2\right)\\
\Rightarrow \frac{\xi^2}{E^2}&=&1-\frac{2h^4a_0^2l'\xi q}{E^2}-\frac{a_0^2h^4q^2}{E^2},
\end{eqnarray*}
in particular 
\begin{eqnarray}
\label{9.20ca}
\frac{\xi_0^2}{E^2}=1+{\cal O}(L^{-1}).
\end{eqnarray}
Using (\ref{9.20ca}) we see that 
\begin{eqnarray}
\label{9.20cb}
\frac{q_0^*}{E}=\frac{q_0}{E}+a\cos\theta_0\frac{\xi_0}{E}=a\cos\theta_0+{\cal O}(L^{-1}).
\end{eqnarray}
Therefore
\begin{eqnarray}
\label{9.20cc}
\frac{\cal K}{E^2}=a^2+{\cal O}(L^{-1}).
\end{eqnarray}
These estimates are uniform for $(\hat{r}_0,\theta_0,\xi_0,q_{\theta0})\in {\cal E}_L\cap{\cal I}_L^{t_0}$.
We note that 
\begin{eqnarray}
\label{C10.37}
|\xi^*|^2=E^2\left(1-\frac{\Delta{\cal K}}{(r^2+a^2)^2E^2}\right),\quad |q^*|^2={\cal K}-a^2E^2\sin^2\theta.
\end{eqnarray}
Using (\ref{9.20cc}) and (\ref{C10.37}) we see that for $L$ sufficiently large 
$\xi^*$ does not change its sign. 
We now claim that there exists a constant $C(t_0)$ s.t. 
\begin{eqnarray}
\label{9.41a}
\dot{\hat{r}}\ge 1-C(t_0)L^{-1/2}
\end{eqnarray}
uniformly in $(\hat{r}_0,\theta_0,\xi_0,q_{\theta 0})\in {\cal I}^{t_0}_L\cap {\cal E}_L.$
We first argue that (\ref{9.41a}) proves the Lemma. By (\ref{star2})
we see that for $t_0$ sufficiently large we have :
\[ \hat{z}(t_0,\theta)<-t_0\quad \forall \theta\in {[}0,\pi{]}. \]
If $t_0$ is fixed in this way, then there exists $\delta>0$ s.t.
\[ \dot{\hat{z}}(\tau,\theta)>-1+\frac{2\delta}{t_0}\quad \forall 0\le \tau\le t_0. \]
We have 
\[ \hat{r}(0)-\hat{z}(t_0,\theta)\ge -L^{-1}-t_0-\hat{z}(t_0,\theta)>0\, \forall \theta\in {[}0,\pi{]} \]
for $L$ sufficiently large and 
\[ \frac{d}{ds} (\hat{r}(s)-\hat{z}(t_0-s,\theta))\ge \frac{2\delta}{t_0}-C(t_0)L^{-1/2}>\frac{\delta}{t_0} \]
for $L$ sufficiently large. It follows :
\[ \hat{r}(s)\ge \hat{z}(t_0-s,\theta)+\delta\quad \forall 0\le s \le t_0. \]
It remains to show (\ref{9.41a}). We have 
\begin{eqnarray*}
\dot{\hat{r}}&=&h^4\sqrt{1-\frac{a^2\Delta}{(r^2+a^2)^2}}\sqrt{1-\frac{{\cal K}\Delta}{(r^2+a^2)^2E^2}}
+a\cos\theta\frac{h^4a_0^2}{E}q^*\\
&=&1+\sqrt{1-\frac{a^2\Delta}{(r^2+a^2)^2}}\left(\sqrt{1-\frac{{\cal K}\Delta}{(r^2+a^2)^2E^2}}
-\sqrt{1-\frac{a^2\Delta}{(r^2+a^2)^2}}\right)h^4\\
&+&a\cos\theta {h^4a_0^2}\left(\frac{q^*}{E}-a\cos\theta\right)\\
&\ge&1-{\cal O}(L^{-1})+a\cos\theta h^4a_0^2\left(\frac{q^*}{E}-a\cos\theta\right),
\end{eqnarray*}
where we have used (\ref{9.20cc}). It is therefore sufficient to show :
\begin{eqnarray}
\label{9.28.0}
\left|a\cos\theta\left(\frac{q^*}{E}-a\cos\theta\right)\right|={\cal O}(L^{-1/2}).
\end{eqnarray}
We distinguish two cases :
\begin{eqnarray}
\label{9.28a}
1. \quad\forall 0\le s \le t_0\quad \cos^2\theta(s)>L^{-1/2}.\\
\label{9.28b}
2. \quad \exists 0\le s_0\le t_0\quad \cos^2\theta(s_0)\le L^{-1/2},\\
\label{9.28c}
s_0\neq 0\Rightarrow \cos^2\theta(s)>L^{-1/2}\, \forall 0\le s<s_0.
\end{eqnarray}
We will treat only the second case. The first case can be considered in some sense as 
a special case of the second one with $s_0=t_0$. We first suppose $s_0>0$. From (\ref{9.20cb}), (\ref{9.28c}) we infer :
\[ sign q_0^*=sign (\cos\theta_0) \]
for $L$ sufficiently large and from (\ref{9.28c}) we infer, using  also (\ref{9.20cc}) and (\ref{C10.37}), that 
\[ |q^*|^2\ge a^2E^2L^{-1/2}(1-{\cal O}(L^{-1/2}))>0 \]
if $L$ sufficiently large. Therefore 
\[ sign(q^*(s))=sign\cos\theta(s)\quad \forall 0\le s\le s_0 \]
because neither $q^*(s)$ nor $\cos\theta(s)$ can change its sign on ${[}0,s_0{]}$.
We find with (\ref{9.20cc}), (\ref{C10.37}) :
\begin{eqnarray*}
\forall 0\le s\le s_0\quad \frac{q^*}{E}&=&sign(\cos\theta)\sqrt{a^2\cos^2\theta+{\cal O}(L^{-1})}\\
\Rightarrow \left|\frac{q^*}{E}-a\cos\theta\right|&=&{\cal O}(L^{-1/2})\quad\forall 0\le s\le s_0.
\end{eqnarray*}
We now claim that :
\begin{eqnarray}
\label{9.28e}
\forall s_0\le s\le t_0\quad \cos^2\theta(s)\le C(t_0)L^{-1/2}.
\end{eqnarray}
Indeed using (\ref{10.33}), (\ref{9.20cc}) and (\ref{C10.37}) we can estimate 
\begin{eqnarray*}
\frac{d}{dt}a^2\cos^2\theta&=&-2a^2\cos\theta\sin\theta\dot{\theta}\le 2aA_1|\cos\theta|\sqrt{\frac{\cal K}{E^2}-a^2\sin^2\theta}\\
&\le&2a^2A_1\cos^2\theta+{\cal O}(L^{-1/2}),
\end{eqnarray*}
where $A_1:=\max_{r\ge r_+}h^4a_0^2a.$ By the Gronwall lemma we obtain :
\begin{eqnarray*}
\cos^2\theta&\le&\int_{s_0}^se^{2A_1(s-\tau)}{\cal O}(L^{-1/2})d\tau+e^{2A_1(s-s_0)}\cos^2\theta(s_0)\\
&\le&\tilde{C}(t_0)L^{-1/2}.
\end{eqnarray*}
(\ref{9.28e}) follows and therefore :
\[ |a\cos\theta|={\cal O}(L^{-1/4}),\,\left|\frac{q^*}{E}\right|={\cal O}(L^{-1/4}),\]
i.e. (\ref{9.28.0}). If $s_0=0$ we can start with (\ref{9.28e}).  
\qed
\section{Proof of Proposition \ref{prop10.1}}

Let us first show the following lemma :
\begin{lemma}
\label{lem10.12}
\begin{eqnarray*}
&\forall N,R>0,\, \exists C_{N,R}>0,t_0>0\,\\
& \forall t_{\epsilon}\ge t_0\,\exists \delta=\delta(t_{\epsilon},N,R)>0,
T_0=T_0(N,R,t_{\epsilon})>0\, \forall 0\le s\le t_{\epsilon},\,T\ge T_0\\
&||(1-\phi_{\delta})(.-\hat{z}(t_{\epsilon}-s,\theta))e^{-is H}F^T_{t_{\epsilon}}||\le C_{N,R} e^{-\kappa_+t_{\epsilon}}.
\end{eqnarray*}
\end{lemma}
{\bf Proof.}

Because of the finite propagation speed we can replace $1-\phi_{\delta}(.-\hat{z}(t_{\epsilon}-s,\theta))$
by 

$\chi_{\epsilon}(1-\phi_{\delta}(.-\hat{z}(t_{\epsilon}-s,\theta)))$ with $\chi_{\epsilon}\in C^{\infty}(\R)$
s.t.
\begin{eqnarray*}
\chi_{\epsilon}(\hat{r})=\left\{\begin{array}{cc} 0 & \hat{r}\le -4t_{\epsilon} \\ 1 & \hat{r}\ge-3t_{\epsilon}. \end{array} \right.
\end{eqnarray*}
By the results of the preceding sections it is sufficient to show :
\[ \chi_{\epsilon}(1-\phi_{\delta})(.-\hat{z}(t_{\epsilon}-s,\theta))e^{-isH_d}Op(N^*_{\leftarrow})F^T_{t_{\epsilon}}\rightarrow\infty,\quad T\rightarrow\infty. \]
Here we have used (\ref{DiagError}) and that 
\[ {[}\chi_{\epsilon}(1-\phi_{\delta})(.-\hat{z}(t_{\epsilon}-s,\theta)),Op(W)M{]}\in \Psi^{-1}(\Sigma_1). \]
Using Lemma \ref{lem10.5} we can replace $Op(N^*_{\leftarrow})F^T_{t_{\epsilon}}$ by $Op\left(N_{\leftarrow}^*\chi_L\left(\left|\frac{\xi}{q}\right|\right)\right)F^T_{t_{\epsilon}}$
\footnote{As in the definition of $W(\hat{r},\theta,\xi,q)$, $\chi(h^2\sqrt{k'^2\xi^2+a_0^2(l'\xi+q)^2})
\chi_L\left(\left|\frac{\xi}{q}\right|\right)$ can be extended to a smooth function.}
with $\chi_L\in C^{\infty}(\R),$
\[ \chi_L(x)=\left\{\begin{array}{cc} 1 & |x|\ge 2L \\ 0 & |x|\le L \end{array}
\right. \]

Let $\tilde{\chi}_{\epsilon}\in C_0^{\infty}(\R)$ and 
\begin{eqnarray*}
\tilde{\chi}_{\epsilon}&=&1\quad\mbox{on}\quad \supp F_{t_{\epsilon}}^T,\\
\tilde{\chi}_{\epsilon}&=&0\quad\mbox{on}\quad (-\infty,-t_{\epsilon}-L^{-1})
\end{eqnarray*}
for all $T\ge T_0$. By \cite[Proposition 18.1.26]{Hoe1} we have for all $f\in L^2(\R\times {[}0,\pi{]})$:
\begin{eqnarray*}
\lefteqn{WF\left(Op\left(j_+\left(\frac{\xi}{|q|}\right)\chi\left(h^2\sqrt{k'^2\xi^2+a_0^2(l'\xi+q)^2}\right)
\chi_L\left(\left|\frac{\xi}{q}\right|\right)\right)\tilde{\chi}_{\epsilon}f\right)}\\
&\subset&\{(\hat{r},\omega;\xi,q);\, \hat{r}\ge -t_{\epsilon}-L^{-1},|\xi|\ge L |q|,\xi>0\},\\
\lefteqn{WF\left(Op\left(j_-\left(\frac{\xi}{|q|}\right)\chi\left(h^2\sqrt{k'^2\xi^2+a_0^2(l'\xi+q)^2}\right)
\chi_L\left(\left|\frac{\xi}{q}\right|\right)\right)\tilde{\chi}_{\epsilon}f\right)}\\
&\subset&\{(\hat{r},\omega;\xi,q);\, \hat{r}\ge -t_{\epsilon}-L^{-1},|\xi|\ge L |q|,\xi<0\},\\
\end{eqnarray*}
where $WF$ denotes the wave front set (see \cite[Chapter VIII]{Hoe1}). Then by the classical results of 
propagation of singularities (see e.g. \cite[Theorem 26.1.1]{Hoe1}), Sobolev embeddings and
Lemma \ref{lem10.7}, we find that the operators 
{\small
\begin{eqnarray*}
&\chi_{\epsilon}(1-\phi_{\delta})(.-\hat{z}(s,\theta))e^{-isH_d}Op\left(j_+\left(\frac{\xi}{|q|}\right)\chi\left(h^2\sqrt{k'^2\xi^2+a_0^2(l'\xi+q)^2}\right)
\chi_L\left(\left|\frac{\xi}{q}\right|\right)\right)\tilde{\chi}_{\epsilon}P_{14},\\
&\chi_{\epsilon}(1-\phi_{\delta})(.-\hat{z}(s,\theta))e^{-isH_d}Op\left(j_-\left(\frac{\xi}{|q|}\right)
\chi\left(h^2\sqrt{k'^2\xi^2+a_0^2(l'\xi+q)^2}\right)\chi_L\left(\left|\frac{\xi}{q}\right|\right)\right)\tilde{\chi}_{\epsilon}P_{23}
\end{eqnarray*}}
are compact for $L,\, t_{\epsilon}$ sufficiently large. The lemma now follows from Corollary \ref{cor10.4}
and the observation that :
\[ (F^T_{t_{\epsilon}})_{2,3}=0,\quad (Op((N^-_{\leftarrow})^*)F^T_{t_{\epsilon}})_{1,4}=0,\quad
\tilde{\chi}_{\epsilon}F^T_{t_{\epsilon}}=F^T_{t_{\epsilon}}. \]
for $T$ sufficiently large.
\qed
\vspace{0.5cm}

{\large \bf Proof of Proposition \ref{prop10.1}}

We first show that:
\begin{eqnarray}
\label{10.47}
||(U(0,t_{\epsilon})-\phi_{\delta}(.-\hat{z}(0,\theta)))e^{-it_{\epsilon}H})F_{t_{\epsilon}}^T||\le C_{N,R} t_{\epsilon}e^{-\kappa_+ t_{\epsilon}}
\end{eqnarray}
uniformly in $T$ large for $\delta$ sufficiently small. 
Let $\phi_{\delta}(.,t)=\phi_{\delta}(.-\hat{z}(t,\theta))$. We define for $g=P_{1,4}g\in {\cal H}:$
\[ v(t,\hat{r},\omega)=(U(t,t_{\epsilon})-\phi_{\delta}(t)e^{i(t-t_{\epsilon})H})g. \]
If $\supp g\subset(\hat{z}(t_{\epsilon},\theta)+\delta,\infty)\times S^2$, then $v$ is a solution of :
\begin{eqnarray*}
\left. \begin{array}{rcl} v(t_{\epsilon},\hat{r},\omega)&=&0,\\
\partial_tv&=&iH_tv+h(t)\quad \hat{r}>\hat{z}(t,\theta),\\
\sum_{\mu \in \{t,\hat{r},\omega\}}{\cal N}_{\mu}\hat{\gamma}^{\mu}v(t,\hat{z}(t,\theta),\omega)&=&-iv(t,\hat{z}(t,\theta),\omega),
\end{array} \right\}
\end{eqnarray*}
where 
\[ h(t)=-\left(\left(\frac{d}{dt}\phi_{\delta}(t)\right)+{[}\phi_{\delta},H{]}\right)e^{i(t-t_{\epsilon})H}
g=B(t)g. \]
Thus :
\[ v(0)=\int_{t_{\epsilon}}^0U(0,s)B(s)gds. \]
We obtain by Lemma \ref{lem10.12} for $\delta$ small enough :
\[ ||(U(0,t_{\epsilon})-\phi_{\delta}(0)e^{-it_{\epsilon}H})F_{\epsilon}^T||\le \int_0^{t_{\epsilon}}
||B(s)F_{t_{\epsilon}}^T||ds\le C_{N,R} t_{\epsilon}e^{-\kappa_+t_{\epsilon}} \]
uniformly in $T$ large. We now estimate for $\epsilon>0$ given :
\begin{eqnarray}
\label{9.30ca}
\lefteqn{||{\cal J}(\hat{r})U(0,T)f-\phi_{\delta}(.-\hat{z}(0,\theta))e^{-it_{\epsilon}H}F^T_{t_{\epsilon}}||}\nonumber\\
&\le&||{\cal J}(\hat{r})U(0,t_{\epsilon})(U(t_{\epsilon},T)f-U_{\leftarrow}(t_{\epsilon},T)\Omega_{\leftarrow}^-f)||\nonumber\\
&+&||{\cal J}(\hat{r})U(0,t_{\epsilon})(U_{\leftarrow}(t_{\epsilon},T)(\Omega_{\leftarrow}^-f)^N_R-F^T_{t_{\epsilon}})||
+||(\Omega_{\leftarrow}^-f)^N_R-\Omega_{\leftarrow}^-f||\nonumber\\
&+&||{\cal J}(\hat{r})(U(0,t_{\epsilon})-\phi_{\delta}(.-\hat{z}(0,\theta))e^{-it_{\epsilon}H})F^T_{t_{\epsilon}}||\nonumber\\
&<&||{\cal J}(\hat{r})U(0,t_{\epsilon})(1-{\cal J}(\hat{r}+t_{\epsilon}))U(t_{\epsilon},T)f||\nonumber\\
&+&||{\cal J}(\hat{r})U(0,t_{\epsilon})(1-{\cal J}(\hat{r}+t_{\epsilon}))U_{\leftarrow}(t_{\epsilon},T)\Omega_{\leftarrow}^-f||
+\epsilon/3
\end{eqnarray}
fixing first $N,R$ and choosing then $t_{\epsilon}, T$ sufficiently large and $\delta=\delta(t_{\epsilon},N,R)$
sufficiently small. Here we have used Proposition \ref{prop9.1}, Lemma
\ref{lem10.3} and (\ref{10.47}). We now claim that for $t_{\epsilon}$ fixed the 
first two terms in (\ref{9.30ca}) go to zero when $T$ goes to infinity. Indeed 
\begin{eqnarray*}
\lefteqn{{\cal J}(\hat{r})U(0,t_{\epsilon})(1-{\cal J}(\hat{r}+t_{\epsilon}))U(t_{\epsilon},T)f}\\
&=&{\cal J}(\hat{r})U(0,t_{\epsilon})(1-{\cal J}(\hat{r}+t_{\epsilon}))e^{i(t_{\epsilon}-T)H}f\\
&=&{\cal J}(\hat{r})U(0,t_{\epsilon})(1-{\cal J}(\hat{r}+t_{\epsilon}))e^{i(t_{\epsilon}-T)H}{\bf 1}_{(-\infty,0{]}}(P^-)f\\
&+&{\cal J}(\hat{r})U(0,t_{\epsilon})(1-{\cal J}(\hat{r}+t_{\epsilon}))e^{i(t_{\epsilon}-T)H}{\bf 1}_{{[}0,\infty)}(P^-)f.
\end{eqnarray*}
We can suppose $f=\chi(H)f,\, \chi\in C_0^{\infty}(\R),\, \supp\chi\subset\{-m+\eta^n,m+\eta^n\}$.
By the minimal velocity estimate we can replace ${\bf 1}_{(-\infty,0)}(P^-)$ by ${\bf 1}_{(-\infty,-\epsilon_{\chi})}(P^-)$
and ${\bf 1}_{{[}0,\infty)}(P^-)$ by ${\bf 1}_{{[}1-\tilde{\epsilon},\infty)}(P^-)$,
where $\tilde{\epsilon}>0$ and $\epsilon_{\chi}$ is given by Lemma \ref{lem7.1a}. But,
\begin{eqnarray*}
\lefteqn{\lim_{T\rightarrow\infty}{\cal J}(\hat{r})U(0,t_{\epsilon})(1-{\cal J}(\hat{r}+t_{\epsilon}))e^{i(t_{\epsilon}-T)H}{\bf 1}_{(-\infty,-\epsilon_{\chi}{]}}(P^-)f}\\
&=&\lim_{T\rightarrow\infty}{\cal J}(\hat{r})U(0,t_{\epsilon})(1-{\cal J}(\hat{r}+t_{\epsilon}))
{\bf 1}_{(-\infty,-\epsilon_{\chi}{]}}\left(\frac{\hat{r}}{t_{\epsilon}-T}\right)e^{i(t_{\epsilon}-T)H}f=0
\end{eqnarray*}
by the finite propagation speed for $U(0,t_{\epsilon}).$ In a similar way we find
\begin{eqnarray*}
\lefteqn{\lim_{T\rightarrow\infty}{\cal J}(\hat{r})U(0,t_{\epsilon})(1-{\cal J}(\hat{r}+t_{\epsilon}))e^{i(t_{\epsilon}-T)H}{\bf 1}_{{[}1-\tilde{\epsilon},\infty)}(P^-)f}\\
&=&\lim_{T\rightarrow\infty}{\cal J}(\hat{r})U(0,t_{\epsilon})(1-{\cal J}(\hat{r}+t_{\epsilon}))
{\bf 1}_{{[}1-\tilde{\epsilon},\infty)}\left(\frac{\hat{r}}{t_{\epsilon}-T}\right)e^{i(t_{\epsilon}-T)H}=0
\end{eqnarray*}
because $(1-{\cal J}(\hat{r}+t_{\epsilon}))
{\bf 1}_{{[}1-\tilde{\epsilon},\infty)}\left(\frac{\hat{r}}{t_{\epsilon}-T}\right)=0$
for $T$ sufficiently large. The second term can be treated in a similar way. 
This concludes the proof of the proposition.
\qed
\chapter{Proof of the main theorem}
\label{mainth}
In this chapter we put the results of the previous chapters together to prove the main theorem.
Recall that we are working with the operators $H^{\nu,n}_{\eta^n}$ etc. and that 
the indices are suppressed. All operators are considered as acting on ${\cal H}^n$.
\section{The energy cut-off}
\label{mainth1}
In this section $N,R>0$ will be fixed.
Let 
\[\Sigma_0^-=\{(\hat{r},\omega);\, \hat{r}\le \hat{z}(0,\theta)\},\,{\cal H}_0^-=(L^2(\Sigma_0^-,d\hat{r}d\omega))^4.\]
On ${\cal H}_0^-$ we define 
\begin{eqnarray*}
H_0^-&=&H,\\
D(H_0^-)&=&\{u\in {\cal H}_0^-;\, H u\in {\cal H}_0^-;\,
\Sigma_{\hat{\mu}\in\{(t,\hat{r},\theta,\varphi)\}}{\cal N}_{\hat{\mu}}\hat{\gamma}^{\hat{\mu}}\Psi(\hat{z}(0,\theta),\omega)=-i\Psi(0,\hat{z}(0,\theta),\omega)\}.
\end{eqnarray*}
We need the following
\begin{lemma}
\label{lem11.1}
$(i)$ Let $\phi\in C^{\infty}(\R\times S^2),\partial_{\hat{r}}\phi\in C_0^{\infty}(\R\times S^2),\phi\equiv 0 $ on $\Sigma_0^-,
\,\chi\in {\cal S}^0(\R)$. Then 
$(\chi(H_0^-\oplus H_0)-\chi(H))\phi$ is compact.

$(ii)$ Let $\phi\in C_0^{\infty}(\Sigma_0),\, \chi\in {\cal S}^{-1}(\R)$. 
Then $\chi(H_0^-\oplus H_0)\phi,\, \chi(H)\phi$ are
compact.
\end{lemma}
{\bf Proof}

$(i)$ We write $\chi(x)=\hat{\chi}(x)(x+i),\,
\hat{\chi}\in {\cal S}^{-1}$. We have as an identity between bounded operators~:
\begin{eqnarray}
\label{11.1}
\chi(H_0^-\oplus H_0)\phi=\hat{\chi}(H_0^-\oplus H_0)(H+i)\phi.
\end{eqnarray}
Therefore :
\begin{eqnarray*}
(\chi(H_0^-\oplus H_0)-\chi(H))\phi
&=&(\hat{\chi}(H_0^-\oplus H_0)-\hat{\chi}(H))(H+i)\phi\\
&=&(\hat{\chi}(H_0^-\oplus H_0)-\hat{\chi}(H))\tilde{\phi}(H+i)\phi\\
&+&(\hat{\chi}(H_0^-\oplus H_0)-
\hat{\chi}(H)){[}H,\tilde{\phi}{]}\phi,
\end{eqnarray*}
where $\tilde{\phi}\in C^{\infty}(\R\times S^2),\, \tilde{\phi}\phi=\phi,\tilde{\phi}\equiv 0$ on 
$\Sigma_0^-$. We only show that the first term is compact, the second term can be treated
in a similar manner. Let $\tilde{\chi}$ 
be an almost analytic extension of $\hat{\chi}$ with 
\begin{eqnarray}
\label{11.2}
\tilde{\chi}|_{\R}=\hat{\chi},\,\forall N\, |\bar{\partial}\tilde{\chi}(z)|\le C_N|\Im z|^N
\langle \Re z \rangle^{-2-N},\, \supp\tilde{\chi}\subset\{x+iy;\, |y|\le C\langle x \rangle \}.
\end{eqnarray}
By the Helffer-Sj\"ostrand formula we can write :
\begin{eqnarray}
\label{10.2w}
\lefteqn{(\hat{\chi}(H_0^-\oplus H_0)-\hat{\chi}(H))\tilde{\phi}(H+i)}\nonumber\\   
&=&\int\bar{\partial}\tilde{\chi}(z)((z-H_0^-\oplus H_0)^{-1}-(z-H)^{-1})\tilde{\phi}(H+i)dz\wedge d\bar{z}. 
\end{eqnarray}
This identity is first understood as an identity between operators $D(H)\rightarrow {\cal H}$.
Now,
\begin{eqnarray*}
\lefteqn{\bar{\partial}\tilde{\chi}(z)((z-H_0^-\oplus H_0)^{-1}-(z-H)^{-1})\tilde{\phi}(H+i)}\\
&=&\bar{\partial}\tilde{\chi}(z)((z-H_0^-\oplus H_0)^{-1}(z-H)\tilde{\phi}(z-H)^{-1}
+(z-H_0^-\oplus H_0)^{-1}{[}H,\tilde{\phi}{]}(z-H)^{-1}\\
&-&\tilde{\phi}(z-H)^{-1}-(z-H)^{-1}{[}H,\tilde{\phi}{]}(z-H)^{-1})(H+i)\\
&=&\bar{\partial}\tilde{\chi}(z)((z-H_0^-\oplus H_0)^{-1}{[}H,\tilde{\phi}{]}(z-H)^{-1}-(z-H)^{-1}{[}H,\tilde{\phi}{]}(z-H)^{-1})(H+i).
\end{eqnarray*}
Both terms are compact and can be estimated by $C\langle x \rangle^{-3}$ according
to (\ref{11.2}). This shows first that the identity (\ref{10.2w}) can be extended 
to an identity between bounded operators and then that the operator is compact.

$(ii)$ The fact that $\chi(H)\phi$ is compact follows from the estimates (\ref{6.21}), (\ref{6.22})
and the same arguments as in \cite[Corollary 4.2]{HN}.
This entails that $\chi(H_0^-\oplus H_0)\phi$ is compact by part $(i)$ of the lemma.
\qed
\\

Let now $\delta=\delta(t_{\epsilon})$ and $\phi_{\delta}$ be as in Proposition \ref{prop10.1}, 
$\tilde{\phi}_{\delta}(\hat{r},\theta)=\phi_{\delta}(\hat{r}-\hat{z}(0,\theta))$
 and $\chi_1,\chi_2\in C^{\infty}(\R)$ s.t. 
\begin{eqnarray*}
\chi_1(x)=\left\{\begin{array}{cc} 1 & |x|\le 1+\eta,\\ 0 & |x|\ge 2+\eta \end{array}
\right., \quad \chi_2(x)=\left\{\begin{array}{cc} 0 & x\le 1+\eta \\ 1 & x\ge 2+\eta \end{array}\right.
\end{eqnarray*}
and s.t. for all $x\ge 0\quad \chi_1(x)+\chi_2(x)=1.$ Here and in the following $\eta=\eta^n$.
Then we have by Lemma \ref{lem11.1} and Corollary \ref{cor10.4} :
\begin{eqnarray}
\label{11.3}
\lim_{T\rightarrow\infty}||\chi_1(H_0)\tilde{\phi}_{\delta}e^{-it_{\epsilon}H}F^T_{t_{\epsilon}}||_0^2=0.
\end{eqnarray}
Indeed we can replace $\tilde{\phi}_{\delta}$ by a compactly supported function
using the support property (\ref{10.2}) of $F^T_{t_{\epsilon}}$.  
Furthermore we have :
\begin{eqnarray}
\label{11.4}
\lefteqn{||\chi_2(H_0)\tilde{\phi}_{\delta}e^{-it_{\epsilon}H}F^T_{t_{\epsilon}}||_0-||\chi_1(H_0)\tilde{\phi}_{\delta}e^{-it_{\epsilon}H}F^T_{t_{\epsilon}}||_0}\nonumber\\ 
&\le& ||{\bf 1}_{{[}\eta,\infty)}(H_0)\tilde{\phi}_{\delta}e^{-it_{\epsilon}H}F^T_{t_{\epsilon}}||_0\nonumber\\
&\le& ||\chi_1(H_0)\tilde{\phi}_{\delta}e^{-it_{\epsilon}H}F^T_{t_{\epsilon}}||_0\nonumber\\
&+&||\chi_2(H_0)\tilde{\phi}_{\delta}e^{-it_{\epsilon}H}F^T_{t_{\epsilon}}||_0.
\end{eqnarray}
Using (\ref{11.3}), (\ref{11.4}) and Lemma \ref{lem11.1} we obtain :
\begin{eqnarray}
\label{SCO}
\lim_{T\rightarrow\infty}||{\bf 1}_{{[}\eta,\infty)}(H_0)\tilde{\phi}_{\delta}e^{-it_{\epsilon}H}F^T_{t_{\epsilon}}||_0
&=&\lim_{T\rightarrow\infty}||\chi_2(H_0)\tilde{\phi}_{\delta}e^{-it_{\epsilon}H}F^T_{t_{\epsilon}}||_0\nonumber\\
&=&\lim_{T\rightarrow\infty}||\chi_2(H_0^-\oplus H_0)\tilde{\phi}_{\delta}e^{-it_{\epsilon}H}(0\oplus F^T_{t_{\epsilon}})||\nonumber\\
&=&\lim_{T\rightarrow\infty}||\chi_2(H)\tilde{\phi}_{\delta}e^{-it_{\epsilon}H}F^T_{t_{\epsilon}}||.
\end{eqnarray}
\begin{lemma}
\label{lem10.2b}
Let $\chi\in {\cal S}^0(\R), \, f\in D(\langle \notD_{S^2} \rangle )$. Let $\tilde{\chi}$ be an 
almost analytic extension of $\chi$ with :
\[ \tilde{\chi}|_{\R}=\chi, \quad \forall N \, |\bar{\partial}\tilde{\chi}(z)|\le C_N |\Im z|^N\langle \Re z \rangle ^{-1-N},\, 
\supp \tilde{\chi}\subset \{ x+iy;|y|\le C\langle \Re z \rangle \}. \]
Then we have :
\begin{eqnarray}
\label{Dcutoff}
(\chi(H)-\chi(H_{\leftarrow}))f=\int \bar{\partial}\tilde{\chi}(z)(z-H)^{-1}(H-H_{\leftarrow})(z-H_{\leftarrow})^{-1}f dz \wedge d\bar{z}.
\end{eqnarray}
\end{lemma}
\begin{remark}
Note that neither 
\begin{eqnarray*}
\int\bar{\partial}\tilde{\chi}(z)(z-H)^{-1}f dz\wedge d\bar{z} \quad \mbox{nor} \quad
\int \bar{\partial}\tilde{\chi}(z)(z-H_{\leftarrow})^{-1} f dz \wedge d\bar{z}
\end{eqnarray*}
is convergent in ${\cal H}$, but the R.H.S. of (\ref{Dcutoff}) is.
\end{remark}
{\bf Proof.}

Let $\chi_0\in C_0^{\infty}(\R)$ with $\chi_0=1$ in a neighborhood of $0$ and 
$\chi_0^mù(x)=\chi_0\left(\frac{x}{m}\right).$ Let $\tilde{\chi}_0\in C_0^{\infty}(\C)$
be an almost analytic extension of $\chi_0$ with :
\[ |\bar{\partial}\tilde{\chi}_0|\le C_N |\Im z|^N\quad \forall N. \]
Then $\tilde{\chi}_0\left(\frac{z}{m}\right)$ is an almost analytic extension of $\chi_0^m$. 
Clearly 
\[ \lim_{m\rightarrow \infty}\chi(H)\chi_0^m(H)f=\chi(H)f,\, \lim_{m\rightarrow\infty}\chi(H_{\leftarrow})\chi_0^m(H_{\leftarrow})f
=\chi(H_{\leftarrow})f. \]
We have 
\begin{eqnarray}
\label{Apcutoff}
\chi(H)\chi_0^m(H)=\int\bar{\partial}\left(\tilde{\chi}(z)\tilde{\chi}_0\left(\frac{z}{m}\right)\right)(z-H)^{-1}dz \wedge d\bar{z}
\end{eqnarray}
and the R.H.S. of (\ref{Apcutoff}) is convergent in norm for all $m$. Therefore :
\begin{eqnarray*}
\lefteqn{(\chi(H)\chi_0^m(H)-\chi(H_{\leftarrow})\chi_0^m(H_{\leftarrow}))f}\\
&=&\int\bar{\partial}\left(\tilde{\chi}(z)\tilde{\chi}_0\left(\frac{z}{m}\right)\right)
(z-H)^{-1}(H-H_{\leftarrow})(z-H_{\leftarrow})^{-1} f dz \wedge d\bar{z}\\
&=&\int(\bar{\partial}\tilde{\chi}(z))\tilde{\chi}_0\left(\frac{z}{m}\right)
(z-H)^{-1}(H-H_{\leftarrow})(z-H_{\leftarrow})^{-1} f dz \wedge d\bar{z}\\
&+&\int\tilde{\chi}(z)(\bar{\partial}\tilde{\chi}_0)\left(\frac{z}{m}\right)\frac{1}{m}
(z-H)^{-1}(H-H_{\leftarrow})(z-H_{\leftarrow})^{-1} f dz \wedge d\bar{z}=:I_1+I_2.\\
\end{eqnarray*}
We have :
\[ ||\bar{\partial}\tilde{\chi}(z)\tilde{\chi}_0\left(\frac{z}{m}\right)(z-H)^{-1}(H-H_{\leftarrow})(z-H_{\leftarrow})^{-1}f||
\lesssim \langle x \rangle ^{-3}||\langle D_{S^2} \rangle f||\]
uniformly in $m$ and for all $z\in \C\setminus(\sigma(H)\cup\sigma(H_{\leftarrow})$:
\begin{eqnarray*}
\lefteqn{(\bar{\partial}\tilde{\chi}(z))\tilde{\chi}_0\left(\frac{z}{m}\right)
(z-H)^{-1}(H-H_{\leftarrow})(z-H_{\leftarrow})^{-1} f}\\
&&\rightarrow 
(\bar{\partial}\tilde{\chi}(z))(z-H)^{-1}(H-H_{\leftarrow})(z-H_{\leftarrow})^{-1} f.
\end{eqnarray*}
Then by the Lebesgue Theorem $I_1$ converges to the R.H.S. of (\ref{Dcutoff}).
The change of coordinates $u=\frac{z}{m}$ gives :
\[ I_2=\int m\tilde{\chi}(um)\bar{\partial}\tilde{\chi}_0(u)(um-H)^{-1}(H-H_{\leftarrow})(um-H_{\leftarrow})^{-1} f du\wedge d\bar{u}. \]
As $\tilde{\chi}(z)$ is bounded we can estimate :
\begin{eqnarray*}
\lefteqn{||m\tilde{\chi}(um)\bar{\partial}\tilde{\chi}_0(u)(um-H)^{-1}(H-H_{\leftarrow})(um-H_{\leftarrow})^{-1}f||}\\
&\lesssim& |\Im u|^2\frac{m}{m^2|\Im u|^2}||\langle D_{S^2} \rangle f||=
\frac{||\langle \notD_{S^2}\rangle f||}{m}\rightarrow 0,\, m\rightarrow \infty. 
\end{eqnarray*}
Thus 
\[I_2\rightarrow 0,\, m\rightarrow \infty. \]
\qed
\begin{lemma}
\label{lem11.2}
We have :
\begin{eqnarray*}
&&\forall \epsilon>0\, \exists t_0>0\, \forall t_{\epsilon}\ge t_0,\, \exists T_0=T_0(t_{\epsilon}),\, \forall T\ge T_0\\
&&|||\chi_2(H)\tilde{\phi}_{\delta(t_{\epsilon})}e^{-it_{\epsilon}H}F^T_{t_{\epsilon}}||-||{\bf 1}_{{[}0,\infty)}(H_{\leftarrow})F^T_{t_0}|||<\epsilon.
\end{eqnarray*}
Here $\delta(t_{\epsilon})$ is chosen as in Proposition \ref{prop10.1}. 
\end{lemma}
{\bf Proof.}
Using Lemma \ref{lem10.12} we see that :
\begin{eqnarray*}
&&\forall \epsilon>0\, \exists t_0>0\, \forall t_{\epsilon}\ge t_0,\, \exists T_0=T_0(t_{\epsilon}),\, \forall T\ge T_0\\
&&|||\chi_2(H)\tilde{\phi}_{\delta(t_{\epsilon})}e^{-it_{\epsilon}H}F^T_{t_{\epsilon}}||-
||\chi_2(H)e^{-it_{\epsilon}H}F^T_{t_{\epsilon}}|||<\epsilon.
\end{eqnarray*}
We then have to show that :
\begin{eqnarray}
\label{11.5}
(\chi_2(H)-\chi_2(H_{\leftarrow}))F^T_{t_{\epsilon}}\rightarrow 0,\, T\rightarrow 0.
\end{eqnarray}
To this purpose let $\tilde{\chi}_2$ be an almost analytic extension of $\chi_2$ with :
\[ \tilde{\chi}_2|_{\R}=\chi,\,\forall N\, |\bar{\partial}\tilde{\chi}_2(z)|\le C_N|\Im z|^N
\langle \Re z \rangle^{-1-N},\, \supp\tilde{\chi}_2\subset\{x+iy;\, |y|\le C\langle x \rangle\}. \]
Then we have by Lemma \ref{lem10.2b} :
\begin{eqnarray*}
\lefteqn{(\chi_2(H)-\chi_2(H_{\leftarrow}))F_{t_{\epsilon}}^T}\\
&=&\int\bar{\partial}\tilde{\chi}_2(z)(z-H)^{-1}
(H-H_{\leftarrow})(z-H_{\leftarrow})^{-1}F_{t_{\epsilon}}^Tdz\wedge d\bar{z}\\
&=&\int\bar{\partial}\tilde{\chi}_2(z)(z-H)^{-1}(P_{\omega}+W)\langle \notD_{S^2} \rangle^{-1}
(z-H_{\leftarrow})^{-1}\langle \notD_{S^2}\rangle F^T_{t_{\epsilon}}dz\wedge d\bar{z}.
\end{eqnarray*}
By Lemma \ref{lem10.5} we find :
\begin{eqnarray}
\label{11.6}
(z-H_{\leftarrow})^{-1}\langle \notD_{S^2} \rangle F^T_{t_{\epsilon}}\rightarrow 0,\, T\rightarrow \infty
\end{eqnarray}
and we have 
\begin{eqnarray}
\label{11.8}
||\bar{\partial}\tilde{\chi}_2(z)(z-H)^{-1}(P_{\omega}+W)(z-H_{\leftarrow})^{-1}F^T_{t_{\epsilon}}||\lesssim\langle x \rangle^{-3}||\langle \notD_{S^2} \rangle F^T_{t_{\epsilon}}||.
\end{eqnarray}
Equations (\ref{11.6}), (\ref{11.8}) give (\ref{11.5}). Therefore we have :
\begin{eqnarray*}
\lim_{T\rightarrow\infty}||\chi_2(H)e^{-it_{\epsilon}H}F^T_{t_{\epsilon}}||
&=&\lim_{T\rightarrow\infty}||\chi_2(H_{\leftarrow})e^{-it_{\epsilon}H_{\leftarrow}}F^T_{t_{\epsilon}}||\\
&=&\lim_{T\rightarrow\infty}||\chi_2(H_{\leftarrow})F^T_0||=\lim_{T\rightarrow\infty}||{\bf 1}_{{[}0,\infty)}(H_{\leftarrow})F^T_0||.
\end{eqnarray*}
In order to replace $\chi_2(H_{\leftarrow})$ by ${\bf 1}_{{[}0,\infty)}(H_{\leftarrow})$ we use the same arguments as before.
\qed

\section{The term near the horizon}
\label{mainth2}
We first compute the radiation explicitly for the asymptotic dynamics :
\begin{lemma}
\label{lem12.1}
We have :
\[ \lim_{T\rightarrow\infty}||{\bf 1}_{{[}0,\infty)}(H_{\leftarrow})F_0^T||^2
=\langle \tilde{f},e^{\frac{2\pi}{\kappa_+}H_{\leftarrow}}\left(1+e^{\frac{2\pi}{\kappa_+}H_{\leftarrow}}\right)^{-1}
\tilde{f} \rangle. \]
\end{lemma}
{\bf Proof.}

The proof is analogous to the proof of \cite[Lemma VI.6]{Ba6}. We repeat it here 
for the convenience of the reader. Let ${\cal F}$ be the Fourier transform with respect to $\hat{r}$. We have~:
{\footnotesize
\begin{eqnarray*}
\lefteqn{||{\bf 1}_{{[}0,\infty)}(H_{\leftarrow})F_0^T||^2}\\
&=&\int_{S^2}\int_0^{\infty}|{\cal F}(F_0^T)(\xi)|^2d\xi d\omega\\
&=&\sum_{j=1}^2\lim_{\epsilon\rightarrow 0^+}\int_{S^2}\hat{A}(\theta)\kappa_+\int_0^{\infty}
\left|\int_{\R}e^{i(\hat{A}(\theta)+i\epsilon)\xi e^{\kappa_+y}}e^{\frac{\kappa_+}{2}y}\tilde{f}_jdy\right|^2d\xi d\omega\\
&=&\sum_{j=1}^2\lim_{\epsilon\rightarrow 0^+}\int_{S^2}\frac{\hat{A}(\theta)\kappa_+}{2}
\int_{\R\times \R}\frac{1}{\epsilon\cosh\left(\frac{\kappa_+}{2}(y_1-y_2)\right)
-i\hat{A}(\theta)\sinh\left(\frac{\kappa_+}{2}(y_1-y_2)\right)}\tilde{f}_j(y_1)\bar{\tilde{f}}_j(y_2)dy_1dy_2\\
&=&\sum_{j=1}^2\lim_{\epsilon\rightarrow 0^+}\int_{S^2}\frac{\hat{A}(\theta)\kappa_+}{4\pi}
\int_{\R}|{\cal F}(\tilde{f}_j)(\xi,\omega)|^2{\cal F}\left(\frac{1}{\epsilon\cosh\left(\frac{\kappa_+x}{2}\right)
-i\hat{A}(\theta)\sinh\left(\frac{\kappa_+x}{2}\right)}\right)\left(-\xi\right)d\xi d\omega.
\end{eqnarray*}}
Now given $\epsilon\neq 0,\xi<0$ and $N,\,M>0,$ we evaluate
\[ \oint l(x)dx,\, l(x):=\frac{e^{-ix\xi}}{\epsilon\cosh\left(\frac{\kappa_+x}{2}\right)
-i\hat{A}(\theta)\sinh\left(\frac{\kappa_+x}{2}\right)}, \]
along the path
\[ \{-N\le \Re x\le N,\, \Im x=0,M\}\cup\{0\le \Im x\le M,\, \Re x=\pm N\}. \]
First we have :
\begin{eqnarray*}
\left|\int_{\pm N}^{\pm N+iM}l(x)dx\right|\lesssim e^{-\frac{\kappa_+}{2}N}\int_0^{\infty}e^{x\xi}dx\rightarrow 0,\, N\rightarrow \infty,\\
\left|\int_{-N+iM}^{N+iM}l(x)dx\right|\lesssim e^{M\xi}\int_{-\infty}^{\infty}e^{-\frac{\kappa_+}{2}|x|}dx\rightarrow0,\, M\rightarrow \infty.\\
\end{eqnarray*}
We deduce that
\[ \int_{-\infty}^{\infty}l(x)dx=2i\pi\sum_{n=1}^{\infty}\rho_n(\epsilon), \]
where $\rho_n(\epsilon)$ are the residues of $l(x)$ at the poles 
$z_n(\epsilon)\in\{z\in \C;\, \Im z>0\}$. We easily check that :
\begin{eqnarray*}
&z_n(\epsilon)=\frac{2i}{\kappa_+}\left(n\pi-arctan\left(\frac{\epsilon}{\hat{A}(\theta)}\right)\right),\\
&\sup_{1\le n}|\rho_n(\epsilon)-\frac{2i}{\hat{A}(\theta)\kappa_+}(-1)^ne^{\frac{2n\pi}{\kappa_+}\xi}|\lesssim \epsilon,
\end{eqnarray*}
hence we get that for $\xi<0$ we have :
\[ \left|{\cal F}\left(\frac{1}{\epsilon\cosh\left(\frac{\kappa_+x}{2}\right)
-i\hat{A}(\theta)\sinh\left(\frac{\kappa_+x}{2}\right)}\right)(\xi)
-\frac{4\pi}{\hat{A}(\theta)\kappa_+}e^{\frac{2\pi}{\kappa_+}\xi}\left(1+e^{\frac{2\pi}{\kappa_+}\xi}\right)^{-1}\right|\lesssim \epsilon. \]
In the same manner, for $\xi>0$, choosing $M<0$ and considering the poles\\
$z_n(\epsilon)\in\{z\in \C;\, \Im z<0\}$ we obtain :
\begin{eqnarray*}
&\int_{-\infty}^{\infty}l(x)dx=2i\pi\sum_{n=0}^{-\infty}\rho_n(\epsilon),\\
&\sup_{n\le 0}|\rho_n(\epsilon)-\frac{2i}{\hat{A}(\theta)\kappa_+}(-1)^ne^{-\frac{2n\pi}{\kappa_+}\xi}|\lesssim \epsilon,\\
&\left|{\cal F}\left(\frac{1}{\epsilon\cosh\left(\frac{\kappa_+x}{2}\right)
-i\hat{A}(\theta)\sinh\left(\frac{\kappa_+x}{2}\right)}\right)(\xi)
-\frac{4\pi}{\hat{A}(\theta)\kappa_+}\left(1+e^{-\frac{2\pi}{\kappa_+}\xi}\right)^{-1}\right|\lesssim \epsilon.
\end{eqnarray*}
Eventually we conclude that
\[ {\cal F}\left(\frac{1}{0-i\hat{A}(\theta)\sinh\left(\frac{\kappa_+}{2}x\right)}\right)(\xi)=\frac{4\pi}{\kappa_+}e^{\frac{2\pi}{\kappa_+}\xi}
\left(1+e^{\frac{2\pi}{\kappa_+}\xi}\right)^{-1}, \]
and 
\begin{eqnarray*}
||{\bf 1}_{{[}0,\infty)}(H_{\leftarrow})F_0^T||^2&=&\int_{S^2}\int_{\R}e^{-\frac{2\pi}{\kappa_+}\xi}\left(1+e^{-\frac{2\pi}{\kappa_+}\xi}\right)^{-1}|{\cal F}(\tilde{f}_j)(\xi,\omega)|^2d\xi d\omega\\
&=&\langle \tilde{f},e^{\frac{2\pi}{\kappa_+}H_{\leftarrow}}\left(1+e^{\frac{2\pi}{\kappa_+}H_{\leftarrow}}\right)^{-1}\tilde{f}\rangle.
\end{eqnarray*}
\qed
\begin{proposition}
\label{prop11.3}
We have :
\[ \lim_{T\rightarrow\infty}||{\bf 1}_{{[}\eta,\infty)}(H_0){\cal J}U(0,T)f||^2_0
=\langle \Omega_{\leftarrow}^-f,e^{\frac{2\pi}{\kappa_+} H_{\leftarrow}}\left(1+e^{\frac{2\pi}{\kappa_+}H_{\leftarrow}}\right)^{-1}\Omega_{\leftarrow}^-f \rangle. \]
\end{proposition}
{\bf Proof.}

For $\epsilon>0$ given we estimate :
\begin{eqnarray*}
\lefteqn{|||{\bf 1}_{{[}\eta,\infty)}(H_0){\cal J}U(0,T)f||_0^2
-\langle \Omega_{\leftarrow}^-f,e^{\sigma H_{\leftarrow}}\left(1+e^{\sigma H_{\leftarrow}}\right)^{-1}\Omega_{\leftarrow}^-f\rangle|}\nonumber\\
&\lesssim&||(\Omega_{\leftarrow}^-f)^N_R-\Omega_{\leftarrow}^-f||||\Omega_{\leftarrow}^-f||\\
&+&|||{\bf 1}_{{[}\eta,\infty)}(H_0){\cal J}U(0,T)f||_0^2-||{\bf 1}_{{[}\eta,\infty)}(H_0)\tilde{\phi}_{\delta}e^{it_{\epsilon}H}F^T_{t_{\epsilon}}||^2|\nonumber\\
&+&|||{\bf 1}_{{[}\eta,\infty)}(H_0)\tilde{\phi}_{\delta}e^{it_{\epsilon}H}F^T_{t_{\epsilon}}||^2-||{\bf 1}_{{[}0,\infty)}(H_{\leftarrow})F^T_0||^2|\\
&+&|||{\bf 1}_{{[}0,\infty)}(H_{\leftarrow})F^T_0||^2
-\langle \tilde{f},e^{\sigma H_{\leftarrow}}\left(1+e^{\sigma H_{\leftarrow}}\right)^{-1}\tilde{f}\rangle |=:I_1+I_2+I_3+I_4.
\end{eqnarray*}
Here $\delta=\delta(t_{\epsilon})$ as in Proposition \ref{prop10.1}.
We first choose $R_0,N_0>0$ s.t.
\[ \forall N\ge N_0, R\ge R_0\quad I_1<\epsilon/4. \]
Using Proposition \ref{prop10.1} we can fix $N\ge N_0,\, R\ge R_0,\,t_{\epsilon}>0,\, \delta=\delta(t_{\epsilon})>0$
and $T_0>0$ s.t. 
\[ \forall T\ge T_0\quad I_2<\epsilon/4. \]
By choosing $T_0$ possibly larger we have :
\[ \forall T\ge T_0\quad I_3<\epsilon/4. \]
Here we have used Lemma \ref{lem11.2} and (\ref{SCO}).
Note that the same $t_{\epsilon}$ can be chosen for the estimate of $I_2,I_3$. Indeed
in both cases we use Lemma \ref{lem10.12}.
The parameters $N,R,t_{\epsilon},\delta$ being fixed in this way we conclude by noting 
that $I_4\rightarrow 0,\, T\rightarrow \infty.$
\qed

\section{Proof of Theorem \ref{thkey2}}
\label{mainth3}
We start with the following identity :
\begin{eqnarray}
\label{12.2}
\lefteqn{||{\bf 1}_{{[}\eta,\infty)}(H_0)U(0,T)f||_0^2}\nonumber\\
&=&||{\bf 1}_{{[}\eta,\infty)}(H_0){\cal J}U(0,T)f||_0^2
+||{\bf 1}_{{[}\eta,\infty)}(H_0)(1-{\cal J})U(0,T)f||_0^2\nonumber\\
&+&2\Re \langle {\bf 1}_{{[}\eta,\infty)}(H_0)(1-{\cal J})U(0,T)f,{\bf 1}_{{[}\eta,\infty)}(H_0){\cal J}U(0,T)f\rangle.
\end{eqnarray}
In order to prove that 
\[ \langle {\bf 1}_{{[}\eta,\infty)}(H_0)(1-{\cal J})U(0,T)f,{\bf 1}_{{[}\eta,\infty)}(H_0){\cal J}U(0,T)f\rangle
\rightarrow 0,\, T\rightarrow \infty \]
it is enough, by the results of the preceding chapters, to prove 
\[ \langle {\bf 1}_{{[}\eta,\infty)}(H_0)(1-{\cal J})U(0,T)f,{\bf 1}_{{[}\eta,\infty)}(H_0)\tilde{\phi}_{\delta} e^{-it_{\epsilon}H}F_{t_{\epsilon}}^T\rangle
\rightarrow 0,\, T\rightarrow \infty. \]
Note that 
\[ {\bf 1}_{{[}\eta,\infty)}(H_0)(H_0+i)^{-1}\tilde{\phi}_{\delta}e^{-it_{\epsilon}H}F^T_{t_{\epsilon}}\rightarrow 0,\, T\rightarrow\infty. \]
Indeed
\begin{eqnarray*}
\lefteqn{{\bf 1}_{{[}\eta,\infty)}(H_0)(H_0+i)^{-1}\tilde{\phi}_{\delta}e^{-it_{\epsilon}H}F^T_{t_{\epsilon}}}\\
&=&{\bf 1}_{{[}\eta,\infty)}(H_0)(H_0+i)^{-1}\tilde{\phi}_{\delta}(H+i)e^{-it_{\epsilon}H}
(H+i)^{-1}F^T_{t_{\epsilon}}\rightarrow 0, T\rightarrow\infty. 
\end{eqnarray*}
It is therefore sufficient to show that for all $f\in {\cal H}^1$ the following limit 
exists :
\[ \lim_{T\rightarrow\infty}(H_0+i)e^{iTH_0}(1-{\cal J})e^{-iTH}f. \]
But,
\begin{eqnarray*}
\lefteqn{(H_0+i)e^{iTH_0}(1-{\cal J})e^{-iTH}f}\\
&=&-e^{-iTH_0}{[}H,{\cal J}{]}e^{-iTH}f+e^{iTH_0}(1-{\cal J})e^{-iTH}(H+i)f\\
&=&-e^{iTH_0}\Gamma^1{\cal J}'e^{-iTH}f+e^{iTH_0}(1-{\cal J})e^{-iTH}(H+i)f.
\end{eqnarray*}
The first term goes to zero because $\sigma_{sc}(H)=\emptyset$ and the second term possesses
a limit by Lemma \ref{lem7.4}. By the same lemma we find :
\begin{eqnarray*}
\lim_{T\rightarrow\infty}||{\bf 1}_{{[}\eta,\infty)}(H_0)(1-{\cal J})U(0,T)f||
&=&||{\bf 1}_{{[}\eta,\infty)}(H_0)W_0^-f||=||W_0^-{\bf 1}_{{[}\eta,\infty)}(H)f||\\
&=&||{\bf 1}_{\R^-}(P^-){\bf 1}_{{[}\eta,\infty)}(H)f||\\
&=&||{\bf 1}_{{[}\eta,\infty)}(H){\bf 1}_{\R^-}(P^-)f||.
\end{eqnarray*}
This concludes the proof of the theorem.
\qed
\begin{appendix}
\chapter{Proof of Proposition \ref{prop6.3}}
\label{AppA}
We will work with the $(r_*,\omega)$ coordinate system. For technical reasons we need to
fix the angular momentum $D_{\varphi}=n$. 
Let $d\tilde{\omega}=\sin\theta d\theta,$
\begin{eqnarray*}
\tilde{\cal H}^n&=&\left\{u=e^{in\varphi}v;\, v\in (L^2(\R\times {[}0,\pi{]},dr_*d\tilde{\omega}))^4\right\},\\
\tilde{\cal H}_t^n&=&\left\{u=e^{in\varphi}v;\, v\in (L^2(\tilde{\Sigma}_t^{col,\varphi},dr_*d\tilde{\omega}))^4\right\},\\
\tilde{\Sigma}_t^{col,\varphi}&=&\left\{(r_*,\theta)\in \R\times {[}0,\pi{]};\, r_*\ge z(t,\theta)\right\}.
\end{eqnarray*}
We define :
\begin{eqnarray}
\label{A.1a}
\tilde{\notD}^{\nu,n}&:=&h\tilde{\notD}_{\gs}^{\nu,n}h+V^n,\\
\label{A.1b}
\tilde{\notD}^{\nu,n}_{\gs}&:=&\Gamma^1D_{r_*}+a_0\Gamma^2(D_{\theta}+\frac{\cot\theta}{2})+a_0\Gamma^3\frac{n}{\sin\theta}+b_0\Gamma^{\nu}+c^n,\\
\tilde{\cal H}^{n1}&=&D(\tilde{\notD}^{\nu,n})=\{u\in\tilde{\cal H}^n;\tilde{\notD}^{\nu,n}u\in\tilde{\cal H}^n\}\nonumber
\end{eqnarray}
Recall that the singularity $\sin\theta=0$ appearing in the expression for $\tilde{\notD}^{\nu,n}_{\gs}$
is a coordinates singularity.
The operator $\tilde{\notD}^{\nu,n}_t$
is the operator acting on $\tilde{\cal H}^n_t$ with formal expression $\tilde{\notD}^{\nu,n}$ and
domain:
\begin{eqnarray*}
\tilde{\cal H}^{n1}_t&:=&D(\tilde{\notD}_t^{\nu,n})\\
&=&\left\{u\in \tilde{\cal H}_t^n;\, \tilde{\notD}_t^{\nu,n}u\in \tilde{\cal H}^n_t,\,
\sum_{\mu\in \{t,r_*,\theta,\varphi\}}{\cal N}_{\mu}\gamma^{\mu}u(z(t,\theta),\theta)
=-iu(z(t,\theta),\theta)\right\}.
\end{eqnarray*}
The extension ${[}..{]}^*_H$ is defined in an analogous way to Section \ref{sec3.5}
as an extension from $\tilde{\cal H}^{n1}_t$ to $\tilde{\cal H}^{n1}$.
We consider the following problem :
\begin{eqnarray}
\label{A.2}
\left. \begin{array}{rcl} \partial_t\Psi&=&i\tilde{\notD}^{\nu,n}_t\Psi,\quad r_*>z(t,\theta) \\
\sum_{\mu\in \{t,r_*,\theta,\varphi\}}{\cal N}_{\mu}\gamma^{\mu}\Psi(t,z(t,\theta),\theta)
&=&-i\Psi(t,z(t,\theta),\theta),\\
\Psi(t=s,.)&=&\Psi_s(.)\in D(\tilde{\notD}^{\nu,n}_s). \end{array} \right\}
\end{eqnarray} 
Proposition \ref{prop6.3} follows from 
\begin{proposition}
\label{propA.1}
Let $\Psi_s\in D(\tilde{\notD}_s^{\nu,n}).$ Then there exists a unique solution 
\[{[}\Psi(.){]}^*_H={[}\tilde{U}^{\nu,n}(.,s)\Psi_s{]}^*_H\in C^1(\R;\tilde{\cal H}^n)\cap C(\R;\tilde{\cal H}^{n1})\] 
of (\ref{A.2}) s.t. for all $t\in \R$ $\Psi(t)\in D(\tilde{\notD}^{\nu,n}_t)$. Furthermore we have
$||\Psi(t)||=||\Psi_s||$ and $\tilde{U}^{\nu,n}(t,s)$ possesses an extension to an isometric
and strongly continuous propagator from $\tilde{\cal H}^n_s$ to $\tilde{\cal H}^n_t$
s.t. for all $\Psi_s\in D(\tilde{\notD}_s^{\nu,n})$ we have :
\[  \frac{d}{dt}\tilde{U}^{\nu,n}(t,s)\Psi_s=i \tilde{\notD}^{\nu,n}_t\tilde{U}^{\nu,n}(t,s)\Psi_s.
\]
\end{proposition}
Note that in Proposition \ref{prop6.3}, there is an additional statement about the finite
propagation speed which follows from ($\Psi$ is supposed to be a solution of (\ref{6.13})) 
\begin{eqnarray*}
\frac{d}{dt}\int_{S^2}\int_{R+|t-s|}^{\infty}|\Psi(t,\hat{r},\omega)|^2d\hat{r}d\omega
&=&2{\bf 1}_{\R^+}(s-t)\int_{S^2}(|\Psi_2|^2+|\Psi_3|^2)d\omega\\
&-&2{\bf 1}_{\R^+}(t-s)\int_{S^2}(|\Psi_1|^2+|\Psi_4|^2)d\omega.
\end{eqnarray*}
\\

\underline{\large Proof of Proposition \ref{propA.1}}

We will drop the indices $\nu,n$ in what follows. 
\begin{itemize}
\item 
Let us first show uniqueness. If $\Psi$ is solution of (\ref{A.2}), then we have :
\begin{eqnarray*}
\frac{d}{dt}\int_0^{\pi}\int_{z(t,\theta)}^{\infty}|\Psi|^2dr_*d\tilde{\omega}&=&-\int_0^{\pi}\dot{z}(t,\theta)|\Psi|^2(z(t,\theta),\theta)d\omega\\
&+&2\Re\langle i\tilde{\notD}_t\Psi,\Psi\rangle=2\Re\langle i\check{\notD}_t\Psi,\Psi\rangle=0.
\end{eqnarray*}
Here $\check{\notD}_t:=\dot{z}D_{r_*}+\tilde{\notD}_t$ is selfadjoint with domain 
$\tilde{\cal H}_t^{n1}$ (see Lemma \ref{Selfadj}).
\item
Let us now prove existence. We introduce the operators :
\begin{eqnarray*}
R(t)&=&(N^2+Z_1^2-Z_2^2)^{-1/2}\left(\begin{array}{cccc} N & 0 & -Z_1 & -Z_2 \\
0 & N & -Z_2 & -Z_1 \\ Z_1 & -Z_2 & N & 0\\ -Z_2 & Z_1 & 0 & N \end{array} \right),\\
N&=&w^{-1},\, Z_1=h^2+\dot{z},\, Z_2=a_0h^2\partial_{\theta}z,\\
{\cal T}&:&{\cal H}_0\ni f\mapsto {\cal T}(t)f\in {\cal H}_t,\\
{[}{\cal T}(t)f{]}(r_*,\omega)&=&f(r_*-z(t,\theta)+z(0,\theta),\omega).
\end{eqnarray*}
We remark that 
\begin{eqnarray*}
R^{-1}(t)&=&(N^2+Z_1^2-Z_2^2)^{-1/2}\left(\begin{array}{cccc} N & 0 & Z_1 & Z_2 \\
0 & N & Z_2 & Z_1 \\ -Z_1 & Z_2 & N & 0\\ Z_2 & -Z_1 & 0 & N \end{array} \right).
\end{eqnarray*}
Furthermore we notice that 
\begin{eqnarray*}
{\cal T}&\in&C^1(\R_t;{\cal L}((C_0^1(\R_{r_*}\times {[}0,\pi{]}))^4,(C_0^0(\R_{r_*}\times {[}0,\pi{]}))^4),\\
\dot{\cal T}(t)&=&-\dot{z}(t,\theta){\cal T}(t)\partial_{r_*}.
\end{eqnarray*}
Then $u$ is a solution of our problem iff $w(t)=R^{-1}(t){\cal T}^{-1}(t)u$ is solution of
\begin{eqnarray}
\label{A.1}
\left. \begin{array}{rcl} \partial_tw&=&iA(t)w,\quad r_*>z(0,\theta)\\
w_2(t,z(0,\theta),\omega)&=&w_3(t,z(0,\theta),\omega)=0, \end{array} \right\}
\end{eqnarray}
where $A(t)=R^{-1}(t){\cal T}^{-1}(t)(\notD+\dot{z}D_{r_*})R(t){\cal T}(t)-R^{-1}(t)\dot{R}(t)$.
We first need to analyze $A(t)$ for fixed $t$. We put
\[ D(A(t))=\{u\in \tilde{\cal H}_0;\, A(t)u\in \tilde{\cal H}_0,\,u_2(z(0,\theta),\omega)=u_3(z(0,\theta),\omega)=0\}. \]
We equip $D(A(t))$ with the graph norm of $A(t)$. The operator $(A(t),D(A(t))$
is selfadjoint (see Lemma \ref{Selfadj}). Let 
\[ d_{\theta}:=D_{\theta}+\frac{\cot\theta}{2i}. \]
We first note that the domain is independent of $t$ :
\begin{eqnarray*}
\lefteqn{\forall t\, D(A(t))=D}\\
&:=&\{u\in \tilde{\cal H}_0;\, D_{r_*}u\in \tilde{\cal H}_0,\, d_{\theta}u\in \tilde{\cal H}_0, u_2(z(0,\theta),\theta)
=u_3(z(0,\theta),\theta)=0\}
\end{eqnarray*}
and the graph norm of $A(t)$ is equivalent to the norm 
\[ ||u||_D=||D_{r_*}u||+||a_0d_{\theta}u||+||u||. \]
This follows from the following estimates :
\begin{eqnarray}
\label{A.3}
\forall u\in D(A(t))\quad ||D_{r_*}u||&\lesssim& ||A(t)u||+||u||,\\
\label{A.4}
\forall u\in D(A(t))\quad ||a_0d_{\theta}u||&\lesssim& ||A(t)u||+||u||.
\end{eqnarray}
To show (\ref{A.3}), (\ref{A.4}) we introduce the operator $\tilde{\notD}^t$ which is
obtained from $\tilde{\notD}$ by evaluating all functions in the definition of $\tilde{\notD}$
at $(r_*+z(t,\theta)-z(0,\theta),\theta)$. Then we have 
\[ A(t)=R^{-1}(t)(\tilde{\notD}^t+\dot{z}D_{r_*})R(t)-R^{-1}(t)\dot{R}(t). \]
For $u$ regular enough we estimate 
\begin{eqnarray*}
||A(t)u||&\gtrsim& ||\tilde{\notD}^tR(t)u||-\max_{\theta}|\dot{z}(t,\theta)|||D_{r_*}R(t)u||-C||u||\\
&\ge& ||\tilde{\notD}^tR(t)u||-(1-\delta)||D_{r_*}R(t)u||-C||u||
\end{eqnarray*}
for some $0<\delta<1$. We henceforth drop the subscript $t$. Let $v:=R(t)u$. The first term on the R.H.S. in the above inequality
can be estimated in the following way :
\[ ||\tilde{\notD}^tv||^2\ge ||h^2D_{r_*}v||^2+||h^2a_0(D_{\theta}+\frac{\cot\theta}{2})v||^2+2\Re I-C||v||^2 \]
with 
\[ I=\int_0^{\pi}\int_{z(t,\theta)}^{\infty}\langle ha_0\Gamma^2(D_{\theta}+\frac{\cot\theta}{2})hv,h\Gamma^1D_{r_*}hv\rangle dr_*d\tilde{\omega}. \]
We have :
\begin{eqnarray*}
2I&=&-\int_0^{\pi}\langle ha_0\Gamma^2\left(\partial_{\theta}+\frac{\cot\theta}{2}\right)hv,\Gamma^1h^2v\rangle(z(t,\theta),\theta)d\tilde{\omega}\\
&-&\int_0^{\pi}\int_{z(t,\theta)}^{\infty}\langle \Gamma^1\partial_{r_*}h^2a_0\Gamma^2\left(\partial_{\theta}+\frac{\cot\theta}{2}\right)hv,hv\rangle dr_*d\tilde{\omega}\\
&+&\int_0^{\pi}(\partial_{\theta}z)\langle hv,ha_0\Gamma^2h\Gamma^1\partial_{r_*}hv\rangle (z(t,\theta),\theta) d\tilde{\omega}\\
&-&\int_0^{\pi}\int_{z(t,\theta)}^{\infty}\langle hv,a_0\Gamma^2\left(\partial_{\theta}+\frac{\cot\theta}{2}\right)h^2\Gamma^1\partial_{r_*}hv\rangle dr_*d\tilde{\omega}.
\end{eqnarray*}
Therefore :
\begin{eqnarray*}
2 \Re I&=&-\Re\int_0^{\pi}\langle h a_0\Gamma^2\left(\partial_{\theta}+\frac{\cot\theta}{2}\right)hv,\Gamma^1h^2v\rangle (z(t,\theta),\theta)d\tilde{\omega}\\
&+&\Re \int_0^{\pi}(\partial_{\theta}z)\langle hv,h^2a_0\Gamma^2\Gamma^1\partial_{r_*}hv\rangle (z(t,\theta),\theta)d\tilde{\omega}\\
&-&2\Re \int_0^{\pi}\int_{z(t,\theta)}^{\infty}\left\langle \left\{h\Gamma^1D_{r_*}h,ha_0\Gamma^2\left(D_{\theta}+\frac{\cot\theta}{2i}\right)h\right\}v,hv\right\rangle dr_*d\tilde{\omega}\\
&=:&B+AC.
\end{eqnarray*}
We first estimate the boundary term. We have :
\[ B= \Re \int_0^{\pi}\langle ha_0\Gamma^2 hv,\Gamma^1 \tilde{\partial}_{\theta}((hv)(z(t,\theta),\theta))\rangle (z(t,\theta),\theta)d\tilde{\omega}, \]
where $\tilde{\partial}_{\theta}=\partial_{\theta}+\frac{\cot\theta}{2}.$
Recalling that 
\[ hv(z(t,\theta),\theta)=\left(\tilde{N}h\left(\begin{array}{c} Nu_1-Z_2u_4 \\
-Z_1u_4 \\ Z_1 u_1 \\ -Z_2 u_1+N u_4 \end{array} \right) \right)(z(t,\theta),\theta)\] 
with $\tilde{N}=(N^2+Z_1^2-Z_2^2)^{-1/2}$ we find :
{\footnotesize
\begin{eqnarray*}
B&\ge&-C||u||_{(L^2({[}0,\pi{]})^4}\\
&+&\int_0^{\pi}\tilde{N}^2\left\langle h
\left(\begin{array}{c} Nu_1-Z_2u_4 \\ -Z_1u_4 \\ Z_1 u_1 \\ -Z_2u_1+Nu_4 \end{array} \right),
h^3a_0\Gamma^2\Gamma^1\left(\begin{array}{c} N\tilde{\partial}_{\theta}(u_1)-Z_2\tilde{\partial}_{\theta}(u_4) \\
-Z_1\tilde{\partial}_{\theta}(u_4) \\ Z_1 \tilde{\partial}_{\theta}(u_1) \\ -Z_2\tilde{\partial}_{\theta}(u_1)+N\tilde{\partial}_{\theta}(u_4) \end{array} \right) \right\rangle
\sin\theta d\theta.
\end{eqnarray*}}
Here $\tilde{\partial}_{\theta}(u_j)$ stands for $\tilde{\partial}_{\theta}(u_j(z(t,\theta),\theta))$.
The term under the second integral in the above inequality equals 

\begin{eqnarray*}
\lefteqn{\tilde{N}^2\left\langle h^4a_0\left(\begin{array}{c} Nu_1-Z_2u_4 \\ -Z_1u_4 \\ Z_1 u_1 \\ -Z_2u_1+Nu_4\end{array} \right),
\left(\begin{array}{c} Z_1\tilde{\partial}_{\theta}(u_4) \\ N\tilde{\partial}_{\theta}(u_1)-Z_2\tilde{\partial}_{\theta}(u_4) \\
Z_2\tilde{\partial}_{\theta}(u_1)-N\tilde{\partial}_{\theta}(u_4) \\ Z_1\tilde{\partial}_{\theta}(u_1) \end{array} \right) \right\rangle}\\
&=&\tilde{N}^2(h^4a_0(NZ_1u_1\tilde{\partial}_{\theta}(\bar{u}_4)-Z_2u_4Z_1\tilde{\partial}_{\theta}(\bar{u}_4)-Z_1u_4N\tilde{\partial}_{\theta}(\bar{u}_1)
+Z_1Z_2u_4\tilde{\partial}_{\theta}(\bar{u}_4)\\
&+&Z_1Z_2u_1\tilde{\partial}_{\theta}(\bar{u}_1)-Z_1u_1N\tilde{\partial}_{\theta}(\bar{u}_4)
-Z_2Z_1u_1\tilde{\partial}_{\theta}(\bar{u}_1)+NZ_1u_4\tilde{\partial}_{\theta}(\bar{u}_1)))=0.
\end{eqnarray*}
Thus 
\[ B\ge -C||u(z(t,\theta),\theta)||_{(L^2({[}0,\pi{]}))^4}. \]
By the usual trace theorems we find :
\[ B\ge -C||u||_{H^{1/2}(\Sigma_t^{col,\varphi})}\ge -\epsilon(||D_{r_*}R(t)u||+||a_0d_{\theta}R(t)u||)-C_{\epsilon}||u||. \]
Let us now consider the anticommutator:
\[ \{h\Gamma^1 D_{r_*}h,ha_0\Gamma^2d_{\theta}h\}=h({[}D_{r_*},h^2a_0{]}\Gamma^1\Gamma^2d_{\theta}+a_0{[}d_{\theta},h^2{]}\Gamma^2\Gamma^1D_{r_*})h. \]
Thus :
\[ AC\ge -\epsilon(||a_0d_{\theta}v||+||D_{r_*}v||)-C_{\epsilon}||v||,\quad \epsilon>0. \]
Gathering everything together and using $h^2\ge 1$ we find :
\begin{eqnarray*}
||A(t)u||&\ge&(1-\tilde{\epsilon})(||D_{r_*}R(t)u||+||a_0d_{\theta}R(t)u||)\\
&-&(1-\delta)
||D_{r_*}R(t)u||-C_{\epsilon}||u||\\
&\ge&(\delta-\tilde{\epsilon})||D_{r_*}u||+(1-\epsilon)||a_0d_{\theta}u||-C_{\epsilon}||u||,\, \tilde{\epsilon}<\delta.
\end{eqnarray*} 
Recall that we have dropped the index $t$. But $a_0\lesssim a_0^t\lesssim a_0.$
The above inequality proves (\ref{A.3}), (\ref{A.4}).
The operators $(A(t),D)$ are selfadjoint and the family $\{(iA(t),D)\}_{t\in \R}$ is a stable family
in the sense of \cite[Definition 5.2.1]{Pa}. We want to check that for $v\in D$, $t\mapsto A(t)v$ is continuously
differentiable. We have 
\begin{eqnarray*}
\lefteqn{\frac{A(t+\delta)-A(t)}{\delta}v}\\
&=&\frac{1}{\delta}(R^{-1}(t+\delta)(\tilde{\notD}^{t+\delta}
+\dot{z}D_{r_*})R(t+\delta)-R^{-1}(t)(\tilde{\notD}^t+\dot{z}D_{r_*})R(t))v\\
&-&(R^{-1}(t)\dot{R}(t))'v+o(\delta)\\
&=&\frac{R^{-1}(t+\delta)-R^{-1}(t)}{\delta}(\tilde{\notD}^{t+\delta}+\dot{z}D_{r_*})R(t+\delta)v\\
&+&\frac{1}{\delta}R^{-1}(t)(\tilde{\notD}^{t+\delta}+\dot{z}D_{r_*})R(t+\delta)v\\
&-&\frac{1}{\delta}R^{-1}(t)(\tilde{\notD}^t+\dot{z}D_{r_*})R(t)v-(R^{-1}(t)\dot{R}(t))'v+o(\delta)\\
&=&\frac{R^{-1}(t+\delta)-R^{-1}(t)}{\delta}(\tilde{\notD}^{t+\delta}+\dot{z}D_{r_*})R(t+\delta)v\\
&+&R^{-1}(t)\frac{\tilde{\notD}^{t+\delta}-\tilde{\notD}^t}{\delta}R(t+\delta)v\\
&+&R^{-1}(t)(\tilde{\notD}^t+\dot{z}D_{r_*})\frac{R(t+\delta)-R(t)}{\delta}v-(R^{-1}(t)\dot{R}(t))'v+o(\delta)\\
&=&:I_1^{\delta}+I_2^{\delta}+I_3^{\delta}-(R^{-1}(t)\dot{R}(t))'v+o(\delta).
\end{eqnarray*}
Here and henceforth a prime denotes a derivative with respect to $t$.
Let 
\begin{eqnarray*}
I_1&:=&(R^{-1}(t))'(\tilde{\notD}^t+\dot{z}D_{r_*})R(t)v,\\
I_2&:=&R^{-1}(t)(\tilde{\notD}^t)'R(t)v,\\
I_3&:=&R^{-1}(t)(\tilde{\notD}^t+\dot{z}D_{r_*})(R)'(t)v.
\end{eqnarray*}
The operator $(\tilde{\notD}^t)'$ is defined by differentiating $\tilde{\notD}^t$ formally with
respect to $t$. Let us first consider the second term. We have :
\begin{eqnarray*}
\lefteqn{||(R^{-1}(t)\frac{\tilde{\notD}^{t+\delta}-\tilde{\notD}^t}{\delta}R(t+\delta)-R^{-1}(t)(\tilde{\notD}^t)'R(t))v||}\\
&\le&||R^{-1}(t)(\tilde{\notD}^t)'(R(t+\delta)-R(t))v||\\
&+&||R^{-1}(t)\left(\frac{{\tilde{\notD}}^{t+\delta}-\tilde{\notD}^t}{\delta}-(\tilde{\notD}^t)'\right)R(t+\delta)v||\\
&\le&||D_{r_*}(R(t+\delta)-R(t))v||+||a_0'd_{\theta}(R(t+\delta)-R(t))v||\\
&+&o(\delta)||R^{-1}(t)d_{\theta}R(t+\delta)v||+o(\delta)||R^{-1}(t)D_{r_*}R(t+\delta)v||\\
&\le& o(\delta)||v||_D.
\end{eqnarray*}
Here we have used that $a_0\ge \epsilon$ on ${[}0,\infty)\times {[}0,\pi{]}$ as well as the estimates 
(\ref{A.3}), (\ref{A.4}). The other terms can be treated in a similar way. We find :
\[ ||I_j^{\delta}-I_j||\le o(\delta)||v||_D,\, 1\le j\le 3. \]
This shows that $t\mapsto A(t)v$ is continuously differentiable. We can therefore 
use \cite[Theorem 5.4.8]{Pa} to find a strongly continuous propagator $S(s,t)$ on 
$\tilde{\cal H}^n_0$ s.t. for $f\in D,\, S(t,s)f\in D$ is a strongly continuous 
differentiable map from $\R_t\times\R_s$ to $\tilde{\cal H}^n_0$ satisfying :
\[ \frac{d}{dt}S(t,s)f=A(t)S(t,s)f,\quad \frac{d}{ds}S(t,s)f=-S(t,s)A(s)f. \]
The propagator 
\[ \tilde{U}(t,s)=R(t){\cal T}(t)S(t,s)R^{-1}(s){\cal T}^{-1}(s) \]
has the required properties.
\end{itemize} 
\qed
\chapter{Penrose Compactification of block $I$}
\label{AppB}
The Penrose compactification is usually obtained by a construction based on the 
PNG's. We present here an analogous construction based on the SNG's. In order to emphasize
the analogy with the PNG construction we name the different coordinate systems that we
introduce as in the PNG case. We will suppose $a>0$. The 
PNG construction is explained in detail in \cite{ONe}.
\section{Kerr-star and Star-Kerr coordinates}
\label{AppB.1}
A part of the construction of this section can be found in \cite{FL} for the Kerr case. 
We will suppose $E=1$. Let $\hat{r}$ be as in Chapter
\ref{sec2}. 
The star-Kerr coordinate system ($\st ,r,\stheta ,\sphi$) is based on
outgoing simple null geodesics. The new coordinates
$\st$, $\stheta$ and $\sphi$ are of the form
\[ \st = t - \hat{r}(r,\theta) \, ,~ \sphi = \varphi - \Lambda (r) \, ,~\stheta=\theta^\sharp(r,\theta)\, , \]
where the function $\Lambda$ is required to satisfy
\[ \Lambda'(r)=\frac{a(2Mr-Q^2)}{B(r)\Delta}, \, B(r)=(r^2+a^2)k'(r_*(r)). \]
The function $\theta^\sharp$ is defined in the following way. Let 
\[ \alpha(r)=-\int_r^{\infty}\frac{a}{B(\nu)}d\nu. \]
For later simplicity we define :
\[ F=1+\tanh(\alpha)\sin \theta^\sharp,\quad G=\tanh(\alpha)+\sin\theta^\sharp. \]
Then the function $\theta^\sharp$ is defined by
\begin{eqnarray}
\label{B.1}
\sin\theta=\frac{G}{F}.
\end{eqnarray}
We note that
\begin{eqnarray}
\cos^2\theta=\frac{\cos^2\stheta}{\cosh^2\alpha F^2}.
\end{eqnarray}
The important property of the coordinate system $(\st,r,\stheta,\sphi)$ is that
\[ \dot{\st}=\dot{\sphi}=\dot{\stheta}=0 \]
along outgoing SNG's with the correct sign of $\theta_0'$. Therefore we have :
\[ N^{a,+}=\frac{r^2+a^2}{\rho^2}k'\partial_r \]
in this coordinate system. We find
\begin{eqnarray}
\label{B.3}
\left. \begin{array}{rcl} d\theta&=&\alpha'\cos\theta dr+\frac{\cos\theta}{\cos\stheta}d\stheta, \\
dt&=&d\st+\frac{\sigma^2}{\Delta B}dr+\frac{a\cos^2\theta}{\cos\stheta}d\stheta,\\
d\varphi&=&d\sphi+\Lambda'(r)dr. \end{array}\right\}
\end{eqnarray}
Using (\ref{B.3}) we can rewrite the metric as 
\begin{eqnarray}
\label{B.4b}
g&=&\left(1+\frac{Q^2-2Mr}{\rho^2}\right)d\st^2+2\frac{\rho^2}{B}d\st dr+2\left(1+\frac{Q^2-2Mr}{\rho^2}\right)
\frac{a\cos\stheta}{F^2\cosh^2\alpha}d\st d\stheta\nonumber\\
&+&2a\frac{(2Mr-Q^2)G^2}{\rho^2F^2}d\st d\sphi-\frac{r^2\rho^2F^2\cosh^2\alpha+a^2(2Mr-Q^2)
\cos^2\stheta}{\rho^2F^4\cosh^4\alpha}d\stheta^2\nonumber\\
&+&2\frac{a^2(2Mr-Q^2)G^2\cos^*\theta}{\rho^2F^4\cosh^2\alpha}d\stheta d\sphi-\frac{\sigma^2G^2}{\rho^2F^2}d\sphi^2.
\end{eqnarray}
The expression (\ref{B.4b}) shows that $g$ can be extended smoothly
across the horizon $\{ r=r_+ \}$. Besides, it does not degenerate
there since its determinant is given by
\[ \det (g) = -\frac{\rho^4 \sin^2 \theta}{F^2\cosh^2\alpha} \]
and does not vanish for $r=r_+$ \footnote{Note however that there is the usual coordinate
singularity at $\sin\theta=0$.}. Thus, we can add the horizon to block
$I$ as a smooth boundary. It is called the past event horizon and given by
\[ \hor^- := \R_{\st} \times \left\{ r=r_+ \right\}_{r} \times
S^2_{\stheta , {\sphi} }. \]
The metric induced by $g$ on hypersurfaces of constant $r$, $g_r$, has determinant
\[ det(g_r)=-\frac{\rho^2\sin^2\theta\Delta}{\cosh^2\alpha F^2} \]
and thus degenerates for $\Delta=0$, i.e. at $\hor^-$. Since $g$ does not degenerate,
it follows that one of the generators of $\hor^-$ is zero, i.e. $\hor^-$ is a null 
hypersurface. Kerr-star coordinates $(\ts,r,\thetas,\phis)$ are constructed using the incoming 
SNG's :
\[ \ts=t+\hat{r}(r,\theta),\, \phis=\varphi+\Lambda(r),\, \thetas=\theta^\sharp(r,\theta). \]
This coordinate system allows to add the future event horizon 
\[ \hor^+=\R_{\ts}\times\{r=r_+\}\times S^2_{\thetas,\phis} \]
as a smooth null boundary to block $I$.
\section{Kruskal-Boyer-Lindquist coordinates}
\label{AppB.2} 
The Kruskal-Boyer-Lindquist coordinate system is a combination of the two Kerr
coordinate systems, modified in such a way that it is regular on both
the future and the past horizons. The time and radial variables are replaced by
\begin{equation}
U = e^{-\kappa_+ \st} \, ,~ V = e^{\kappa_+ \ts} \, ,
\end{equation}
where $\kappa_+$ is the surface gravity at the outer horizon, see (\ref{surfgr}). 
The coordinate $\theta^\sharp=\stheta=\thetas$ is kept unchanged. The longitude function is
defined by
\begin{equation} \label{PhiSharp}
\varphi^\sharp = \varphi - \frac{a}{r_+^2
  +a^2} \, t \, .
\end{equation}
The functions ($U,V,\theta ,
\varphi^\sharp $) form an analytic coordinate system on $\ExtB \cup
\hor^+ \cup \hor^- - \mathit{(axes)}$. In this coordinate
system, we have
\begin{gather}
\ExtB = ]0,+\infty [_U \times ]0,+\infty [_V \times S^2_{\theta^\sharp ,
  \varphi^\sharp} \, ,\nonumber \\
\hor^+ = \{ 0\}_U \times [0,+\infty [_V \times S^2_{\theta^\sharp ,
  \varphi^\sharp} \, , ~ \hor^- = [0,+\infty [_U \times \{ 0\}_V
  \times S^2_{\theta^\sharp ,\varphi^\sharp} \, , \nonumber
\end{gather}
simply because $\ts$ (resp. $\st$) is regular at $\hor^+$
(resp. $\hor^-$), takes all real values on $\hor^+$
(resp. $\hor^-$), and tends to $-\infty$ (resp. $+\infty$) at $\hor^-$
(resp. $\hor^+$). We want to build the crossing sphere $S^2_c=\{U=V=0\}$. To 
this purpose we introduce the function 
\[ L:=\frac{r-r_+}{UV}. \]
\begin{lemma}
\label{studyL}
$(i)$ The functions $r$ and $\sin\theta$ are well defined analytic functions on ${[}0,\infty)_U
\times {[}0,\infty)_V\times S^2_{\theta^{\sharp},\varphi^{\sharp}}.$ 

$(ii)$ The function $L$ extends to a nonvanishing analytic function on ${[}0,\infty)_U\times {[}0,\infty)_V\times
S^2_{\theta^{\sharp},\varphi^{\sharp}}.$
\end{lemma}
{\bf Proof.}
$(i)$ We have 
\[ \left(\begin{array}{c} UV \\ \sin\theta^{\sharp} \end{array} \right)
=\left(\begin{array}{c} e^{2\kappa_+\hat{r}(r,\theta)} \\ \frac{\tanh(\alpha)-\sin\theta}{\tanh(\alpha)\sin\theta-1} \end{array}
\right)=H(r,\sin\theta)=\left(\begin{array}{c} H_{UV}(r,\sin\theta) \\ H_{\theta^{\sharp}}(r,\sin\theta) \end{array}
\right). \]
Using 
\[ e^{2\kappa_+\hat{r}}=(r-r_+)(r-r_-)^{\kappa_+/\kappa_-}e^{2\kappa_+\int_{r_+}^r
\left(\sqrt{1-\frac{a^2\Delta}{(r^2+a^2)^2}}-1\right)\frac{\tau^2+a^2}{\Delta}d\tau}e^{2\kappa_+a\sin\theta} \]
we easily check that 
\[ \left(\begin{array}{cc} \partial_1H_{UV} & \partial_2 H_{UV} \\
\partial_1H_{\theta^{\sharp}} & \partial_2H_{\theta^{\sharp}} \end{array} \right) (r_+,\sin\theta) \]
is invertible. Thus $(r,\sin\theta)=H^{-1}(UV,\sin\theta^{\sharp})$ is well defined in a 
neighborhood of the crossing sphere. 

$(ii)$ We have 
\[ L=(r-r_-)^{-\kappa_+/\kappa_-}e^{-2\kappa_+\int_{r_+}^r
\left(\sqrt{1-\frac{a^2\Delta}{(r^2+a^2)^2}}-1\right)\frac{\tau^2+a^2}{\Delta}d\tau}e^{-2\kappa_+a\sin\theta}. \]
This and the result of $(i)$ give $(ii)$.
\qed

We have :
\begin{eqnarray}
\label{B.4}
\left. \begin{array}{rcl} 
dt&=&\frac{L}{2\kappa_+(r-r_+)}(UdV-VdU),\\ 
dr&=&\frac{(r-r_-)LB}{2\kappa_+\sigma^2}(UdV+VdU)-\frac{B\Delta}{\sigma^2}\frac{a\cos^2\theta}{\cos\theta^{\sharp}}d\theta^\sharp,\\
d\varphi&=&d\varphi^\sharp+\frac{a}{r_+^2+a^2}\frac{L}{2\kappa_+(r-r_+)}(UdV-VdU),\\
d\theta&=&\frac{a(r-r_-)\cos\theta L}{2\kappa_+\sigma^2}(UdV+VdU)+\frac{\cos\theta B^2}{\cos\theta^{\sharp} \sigma^2}
d\theta^\sharp. \end{array} \right\}
\end{eqnarray}
Using (\ref{B.4}) we find that the Kerr-Newman metric in these coordinates takes the form
\begin{eqnarray}
\label{B.5}
g&=&-\frac{2L(r-r_-)}{4\kappa_+^2\rho^2}\left(\frac{\rho_+^4}{(r_+^2+a^2)^2}+\frac{\rho^4B^2}{\sigma^4}\right)dUdV\nonumber\\
&+&\frac{L^2(r-r_-)}{4\kappa_+^2\rho^2(r-r_+)}\left(\frac{\rho_+^4}{(r_+^2+a^2)^2}-\frac{\rho^4B^2}{\sigma^4}\right)(V^2dU^2+U^2dV^2)\nonumber\\
&-&\rho^2\frac{a^2(r-r_-)^2\cos^2\theta L^2}{4\kappa_+^2\sigma^4}(UdV+VdU)^2\nonumber\\
&-&\frac{\sin^2\theta}{\rho^2}\frac{(r+r_+)^2a^2L^2}{(r_+-r_-)^2}(UdV-VdU)^2\nonumber\\
&-&\frac{\sin^2\theta a L}{\kappa_+\rho^2(r_+^2+a^2)}((r-r_-)\rho_+^2+(r^2+a^2)(r+r_+))(UdV-VdU)d\varphi^\sharp\nonumber\\
&-&\frac{\rho^2B^2}{\sigma^2\cosh^2\alpha F^2}(d\theta^\sharp)^2-\frac{\sigma^2\sin^2\theta}{\rho^2}(d\varphi^\sharp)^2,
\end{eqnarray}
where $\rho_+^2=r_+^2+a^2$.
Clearly $\frac{1}{r-r_+}\left(\frac{\rho_+^4}{(r_+^2+a^2)^2}-\frac{\rho^4B^2}{\sigma^4}\right)$
extends to an analytic function on ${[}0,\infty)_U\times{[}0,\infty)_V\times S^2_{\theta^\sharp,\varphi^\sharp}.$
The expression (\ref{B.5}) then shows that $g$ is smooth on $\ExtB \cup \hor^+ \cup \hor^-$ and can be
extended smoothly on $[0,+\infty [_U \times [0,+\infty [_V \times
S^2_{\theta^\sharp ,\varphi^\sharp }$. The crossing sphere $S^2_c=\{U=V=0\}$ is a regular surface in the extended
space-time
\[ \left( \KB := [0,+\infty [_U \times [0,+\infty [_V \times
  S^2_{\theta^\sharp , \varphi^\sharp} \, ,~g \right).
\]
Hence, the Kruskal-Boyer-Lindquist coordinates give us a
global description of the horizon
\[ \hor = \hor^- \cup S^2_c \cup \hor^+ = \left( [0,+\infty [_U
\times \{ 0\}_V \times S^2_{\theta^\sharp ,\varphi^\sharp} \right) \cup
\left( \{ 0\}_U \times [0,+\infty [_V \times S^2_{\theta^\sharp
,\varphi^\sharp} \right) \]
as a union of two smooth null boundaries $\hor^+ \cup
S^2_c$ and $S^2_c \cup \hor^-$.

\section{Penrose compactification of Block $I$}
\label{AppB.3}
The Penrose compactification of the exterior of a Kerr-Newman black hole is
performed using two independent and symmetric constructions, one based on
star-Kerr, the other on Kerr-star coordinates. We describe explicitly
only the first of these two constructions.

Future null infinity is defined as the set of limit points of outgoing
simple null geodesics as $r\rightarrow +\infty$. This rather
abstract definition of a 3-surface, describing the congruence of
outgoing simple null geodesics, can be given a precise meaning
using star-Kerr coordinates. We consider the expression (\ref{B.4b}) of
the Kerr-Newman metric in star-Kerr coordinates and replace the variable $r$
by $w =1/r$. In these new variables, the exterior of the black hole is
described as
\[ \ExtB = \R_{\st} \times \left] 0 , \frac{1}{r_+} \right[_w \times
S^2_{\stheta ,\sphi } \, .\]
The conformally rescaled metric
\begin{equation} \label{RescaledMetric}
\hat{g} = \Omega^2 g \, ,~\Omega =w = \frac{1}{r}
\end{equation}
takes the form
\begin{eqnarray*}
\hat{g}& = & \left( w^2 + \frac{Q^2w^4-2Mw^3}{1+a^2 w^2 \cos^2 \theta}
\right) \d \st^2\\
&-&2\frac{1+a^2w^2\cos^2\theta}{\sqrt{(1+a^2w^2)^2-a^2w^2(1+(a^2+Q^2)w^2-2Mw)}}d\st dw\\
&+&2w^2\left(1+\frac{Q^2w^2-2Mw}{1+a^2w^2\cos^2\theta}\right)\frac{a\cos\stheta}{F^2 \cosh^2\alpha}d\st d\stheta\\
&+&2aw^3\frac{(2M-Q^2w)G^2}{(1+a^2\cos^2\theta w^2)F^2}d\st d\sphi\\
&-&\frac{(1+a^2\cos^2\theta w^2)F^2\cosh^2\alpha+a^2(2Mw^3-Q^2w^4)\cos^2\stheta}{(1+a^2\cos^2\theta w^2)F^4 \cosh^4\alpha}d\stheta^2\\
&+&2aw^3\frac{(2M-Q^2w)G^2\cos\stheta}{(1+a^2\cos^2\theta w^2)\cosh^2\alpha F^4}d\stheta d\sphi\\
&-&\frac{((1+a^2w^2)^2-a^2w^2(1+(a^2+Q^2)w^2-2Mw)\sin^2\theta)G^2}{(1+a^2\cos^2\theta w^2)F^2}\sin^2\theta d\sphi^2.
\end{eqnarray*}
The functions $\sin\theta,\, \cos\theta, \, F, \, G, \, \alpha$ have to be understood
as functions of $\frac{1}{w},\stheta$. It is clear from the formulas in Section \ref{AppB.1}
that they possess analytic extensions to ${[}0,\frac{1}{r_+}{]}\times S^2_{\stheta}$. The 
expression above shows that $\hat{g}$ can be extended smoothly on the domain
\[ \R_{\st} \times \left[ 0 , \frac{1}{r_+} \right]_w \times
S^2_{\stheta ,\sphi } \, .\]
The hypersurface
\[ \scri^+ := \R_{\st} \times \left\{ w=0 \right\} \times S^2_{\stheta ,
  \sphi } \]
can thus be added to the rescaled space-time as a smooth hypersurface,
describing future null infinity as defined above. This hypersurface is
indeed null since
\[ \hat{g}_{|_{w=0}} = -\d \stheta^2 -\sin^2 \theta \, \d \sphi^2 \]
is degenerate (recall that $\scri^+$ is a 3-surface) and
\[ \det \left( \hat{g} \right) = -\frac{(
  1 + a^2 w^2 \cos^2 \theta)^2}{F^2\cosh^2\alpha} \sin^2 \theta \]
does not vanish for $w=0$.

Similarly, using Kerr-star instead of star-Kerr coordinates, we
define past null infinity, the set of limit points as
$r\rightarrow +\infty$ of incoming simple null geodesics, as
\[ \scri^- := \R_{\ts} \times \left\{ w=0 \right\} \times S^2_{\thetas
  , \phis } \, . \]
\begin{figure}
\centering\epsfig{figure=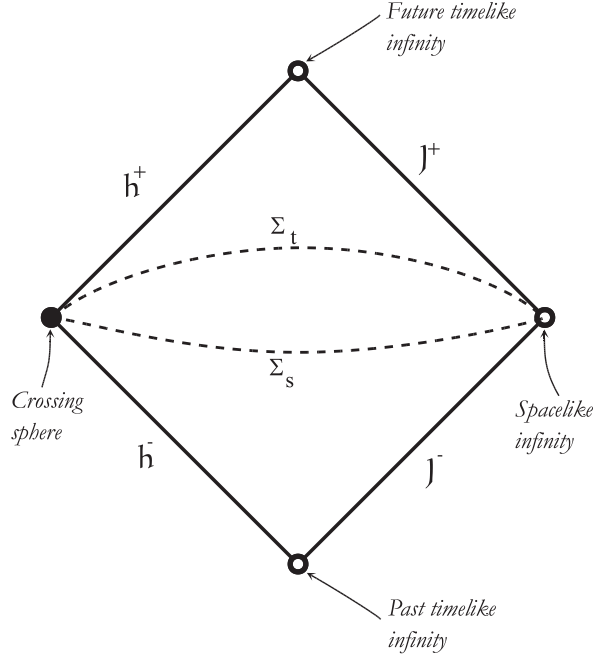,width=8cm}
\caption{The Penrose compactification of block $I$, with two
  hypersurfaces $\Sigma_s$ and $\Sigma_t$, $t>s$.}
\label{PenCompFigure}
\end{figure}
The Penrose compactification of block $I$ is then the space-time
\[ \left( \overline{\ExtB} \, ,~ \hat{g} \right) \,
,~\overline{\ExtB} = \ExtB \cup \hor^+ \cup S^2_\mathrm{c} \cup \hor^-
\cup \scri^+ \cup \scri^- \, ,\]
$\hat{g}$ being defined by (\ref{RescaledMetric}). In spite of the
terminology used, the compactified space-time is not compact. There
are three ``points'' missing to the boundary~: $i_+$, or future
timelike infinity, defined as the limit point of uniformly timelike
curves as $t\rightarrow +\infty$, $i_-$, past timelike infinity,
symmetric of $i_+$ in the distant past, and $i_0$, spacelike infinity,
the limit point of uniformly spacelike curves as $r\rightarrow
+\infty$. These ``points'' are singularities of the rescaled metric. See
Figure \ref{PenCompFigure} for a representation of the compactified
block $I$.
\\

Writing 
\[ U=e^{-\kappa_+(t-\hat{z}(t,\theta))},\, V=e^{\kappa_+(t+\hat{z}(t,\theta))}\]
and suppressing two dimensions gives the picture of the collapse of Figure \ref{DomInt}
(see Section \ref{sec3.5}).
\end{appendix}

\end{document}